%% file: OSconj.tex
\documentclass[11pt]{article}
\usepackage{amsmath,amssymb,graphicx,graphics}
\newtheorem{lem}{Lemma}
\newtheorem{conj}{Conjecture}
\newtheorem{defi}{Definition}
\newtheorem{thm}{Theorem}
\newtheorem{cor}{Corollary}

\newtheorem{prop}{Proposition}

\def\odots{\reflectbox{\text{$\ddots$}}}

\begin{document}
\title{ A complete list of lens spaces constructed by Dehn surgery I}
\author{Motoo Tange}
\date{}
\maketitle
\abstract
Berge in \cite{[B]} defined doubly primitive knots, which yield lens spaces by Dehn surgery.
At the same paper he listed the knots into several types.
In this paper we will prove the list is complete when $\tau>1$.
The invariant $\tau$ is a quantity with regard to lens space surgery, which is defined in this paper.
Furthermore at the same time we will also prove that Table~6 in \cite{[tan2]} is complete
as Poincar\'e homology sphere surgery when $\tau>1$.

\section{Introduction}
\label{intro}

\subsection{Several necessary conditions for lens space surgery}
Let $K$ be a knot in a 3-manifold $M$.
Removing the open tubular neighborhood $\text{nbd}(K)$ of $K$, and gluing the solid torus $V$
by some map between the boundary of $V$ and $M-\text{nbd}(K)$ we 
obtain a new 3-manifold $M'$. 
By iterating some Dehn surgeries we can obtain all 3-manifolds.
However we can not obtain all 3-manifolds by single Dehn surgeries of knots in $S^3$.
When a lens space is given by a Dehn surgery of a knot,
certain restrictions are imposed on the lens space and the knot.
In this section we share with us some of the restrictions.

Let $K$ be a knot in a homology sphere $Y$.
The manifold $Y_p(K)$ stands for Dehn surgery of $K$ with slope $p$.
We define lens space $L(p,q)$ to be $S^3_{p/q}(\text{unknot})$.
The problem of when $S^3_p(K)$ is lens space and when lens space $L(p,q)$ is
obtained from Dehn surgery of $K$ is fundamental.

Kadokami and Yamada proved the following by using torsion invariant.
\begin{thm}[\cite{[KaYa]}]
\label{ky}
Let $Y$ be a homology sphere.
If $K\subset Y$ is a knot and $Y_p(K)$ is lens space $L(p,q)$, then
the Alexander polynomial is the following form:
$$\Delta_K(t)=\frac{(t^{hg}-1)(t-1)}{(t^{g}-1)(t^h-1)}\ \ \ (t^p-1)$$
up to multiplication of $\pm t^{\pm1}$,
where $h,g$ are coprime integers and satisfy $hg=\pm1(p),$ and $h^2=\pm q\ (p)$.
\end{thm}
We can see that in the case where there exist $h$ and $g$ with $\gcd(h,g)=1$ and $hg+1\le p$, the Alexander polynomial is the same as the polynomial of $(h,g)$-torus knot $T_{h,g}$.
This type corresponds to type (I),(II) in Table~\ref{rs}.
We have to notice that any coefficient in Theorem~\ref{ky} is equivalent to the formula by using
correction term in \cite{[OS2]}.

Ozsv\'ath and Szab\'o \cite{[OS2]}, and Greene \cite{[Gr]} showed that the Seifert genus $g(K)$ of any knot $K$ yielding lens spaces has
an inequality $2g(K)-1\le p$ (Ozsv\'ath-Szab\'o) and $2g(K)-1\le p-2\sqrt{\frac{4p+1}{5}}$ (Greene),
where $p$ is the order of $H_1$ of the lens space.
The similar inequality holds for $\Sigma(2,3,5)$ (see \cite{[tan2]}).

Let $Y$ be a homology sphere.
A Dehn surgery $Y_p(K)$ has the core circle of the surgery which is called the circle {\it dual knot}.
In the case where $Y_p(K)$ is lens space $L(p,q)$,
the first homology class of the dual knot $\tilde{K}$ assigns an integer $k$.
In fact since there is a core circle by genus $1$ Heegaard splitting of the lens space,
we have only to take $k$ as the difference $[\tilde{K}]=k[c]$, where $c$ is the core circle
so that $q=k^2(p)$ holds.
However the assignment has two ambiguities: the choice of the non-trivial core circle in the two solid tori
and the choice of the orientation of the knot.
Thus we consider the integer $k$ as the set ${\mathcal D}(p,K):= \{k,-k,k^{-1},-k^{-1}\}\subset{\mathbb Z}/p{\mathbb Z}$.
Here we call it {\it the dual class invariant}.
The integers $h,g$ in Theorem~\ref{ky} are elements of ${\mathcal D}(p,K)$.
Abstractly we define ${\mathcal D}(p,k)$ to be the set $\{k,-k,k^{-1},-k^{-1}\}$.

Here we define a Laurent polynomial as follows.
Let $T_{r,s}$ be the $(r,s)$-torus knot.
\begin{defi}
Let $(p,k)$ be a coprime integers pair, and $h,g$ two coprime representations in ${\mathcal D}(p,k)$ with $hg=\pm1\ (p)$.
Adding $\pm t^n(t^p-1)$ several times to the symmetrized Alexander polynomial $\Delta_{T_{h,g}}(t)$,
we can make any exponent of the terms of the polynomial change a term between $-\lfloor\frac{p}{2}\rfloor+1$ and $\lfloor\frac{p}{2}\rfloor$.
Then we get a Laurent polynomial $\tilde{\Delta}_{p,k}(t)$.
If $p$ is even and the $p/2$-th term $a_{\frac{p}{2}}$ is not zero, then we modify $\tilde{\Delta}_{p,k}(t)$ into
$$\tilde{\Delta}_{p,k}(t)+\frac{a_{\frac{p}{2}}}{2}t^{-\frac{p}{2}}-\frac{a_{\frac{p}{2}}}{2}t^{\frac{p}{2}}$$
In this way we get a symmetric Laurent polynomial and denote it by $\Delta_{p,k}(t)$.
\end{defi}
We define $g(p,k_1)$ to be the degree of $\Delta_{p,k_1}(t)$.

By Fintushel and Stern's work in \cite{[fs]} conversely if there exists a coprime integer pair $(p,k)$ satisfying
$q=k^2(p)$, one can realize a Dehn surgery $L(p,q)=Y_p(K)$ on a homology sphere $Y$ with $k\in{\mathcal D}(p,K)$.
Hence one of lens space surgery problems is to consider when for a coprime pair $(p,k)$
there exists a knot $K$ in $S^3$ (in general in a fixed homology sphere) with $k\in {\mathcal D}(p,K)$.

Our starting point is a coprime pair $(p,k)$ then we call $(p,k)$  {\it the initial data} for 
lens space surgery by taking the minimal $k$ in ${\mathcal D}(p,k)$.
Here any element in ${\mathcal D}(p,k)$ is reduced between $0$ and $p-1$.
For an initial data $(p,k_1)$ if there exists a lens space surgery on a homology sphere $Y$
such that $k\in{\mathcal D}(p,K)$,
where $K$ is the knot in $Y$, we call $(p,k)$ the initial data {\it realized by the lens space surgery}.

We introduce an explicit coefficient formula of $\Delta_{p,k}(t)$, which is computed by Theorem~\ref{ky}, that 
was proved in \cite{[tan1]}.
Let $(p,k)$ be a coprime integer satisfying $0<k<p$ and $k'$ the inverse of $k\mod p$ satisfying $0<k'<p$.
We use this notation in other places as long as we do not indicate.
We put $m:=\frac{kk'-1}{p}$, $q=k^2 (p)$, $c:=\frac{(k+1-p)(k-1)}{2}$, and 
$$\Phi_{p,q}^l(k)=\#\{j\in\{1,2,\cdots k'\}|[qj-l]_p\in \{1,2,\cdots,k\}\},$$
where $[\alpha]_p$ stands for the reduction of any integer $\alpha$ to $0\le [\alpha]_p<p$.
\begin{thm}[\cite{[tan2]}]
Let $(p,k)$ be a coprime integer pair and $a_i$ the $i$-th coefficient of the Laurent polynomial $\tilde{\Delta}_{p,k}(t)$.
Then $a_i\ (-\lfloor\frac{p}{2}\rfloor+1\le i\le \lfloor\frac{p}{2}\rfloor)$ has the following explicit formula
\begin{equation}
a_i=-m+\Phi_{p,q}^{ki+c}(k).
\label{myf}
\end{equation}
\end{thm}
Hence if a knot $K\subset S^3$ yields lens space $L(p,q)$, the symmetrized Alexander polynomial
$a_i(K)$ is computed by the formula.
From the genus bound $2g(K)-1\le p$ by Ozsv\'ath and Szab\'o we can get the following.
\begin{cor}
For a lens space surgery $L(p,q)=S^3_p(K)$ with initial data $(p,k)$ we have
$$\Delta_K(t)=\Delta_{p,k}(t).$$
In particular $g(K)=g(p,k_1)$ holds.
\end{cor}
When we replace $S^3$ with $\Sigma(2,3,5)$, we get the same equality.
In this paper we assume that the genus of any knot yielding lens space surgery satisfies $2g(K)\le p$.
Ozsv\'ath and Szab\'o's result does not deny the $2g(K)-1=p$ case.
The author was informed that Dr. Greene showed that at least the $2g(K)-1=p$ case can be ruled out as long as lens space surgery over $S^3$.

Ozsv\'ath and Szab\'o proved the new restrictions in \cite{[OS1],[OS2]}.
\begin{thm}[\cite{[OS1]}]
\label{alter}
If $K\subset S^3$ is a knot and $S^3_p(K)$ is lens space, then
the symmetrized Alexander polynomial $\Delta_K(t)$ is the following form:
$$\Delta_K(t)=(-1)^m+\sum_{j=1}^m(-1)^k(t^{n_j}+t^{-n_j}),$$
where the sequence $\{n_j\}$ is increasing natural numbers: $n_1<n_2<\cdots<n_m$.
\end{thm}
The same assertion is also satisfied for $\Sigma(2,3,5)$ in place of $S^3$.
\begin{thm}[\cite{[OS2]}]
\label{pos}
If $K\subset S^3$ is a knot and $S^3_p(K)$ is a lens space $L(p,q)$, then
the following quantity
$$t_i:=\begin{cases}
d(L(p,q),ki+c)-d(L(p,1),i) &\text{if } 2|i|\le p\\
0 & \text{otherwise.}
\end{cases}$$
is non-negative integer and is
coincident with the torsion invariant $t_i(K)$ for any integer $i$.
Here $t_i(K)$ is the $i$-th Turaev torsion of $S^3_0(K)$.
\end{thm}
Some of lens space surgeries on $\Sigma(2,3,5)$ fail to this theorem, 
for example $L(22,3)$ is obtained as a lens space surgery on $\Sigma(2,3,5)$ and $t_{11}=-2$ holds.
 
By using these restrictions, we can make a sharp distinction as in \cite{[OS1]}.
However either of Theorem~\ref{alter} and \ref{pos} is not perfect so as to distinguish lens space surgery.

\subsection{Berge's examples}
Berge in \cite{[B]} defined a class of knots in $S^3$, which is most important to lens space surgery so far.
\begin{defi}
Let $K$ be a knot in $S^3$.
We call $K$ a doubly primitive knot if $K$ isotopic to a knot $L$ in a genus $2$ Heegaard surface in $S^3$
and both classes induced in $\pi_1(V_i)$ are primitive elements, where $V_i\ (i=1,2)$ are
the genus $2$ handlebodies.
\end{defi}
This definition can be easily extended to any knot in a homology sphere with genus $2$ Heegaard decomposition,
which satisfies the same conditions.

He proved that all doubly primitive knots yield lens spaces by Dehn surgery with an integer slope.
In the other words this condition is a sufficient condition for lens space surgery.
He conjectures that doubly primitive knots are all knots yielding lens spaces.
This conjecture is still open.

He also listed the doubly primitive knots in \cite{[B]}.
Thus any of the list is realized by a lens space surgery on $S^3$, however
it is open question that this list is complete.

In \cite{[Ras]} J. Rasmussen rearranged the list to the equivalent one below (Table~\ref{rs}).
\begin{table}[thbp]
$$\begin{array}{|c|c|c|}\hline
\text{I,II} & p=ik\pm 1\ (k^2)& \gcd(i,k)=1,2\\\hline
\text{III}_\pm & p=(2k\mp 1)J\ (k^2) & J\in{\mathbb Z}; J|k\pm1;\frac{k\pm1}{J}:\text{ odd}\\\hline
\text{IV}_\pm & p=(k\pm1)J\ (k^2) & J:\text{odd}; J|2k\mp1\\\hline
\text{V}_\pm & p=(k\pm1)J\ (k^2) & J\in{\mathbb Z}; J|k\pm1; J:\text{ odd}\\\hline
\text{VII, VIII} & k^2\pm k\pm1=0(p) & \\\hline
\text{IX} & p=22J^2+9J+1& k=11J+2, J\in{\mathbb Z}\\\hline
\text{X} & p=22J+13J+2 & k=11J+3, J\in{\mathbb Z}\\\hline
\end{array}$$
\caption{Berge's list}
\label{rs}
\end{table}

\subsection{Poincar\'e homology sphere version of Berge's examples}
Next we introduce the result \cite{[tan2]}, which is a sufficient condition for lens space surgery on Poincar\'e homology sphere.
\begin{thm}[\cite{[tan2]}]
Let $(p,k)$ be any of Table~\ref{po}.
Then $(p,k)$ is realized by a lens space surgery on $K\subset \Sigma(2,3,5)$.
Furthermore we can take a doubly primitive knot as a knot for the surgery.
\end{thm}
\begin{table}[tbp]
$$\begin{array}{|c|c|c|c|}\hline
\text{type} & p & k\\\hline
\text{A}_1 & 14J^2++7J+1 & 7J+2\\\hline
\text{A}_2 & 20J^2+15J+3 & 5J+2\\\hline
\text{B} & 30J^2+9J+1 & 6J+1\\\hline
\text{C}_1 & 42J^2+23J+3 & 7J+2\\\hline
\text{C}_2 & 42J^2+47J+13 & 7J+4\\\hline
\text{D}_1 & 52J^2+15J+1 & 13J+2\\\hline
\text{D}_2 & 52J^2+63J+19 & 13J+8\\\hline
\text{E}_1 & 54J^2+15J+1 & 27J+4\\\hline
\text{E}_2 & 54J^2+39J+7 & 27J+10\\\hline
\text{F}_1 & 69J^2+17J+1 & 23J+3\\\hline
\text{F}_2 & 69J^2+29J+3 & 23J+5\\\hline
\text{G}_1 & 85J^2+19J+1 & 17J+2\\\hline
\text{G}_2 & 85J^2+49J+7 & 17J+5\\\hline
\text{H}_1 & 99J^2+35J+3 & 11J+2\\\hline
\text{H}_2 & 99J^2+53J+7 & 11J+3\\\hline
\text{I}_1 & 120J^2+16J+1 & 12J+1\\\hline
\text{I}_2 & 120J^2+20J+1 & 20J+2\\\hline
\text{I}_3 & 120J^2+36J+3 & 12J+2\\\hline
\text{J} &  120J^2+104J+22 & 12J+5\\\hline
\text{K} & 191 & 15\\\hline
\end{array}$$
\caption{The Poincar\'e homology sphere version of Berge's examples.}
\label{po}
\end{table}

Theorem~\ref{ky} and Theorem~\ref{alter} also hold for lens space surgery on the Poincar\'e homology sphere.
However Theorem~\ref{pos} does not hold because the initial data $(22,5)$ does not satisfy Theorem~\ref{pos} but
admits lens space surgery $\Sigma(2,3,5)_{22}(K)$.

Then the following is conjectured.
\begin{conj}
\label{s3conj}
Let $(p,k)$ be a coprime integer.
Suppose that the polynomial $\Delta_{k,p}(t)$ satisfies the alternating condition in Theorem~\ref{alter}, 
and the pair satisfies the condition in Theorem~\ref{pos}.
Then $(p,k)$ is realized by a lens space surgery on $S^3$.
\end{conj}
This conjecture is equivalent to Conjecture~1.12 in \cite{[OS2]}.
\begin{conj}
\label{poconj}
Let $(p,k)$ be a coprime integer.
Suppose that the polynomial $\Delta_{p,k}(t)$ satisfies the alternating condition in Theorem~\ref{alter}, 
Then $(p,k)$ is realized by a lens space surgery on $S^3$ or $\Sigma(2,3,5)$.
\end{conj}
If Conjecture~\ref{s3conj} is true, then it means that Table~\ref{rs} is the complete list of doubly primitive knots.
If Conjecture~\ref{poconj} is true, then it means that Table~\ref{rs} and \ref{po} are the complete list of 
doubly primitive knots in $S^3$ and $\Sigma(2,3,5)$.
Even if these conjectures are proven, the problem of whether doubly primitive knots are all the knots yielding lens space or not
remains open.
The aim of this paper is a partial contribution to Conjecture~\ref{s3conj} and \ref{poconj}.
The meaning of ``partiality" is described  later.


\section{Preliminaries}
\subsection{Quadratic relations for lens space surgery}
\label{qeq}
In this section we shall define a quadratic relation in ${\mathbb Z}/p{\mathbb Z}$.
Let $(p,k_1)$ be a coprime integer pair with $p$ positive and $k_2$ the reduced element with
$k_1\le k_2<\frac{p}{2}$ in $\{k_1,[-k_1]_p,[k_1^{-1}\bmod p]_p,[-k_1^{-1}\bmod p]_p\}$.

We denote by $[[\alpha]]_p$ the reduction between $-\lfloor\frac{p}{2}\rfloor+1$ and $\lfloor\frac{p}{2}\rfloor$.
We define $q_2$ to be $[[k_2^2]]_p$ and $a$ to be $||q_2|-k_2|$.
Hence we can find a relation on ${\mathbb Z}/p{\mathbb Z}$
\begin{equation}
ak_1^2+\epsilon_1k_1+\epsilon_2=0\ (p),
\label{asseq}
\end{equation}
where $\epsilon_i=\pm1$.
This relation is called {\it the associated (quadratic) relation} in this paper.
Note that relations of $(p,k_1)$ in the form of (\ref{asseq}) are not always unique for any data.
For example $(43,12)$ has two relations $2k_1^2+k_1+1=0\ (43)$ and $5k_1^2+k_1-1=0\ (43)$.

\begin{prop}
The $a=0$ case is equivalent to $k_1=k_2=1$ then this initial data can be
realized by $p$-surgery of the unknot.
\end{prop}
{\bf Proof.} \
Therefore $k_1=\pm1\ (p)$ holds.
\hfill$\Box$\\

The $a=1$ case is (VII) and (VIII) in Table~\ref{rs}.
Then we assume that $a\ge 2$.

\begin{lem} 
Let $(p,k_1)$ be an initial data with $k_1>1$.
The associated quadratic relation of $(p,k_1)$ is the relation which
the term $a$, which is the coefficient of degree $2$ of the relation, is minimal 
among the relations having form of (\ref{asseq}).

Any relation of the form of (\ref{asseq}) with the minimal degree $2$ term is the associated relation.
\end{lem}
{\bf Proof}
The former of the assertion is obvious from the definition of $a$.
We show the latter part.
Suppose that $(p,k_1)$ satisfies $ak_1^2+\epsilon_1k_1+\epsilon_2=0\ (p)$ and $ak_1^2+\epsilon'_1k_1+\epsilon_2'=0\ (p)$,
and $\epsilon_1<\epsilon_1'$.
Then we have $(\epsilon_1'-\epsilon_1)k_1=\epsilon_2-\epsilon_2'\ (p)$.
Since we have $0<2k_1<p$, $2k_1=2$ holds.
Then we have $k_1=1$.
If $\epsilon_1=\epsilon_1'$, then this relation is the same.\hfill$\Box$

Now let $(p,k_1)$ be a coprime pair as above.
We take a relation $ak_1^2+\epsilon_1k_1+\epsilon_2=0\ (p)$ and 
define an integer $n$ to be $ak_1^2+\epsilon_1k_1+\epsilon_2-np=0$.
Solving the quadratic equation, we get 
$$k_1=\frac{-\epsilon_1+X}{2a}$$
where $X$ is positive integer and we put $X^2=1-4a(\epsilon_2-np)$.
Hence we have $p=\frac{X^2-D}{4an}$, where we put $D=1-4a\epsilon_2$.
Here let $\tau,\gamma$ be integers satisfying $X=2an\tau+\gamma$ and $0\le\gamma<2an$.
Hence we have   
\begin{equation}
k_1=n\tau+\frac{\gamma-\epsilon_1}{2a}=n\tau+\gamma',
\label{kfor}
\end{equation}
where we put $\gamma'=\frac{\gamma-\epsilon_1}{2a}$\ ($-\frac{\epsilon_1}{2a}\le \gamma'<n-\frac{\epsilon_1}{2a}$),
thus $\gamma=\epsilon_1\ (2a)$ holds.
Then $p$ is described as follows:
\begin{equation}
p=an\tau^2+\gamma \tau+\frac{\gamma^2-D}{4an}.
\label{peq}
\end{equation}
Since $\gamma^2-4an\frac{\gamma^2-D}{4an}=D=1-4a\epsilon_2$, if $\epsilon_2=1$, then the quadratic function has a positive value for any $\tau$
as long as $a>0$.
If $\epsilon_2=-1$, then considering $p=an\left(\tau+\frac{\gamma}{2an}\right)^2-\frac{D}{4an}$, we have $an\left(0+\frac{\gamma}{2an}\right)^2-\frac{D}{4an}=\frac{a(\gamma')^2+\epsilon_1\gamma'-1}{n}\ge 0$ as long as $a>0$ and $\gamma'>1$.
Thus we have $an\left(-2+\frac{\gamma}{2an}\right)^2-\frac{D}{4an}=4an-2\gamma+\frac{a(\gamma')^2+\epsilon_1\gamma'-1}{n}\ge 0$
\begin{defi}
Let $(p,k_1)$ be an initial data.
We denote by $\tau$ the parameter computed as above by using the associated quadratic relation.
If the parameter $\tau\ge 2$ holds we call $(p,k_1)$ a stable (initial) data.
\end{defi}

Here we list the associated relation for lens surgeries in Table~\ref{rs} and (\ref{po}),
except (I), (II), ($\text{V}_{\pm}$), (VII), and (VIII).
\begin{table}[htbp]
$$\begin{array}{|c|c|c|c|c|}\hline
\text{type} & J & \text{the associated relation} &n\\\hline
\text{III}_{\pm} &|J|\ge 2&  ak_1^2\pm k_1-1=0\ \ n|a-2& \\\hline
\text{IV}_{\pm} &|J|\ge 2& ak_1^2\pm k_1-1=0\ \ n|a-2&  \\\hline
\text{IX,X} & |J|\ge 2&2k_1^2\pm k_1+1=0&11\\\hline
\text{A}_1 &|J|\ge 2& 2k_1^2-\delta k_1+1=0\ &7 \\\hline
\text{A}_2 & |J|\ge 2,J=1&4k_1^2-\delta k_1+1=0\ &5\\\hline
\text{E}_1,\text{E}_2 &|J|\ge 2& 2k_1^2-\delta k_1-1=0\ &27 \\\hline
\text{F}_1,\text{F}_2 &|J|\ge 2& 3k_1^2-\delta k_1-1=0\ &23\\\hline
\text{D}_1 &|J|\ge 2& 4k_1^2-\delta k_1-1=0\ &13\\\hline
\text{D}_2 &|J|\ge 2,J=1& 4k_1^2-\delta k_1-1=0\ &13\\\hline
\text{G}_1,\text{G}_2 &|J|\ge 2,J=1&5k_1^2-\delta k_1-1=0&17\\\hline
\text{C}_1,\text{C}_2 &|J|\ge 2,J=1 &6k_1^2-\delta k_1-1=0 &7\\\hline
\text{H}_1,\text{H}_2 &|J|\ge 2,J=1& 9k_1^2-\delta k_1-1=0&11\\\hline
\text{B} &J\ge 2& (15J^2-18J-5)k_1^2-k_1-1=0&18J^2-21J-7\\\hline
\text{B} & J\le -2&(15J^2+27J+6)k_1^2-k_1+1=0&18J^2+33J+8\\\hline
\text{B} &J=1& 8k_1^2+k_1+1=0&10 \\\hline
\text{I}_1 & J\ge 2&(40J^2-28J-3)k_1^2-k_1-1=0&48J^2-32J-5 \\\hline
\text{I}_1 &J\le -1 &(40J^2+12J-1)k_1^2-k_1-1=0&48J^2+16J-1 \\\hline
\text{I}_2 & J\ge 2 &(60J^2-50J-5)k_1^2-k_1-1=0&200J^2-160J-23\\\hline
\text{I}_2 &J\le -2& (60J^2+70J+6)k_1^2-k_1+1=0&200J^2+240J+27 \\\hline
\text{I}_3 &J\ge 2& (60J^2-42J-9)k_1^2-k_1-1=0&72J^2-48J-13 \\\hline
\text{I}_3 &J\le -2& (60J^2+78J+12)k_1^2-k_1+1=0&72J^2+96J+17\\\hline
\text{J} &J\ge 1& (40J^2+28J+6)k_1^2+k_1-1=0&48J^2+32J+7 \\\hline
\text{J} &J\le -2& (40J^2-12J-4)k_1^2+k_1-1=0&48J^2+80J+27\\\hline
\text{K} && 22k_1^2+k_1+1=0&26\\\hline
\end{array}$$
\caption{The associated relations. ($\delta=\text{sgn}(J)$).
}
\end{table}

\begin{table}[htbp]
$$\begin{array}{|c|c|c|}\hline
(p,k_1) & (\text{type},J)& \text{the associated relation} \\\hline
(43,12) & (\text{A}_1,-2), (\text{F}_2,-1),(\text{G}_2,-1) & 2k_1^2+k_1+1=0\\\hline
(8,3) & (\text{A}_1,-1), (\text{A}_2,-1),(\text{C}_2,-1),(\text{D}_2,-1)& 2k_1^2-k_1+1=0\\\hline
(38,7) & (\text{A}_2,1), (\text{D}_1,-1), (\text{J},-1) & 4k_1^2-k_1+1=0\\\hline
(53,8) & (\text{A}_2,-2), (\text{F}_1,-1),(\text{H}_2,-1) & 4k_1^2+k_1+1=0\\\hline
(68,9) & (\text{C}_1,1), (\text{D}_1,1) & 6k_1^2-k_1-1=0\\\hline
(125,12) & (\text{C}_1,-2) & 6k_1^2+k_1-1=0\\\hline
(102,11) & (\text{C}_2,1) & 6k_1^2-k_1-1=0\\\hline
(87,10) & (\text{C}_2,-2), (\text{I}_2,-1),(\text{F}_1,1) & 6k_1^2+k_1-1=0\\\hline
(179,24) & (\text{D}_1,-2) & 4k_1^2+k_1-1=0\\\hline
(134,21) & (\text{D}_2,1) & 4k_1^2-k_1-1=0\\\hline
(101,18) & (\text{D}_2,-2), (\text{I}_2,-1),(\text{F}_2,1) & 4k_1^2+k_1-1=0\\\hline
(187,50) & (\text{E}_1,-2) & 2k_1^2+k_1-1=0\\\hline
(145,44) & (\text{E}_2,-2)& 2k_1^2+k_1-1=0\\\hline
(243,43) & (\text{F}_1,-2) & 3k_1^2+k_1-1=0\\\hline
(221,41) & (\text{F}_2,-2) & 3k_1^2+k_1-1=0\\\hline
(303,32) & (\text{G}_1,-2) & 5k_1^2+k_1-1=0\\\hline
(249,29) & (\text{G}_2,-2) & 5k_1^2+k_1-1=0\\\hline
(141,22) & (\text{G}_2,1),(\text{I}_2,1) & 5k_1^2-k_1-1=0\\\hline
(329,20) & (\text{H}_1,-2) & 9k_1^2+k_1-1=0\\\hline
(137,13) & (\text{H}_1,1), (\text{I}_1,1) & 9k_1^2-k_1-1=0\\\hline
(297,19) & (\text{H}_2,-2) & 9k_1^2+k_1-1=0\\\hline
(159,14) & (\text{H}_2,1), (\text{I}_2,1) & 9k_1^2-k_1-1=0\\\hline
\end{array}$$
\caption{The associated relations having lens space surgeries on $\Sigma(2,3,5)$ with $\tau=1$.}
\label{tau1d}
\end{table}


Here we state Main theorem in this paper.
\begin{thm}[Main theorem]
Suppose that $(p,k_1)$ is a stable data.
If $\Delta_{p,k}(t)$ satisfies the alternating condition in Theorem~\ref{ky} and Theorem~\ref{alter}, then
$(p,k)$ is realized by a lens space surgery of types below:
$$I,II,III_\pm,IV_\pm,VII,VIII,IX,X,A_i,C_i,D_i,E_i,F_i,G_i,H_i,$$
where $i=1$ or $2$.
\end{thm}
The partiality in Section~\ref{intro} means the stability.

\begin{lem}
If $\tau\ge 2$ then, we have $k_2=ak_1+\epsilon_1$.
\end{lem}
{\bf Proof.}
Putting $\kappa=ak_1+\epsilon_1$, we have $p=\kappa(\tau+\frac{\gamma'}{n})+\frac{\epsilon_2}{n}$ and 
$$-\frac{\epsilon_1}{2an}<\frac{\gamma'}{n}<1-\frac{\epsilon_1}{2an}.$$
Hence $\kappa\le \frac{p}{2}$ and $\kappa k_1=\pm1\ (p)$, therefore we have $\kappa=k_2$.\hfill$\Box$\\

In this section we shall concentrate on seeing the coefficients of the polynomial $\tilde{\Delta}_{p,k_1}(t)$ 
for the stable data $(p,k_1)$.

Table~\ref{tau1d} is the list of lens spaces with $\tau=1$ in Table~2 in \cite{[tan2]}.
These lens spaces are not dealt with in this paper.
We conjecture that all lens spaces obtained by Dehn surgeries on $\Sigma(2,3,5)$ with $\tau=1$ are 
included in Table~\ref{tau1d}.

\medskip
\subsection{The $\gamma'=0$ or $\gamma'=n$ case.}
First we classify all stable data realized lens space surgery on $S^3$ or $\Sigma(2,3,5)$ in the cases of $(\epsilon_1,\gamma')=(1,0),(-1,n)$.
\begin{thm}
In the case of $(\epsilon_1,\gamma')=(1,0)$
the data is realized by type (I).
\end{thm}
{\bf Proof.} \
In the case of $(\epsilon_1,\gamma')=(1,0)$ we have $\gamma=1$.
From integrality of $p$ and Formula (\ref{peq}), $\frac{1}{n}$ is an integer thus $n=1$.
Therefore $(p,k_1)=(a\tau^2+\tau+\epsilon_2,\tau)$ holds.
Thus we have $k_2=a\tau+1$.
From this and Theorem~12 in \cite{[tan1]}, $\gcd(k_1,k_2)=1$ or $2$.
Therefore this data is realized by type (I).

In the case of $(\epsilon_1,\gamma')=(-1,n)$ we have $\gamma=2an-1$.
In the same way from integrality of $p$, we get $n=1$.
Thus $(p,k)=((\tau+1)(a(\tau+1)-1)+\epsilon_2,\tau+1)$.
Hence we have $k_2=a(\tau+1)-1$.
From Theorem~12 in \cite{[tan1]}, $\gcd(k_1,k_2)=1$ or $2$.
This data is realized by type (I).\hfill$\Box$\\

Then the inverse element $k_1'\bmod p$ is reduced as follows:
$$k_1':=\begin{cases}
k_2& \epsilon_2=-1\\
p-k_2& \epsilon_2=1.
\end{cases}$$

\subsection{The case of $k_1=1,2,3$.}
In the case of $k_1=1$, $k_2$ is also $1$ and the initial data is realized by the unknot surgery $L(p,1)=S^3_p(\text{unknot})$.

In the case of $k_1=2$, $p$ is odd.
Then we put $p$ by $2n+1$.
Hence $k_2=n$ holds.
From this and Theorem~12 in \cite{[tan1]}, $n$ is odd.
In this case each of data is realized by the $(2,n)$-torus knot surgery.

In the case of $k_1=3$, $p$ is $3n\pm1$.
Then $k_2=n$ holds.
From this and Theorem~12 in \cite{[tan1]}, $\gcd(3,n)=1$ holds.
In this case each of data is realized by the $(3,n)$-torus knot surgery.

Therefore now we may assume that $k_1\ge 4$ holds.
\subsection{The $n=1$ case.}
If $n=1$ holds, then we have $ak_1^2+\epsilon_1k_1+\epsilon_2=p$.
This data is realized by $(k_1,ak_1+\epsilon_1)$-torus knot surgery.
This is type (I) in Berge's list.

Therefore we assume that $n\ge 2$.
\subsection{The coefficients of the Alexander polynomial} 
\label{coeffformula}

For a stable data $(p,k_1)$ the reduced square $q_2$ satisfies $q_2=\epsilon_1k_2-\epsilon_2a$ and we have $k_1'=[-\epsilon_2k_2]_p$.
We put $q_1:=[k_1^2]_p$, $c:=\frac{(k_1+1-p)(k_1-1)}{2}$ and $\bar{g}=k'_1(c-q_1)+1=-\epsilon_2k_2(c-q_1)+1\ (p)$.
Then we get
\begin{equation}
-2\bar{g}=k_1+k_1'-2=k_1-\epsilon_2k_2-2\ (p)
\label{barg}
\end{equation}
We define a function $A$ to be 
$$A(n_1,n_2,n_3)=\tilde{a}_{-\bar{g}+n_1k_1+n_2 k_1'+n_3},$$
where the coefficients $\tilde{a}_i=a_{[[i]]_p}$ for any integer $i\in {\mathbb Z}$.
The coefficients $\tilde{a}_i$ are regarded as the function that cyclically lifts
the coefficients of $\Delta_{p,k_1}(t)$ to ${\mathbb Z}$.

Here we prepare the following lemma.
\begin{lem}
\label{differ}
Let $(p,k_1)$ be any initial data.
We have
$$A(m_1,m_2,m_3)-A(m_1,m_2-1,m_3)=E_{k_1'}(m_2q_2+m_3k_1'+m_1+1)-E_{k_1'}(m_2q_2+(m_3-1)k'_1+m_1+1)$$
$$=E_{p-k_1'}(m_2q_2+(m_3-2)k_1'+m_1+1)-E_{p-k_1'}(m_2q_2+(m_3-1)k_1'+m_1+1),$$
where the function $E_x(y)$ is defined to be
$$E_x(y)=\begin{cases}
1 & y\equiv 1,2,\cdots,[x]_p\mod p\\
0 & \text{otherwise}.
\end{cases}$$
\end{lem}
{\bf Proof.}
By definition we have
$$\Phi_{p,q_1}^{k_1i+c}(k_1)=\sum_{j=1}^{k_1}E_{k_1'}(q_2j+k_1'i+q_2c).$$
Now we put $\sigma=-\bar{g}+m_1k_1+m_3$
\begin{eqnarray}
&&A(m_1,m_2,m_3)-A(m_1,m_2-1,m_3)\nonumber\\
&=&\sum_{j=1}^{k_1}E_{k_1'}(q_2j+k_1'(\sigma+m_2k_1')+q_2c)-\sum_{j=1}^{k_1}E_{k_1'}(q_2j+k_1'(\sigma+(m_2-1)k'_1)+q_2c)\nonumber\\
&=&\sum_{j=1}^{k_1}E_{k_1'}(q_2j+k_1'(\sigma+m_2k_1')+q_2c)-\sum_{j=0}^{k_1-1}E_{k_1'}(q_2j+k_1'(\sigma+m_2k_1')+q_2c)\nonumber\\
&=&E_{k_1'}(m_2q_2+(\sigma+1)k_1'+q_2c)-E_{k_1'}(m_2q_2+\sigma k_1'+q_2c).
\label{siki1}
\end{eqnarray}
Here we have
\begin{eqnarray}
\sigma k_1'+q_2c&=&(-\bar{g}+m_1k_1+m_3)k_1'+q_2c\nonumber\\
&=&-(k_1'(c-q_1)+1)k_1'+m_1+m_3k_1'+q_2c=(m_3-1)k_1'+m_1+1\ (p).\nonumber
\end{eqnarray}
Therefore we have
$$(\ref{siki1})=E_{k_1'}(m_2q_2+m_3k_1'+m_1+1)-E_{k_1'}(m_2q_2+(m_3-1)k'_1+m_1+1).$$
Furthermore, by using the obvious relation $E_{k_1'}(x)+E_{p-k_1'}(x-k_1')=1$ we get
$$(\ref{siki1})=E_{p-k_1'}(m_2q_2+(m_3-2)k_1'+m_1+1)-E_{p-k_1'}(m_2q_2+(m_3-1)k'_1+m_1+1).$$
\hfill$\Box$\\

We can describe the formulas of Lemma~\ref{differ} as follows:
\begin{eqnarray}
&&A(m_1,m_2,m_3)-A(m_1,m_2-1,m_3)\nonumber\\
&=&E_{k_2}(m_2q_2-\epsilon_2m_3k_2+m_1+1)-E_{k_2}(m_2q_2-\epsilon_2(m_3-1)k_2+m_1+1).
\label{Aform}
\end{eqnarray}

\subsection{The two reductions: $q_1$ and $q_2$.}
Let $(p,k_1)$ be a stable data.
Suppose that $(p,k_1)$ satisfies a relation
$$ak_1^2+\epsilon_1k_1+\epsilon_2=0.\ (p)$$
Here we recall $a=-\epsilon_2 q_2+\epsilon_2 \epsilon_1k_2$.
Then we have 
\begin{equation}\
k_2=\epsilon_1q_2+\epsilon_1\epsilon_2 a.\label{kq}
\end{equation}
Hence
$$q_2=\epsilon_1ak_1-\epsilon_2a+1\text{ and }|q_2|=\epsilon_1q_2.$$
From 
$$\frac{|q_2|}{a}=k_1-\epsilon_1\epsilon_2+\frac{\epsilon_1}{a}$$
we have $\frac{|q_2|}{a}\ge 4-\epsilon_1\epsilon_2+\frac{1}{a}>2$, as a result
$\frac{a}{|q_2|}<\frac{1}{2}$ holds.
Furthermore from 
$\epsilon_1q_2=k_2-\epsilon_1\epsilon_2a\ge 4a+\epsilon_1-\epsilon_1\epsilon_2a\ge 3a-1>0$,
$$|q_2|=\epsilon_1q_2$$

Let $\eta$ be the integer satisfying $aq_1+\epsilon_1k_1+\epsilon_2=\eta p$ and $0<q_1<p$.
$$q_1=\frac{\eta p-\epsilon_1k_1-\epsilon_2}{a},$$
then $k_1^2-q_1=\frac{n-\eta}{a}p$ implies $\eta=n\ (a)$, namely $\eta=[n]_a$.

\subsection{A symmetry of the stable data}
Over all stable data $\{(p,k_1)\}$ we define an involution as follows:
$$(p,k_1)\to (p-(\gamma-an)(2\tau+1),k_1+n-2\gamma'):=(\underline{p},\underline{k}_1)$$
$$(\gamma',\gamma,X)\to (n-\gamma',2an-\gamma,X+2(an-\gamma)$$
$$(\epsilon_1,\epsilon_2)\to (-\epsilon_1,\epsilon_2)$$
$$(a,n,\tau,D)\to (a,n,\tau,D).$$
In fact we have $\underline{k}_1\ge k_1-n=n(\tau-1)+\gamma'\ge 0$ and $\underline{p}= an(\tau+1)^2-\gamma(\tau+1)+\frac{\gamma^2-D}{4an}$
\begin{eqnarray}
a\underline{k}_1^2-\epsilon_1\underline{k}_1+\epsilon_2&=&ak_1^2+\epsilon_1k_1+\epsilon_2+2ak_1(n-2\gamma')+a(n-2\gamma')^2-\epsilon_1(n-2\gamma')-2\epsilon_1k_1\nonumber\\
&=&n(p+(an-\gamma)(2\tau+1))=n\underline{p}.\nonumber
\end{eqnarray}
Hence $\underline{p}$ and $\underline{k}_1$ are relatively prime and $\underline{p}>0$.
Similarly putting
$$\underline{k}_2:=a\underline{k}_1-\epsilon_1=ak_1+an-2a\gamma'-\epsilon_1=k_2+an-\gamma-\epsilon_1$$
we have
$\frac{1}{n}(\underline{k}_1\underline{k}_2+\epsilon_1)=\underline{k}_2(\tau+1-\frac{\gamma'}{n})+\frac{\epsilon_1}{n}\ge 2\underline{k}_2+\frac{\epsilon_1}{n}$,
hence $\underline{k}_2\le\frac{1}{2}(\underline{p}-\frac{\epsilon_2}{n})$.
Therefore $\{\underline{k}_1,\underline{k}_2,p-\underline{k}_2,p-\underline{k}_1\}$ is the dual class invariant for the data $(\underline{p},\underline{k}_1)$
and $0<\underline{k}_1\le \underline{k}_2<\frac{p}{2}$ is satisfied.

\begin{lem}
\label{asspres}
Let $(p,k_1)$ be a stable data.
Let $(\underline{p},\underline{k}_1)$ be the data from $(p,k_1)$ via the involution.
Then $a\underline{k}_1^2-\epsilon_1\underline{k}_1+\epsilon_2=0(\underline{p})$ is the associated relation.
\end{lem}
{\bf Proof.} \
Since $\underline{q}_2=[[\underline{k}_2^2]]_{\underline{p}}=-\epsilon_1\underline{k}_2-\epsilon_2a$,
the absolute value $|\underline{q}_2|=-\epsilon_1\underline{q}_2=\underline{k}_2+\epsilon_1\epsilon_2a$.
Hence we have $||\underline{q}_2|-k_2|=a$.
Therefore $a\underline{k}_1^2-\epsilon_1\underline{k}_1+\epsilon_1=0(p)$ is the associated relation for $(p,\underline{k}_1)$.
\hfill$\Box$

The correspondence on non-stable data $(p,k_1)$ does not always preserve the associated relation in the same way as Lemma~\ref{asspres}.
Namely the correspondence is not involution.
For example $(p,k_1)=(191,15)$ is gotten the values $a=22$, $n=26$, $\tau=0$, $q_2=-73$.
Then the image of $(191,15)$ is $(102,11)$, but the associated relation of $(102,11)$ is
$6k_1^2-k_1-1=0(102)$ and $\tau=1$.
Moreover the image of $(102,11)$ is $(87,10)$ with the same associated relation as $(102,11)$
and $\tau=1$.
\subsection{The difference of $A(x)$.}
We define $A(x)$ to be 
$$A(\epsilon_1\epsilon_2a-1,\epsilon_1x,1).$$
By using Formula~(\ref{Aform}), we calculate the difference of $A(x)$ in this section.

\subsubsection{The $\epsilon_2=1$ case}
\begin{lem}
Let $p,k_1,a$ be the same parameters as above.
Then the differences
\begin{eqnarray}
&&A(\epsilon_1a-1,\epsilon_1x,1)-A(\epsilon_1a-1,\epsilon_1x-1,1)\nonumber\\
&=&E_{k_2}(\epsilon_1xq_2+k_2+\epsilon_1a)-E_{k_2}(\epsilon_1xq_2+\epsilon_1a)\label{e1}
\end{eqnarray}
are divided into $3$ cases as follows.
\end{lem}
\begin{enumerate}
\item[(a)] The condition (\ref{e1})$=-1$ is equivalent to the following.\\
($\epsilon_1=1$)
There exists an integer $\ell$ such that
$$x=\begin{cases}
\lfloor \frac{p\ell}{|q_2|}\rfloor+1\ \text{or}\\
-q_1j\ (p)\ \ (j=0,\cdots,a-1)\ \ \left(\lfloor \frac{p\ell}{|q_2|}\rfloor\right).\\
\end{cases}$$
($\epsilon_1=-1$)
There exists an integer $\ell$ such that
$$x=\lfloor\frac{p\ell}{|q_2|}\rfloor+1\text{ and }\neq -q_1j\ (p)\ \ (j=1,\cdots,a).$$
\item[(b)] The condition (\ref{e1})$=1$ is equivalent to the following.\\
($\epsilon_1=1$)
There exists an integer $\ell$ such that
$$x=\begin{cases}
\lfloor \frac{p\ell}{|q_2|}\rfloor,\ \ \text{and}\ \ \neq -q_1j\ (p)\ \ (j=0,1,\cdots,a-1),\\
-q_1j-1\ (p)\ \ (j=0,1,\cdots,a-1),\ \  (\lfloor \frac{p\ell}{|q_2|}\rfloor-1)\text{, or}\\
-q_1j+k_1\ (p)\ \ (j=0,1,\cdots,a-1),\ \  (\lfloor \frac{p\ell}{|q_2|}\rfloor-1).
\end{cases}$$
($\epsilon_1=-1$)
There exists an integer $\ell$ such that
$$x=\begin{cases}
-q_1j\ \ (p)\ \ (j=1,\cdots,a)\ \ \left(\lfloor\frac{p\ell}{|q_2|}\rfloor+1\right)\text{, or}\\
\lfloor \frac{p\ell}{|q_2|}\rfloor\text{ and }\neq-q_1j-1\ (p),\ -q_1j-k_1\ \ (j=1,\cdots,a).
\end{cases}$$
\item[(c)] (\ref{e1})$=0$ otherwise.
\end{enumerate}
{\bf Proof.} \
{\bf (a)} \
The condition (\ref{e1})$=-1$ is equivalent to $p\ell<\epsilon_1xq_2+\epsilon_1a\le p\ell+k_2$\\
$$\Leftrightarrow x-1\le\frac{p\ell}{|q_2|}<x+\frac{a}{q_2}.$$
If $\epsilon_1=1$ then this condition is 
$$x=\lfloor\frac{p\ell}{|q_2|}\rfloor+1$$
or 
$$x=-q_1j\ \ (p)\ \ \ (j=0,\cdots,a-1)\ \ \ \ (x=\lfloor\frac{p\ell}{|q_2|}\rfloor).$$
If $\epsilon_1=-1$ then this condition is
$$x=\lfloor\frac{p\ell}{|q_2|}\rfloor+1\text{ and }\neq -q_1j\ \ (j=1,\cdots,a).$$
{\bf (b)} \
The condition (\ref{e1})$=1$ is equivalent to $p\ell<\epsilon_1xq_2+k_2+\epsilon_1a\le p\ell+k_2$
$$\Leftrightarrow x+\frac{a}{q_2}\le \frac{p\ell}{|q_2|}<x+1+\frac{2a}{q_2}.$$
If $\epsilon_1=1$ then this condition is
$$x=\lfloor\frac{p\ell}{|q_2|}\rfloor \text{ and } \neq -q_1j\ \ (p)\ \ (j=0,\cdots a-1)$$
or
$$x=-q_1j-1\ (p)\ (j=0,1,\cdots,2a-1).$$
If $\epsilon_1=-1$ then this condition is
$$x=-q_1j\ (p)\ (j=1,\cdots,a)\ \  \left(x=\lfloor\frac{p\ell}{|q_2|}\rfloor+1\right)$$
or
$$x=\lfloor\frac{p\ell}{|q_2|}\rfloor\text{ and }\neq -q_1j-1\ (p)\ (j=1,\cdots,2a).$$
\hfill$\Box$

\subsubsection{The $\epsilon_2=-1$ case}
\begin{lem}
Let $p,k_1$, and $a$ be the same parameters as above.
Then the difference
\begin{equation}
A(-\epsilon_1a-1,\epsilon_1x,1)-A(-\epsilon_1a-1,\epsilon_1x-1,1)
\label{e2}
\end{equation}
is divided into three cases as follows.
\end{lem}
\begin{enumerate}
\item[(a)] The condition (\ref{e2})$=-1$ is equivalent to the following.\\
($\epsilon_1=1$)
There exists an integer $\ell$ such that
$$x=\lfloor \frac{p\ell}{|q_2|}\rfloor+1\ \text{ and } \neq q_1j\ (p)\ \ (j=1,\cdots,a).$$
($\epsilon_1=-1$)
There exists an integer $\ell$ such that
$$x=\begin{cases}
\lfloor \frac{p\ell}{|q_2|}\rfloor+1\text{, or}\\
q_1j\ (p)\ \ (j=0,\cdots,a-1)\ \ \left(\lfloor\frac{p\ell}{|q_2|}\rfloor\right). 
\end{cases}$$
\item[(b)] The condition (\ref{e2})$=1$ is equivalent to the following.\\
($\epsilon_1=1$)
There exists an integer $\ell$ such that
$$x=\begin{cases}
q_1j\ (p) \ \ (j=1,\cdots,a)\ \ (\lfloor \frac{p\ell}{|q_2|}\rfloor+1)\text{, or}\\
\lfloor \frac{p\ell}{|q_2|}\rfloor\text{ and }\neq q_1j-1,\ \ q_1j-k_1\ (p)\ \ (j=1,2,\cdots,a).
\end{cases}$$
($\epsilon_1=-1$)
There exists an integer $\ell$ such that
$$x=\begin{cases}
\lfloor\frac{p\ell}{|q_2|}\rfloor\text{ and }\neq q_1j\ (p)\ \ (j=0,\cdots,a-1),\\
q_1j-1\ (p)\ \ (j=0,\cdots,a-1)\ \ \left(\lfloor \frac{p\ell}{|q_2|}\rfloor-1\right)\text{, or}\\
q_1j+k_1\ (p)\ \ (j=0,\cdots,a-1)\ \ \left(\lfloor \frac{p\ell}{|q_2|}\rfloor-1\right).
\end{cases}$$
\item[(c)] (\ref{e2})$=0$ otherwise.
\end{enumerate}
{\bf Proof.} \
{\bf (a)} \
The condition (\ref{e1})$=-1$ is equivalent to $p\ell<\epsilon_1xq_2-\epsilon_1a\le p\ell+k_2$\\
$$\Leftrightarrow x-1\le   \frac{p\ell}{|q_2|}< x-\frac{a}{q_2}.$$
If $\epsilon_1=1$ then this condition is
$$x=\lfloor\frac{p\ell}{|q_2|}\rfloor+1\ \text{ and }\neq q_1j\ \ (j=1,\cdots,a)$$ 
If $\epsilon_1=-1$, then this condition is
$$x=\lfloor\frac{p\ell}{|q_2|}\rfloor+1$$
or
$$x=q_1j\ \ (j=0,\cdots,a-1)\ \ \left(x=\lfloor\frac{p\ell}{|q_2|}\rfloor\right).$$
{\bf (b)} \
This condition (\ref{e1})$=1$ is equivalent to $p\ell<\epsilon_1xq_2+k_2-\epsilon_1a\le p\ell+k_2$\\
$$\Leftrightarrow x-\frac{a}{q_2}\le \frac{p\ell}{|q_2|}<x+1-\frac{2a}{q_2}.$$
If $\epsilon_1=1$, then this condition is
$$x=q_1j (p)\ \ (j=1,\cdots, a)\ \ (x=\lfloor\frac{p\ell}{q_2}\rfloor+1)$$
or
$$x=\lfloor\frac{p\ell}{|q_2|}\rfloor\text{ and }\neq q_1j-1\ \  (j=1,\cdots,2a).$$
If $\epsilon=-1$, then this condition is 
$$x=\lfloor\frac{p\ell}{|q_2|}\rfloor\text{ and }\neq q_1j(p)\ \ (j=0,\cdots,a-1)$$
or
$$x=q_1j(p)\ \ (j=0,\cdots,2a-1)\ \ \left(x=\lfloor\frac{p\ell}{|q_2|}\rfloor-1\right).$$
\hfill$\Box$

As a result we have the following.
\begin{lem}
\label{e1e2}
Let $p,k_1,a$ be the same parameters as above.
Then the difference
\begin{equation}
A(\epsilon_1\epsilon_2a-1,\epsilon_1x,1)-A(\epsilon_1\epsilon_2a-1,\epsilon_1x-1,1)
\label{e3}
\end{equation}
is divided into three cases as follows.
\end{lem}
\begin{enumerate}
\item[(a)] The condition (\ref{e3})$=-1$ is equivalent to the following.\\
($\epsilon_1\epsilon_2=1$)
There exists an integer $\ell$ such that
$$x=\begin{cases}
\lfloor \frac{p\ell}{|q_2|}\rfloor+1\text{, or}\\
-\epsilon_2q_1j\ (p)\ \ (j=0,\cdots,a-1)\ \ \left(\lfloor\frac{p\ell}{|q_2|}\rfloor\right). 
\end{cases}$$
($\epsilon_1\epsilon_2=-1$)
There exists an integer $\ell$ such that
$$x=\lfloor \frac{p\ell}{|q_2|}\rfloor+1\ \text{ and } \neq -\epsilon_2q_1j\ (p)\ \ (j=1,\cdots,a).$$
\item[(b)] The condition (\ref{e3})$=1$ is equivalent to the following.\\
($\epsilon_1\epsilon_2=1$)
There exists an integer $\ell$ such that
$$x=\begin{cases}
\lfloor\frac{p\ell}{|q_2|}\rfloor\text{ and }\neq -\epsilon_2q_1j\ (p)\ \ (j=0,\cdots,a-1),\\
-\epsilon_2q_1j-1\ (p)\ \ (j=0,\cdots,a-1)\ \ \left(\lfloor \frac{p\ell}{|q_2|}\rfloor-1\right)\text{, or}\\
-\epsilon_2q_1j+k_1\ (p)\ \ (j=0,\cdots,a-1)\ \ \left(\lfloor \frac{p\ell}{|q_2|}\rfloor-1\right).
\end{cases}$$
($\epsilon_1\epsilon_2=-1$)
There exists an integer $\ell$ such that
$$x=\begin{cases}
-\epsilon_2q_1j\ (p) \ \ (j=1,\cdots,a)\ \ (\lfloor \frac{p\ell}{|q_2|}\rfloor+1)\text{, or}\\
\lfloor \frac{p\ell}{|q_2|}\rfloor\text{ and }\neq -\epsilon_2q_1j-1,\ \ -\epsilon_2q_1j-k_1\ (p)\ \ (j=1,2,\cdots,a).
\end{cases}$$
\item[(c)] (\ref{e3})$=0$ otherwise.
\end{enumerate}

\section{The global view of the function $A(i+jk_1)$}
\subsection{The local behavior of $A(x)$} 
From Lemma~\ref{e1e2}, $x$ in which $A(x)\neq A(x-\epsilon_1)$ is
$x=\lfloor\frac{p\ell}{|q_2|}\rfloor-1,\lfloor\frac{p\ell}{|q_2|}\rfloor,\lfloor\frac{p\ell}{|q_2|}\rfloor+1$.
We consider the quotient
\begin{equation}
\frac{p}{|q_2|}=\tau+\frac{\gamma'+\epsilon_1\epsilon_2}{n}+\frac{a}{n|q_2|}.
\label{p/q2}
\end{equation}
Here we put $\alpha=\gamma'+\epsilon_1\epsilon_2+\frac{a}{|q_2|}$.
In particular we have
$$\lfloor\frac{np}{|q_2|}\rfloor=n\tau+\gamma'+\epsilon_1\epsilon_2+\lfloor\frac{a}{|q_2|}\rfloor=k_1+\epsilon_1\epsilon_2.$$

\subsubsection{The $\epsilon_1\epsilon_2=1$ case.} 
Since by Theorem~\ref{alter} the values of $\tilde{a}_i$ are $0$, $\pm1$, or $2$,
the possibilities of the local behavior of $A(x)$ around $x=\lfloor\frac{p\ell}{|q_2|}\rfloor$ are Figure~\ref{noteq}, \ref{eqqq}, and \ref{1k1} in the case of $\epsilon_1=1$, and 
Figure~\ref{noteqp}, \ref{eqqqp}, and \ref{1k1p} in the case of $\epsilon_1=-1$.

\begin{figure}[htbp]
\begin{center}
\input{ee=1not.tex}
\caption{$A(x)$: $(\epsilon_1,\epsilon_2)=(1,1)$, $\lfloor\frac{p\ell}{|q_2|}\rfloor\neq -\epsilon_2q_1j$}
\label{noteq}
\end{center}
\begin{center}
\input{ee=1eq.tex}
\end{center}
\caption{$A(x)$: $(\epsilon_1,\epsilon_2)=(1,1), j=0,\cdots,a-1$)}
\label{eqqq}
\begin{center}
\input{ee=1k1.tex}
\caption{$A(x)$: $(\epsilon_1,\epsilon_2)=(1,1), j=0,\cdots,a-1$)}
\label{1k1}
\end{center}
\end{figure}

Now suppose that $\epsilon_1=1$.
If $A(x_0)=2$ holds for some $x_0$, then we cannot find any function $A(x)$ by connecting Figure~\ref{noteq}, \ref{eqqq}, or \ref{1k1}.
For example if the type (+,iii) appears, then the first possible type which is non-constant with $x>x_0$ is (+,v) or (+,iii).
Since the integer $x$ with $A(x)=2$ is unique in ${\mathbb Z}/p{\mathbb Z}$, $\Delta_K(t)=x^{-1}+1+x$ holds.
This polynomial cannot be any Alexander polynomial of a knot in a homology sphere.
Hence we may assume that $|A(x)|\le 1$ holds and ($+$,i), ($+$,ii), ($+$,iv), or ($+$,vi) are 
applicable non-constant local behaviors of $A(x)$.

We put $\alpha(j)=-\epsilon_2q_1j$ and $\beta(j)=-\epsilon_2q_1j+k_1$.
In the case of $\epsilon_1=1$, we can get Figure~\ref{valueab11}, in which there are the values of $A(x)$ from $\alpha(j)-1$ to $\beta(j)+1$,
fitting together of local behaviors for the values to be connected.
\clearpage
\begin{figure}[htbp]
$$\begin{array}{|c|c|c|c|c|c|c|c|c|c|c|c|c|}\hline
x&  \alpha(j)-1 & \alpha(j) & \alpha(j)+1 & \cdots & \lfloor\frac{p\ell}{|q_2|}\rfloor-1 & \lfloor\frac{p\ell}{|q_2|}\rfloor & \lfloor\frac{p\ell}{|q_2|}\rfloor+1 & \cdots\\\hline
A(x)& 1& 0 & -1  & \cdots & -1 & 0 & -1 & \cdots \\\hline
\end{array}$$
$$\begin{array}{|c|c|c|c|c|c|c|c|c|c|c|c|c|}\hline
\cdots & \lfloor\frac{p(\ell+1)}{|q_2|}\rfloor-1 & \lfloor\frac{p(\ell+1)}{|q_2|}\rfloor& \lfloor\frac{p(\ell+1)}{|q_2|}\rfloor+1 & \cdots & \beta(j)-1 & \beta(j) & \beta(j)+1\\\hline
\cdots & -1 & 0 & -1 & \cdots & -1 & 0 & 1 \\\hline
\end{array}$$
\caption{The values $A(x)$ from $\alpha(j)-1$ to $\beta(j)+1$ in the case of $(\epsilon_1,\epsilon_2)=(1,1)$.}
\label{valueab11}
\end{figure}
In the region $\alpha(j)+1\le x\le \beta(j)-1$ there are $0$ and $-1$.
We call the integer points in the region {\it negative region}.

On the other hand Figure~\ref{valueba11} is the values of $A(x)$ between $\beta(j)-1$ and $\alpha(k)+1$.
\begin{figure}[htbp]
$$\begin{array}{|c|c|c|c|c|c|c|c|c|c|c|c|c|}\hline
x& \beta(j)-1 & \beta(j) & \beta(j)+1 &\beta(j)+2 &  \cdots & \lfloor\frac{p\ell'}{|q_2|}\rfloor-1 & \lfloor\frac{p\ell'}{|q_2|}\rfloor & \lfloor\frac{p\ell'}{|q_2|}\rfloor+1 & \cdots\\\hline
A(x)& -1& 0 & 1  & 0 & \cdots & 0 & 1 & 0 & \cdots \\\hline
\end{array}$$
$$\begin{array}{|c|c|c|c|c|c|c|c|c|c|c|c|c|}\hline
\cdots & \lfloor\frac{p(\ell'+1)}{|q_2|}\rfloor-1 & \lfloor\frac{p(\ell'+1)}{|q_2|}\rfloor& \lfloor\frac{p(\ell'+1)}{|q_2|}\rfloor+1 & \cdots & \alpha(k)-2 & \alpha(k)-1 & \alpha(k) & \alpha(k)+1\\\hline
\cdots & 0 & 1 & 0 & \cdots & 0 & 1 & 0 & -1 \\\hline
\end{array}$$
\caption{The values $A(x)$ from $\beta(j)-1$ to $\alpha(j)+1$ in the case of $(\epsilon_1,\epsilon_2)=(1,1)$.}
\label{valueba11}
\end{figure}

Then the values $A(x)$ between $\beta(j)+1$ and $\alpha(k)-1$ are $0$ or $1$.
We call the integer points in this region {\it positive region}.

Thus if $(p,k_1)$ is a data realized by a lens space surgery over $S^3$ or $\Sigma(2,3,5)$, 
then we can decide the values $A(x)$ completely.
Hence we figure out that each positive region and negative region appear alternatively and 
do not overlap each other.

In the same way all the local behaviors with $\epsilon_1=-1$ are from Figure~\ref{noteqp}, \ref{eqqqp}, and \ref{1k1p}.
The case satisfying $A(x)=2$ for some integer $x$ is classified in Section~\ref{an}.

Thus 
we can decide the values $A(x)$ for any $x$.
\begin{figure}[htbp]
\begin{center}
\input{ee=1notp.tex}
\caption{$A(x)$: $(\epsilon_1,\epsilon_2)=(-1,-1),\lfloor\frac{p\ell}{|q_2|}\rfloor\neq -\epsilon_2q_1j$}
\label{noteqp}
\end{center}
\begin{center}
\input{ee=1eqp.tex}
\end{center}
\caption{$A(x)$: $(\epsilon_1,\epsilon_2)=(-1,-1)$, $j=0,\cdots,a-1$)}
\label{eqqqp}
\begin{center}
\input{ee=1k1p.tex}
\caption{$A(x)$: $(\epsilon_1,\epsilon_2)=(-1,-1)$,($j=0,\cdots,a-1$)}
\label{1k1p}
\end{center}
\end{figure}
Thus we call the integer points between $\alpha(j)+2$ and $\beta(j)$ and
integer points between $\beta(j)+3$ and $\alpha(k)$ {\it positive region} and {\it negative region}
respectively.

\begin{figure}[htbp]
\begin{center}
\input{stair.tex}   
\caption{The stair-like sequences in positive region.\ ($\epsilon_1=1$)}
\label{stair}
\end{center}
\end{figure}
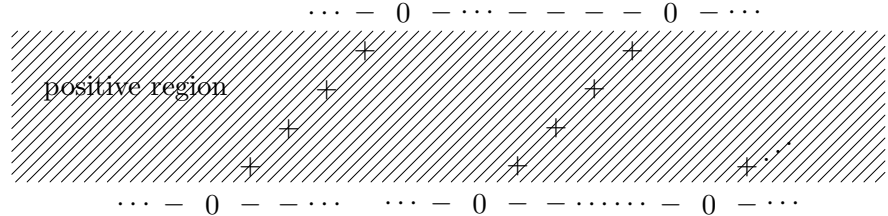

Next we put the values $A(i+jk_1)$ on $i$-$j$ plain.
We call this distribution {\it global view} of the Alexander polynomial.
Since we have
$$\lfloor\frac{mnp}{|q_2|}\rfloor=mn\tau +m(\gamma'+\epsilon_1\epsilon_2)+\lfloor\frac{m a}{|q_2|}\rfloor=m(k_1+\epsilon_1\epsilon_2)+\lfloor\frac{m a}{|q_2|}\rfloor.$$
Then the following Lemma is true.
\begin{lem}
\begin{eqnarray}
\lfloor\frac{mnp}{|q_2|}\rfloor-\lfloor\frac{(m-1)np}{|q_2|}\rfloor-(k_1+\epsilon_1\epsilon_2)&=&\lfloor\frac{ma}{|q_2|}\rfloor-\lfloor\frac{(m-1)a}{|q_2|}\rfloor\label{mnpq2}\\
&=&\begin{cases}1&0\le [ma]_{|q_2|}<a\\0&\text{otherwise}\end{cases}\nonumber
\end{eqnarray}
\end{lem}
Since we have $a'=a^{-1}=-\epsilon_1(k_1-\epsilon_1\epsilon_2)\ \ (|q_2|)$,
when $m=a'i\ \ (i=0,\cdots,a-1)$, we have (\ref{mnpq2}) $=1$.
Hence we have
$$\lfloor\frac{a'inp}{|q_2|}\rfloor=-\epsilon_1iq_1\ \ (p).$$
This implies that if $\epsilon_1=1$ (or $\epsilon_1=-1$), then all $x$ that gives $A(x)=1$ (or $A(x)=-1$) lines like stairs in positive region (or negative region) on $i$-$j$ plain.
See Figure~\ref{stair}.
Here each of places which put $+$ or $-$ represents $A(x)=1$ or $-1$ respectively and the blanks are all $0$.
The starting point of a stair-like sequence is $(x,y)$, where $(i,j)=(x-1,y-1)$ is a point in negative (or positive) region that $A(i+jk_1)$ is $0$,
and the end point is $(z,w)$, where $(i,j)=(z+1,w+1)$ is a point in negative (or positive) region that $A(i+jk_1)$ is $0$.

Combining the results above, we plot $A(i+jk_1)$ on $i$-$j$-plane to get Figure~\ref{g11},~\ref{g-1-1} in the case of $\epsilon_1=1$, $-1$ respectively.
\begin{figure}[htbp]
\begin{center}
\input{tableee=1.tex}
\caption{The global view of $A(i+jk_1)$ in the case of $(\epsilon_1,\epsilon_2)=(1,1)$}
\label{g11}
\end{center}
\end{figure}
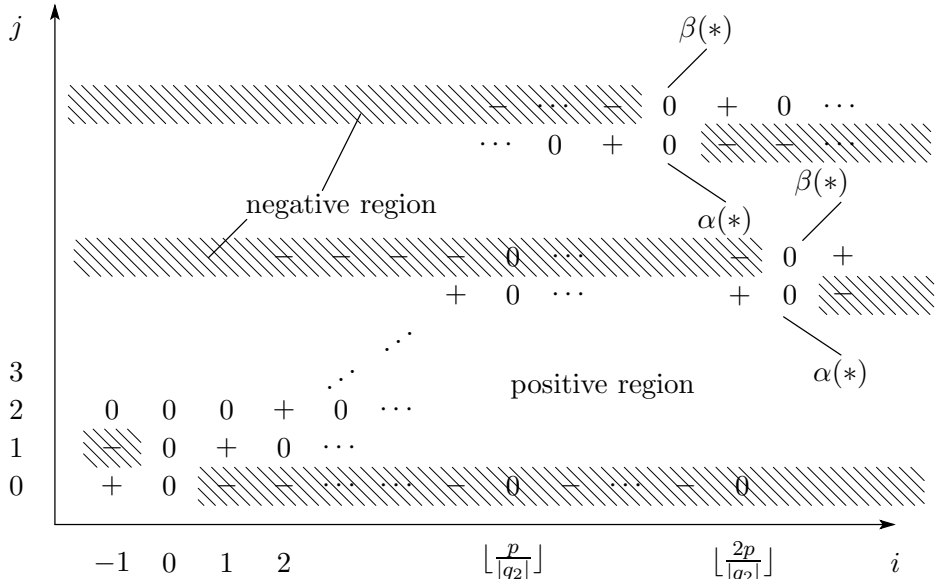
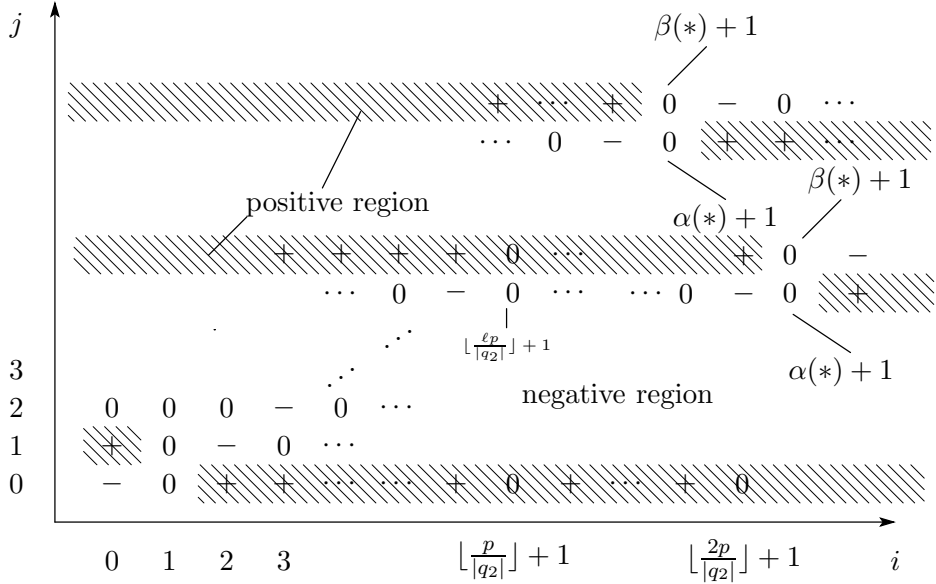
\begin{figure}[htbp]
\begin{center}
\input{tableee=1n.tex}
\caption{The global view of $A(i+jk_1)$ in the case of $(\epsilon_1,\epsilon_2)=(-1,-1)$}
\label{g-1-1}
\end{center}
\end{figure}

\subsubsection{The $\epsilon_1\epsilon_2=-1$ case.} 
In this section we give the global view of $A(i+jk_1)$ by investigating 
the local behavior of $A(x)$ in the case of $\epsilon_1\epsilon_2=-1$.

In the case of $\epsilon_1=1$, using Lemma~\ref{e1e2}, we can easily show
that the possible non-constant local behaviors of the function $A(x)$ 
are Figure~\ref{noteq2p}, \ref{eqqq2p}, \ref{1k12p}.
In the case of $\epsilon_1=-1$ the non-constant local behaviors are Figure~\ref{noteq2m}, \ref{eqqq2m}, \ref{1k12m}.
Here $\gamma(j)$ stands for $-\epsilon_2q_1j-k_1$.
\begin{figure}[htbp]
\begin{center}
\input{ee=-1notp.tex}
\caption{$A(x)$: $(\epsilon_1,\epsilon_2)=(1,-1),\ \lfloor\frac{p\ell}{|q_2|}\rfloor\neq -\epsilon_2q_1j\ (0\le j<a)$}
\label{noteq2p}
\end{center}
\begin{center}
\input{ee=-1eqp.tex}
\caption{$A(x)$: $(\epsilon_1,\epsilon_2)=(1,-1),\ (j=1,\cdots,a)$}
\label{eqqq2p}
\end{center}
\begin{center}
\input{ee=-1k1p.tex}
\caption{$A(x)$: $(\epsilon_1,\epsilon_2)=(1,-1),\ (j=1,\cdots,a)$}
\label{1k12p}
\end{center}
\end{figure}
The case where $A(\alpha(j))=2$ for some integer $j$ will be dealt with in Section~\ref{2ggg}.
We assume that $|A(x)|=1$ holds for any $x$ this case will be dealt in Section~\ref{gvito}.

\subsubsection{The great common divisor $\gcd(a,n)$.}
\label{an}
We put $d=(\eta,a)=(n,a)$, $\eta=\eta_1d$, $a=a_1d$.
We put $\tilde{\eta}:=[\eta_1^{-1}]_{a_1}$.
Let $j$ be $j_1a_1+[-\epsilon_2\tilde{\eta}j_2]_{a_1}$ with $0\le j_1<d$, $0\le j_2<a_1$.
Then we have $-\epsilon_2jq_1=-\epsilon_2\frac{j\eta p-\epsilon_1jk_1-\epsilon_2j}{a}\equiv\frac{j_2}{a_1}p+\frac{j(\epsilon_1\epsilon_2k_1+1)}{a}\ (p)$.
We define intervals to be
$$I_{j_1,j_2}:=[-\epsilon_2(\frac{j_2}{a_1}p+\frac{j(\epsilon_1\epsilon_2k_1+1)}{a}),-\epsilon_2(\frac{j_2}{a_1}p+\frac{j(\epsilon_1\epsilon_2k_1+1)}{a}+k_1)]$$
if $\epsilon_1\epsilon_2=1$, and 
$$I_{j_1,j_2}:=[-\epsilon_2(\frac{j_2}{a_1}p+\frac{j(\epsilon_1\epsilon_2k_1+1)}{a}-k_1),-\epsilon_2(\frac{j_2}{a_1}p+\frac{j(\epsilon_1\epsilon_2k_1+1)}{a})]$$
if $\epsilon_1\epsilon_2=-1$.

For a fixed $j_2$ the intervals $I_{j_1,j_2}\ \ (j_1=0,\cdots,d-1)$ are overlapped each other.
Furthermore $I_{j_1,j_2}\cap I_{j_1,j_3}=\emptyset$ holds unless $j_2=j_3$.
On the other hand if $d\ge 2$, then we have $I_{j_1,j_2}\cap I_{j_3,j_4}=\emptyset$ unless $j_1=j_3$,
since $\frac{p}{a_1}\ge \frac{2p}{a}\ge \frac{2(ak_1^2+\epsilon_1k_1+\epsilon_2)}{an}\ge 2k_1(\tau+\frac{\gamma'}{n})-\frac{2(k_1+1)}{an}>4k_1-(k_1+1)=3k_1-1$.



Thus $A(x)$ increases by $-\epsilon_2$ around $x=-\epsilon_2jq_1$ as shown in Figure~\ref{eqqq}, \ref{eqqqp}, \ref{eqqq2p}, and \ref{eqqq2m}.
In the same way $A(x)$ increases by $\epsilon_2$ around $x=-\epsilon_2jq_1+\epsilon_1\epsilon_2k_1$ as it is shown in Figure~\ref{1k1}, \ref{1k1p}, \ref{1k12p}, and \ref{1k12m}.
Therefore $d\le 2$ holds by Theorem~\ref{alter}.
Thus we showed the following.
\begin{lem}
Let $(p,k_1)$ be a data realized by lens space surgery $S^3$ or $\Sigma(2,3,5)$.
Then the parameters $a,n$ have $\gcd(a,n)=1$ or $2$.
\end{lem}

\subsubsection{The $2g(p,k_1)=p$ case}
\label{2ggg}
For the case where $2g=p$ is satisfied
the case of $(\epsilon_1,\epsilon_2)=(-1,-1)$ and $2g=p$ is remained.
If $A(x_0)=2$ for some integer $x$, then using Figure~\ref{eqqqp} and \ref{1k1p} we can get the following for some integer $j$.
$$\begin{array}{|c|c|c|c|c|c|c|}\hline
-\epsilon_2q_1j-1&-\epsilon_2q_1j&-\epsilon_2q_1j+1& -\epsilon_2q_1j+2& -\epsilon_2q_1j+3& -\epsilon_2q_1j+4& -\epsilon_2q_1j+5\\\hline
1&0&1&2&1&0&1\\\hline
\end{array}$$
Namely $-\epsilon_2q_1j+2=-\epsilon_2q_1j'+k_1$ holds.

\begin{figure}[htbp]
\begin{center}
\input{an2g=p.tex}
\caption{The local view of $x$ with $A(x)=2$.}
\label{lv}
\end{center}
\end{figure}
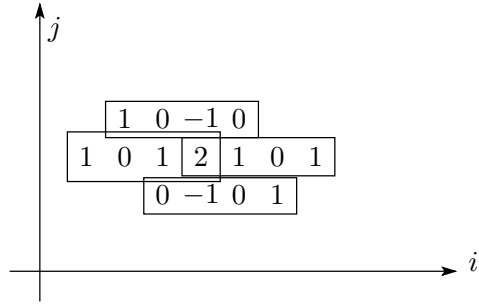

\begin{figure}[htbp]
\begin{center}
\input{int.tex}
\caption{The overlapped interval.}
\label{int}
\end{center}
\end{figure}
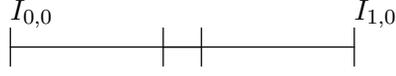

In this case we have $d=2$ and the local view around $(i,j)$ satisfying $x=i+jk_1$ and $A(x)=2$ is Figure~\ref{lv}.
Since this pattern does not appear in other places up to the moving by $p{\mathbb Z}$, $a_1=1$ holds.
Thus $I_{0,0}\cap I_{1,0}$ is overlapped as in Figure~\ref{int}.
Thus we have $-\epsilon_2q_1j+2=q_1+2=k_1$.
From the associated relation $2k_1^2-k_1-1$, we get $k_1=5$ and $q_1=3$.
Thus we have $p|22$.
If $p=11$, we cannot find lens space surgery from the list in \cite{[tan2]}.
Therefore we have $p=22$.
This case is realized by lens space surgery of type $\text{E}_2$ and $J=-1$ in Table~\ref{po}

In the case $\epsilon_1\epsilon_2=-1$ we can show in the same way as above
that the $2g=p$ case does not exist.
\subsubsection{The global view in the case of $\epsilon_1\epsilon_2=-1$}
\label{gvito}
We get back to the description of the global view.
Here we suppose that $d=2$ holds.
In particular $a$ is even number.
Then we get the values as in Figure~\ref{d=2}.
\begin{figure}[htbp]
$$\begin{array}{|c|c|c|c|c|c|c|c|c|c|}\hline
x & \gamma(a_1+j_2) & \gamma(a_1+j_2)+1 & \cdots & \gamma(j_2)-1 & \gamma(j_2) & \gamma(j_2)+1 & \cdots & \alpha(a_1+j_2)-2\\\hline
A(x) & 1 & 0 & \cdots & 0 & 0 &  -1& \cdots & -1\\\hline
\end{array}
$$
$$\begin{array}{|c|c|c|c|c|c|c|c|c|c|}\hline
 \alpha(a_1+j_2)-1 & \alpha(a_1+j_2) & \cdots & \alpha(j_2)-2 & \alpha(j_2)-1 & \alpha(j_2) & \cdots & \gamma(a_1+j_2+1)\\\hline
  -1 & 0 & \cdots & 0 & 0 &  1& \cdots & 1\\\hline
\end{array}
$$
\caption{The values $A(x)$ from $\gamma(a_1+j_2)$ to $\gamma(a_1+j_2+1)$ in the case of $d=2$ and $\epsilon_1\epsilon_2=-1$}
\label{d=2}
\end{figure}
Note that if $\lfloor\frac{p\ell}{|q_2|}\rfloor$ appears between $\alpha(j_2)$ and $\gamma(a_1+j_2+1)$,
it is at most one time.
Then we can reduce the following inequality.
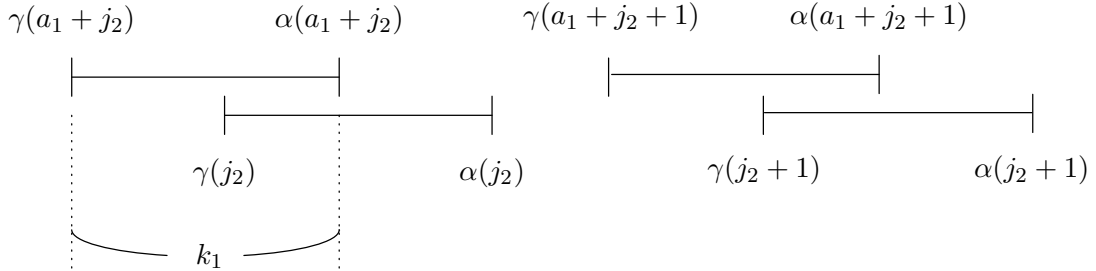
\begin{figure}[htbp]
\begin{center}
\input{interval.tex}
\caption{The $d=2$ case.}
\end{center}
\end{figure}
$$ak_1+(\tau+2)\frac{a}{2}+\tau>p=\frac{ak_1^2+\epsilon_1k_1+\epsilon_2}{n}$$
Thus
$$ak_1(\tau-1+\frac{\gamma'}{n})<\frac{k_1+1}{n}+\frac{n\tau+\gamma'}{n}\frac{a}{2}+a+\tau<\frac{k_1}{n}(2+\frac{a}{2})+\frac{1}{n}+a$$
From $\frac{2}{an}\le \frac{1}{n}$ and $\frac{1}{k_1}<\frac1{2n}$,
$$<ak_1(\frac{2}{an}+\frac{1}{2n}+\frac{1}{k_1})+\frac{1}{n}<ak_1\frac{2}{n}+\frac{1}{n}$$
Therefore
$$ak_1(\tau-1+\frac{\gamma'-2}{n})<\frac{1}{n}$$
Since the left hand side is more than one, this inequality is inconsistent.

Therefore we must have $d=1$, namely $(a,n)=1$.
If $A(\gamma(j))=1$ holds for some $j$, then $A(x)\le 0$ holds for all $x$.
This is equivalent to $\Delta_{p,k_1}(t)=1$ and the data is realized by the unknot surgery.
Therefore the non-constant local possibilities are (-,i), (-,ii), (-,iv) and (-,iix).

In the case of $\epsilon_1=-1$, we can deduce that the non-constant local behaviors are (-,ix), (-,x), (-,xi), (-,xiv)
in the same argument as in the $\epsilon_1=1$ case.
\begin{figure}[htbp]
\begin{center}
\input{ee=-1notm.tex}
\caption{$A(x)$: $(\epsilon_1,\epsilon_2)=(-1,1),\ \lfloor\frac{p\ell}{|q_2|}\rfloor\neq -\epsilon_2q_1j\ (0\le j<a)$}
\label{noteq2m}
\end{center}
\begin{center}
\input{ee=-1eqm.tex}
\end{center}
\caption{$A(x)$: $(\epsilon_1,\epsilon_2)=(-1,1),\ (j=1,\cdots,a)$}
\label{eqqq2m}
\begin{center}
\input{ee=-1k1m.tex}
\caption{$A(x)$: $(\epsilon_1,\epsilon_2)=(-1,1),\ (j=1,\cdots,a)$}
\label{1k12m}
\end{center}
\end{figure}

Combining these results, we plot $A(i+jk_1)$ on $i$-$j$-plane to get Figure~\ref{g1-1},~\ref{g-11} in the case of $\epsilon_1=1$, $-1$ respectively.
\begin{figure}[htbp]
\begin{center}
\input{tableee=-1.tex}
\caption{The global view of $A(i+jk_1)$ in the case of $(\epsilon_1,\epsilon_2)=(1,-1)$}
\label{g1-1}
\end{center}
\end{figure}
\begin{figure}[htbp]
\begin{center}
\input{tableee=-1n.tex}
\caption{The global view of $A(i+jk_1)$ in the case of $(\epsilon_1,\epsilon_2)=(-1,1)$}
\label{g-11}
\end{center}
\end{figure}

\section{A certain integrality condition}


Here for any $1\le i\le a-1$ we define $n_i$ to be $[-\epsilon_2\eta^{-1}i]_a$, so that
we have 
\begin{equation}
-\epsilon_2n_iq_1=\frac{ip+n_i(\epsilon_1\epsilon_2k_1+1)}{a}\ \ (p).
\label{e2niq1}
\end{equation}
In the case of $\epsilon_1\epsilon_2=1$ we put $n_0=0$ then $-\epsilon_2n_0q_1=0$.
In the case of $\epsilon_1\epsilon_2=-1$ we put $n_a=a$ then $-\epsilon_2n_aq_1=\epsilon_1\epsilon_2k_1+1=-k_1+1$.

We define $n'$ to be the integer satisfying $an'=nn_{a-1}-\epsilon_2$.
Then $a(n-n')=nn_1+\epsilon_2$ holds.
We define $n'_i$ to be the integer satisfying $an_i'=nn_i+\epsilon_2i$ ($1\le i\le a-1$).
When $\epsilon_1\epsilon_2=1$, we have $n_0'=0$.
When $\epsilon_1\epsilon_2=-1$, we have $n_a'=n+\epsilon_2$.
Then $n_1'=n-n'$ and $n'_{a-1}=n'+\epsilon_2$ holds.
We shall show the integral condition.
\begin{lem}[integrality condition]
\label{intcond}
Let $\gamma'$, $n_1$ be as above.
Then (\ref{e2niq1}) is integer for any $i=1,2,\cdots,a-1$ if and only if
$$(\gamma')^2-\epsilon_1\epsilon_2n'\gamma'-n'=0\ \ (n)$$
holds.
\end{lem}
{\bf Proof.} \
From the integrality of (\ref{e2niq1}) we have
\begin{eqnarray}
(\ref{e2niq1})&=&in\tau^2+2i\gamma'\tau+\frac{i(\gamma')^2+n_i'(\epsilon_2k_1+1)}{n}\nonumber\\
&=&in\tau^2+(2i\gamma'+\epsilon_1\epsilon_2n'_i)\tau+\frac{i(\gamma')^2+\epsilon_1\epsilon_2n_i'\gamma'+n_i'}{n}.\nonumber
\end{eqnarray}
Thus from $n'_i=in_1'\ (n)$ for any $i$ we have 
$$i(\gamma')^2+\epsilon_1\epsilon_2n_i'\gamma'+n_i'=i((\gamma')^2+\epsilon_1\epsilon_2n_1'\gamma'+n_1')=0\ \ (n).$$
Therefore the integrality of (\ref{e2niq1}) is equivalent to
$$(\gamma')^2+\epsilon_1\epsilon_2n_1'\gamma'+n_1'=0\ \ (n).$$
\hfill$\Box$




\section{Classification}
In this section we classify all stable data realized lens space surgery on $S^3$ or $\Sigma(2,3,5)$ under the conditions
$a\ge 2$, $0<\gamma'< n$, $k_1\ge 4$, $n\ge 2$ and $g(p,k_1)<\frac{p}{2}$.
\begin{prop}
\label{nonhyper}
If the pair $(p,k_1)$ is realized by a surgery over $S^3$ or $\Sigma(2,3,5)$ and $n=1$ holds, then 
$(k_1,k_2)=1$ or $2$ and the surgery is realized by $(I),(II)$.
\end{prop}
{\bf Proof} \
From Theorem~12 in \cite{[tan1]} $(k_1,k_2)=1$ or $2$ holds.
If $(k_1,k_2)=1$, then $(p,k_1)$ is realized by (I) and if $(k_1,k_2)=2$,
then $(p,k_1)$ is realized by (II).
\hfill$\Box$
\subsection{The case of $\gamma'=n-1$ or $1$.}
Hence for any integer $m$ we have
\begin{equation}
\lfloor\frac{mp}{|q_2|}\rfloor-\lfloor\frac{(m-1)p}{|q_2|}\rfloor=\tau+\lfloor\frac{\alpha m}{n}\rfloor-\lfloor\frac{\alpha (m-1)}{n}\rfloor.
\label{defloor}
\end{equation}
Since $0<\frac{\alpha}{n}<2$, (\ref{defloor}) is $\tau$, $\tau+1$ or $\tau+2$.
If for some integer $m$ the value of (\ref{defloor}) is $\tau+2$, then $\alpha>n$, $\gamma'\ge n-1$, and $\epsilon_1\epsilon_2=1$ hold.

\begin{lem}
If (\ref{defloor}) is $\tau+2$, then the stable data $(p,k_1)$ is realized by lens space surgery type (I) or (II).
\end{lem}
{\bf Proof.} \
If (\ref{defloor}) is $\tau+2$, we have $\alpha=\gamma'+\epsilon_1\epsilon_2+\frac{a}{|q_2|}>n$.
From $\gamma'<n$, we have $\gamma'=n-1$ and $\epsilon_1\epsilon_2=1$.
From integral condition $1-n'+n'=1=0(n)$, $n=1$ holds.
This data $(p,k_1)$ is realized by Dehn surgeries of (I) or (II) by Proposition~\ref{nonhyper}.
\hfill$\Box$
\medskip

If $\gamma'=n-1$ and $\epsilon_1\epsilon_2=-1$, then from the integrality condition $1-n'-n'=1-2n'=0(n)$
we have $a=-2\epsilon_2=2\epsilon_1\ (n)$.
Thus we can put $a=un+2\epsilon_1$ for some integer $u$.
Therefore
$k_1=n\tau+n-1$ and 
\begin{eqnarray}
p&=&(un+2\epsilon_1)n\tau^2+(2a(n-1)+\epsilon_1)\tau+an-2na+\epsilon_1+u\nonumber\\
&=&u(n\tau+n-1)^2+\epsilon_1(2n\tau^2+(4n-3)\tau+2n-3)\nonumber\\
&=&uk_1^2+\epsilon_1(2k_1-1)(\tau+1)\nonumber
\end{eqnarray}
Since $(n,a)=1$ holds, $n$ is odd number.
Each of the data is realized by type $(\text{III}_+)$ by Berge's list.
Actually put $J=\epsilon_1(\tau+1)$.
\medskip

In the case of $\gamma'=1$ in the same way,
if we have $\epsilon_1\epsilon_2=1$, we can get 
$$(p,k_1)=(uk_1^2-\epsilon_1(2k_1+1)\tau,n\tau+1)$$
for some integer $u$, where $n$ is odd.
If $\epsilon_1\epsilon_2=-1$ then $n=1$.

The data of the former part is realized by type $(\text{III}_-)$. (Put $J=-\epsilon_1\tau$.
The data of the latter part is realized by type (I) or (II).

\subsection{Forbidden and admitted sequences}
\begin{defi}
Let $m_i (i\in {\mathbb Z})$ be an integer-valued sequence.
We call a subsequence $m_{i_0},m_{i_0+1},\cdots, m_{i_1}$ satisfying
$m_{i_0}=m_{i_1}=-1$ and $m_{i_0+j}=0\ (1\le j\le i_1-i_0-1)$ {\it forbidden subsequence}.
We call a subsequence $m_{i_0},m_{i_0+1},\cdots, m_{i_1}$ satisfying
$m_{i_0}=m_{i_1}=1$ and $m_{i_0+j}=0\ (1\le j\le i_1-i_0-1)$ {\it admitted subsequence}.
\end{defi}
If the sequence $\{\tilde{a}_i\}$ includes a forbidden subsequence,
this data is not realized by any lens space surgery over $S^3$ or $\Sigma(2,3,5)$.
If the sequence $\{\tilde{a}_i\}$ includes an admitted subsequence
and the data is realized by a lens space surgery over $S^3$ or $\Sigma(2,3,5)$,
then this subsequence $\tilde{a}_{i_0},\cdots,\tilde{a}_{i_1}$ satisfies $i_0=g(p,k_1)\ (p)$, $i_1=p-g(p,k_1)\ (p)$ and
$i_1-i_0=p-2g(p,k_1)$.

\subsection{Blocks}
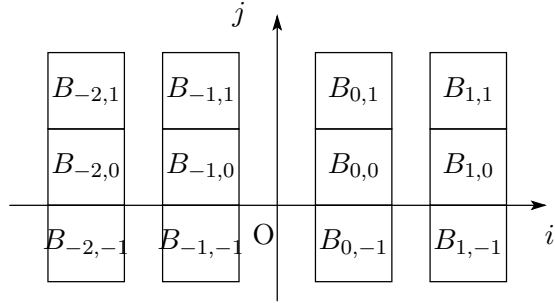
\begin{figure}[htbp]
\begin{center}
\input{Bij.tex}
\caption{$B_{i,j}$:$(i,j)$ component of blocks}
\label{BIJ}
\end{center}
\end{figure}

We consider polygons surrounded by a bold dotted line in the global view on $i$-$j$-plane
as in Figure~\ref{block}.
This polygon is composed of 6 segments with lattice point ${\mathbb Z}^2$ as the ends of the segments.
The top and bottom segments are the maximal one that all the lattice points $x$ on it have $A(x)=\pm1$.
Let $\{(x,y)|y=s_0,t_0\le x\le t_1\}$ and $\{(x,y)|y=s_1,b_0\le x\le b_1\}$ be the top and bottom line
respectively.
The two of the remaining 4 lines are
$$\{(x,y)|y=t_1,\ s_1-1\le y\le s_1 \},\ \{(x,y)|x=b_0,\ s_0\le y\le s_0+1 \},$$
or
$$\{(x,y)|x=t_0,s_1-1\le y\le s_1 \},\ \{(x,y)|x=b_1,\ s_0\le y\le s_0+1\}.$$
for $\epsilon_1\epsilon_2=1$ or $-1$ respectively.
These make two components.
The more two segments are two lines constructed 
by connecting each of right end points of the components and
each of left end points of that.
Furthermore $t_0,t_1,b_0,b_1,s_0$, and $s_1$ satisfy
$$t_1-b_0=s_1-s_0-2$$
or
$$b_1-t_0=s_1-s_0-2.$$
The former appears in global view of $\epsilon_1\epsilon_2=1$
and the latter in global view of $\epsilon_1\epsilon_2=-1$.
Note that here are some $\pm1$ on line connecting $(t_1,s_0-1)$ and $(b_0,s_1)$ or $(t_0,s_0-1)$ and $(b_1,s_1+1)$,
where each $\pm1$ on the line has opposite sign on top and bottom line.
\begin{figure}[htbp]
\begin{center}
\input{block.tex}
\caption{Blocks}
\label{block}
\end{center}
\end{figure}
We call this polygon {\it a block}.
Let denote by $B_{i,j}$ a block with the $(i,j)$ component
as in Figure~\ref{BIJ}.
However actually in the case of $\epsilon_1\epsilon_2=-1$ the bottom line of $B_{0,0}$ is not
lies on the $i$-axis.
One must note that the accurate location of blocks in the case later.
The left one in Figure~\ref{block} is called (+)-block and the right one (-)-block.
Here there are $1$ or $-1$ in the top or bottom row of the block.
The points except the top, bottom and left (or right) sides of (+)-block (or (-)-block) have all $0$.

We define {\it top width} $\partial_+ B_{i,j}$ for a block to be $t_1-t_0+1$
and {\it bottom width} $\partial_-B_{i,j}$ to be $b_1-b_0-1$, and 
define {\it height} $H(B_{i,j})$ to be $s_1-s_0+1$.
The right hand side of (\ref{defloor}) calculates the top (or bottom) width of the block.

\begin{figure}[htbp]
\begin{center}
\input{blockboun.tex}
\caption{The vertices of a block $B_{i,j}.\ \ (\epsilon_1,\epsilon_2)=(-1,1)$}
\label{bb}
\end{center}
\end{figure}
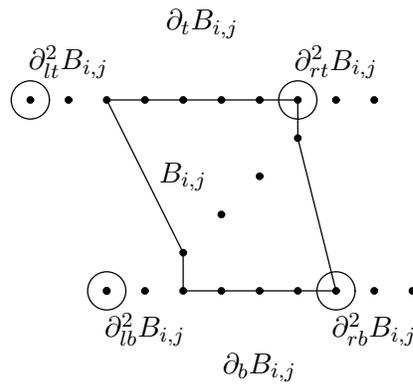

There exists a point $(i_1,j_1)$ satisfying $i_1+j_1k=\lfloor\frac{p\ell}{|q_2|}\rfloor$ around the vertex of the right (or left) top (or bottom) in a block $B_{i,j}$.
We denote the points by $\partial^2_{lt}B_{i,j}$, $\partial^2_{rt}B_{i,j}$, $\partial^2_{lb}B_{i,j}$, and $\partial^2_{rb}B_{i,j}$
corresponding to the left top, right top, left bottom, and right bottom vertices respectively as in Figure~\ref{bb},
in which the center points of the circles are the points $\partial^2_{\ast}B_{i,j}$ for $\ast=lt,rt,lb,$ and $rb$.

Figure~\ref{bll} presents examples of the blocks with the height $\tau-1$, $\tau$, or $\tau+1$.
top width and bottom width is $\tau$ or $\tau+1$.
\begin{figure}[htbp]
\begin{center}
\input{blocks.tex}
\caption{blocks}
\label{bll}
\end{center}
\end{figure}

As a result the next lemma easily follows.
\begin{lem}
The top (or bottom) width of any block is $\tau$ or $\tau+1$.
Furthermore if top (or bottom) width of $B_{i_0,j}$ and $B_{i_0+m,j}$ are $\tau+1$ and 
the top (or bottom) width of $B_{i_0+i_1,j}$ for any $1\le i_1\le m-1$ is $\tau$, then
$m=\lfloor\frac{n}{\alpha}\rfloor$ or $1+\lfloor\frac{n}{\alpha}\rfloor$.
These points are called the 2nd boundary of the block.
In particular, the top width of $B_{-1,0}$ is $\tau+1$
\end{lem}
{\bf Proof.}
We show just the last assertion.
Since $[\alpha0]_n=0$ holds, Lemma~\ref{ifif} implies $B_{-1,0}$ has $\tau+1$ top width.
\hfill$\Box$


\subsection{The $2\le \gamma'\le n-2$ case (i).}
From this section we assume that $2\le \gamma'\le n-2$ holds.
In particular $n\ge 4$, $\alpha<n$, and
$$\lfloor\frac{p}{|q_2|}\rfloor=\tau+\lfloor\frac{\gamma'+\epsilon_1\epsilon_2}{n}\rfloor=\tau$$
hold.
\begin{lem}
\label{ifif}
(\ref{defloor}) is $\tau+1$, if and only if we have $0\le[\alpha m]_n< \alpha$.
\end{lem}
{\bf Proof.} \
The equality
$\lfloor\frac{\alpha m}{n}\rfloor-\lfloor\frac{\alpha (m-1)}{n}\rfloor=1$
implies that there exists an integer $n_0$ satisfying
$\frac{\alpha (m-1)}{n}<n_0\le \frac{\alpha m}{n}$.
Therefore equivalently $0\le [\alpha m]_n<\alpha$ holds.\hfill$\Box$

Since we have
\begin{eqnarray}
\frac{[-\epsilon_2n_1q_1]_p}{k_1}&=&\frac{p+(\epsilon_1\epsilon_2k+1)n_1}{ak_1}\nonumber\\
&=&\frac{ak_1^2+(\epsilon_1\epsilon_2k_1+1)(nn_1+\epsilon_2)}{ank_1}\nonumber\\
&=&\tau+\epsilon_1\epsilon_2+\frac{\gamma'-\epsilon_1\epsilon_2n'}{n}+\frac{n-n'}{nk_1},\label{tauone}
\end{eqnarray}
$\frac{n-n'}{nk_1}<\frac{1}{k_1}\le \frac{1}{4n}$, and
$\gamma'-\epsilon_1\epsilon_2n'=(\gamma')^{-1}n'\neq 0\ (n)$,
we have
$$\lfloor\frac{[-\epsilon_2n_1q_1]_p}{k_1}\rfloor=\tau+\epsilon_1\epsilon_2+\lfloor\frac{\gamma'-\epsilon_1\epsilon_2n'}{n}\rfloor\le \tau+1.$$
\hfill$\Box$

\begin{lem}
\label{taufloorfunc}
For any $0<i<a$ we have 
$$\frac{[-\epsilon_2n_iq_1]_p}{k_1}=i\tau+\frac{i\gamma'+\epsilon_1\epsilon_2n'_i}{n}+\frac{n_i'}{nk_1}$$
and
$$\lfloor\frac{[-\epsilon_2n_iq_1]_p}{k_1}\rfloor=i\tau+\epsilon_1\epsilon_2+\lfloor\frac{i\gamma'-\epsilon_1\epsilon_2[n'i]_n}{n}\rfloor.$$
\end{lem}
{\bf Proof}
\begin{eqnarray}
\frac{[-\epsilon_2n_iq_1]_p}{k_1}&=&\frac{ip+(\epsilon_1\epsilon_2k_1+1)n_i}{ak_1}\nonumber\\
&=&\frac{aik_1^2+(\epsilon_1\epsilon_2k_1+1)(nn_i+i\epsilon_2)}{ank_1}\nonumber\\
&=&i\tau+\frac{i\gamma'+\epsilon_1\epsilon_2n'_i}{n}=i\tau+\frac{i\gamma'+\epsilon_1\epsilon_2n'_i}{n}+\frac{n_i'}{nk_1}\nonumber
\end{eqnarray}
Here $\frac{n_i'}{nk_1}\le \frac{n}{nk_1}\le\frac{1}{2n}$ and $i\gamma'+\epsilon_1\epsilon_2n'_i=ia^{-1}(\gamma')^{-1}(a(\gamma')^2+\epsilon_1\gamma')=ia^{-1}(\gamma')^{-1}(-\epsilon_2)\neq 0(n)$ hold.
Hence we have $\lfloor\frac{i\gamma'+\epsilon_1\epsilon_2n'_i}{n}+\frac{n_i'}{nk_1}\rfloor=\lfloor\frac{i\gamma'+\epsilon_1\epsilon_2n'_i}{n}\rfloor$.
Since $n'_i=[n'_1i]_n=n-[n'i]_n$ holds, we get the formula as above.
\hfill$\Box$

\begin{lem}
\label{liformula}
For any $0<i<a$ the integer $l_i$ satisfying 
$$\lfloor\frac{pl_i}{|q_2|}\rfloor=[-\epsilon_2n_iq_1]_p+\begin{cases}0&\epsilon_1\epsilon_2=1\\-1&\epsilon_1\epsilon_2=-1\end{cases}$$
is $l_i=i(n\tau+\gamma'-\epsilon_1\epsilon_2)+\epsilon_1\epsilon_2n'_i$.
\end{lem}

\subsubsection{The $\gamma'+\epsilon_1\epsilon_2(n-n')\ge n$ case.}
Suppose that $\gamma'+\epsilon_1\epsilon_2(n-n')\ge n$ holds.
Since $\gamma'\le n-2$, we have $\epsilon_1\epsilon_2=1$.
Therefore $\gamma'>n'$ and $\lfloor\frac{[-\epsilon_2n_1q_1]_p}{k_1}\rfloor=\tau+1$ hold.

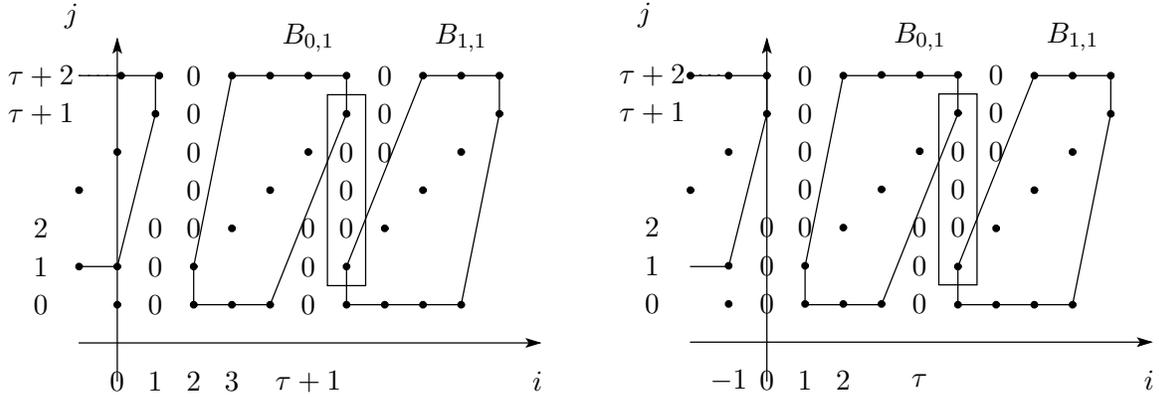
\begin{figure}[htbp]
\begin{center}
\input{tau+1.tex}
\caption{The forbidden and admitted subsequences in the case of $\lfloor\frac{[-\epsilon_2n_1q_1]_p}{k_1}\rfloor=\tau+1$}
\label{tau1}
\end{center}
\end{figure}

If $\epsilon_1=-1$, then we can find a forbidden subsequence as in the left of Figure~\ref{tau1},
which this figure is a part of the global view of $A(i+jk_1)$.
Hence none of data is realized by lens space surgery over $S^3$ or $\Sigma(2,3,5)$.
Thus we assume that $\epsilon_1=1$.
In this case we can get an admitted subsequence as the right in Figure~\ref{tau1}.

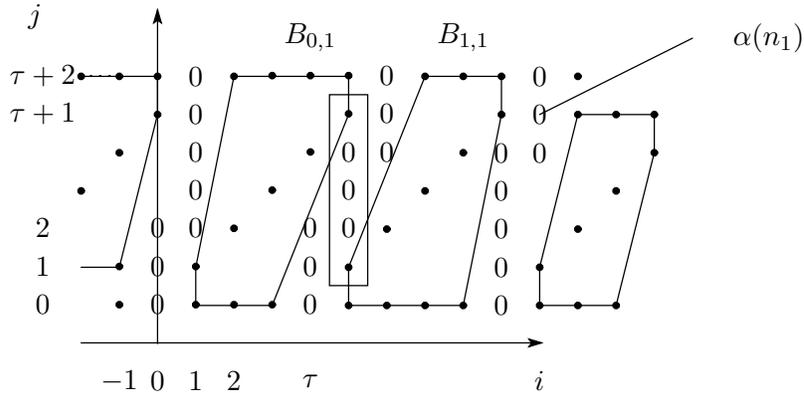
\begin{figure}[htbp]
\begin{center}
\input{tau+1sym.tex}
\caption{The admitted subsequences $\epsilon_2=1$}
\label{tau1sym}
\end{center}
\end{figure}

Note that $A(x)=\tilde{a}_{\iota+\epsilon_1x}$ for a constant $\iota$ independent of $n_2$.
From Lemma~\ref{liformula} and the symmetry of Alexander polynomial we get Figure~\ref{tau1sym}.
In the case of $\epsilon_1\epsilon_2=1$, we have $l_1=n(\tau+1)-1+\gamma'-n'$.
The block that $(i,j)=\partial^2_{lt}B$ satisfies $(i+jk_1)=[-\epsilon_2n_1q_1]_p$
is $B_{\gamma'-n',0}$.
Thus $\gamma'-n'=2$ holds.
From Lemma~\ref{intcond}, we have
$$(n'+2)^2-n'(n'+2)-n'=0\ (n),$$
thus we have $n'=n-4$.
Regarding $a(n-n')=nn_1+1$ as $\bmod n$, 
we have $4a-1=0\ (n)$.
As a result we can get
$$
\begin{cases}
a=\frac{un+1}{4} &  n_1=u\\
n'=n-4, & \gamma'=n-2
\end{cases}
$$
for some integer $u$.

On the other hand we have
$$-(-\bar{g}+(a-1)k_1+1)=-\bar{g}+(a-1)k_1-([-\epsilon_2n_1q_1]_p+k_1)k_2+1.$$
Hence $n_1=a-1$ and $n_{a-1}=1$.

Thus $u=a-1$ and we have
$$u(n-4)=3.$$
Solving this equation, we get $(a,n,\gamma')=(2,7,5),(4,5,3)$,
that is, 
$$(p,k_1)=(14\tau^2+21\tau+8,7\tau+5),(20\tau^2+25\tau+8,5\tau+3).$$
Therefore these families are realized by $\text{A}_1$ and $\text{A}_2$ of Table~\ref{po}.

\subsubsection{The $\gamma'+\epsilon_1\epsilon_2(n-n')<0$ case.}
Suppose that $\gamma'+\epsilon_1\epsilon_2(n-n')<0$ holds.
Since $\gamma'\le n-2$, we have $\epsilon_1\epsilon_2=-1$.
Then
$$\lfloor\frac{[-\epsilon_2n_1q_1]_p}{k_1}\rfloor=\tau-1$$
holds.

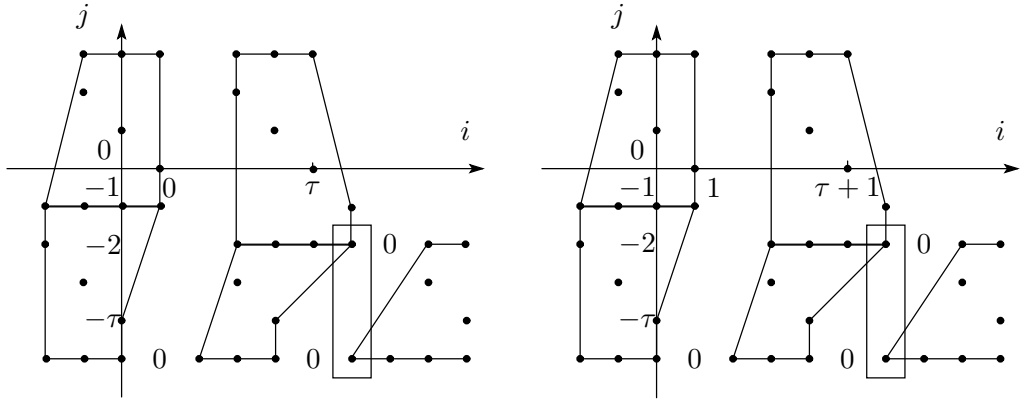
\begin{figure}[htbp]
\begin{center}
\input{tau+1-1.tex}
\caption{The forbidden and admitted subsequence in the case of $\gamma'+\epsilon_1\epsilon_2(n-n')<0$.}
\label{tau+1-1}
\end{center}
\end{figure}

In the case of $\epsilon_1=1$ we get the left of Figure~\ref{tau+1-1} in which 
there is a forbidden subsequence as indicated.
Then if we are to find lens space surgery on $S^3$ or $\Sigma(2,3,5)$, then $\epsilon_1$ must be $-1$.

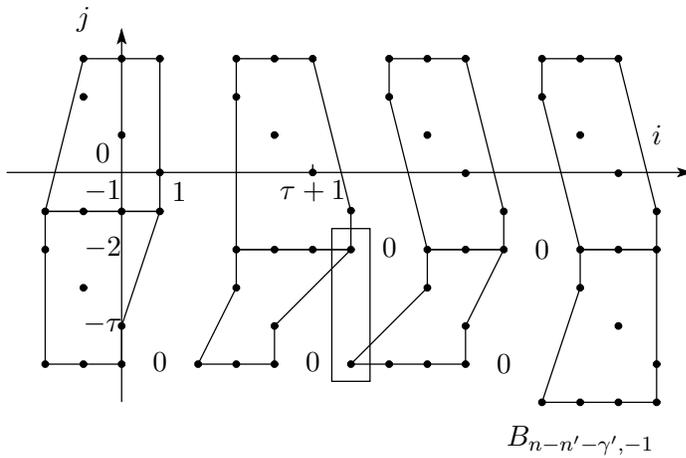
\begin{figure}[htbp]
\begin{center}
\input{tau+1-1sym.tex}
\caption{The admitted subsequence in the case of $\epsilon_1=-1$.}
\label{tau+1-1sym}
\end{center}
\end{figure}

In the case of $\epsilon_1=-1$ we get the right of Figure~\ref{tau+1-1},
so that we can find the admitted subsequence in the box indicated.
From the symmetry of Alexander polynomial $B_{-1,0}$ and $B_{0,0}$ are symmetrical about $(\tau+1,-\frac{\tau+3}{2})$
as in Figure~\ref{tau+1-1sym}.
In the other words, if one rotate the global view by $180^{\circ}$ about $(\tau+1,-\frac{\tau+3}{2})$, then one get the same global view.
From Lemma~\ref{liformula} the block $B_{0,\star}$ that the second boundary $(i,j)=\partial^2_{tl}B_{0,\star}$ satisfies $i+jk_1=[-\epsilon_2n_{a-1}q_1]_p$
is $B_{n-n'-\gamma',1}$.

Hence we have $n-n'-\gamma'=2$.
From Lemma~\ref{intcond}, in the same way 
$$\begin{cases}a=\frac{un+1}{4}& n_1=u\\n'=n-4&\gamma'=2.\end{cases}$$
On the other hand we have
$$-(-\bar{g}+(-a-1)k_1+1)=-\bar{g}+(-a-1)k_1-([-\epsilon_2n_{a-1}q_1]+3-2k_1)k_2+1.$$
Thus we have $n_{a-1}=1$ and $a_1=a-1$.
Hence we have $(n-4)u=3$.
Then $(u,n)=(1,7),(3,5)$.
In the former case
$$(p,k_1)=(14\tau^2+7\tau+1,7\tau+2)$$
In the latter case
$$(p,k_1)=(20\tau^2+15\tau+3,5\tau+2)$$
In this next section we treat the case of $0\le \gamma'+\epsilon_1\epsilon_2(n-n')<n$.

\subsection{The $2\le \gamma'\le n-2$ case (ii).}
\label{part2}
Suppose that $0\le\gamma'+\epsilon_1\epsilon_2(n-n')<n$ holds.
We claim the following Lemmas.

\begin{lem}
\label{pa-i}
For $i$ with $1\le i\le a-1$ we have
$$\lfloor\frac{p-[-\epsilon_2n_{a-i}q_1]_p}{k_1}\rfloor=i\tau+\lfloor\frac{\gamma'i-\epsilon_1\epsilon_2(n-n_i')}{n}\rfloor=i\tau+\lfloor\frac{\gamma'i-\epsilon_1\epsilon_2[n'i]_n}{n}\rfloor.$$
\end{lem}
{\bf Proof} \
Calculate as follows:
\begin{eqnarray}
\frac{p-[-\epsilon_2n_{a-i}q_1]_p}{k_1}&=&\frac{ik_1}{n}+\frac{-\epsilon_1\epsilon_2(n_{a-i}'-\epsilon_2)}{n}+\frac{-n_{a-i}'+\epsilon_2}{nk_1}\nonumber\\
&=&\frac{ik_1}{n}+\frac{-\epsilon_1\epsilon_2(n-n'_i)}{n}-\frac{n-n'_i}{nk_1}\nonumber\\
&=&i\tau+\frac{i\gamma'-\epsilon_1\epsilon_2(n-n'_i)}{n}-\frac{n-n'_i}{nk_1}.\nonumber
\end{eqnarray}
Here since we have
$$i\gamma'-\epsilon_1\epsilon_2(n-n'_i)=i\gamma'+\epsilon_1\epsilon_2n_i'=a^{-1}i(a\gamma'+\epsilon_1)=a^{-1}i(-\epsilon_2(\gamma')^{-1})\neq0\ \ (n)$$
Thus since $\frac{i\gamma'-\epsilon_1\epsilon_2(n-n'_i)}{n}\not\in{\mathbb Z}$ and $\frac{n-n'_{i}}{nk_1}<\frac{1}{k_1}<\frac1{2n}$, we have
$$\lfloor\frac{i\gamma'-\epsilon_1\epsilon_2(n-n'_i)}{n}-\frac{n-n'_i}{nk_1}\rfloor=\lfloor\frac{i\gamma'-\epsilon_1\epsilon_2(n-n'_i)}{n}\rfloor.$$
\hfill$\Box$

Here we put $n_0=0$.


\begin{lem}
\label{ldecide1}
Let $l'_i$ be the integer satisfying 
$$\lfloor\frac{(|q_2|-l'_i)p}{|q_2|}\rfloor=[-\epsilon_2n_{a-i}q_1]_p.$$
Then $l'_i$ is
$$in\tau+i\gamma'-\epsilon_1\epsilon_2(n-n'_i+i)=ik_1-\epsilon_1\epsilon_2(n-n_i'+i).$$
\end{lem}
{\bf Proof.}
Since 
$$\lfloor\frac{(|q_2|-l'_i)p}{|q_2|}\rfloor=p-\lfloor\frac{l'_ip}{|q_2|}\rfloor-1,$$
by easy calculation we get
\begin{eqnarray}
\frac{(ik_1-\epsilon_1\epsilon_2(n-n'_i+i))p}{|q_2|}&=&p-\frac{(a-i)p+(\epsilon_1\epsilon_2k_1+1)n_{a-i}}{a}-\frac{\epsilon_1\epsilon_2n_{a-i}}{|q_2|}\nonumber\\
&=&p-[-\epsilon_2n_{a-i}q_1]_p-\frac{\epsilon_1\epsilon_2n_{a-i}}{|q_2|}.\nonumber
\end{eqnarray}
Thus putting $l'_i=ik_1-\epsilon_1\epsilon_2(n-n'_i+i)$, we get
$$\lfloor\frac{(|q_2|-l'_i)p}{|q_2|}\rfloor=[-\epsilon_2n_{a-i}q_1]_p+\begin{cases}0& \epsilon_1\epsilon_2=1\\-1& \epsilon_1\epsilon_2=-1.\end{cases}$$
\hfill$\Box$

We put $n'=[-\epsilon_2a^{-1}]_n$.
Then we have
$$l'_i=i(n\tau+\gamma'-\epsilon_1\epsilon_2)-\epsilon_1\epsilon_2[n'i]_n,$$
for $n-n_i'=[-\epsilon_2a^{-1}i]_n$ and $n'=[-\epsilon_2a^{-1}]_n$.

\subsubsection{The $a=2$ case} 
\begin{lem}
Let $(p,k_1)$ be a data realized by lens space surgery of $S^3$ or $\Sigma(2,3,5)$.
In the $a=2$ case we get $\text{IX}$, $\text{X}$ or $\text{E}_i$.
\end{lem}

Here we will classify the case where $n_{1}=n_{a-1}\Leftrightarrow a=2$ holds.
We may assume that $n'\ge\gamma'+2$ holds.

Thus we have $n_1=n_{a-1}=1,\ \gamma'=\frac{\gamma-\epsilon_2}{4},\ n'=\frac{n-\epsilon_2}{2}$. 
Hence from $\gamma'\le n'-2$,
$$\gamma\le 2n-\epsilon_2-8$$
holds.

Suppose that $\epsilon_2=1$.
If $\gamma=2n-9$, then from the integrality of $p$,
we get $(\gamma,n)=(13,11)$.
This case $(p,k)=(22\tau^2+13\tau+2,11\tau+3)$ is realized by (X).

Suppose that $\epsilon_2=-1$.
If $\gamma=2n-7$, then the from integrality of $p$,
we get $(\gamma,n)=(3,5)$, hence $\gamma'=1$.
This family before we classified.

We may assume that $\gamma< 2n-\epsilon_2-12$ holds.

Now since $\lfloor\frac{n}{\alpha}\rfloor\ge n'-\gamma'$ holds,
using $n'-\gamma'=\frac{2n-\gamma-\epsilon_2}{4}$, we have
$$\frac{4n}{\gamma-\epsilon_2+4}>\frac{2n-\gamma-1}{4}\Leftrightarrow n<\frac{(\gamma+1)(\gamma-\epsilon_2+4)}{2(\gamma-\epsilon_2-4)}.$$
Hence 
$$\frac{\gamma+12+\epsilon_2}{2}<\frac{(\gamma+1)(\gamma-\epsilon_2+4)}{2(\gamma-\epsilon_2-4)}$$
Thus we have
$$\gamma<\frac{53+15\epsilon_2}{3+\epsilon_2}.$$
Hence if $\epsilon_2=1$, then $\gamma\le 16 $, and if $\epsilon_2=-1$, then $\gamma\le 18$.

If $\epsilon_2=1$, then by easy calculation, 
the possibility is $(\gamma,n)=(9,11)$ only.
Thus $(p,k_1)=(22\tau^2+9\tau+1,11\tau+2)$, and each of this family is realized by (IX).

If $\epsilon_2=-1$, then by easy calculation, 
the possibility is $(\gamma,n)=(15,27)$ only.
This case does not satisfy $\lfloor\frac{n}{\alpha}\rfloor \ge n'-\gamma'$.

Now suppose that $\lfloor\frac{n}{\alpha}\rfloor\le n'-\gamma'$ holds.
Thus $\epsilon_2=-1$ holds and from the symmetry of Alexander polynomial
$$n'-\gamma'=2\lfloor\frac{n}{\alpha}\rfloor\Leftrightarrow \frac{n'-\gamma'}{2}\le \frac{n}{\gamma'+1}<\frac{n'-\gamma'}{2}+1.$$
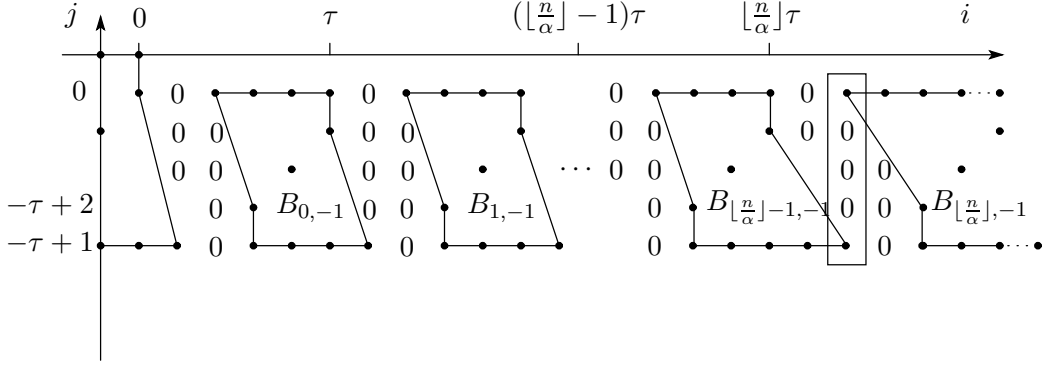
\begin{figure}[htbp]
\begin{center}
\input{tauIIppad.tex}
\caption{An admitted subsequence in the case of $a=2$}
\label{tauIIPPad}
\end{center}
\end{figure}
If $\gamma>3$, then we have
$$\gamma'+2+\frac{2}{\gamma'-3}<n'\le \gamma'+4+\frac{10}{\gamma'-3},$$
in particular $4\le n'-\gamma'$.

We assume that $n'=\gamma'+4$.
Then from Lemma~\ref{intcond} we get $n|27$, namely $n=3,9,27$.
The possibilities of them are $n=27,9$.
If $n=27$, then $n'=14,\ \gamma'=10$, and $\gamma=39$.
Therefore the initial data is 
$$(p,k_1)=(54\tau^2+39\tau+7,27\tau+10),$$
and this family is $\text{E}_2$.
If $n=9$, then $n'=5$, $\gamma'=1$.
This case is classified before.

We assume that $n'\ge \gamma'+6$.
In the inequality above $3\le \gamma'\le 8$ holds.
All the possibilities are $(\gamma',n,n')=(8,27,14),(7,25,13),(6,23,12),(6,25,13),(5,21,11),(5,23,12),$
$(5,25,13),(5,27,14),
(4,21,11),(4,23,12),(4,25,13),(4,27,14),(4,29,15),(4,31,16),$\\
$(4,33,17)$ or $(4,35,18)$.
The data which satisfies Lemma~\ref{intcond} is $(4,27,14)$ only.
This case 
$$(p,k_1)=(54\tau^2+15\tau+1,27\tau+4)$$
is class $\text{E}_{2}$.

If $\gamma'=3$, then $9-3n'-n'=9-2(n+1)=-2n+7=7=0\ (n)$, then $n=7$ holds.
Thus $n'=4,\gamma'=2$ holds.
This does not satisfy $an=nn_{a-1}+1$.
\hfill$\Box$


Now we may assume that $a\ge 3$.
\subsubsection{The $\epsilon_1\epsilon_2=1$ case.}
First we suppose that $n'-\gamma'=1$ holds.
From Lemma~\ref{intcond} $n|2\gamma'+1$ hence
\begin{equation}
[-\epsilon_2n_1q_1]_p=k_1(\tau+1)+1-\tau+\frac{2\gamma'+1}{n}.
\end{equation}
Therefore $n|2\gamma'+1$ holds.
Now since
$$2\gamma'+1\le 2(n-2)+1=2n-3,$$
$2\gamma'+1=n$ must be true as ${\mathbb Z}$.
Thus we have
$$\gamma'=\frac{n-1}{2},\ n'=\frac{n+1}{2}.$$
Taking $(n+1)a=(nn_{a-1}-\epsilon_2)2$ as $\bmod n$,
we get $n|a+2\epsilon_2$.
Let $u$ be the integer with $a=un-2\epsilon_2$.
Hence
$$n_{a-1}=\frac{a+u}{2},\ u_{1}=\frac{a-u}{2}.$$
Therefore
$$(p,k_1)=(uk_1^2-\epsilon_2\tau'(k_1-1),\frac{n\tau'-1}{2}),$$
where $\tau'=2\tau+1$.
These families are realized by $(\text{IV}_-)$. (Put $-2\epsilon_2\tau'=J$.)

Hence we may assume that $n'\ge \gamma'+2$.

\begin{figure}[htbp]
\begin{center}
\input{tauIIpp.tex}
\caption{A forbidden or admitted subsequence in the case of $\tau=5$}
\label{tauIIPP}
\end{center}
\end{figure}
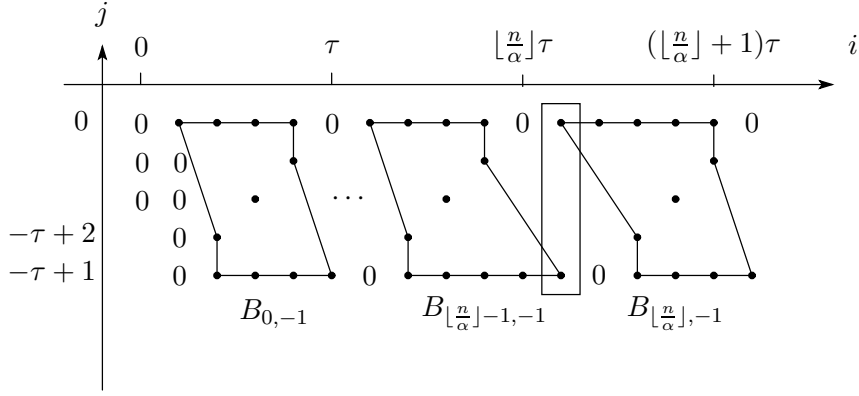

\begin{lem}
Let $(p,k_1)$ be a data with $\epsilon_1\epsilon_2=1$.
Suppose that $0<n'-\gamma'<n$ and $n'\ge \gamma'+2$.
\begin{enumerate}
\item If $\epsilon_1=1$, then $\lfloor \frac{n}{\alpha}\rfloor=n'-\gamma'$ holds.
\item If $\epsilon_1=-1$, then $\lfloor\frac{n}{\alpha}\rfloor=n'-\gamma'$ or $\lfloor\frac{2n}{\alpha}\rfloor=n'-\gamma'$ holds.
\end{enumerate}
\end{lem}
{\bf Proof} \
Suppose that $(\epsilon_1,\epsilon_2)=(1,1)$.
If $\lfloor\frac{n}{\alpha}\rfloor>n'-\gamma'$, then from the global view Figure~\ref{tauIIPP}, the indicated subsequence is 
the forbidden subsequence.
Then we deduce that $\lfloor\frac{n}{\alpha}\rfloor\le n'-\gamma'$, namely $\lfloor\frac{n}{\alpha}\rfloor=n'-\gamma'$ since
$B(\lfloor\frac{n}{\alpha}\rfloor,-1)$ is the leftmost block with the top width $\tau+1$ in $B(\star,-1)$ satisfying $\star\ge 0$.

Suppose that $(\epsilon_1,\epsilon_2)=(-1,-1)$.
If $\lfloor\frac{n}{\alpha}\rfloor>n'-\gamma'$, in the same
way we can find an admitted subsequence as in Figure~\ref{tauIIPP}.
If $\lfloor\frac{2n}{\alpha}\rfloor>n'-\gamma'$, then we can find two admitted subsequences.
Hence we have $\lfloor\frac{2n}{\alpha}\rfloor\le n'-\gamma'$.
Namely we have $\lfloor\frac{n}{\alpha}\rfloor=n'-\gamma'$ or $\lfloor\frac{2n}{\alpha}\rfloor=n'-\gamma'$
\hfill$\Box$



\begin{lem}
\label{CDE}
Let $(p,k_1)$ be a data realized by $S^3$ or $\Sigma(2,3,5)$.
Suppose that $\epsilon_1=-1$ and $\lfloor\frac{2n}{\alpha}\rfloor= n'-\gamma'$, and $a\ge3$ holds.
Then $2\lfloor\frac{n}{\alpha}\rfloor= n'-\gamma'$ holds and $(p,k_1)$ is realized by
type $\text{C}_i$ or $\text{D}_i$ for $i=1,2$.
\end{lem}
{\bf Proof.}
From Section~\ref{part2}, $0< n'-\gamma'\le n$ and Lemma~\ref{pa-i} we have
$\lfloor\frac{p-[-\epsilon_2n_{a-1}q_1]_p}{k_1}\rfloor=\tau+\lfloor\frac{\gamma'-n'}{n}\rfloor=\tau-1$.
Thus the height of $B_{0,-1},\cdots,B_{n'-\gamma'-1,-1}$ are all $\tau$.
\begin{figure}[htbp]
\begin{center}
\input{tauIIpp2.tex}
\caption{A forbidden or admitted subsequence for Lemma~\ref{CDE} and the symmetry.}
\label{tauIIPP2}
\end{center}
\end{figure}
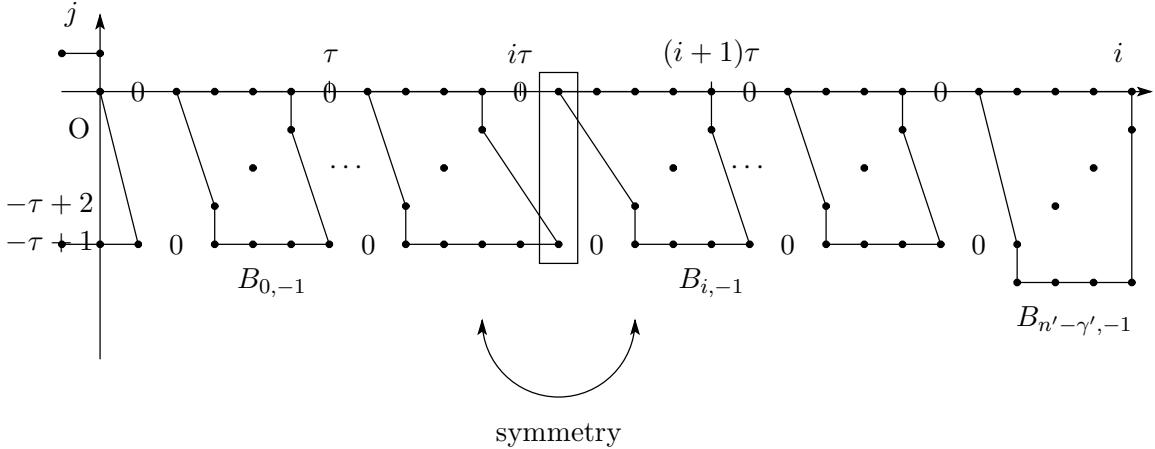
From Lemma~\ref{ifif} and the assumption, $B_{\lfloor\frac{n}{\alpha}\rfloor,-1}$ is top width $\tau+1$.
Then $n'-\gamma'\ge 2$ holds.

We define a positive integer $i$ to be $n'-\gamma'=2i$.
From the top width of $B_{j,-1}\ (0\le j<i)$ is $\tau$ and the top width of $B_{i,-1}$ is $\tau+1$,
we have
$$\alpha\le[\alpha j]_n\ (0\le j<i),\ \ 0\le [\alpha i]_n<\alpha$$
This condition is equivalently
\begin{equation}
\gamma'+1\le[(\gamma'+1)j]_n\ (0\le j<i),\ \ 0\le [(\gamma'+1) i]_n<\gamma'+1,
\label{gammanocond}
\end{equation}
for $0<\frac{aj}{|q_2|}<1$ holds.
The condition (\ref{gammanocond}) implies
$$n<i(\gamma'+1)<n+\gamma'+1<(i+1)(\gamma'+1).$$
Thus the inequality
\begin{equation}
\frac{n-i}{i}<\gamma'<\frac{n+1-i}{i-1}
\label{gammadnocond}
\end{equation}

From the integrality condition and $\gamma'+2i=n'$, we have $(2i+1)\gamma'=-2i$.
Thus for some positive integer $u'$ we have $\gamma'=\frac{u'n-2i}{2i+1}$.
The integer $u'$ satisfying (\ref{gammadnocond})
is $2$ only.
Thus we have
$$\begin{cases}
a=\frac{un+2i+1}{4i^2} & n_1=\frac{(2i-1)a-u}{2i+1}\\
n'=\frac{2n+4i^2}{2i+1} & \gamma'=\frac{2n-2i}{2i+1}.
\end{cases}
$$
The symmetry of Alexander polynomial induces symmetry of the global view as indicated in Figure~\ref{tauIIPP2}.
From the symmetry of the global view we have
$$-(-\bar{g}+(a-1)k_1+1)=-\bar{g}+(a-1)k_1-(-\epsilon_2n_{a-1}q_1+2+k_1)k_2+1,$$
thus we have $n_1=1$.
Thus we have
\begin{equation}
((2i-1)n-4i^2)u=8i^3+1.
\label{in}
\end{equation}

In the case of $i=1$, (\ref{in}) is $u(n-4)=9$, thus $(u,n)=(1,13),(9,5),(3,7)$.
The data satisfying integral condition of $n'$ are 
$$(p,k_1)=(52\tau^2+63\tau+19,13\tau+8),(42\tau^2+47\tau+13,7\tau+4).$$

In the case of $i=2$, (\ref{in}) is $u(3n-16)=65$, thus $(u,n)=(1,27),(13,7)$.
Thus we have
$$(p,k_1)=(54\tau^2+39\tau+7,27\tau+10),(42\tau^2+23\tau+3,7\tau+2).$$
The former case is inconsistent with the assumption $a\ge 3$.

In the case of $i=3$ the equation (\ref{in}) is $u(5n-36)=217$, thus $(u,n)=(1,217),(217,1),(7,31),(31,7)$.
The data satisfying integral condition of $n'$ is $(u,n)=(7,31)$.
The case $(7,31)$ fails to integrality of $a$.

In the case of $i=4$ the equation (\ref{in}) is $u(7n-64)=513$, thus $(u,n)=(19,13)$.
Thus we have
$$(p,k_1)=(52\tau^2+15\tau+1,13\tau+2).$$

In the case of $i=5$ the equation (\ref{in}) is $u(9n-100)=1001$, thus $(u,n)=(19,27)$.
This is inconsistent with the assumption $a\ge 3$.

Suppose that $i\ge 6$.
From
$$n'-\gamma'=2\lfloor\frac{n}{\alpha}\rfloor<2\frac{n}{\gamma'+1}$$
we can get $$n'<\frac{n_{a-1}\gamma'(\gamma'+1)-2}{n_{a-1}(\gamma'+1)-2a}.$$
Hence we have 
$$\gamma'+12<\frac{n_{a-1}\gamma'(\gamma'+1)-2}{n_{a-1}(\gamma'+1)-2a}\Rightarrow\ \ \gamma'<\frac{12a-1}{5n_{a-1}-n_1}-1=\frac{7a+5}{5a-6}\le \frac{26}{9}.$$
Hence $\gamma'=2$ holds.

In the case of $\gamma'=2$, $n=3i+1$ and 
$$u((2i-1)(3i+1)-4i^2)=(2i+1)(4i^2-2i+1)\Leftrightarrow u(2i+1)(i-1)=(i-1)(4i+2)+3.$$
Then $i-1|3$ 
Hence $i=2,4$.
This case is classified before.
Therefore the appearing cases are all realized by type $C_i,D_i$.\hfill$\Box$

Suppose that $n'-\gamma'=\lfloor\frac{n}{\alpha}\rfloor$ holds.

\begin{lem}
\label{sim}
Let $(p,k_1)$ be a data with $\epsilon_1\epsilon_2=1$ and $\lfloor\frac{n}{\alpha}\rfloor=n'-\gamma'\ge2$.
Then we have $0<n'<\frac{n}2$.
\end{lem}
{\bf Proof} \
Suppose that $\frac{n}{2}<n'<n$ holds.
In the case of $n'= \gamma'+i\ \ (i\ge 2)$ we have $(n'-i)^2-n'(n'-i)-n'=0\ (n)$, $(i+1)n'=i^2\ (n)$.
From $\alpha\le [\alpha(n(-\tau+1)+j+1)]_n<n\ \ (j=1,2,\cdots,i-1)$ and $\alpha\le [\alpha(n(-\tau+1)+i+1)]_n<n$,
we have 
$$\frac{n+i(i+1)}{i+1}-1\le n'<\frac{n+i^2}{i}-1.$$
Thus we have $n'=\frac{n+i^2}{i+1}$.
Then we have
$$\begin{cases}
a=\frac{un-\epsilon_2(i+1)}{i^2}& n_1=\frac{ia-u}{i+1}\\
n'=\frac{n+i^2}{i+1} & \gamma'=\frac{n-i}{i+1}.
\end{cases}$$

From $\frac{n}{2}<\frac{n+i^2}{i+1}$ we have $n<2(i+1)+\frac{2}{i-1}$.
If $i\ge 3$, we have $n<2i+3$.
Since $i=\lfloor\frac{n}{\alpha}\rfloor<\frac{n}{3}<\frac{2i+3}{3}$, we have $i<3$.
This is a contradiction.
If $i=2$, then we have $\gamma'=\frac{n-2}{3}$ and $n<8$.
The possibility is $n=5$ only since the integrality of $\gamma'$.
In this case $a=\frac{5u-3\epsilon_2}{4}$, $n_1=\frac{2a-u}{3}$, $n'=3$ and $\gamma'=1$.
This case is already classified before.
Therefore we have $0<n'<\frac{n}{2}$.
\hfill$\Box$\\
Lemma~\ref{sim} holds even in the case of $\epsilon_1\epsilon_2=-1$ replacing $\gamma'$ and $\alpha$ with $n-\gamma'$ and $n-\alpha$.
The proof is skipped because the process is similar to Lemma~\ref{sim}.

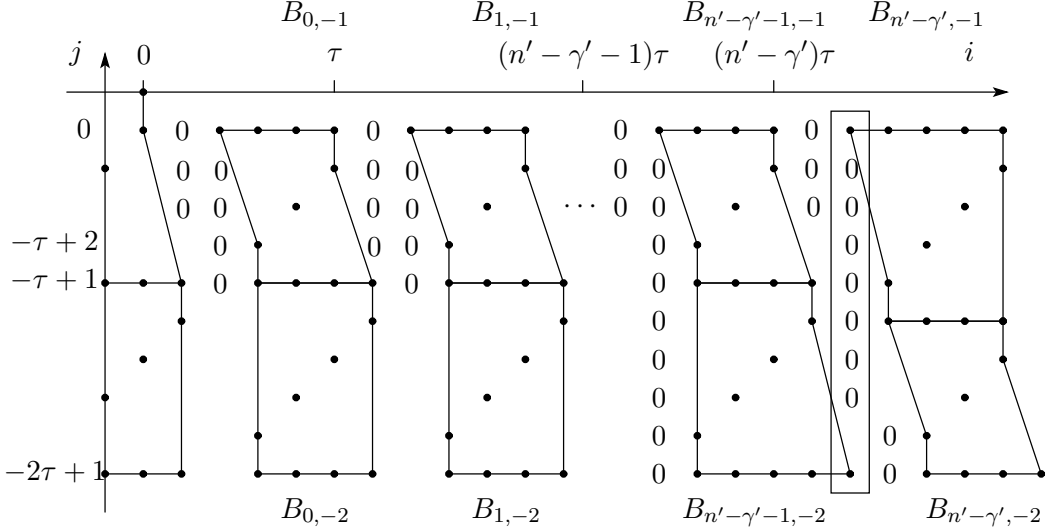
\begin{figure}[htbp]
\begin{center}
\input{height2.tex}
\caption{The place in which the heights of blocks change.($0<n'<\frac{n}{2}$)}
\label{hei2}
\end{center}
\end{figure}

\begin{lem}
Let $(p,k_1)$ be a stable data.
If $a\ge 3$, $0<n'-\gamma'<n$ $\epsilon_1=1$ and $0<n'<\frac{n}{2}$, then none of data satisfying this condition is
realized by $S^3$ or $\Sigma(2,3,5)$.
\end{lem}
{\bf Proof} \
From Lemma~\ref{pa-i}, we get
$$\lfloor\frac{p-[-\epsilon_2n_{a-2}n_1q_1]_p}{k_1}\rfloor=2\tau+\lfloor\frac{2\gamma'-2n'}{n}\rfloor=2\tau-\lfloor\frac{2n'-2\gamma'}{n}\rfloor-1=2\tau-1.$$
Thus we can find a forbidden subsequence from the Figure~\ref{hei2}.
\hfill$\Box$

Then Lemma~\ref{pa-i} $\lfloor\frac{p-[-\epsilon_2n_{a-2}q_1]_p}{k_1}\rfloor=2\tau+\lfloor\frac{2\gamma'-2n'}{n}\rfloor= 2\tau-1$ holds,
because $n'-\gamma'<\frac{n}{2}$.
\begin{lem}
Let $(p,k_1)$ be a data.
Suppose that $\epsilon_1=-1$.
If $n'-\gamma'=\lfloor\frac{n}{\alpha}\rfloor$, then the data $(p,k_1)$ is
either of $\text{F}_1,\text{F}_2,\text{G}_1,\text{G}_2,\text{H}_1$, or $\text{H}_2$.
\end{lem}
{\bf Proof.}
The global view is as in Figure~\ref{hei2} since $0<n'<\frac{n}{2}$ and Lemma~\ref{pa-i}.
From this view we have
$$-(-\bar{g}+(a-1)k_1+1)=-\bar{g}+(a-1)k_1-(-\epsilon_2n_{a-2}q_1+2+k_1)k_2+1,$$
$n_2=1$, hence $n_1=\frac{a+1}{2}$.

Putting $n'-\gamma'=i$, we have 
$$\begin{cases}
a=\frac{un+i+1}{i^2}&n_1=\frac{ia-u}{i+1}\\
n'=\frac{n+i^2}{i+1}& \gamma'=\frac{n-i}{i+1}.
\end{cases}$$
Since $a=\frac{2u+i+1}{i-1}$, 
$$((i-1)n-2i^2)u=(i+1)(i^2-i+1)$$

If $i\le 5$ holds, the possibility is $(i,u,n,a,n',\gamma')=(2,1,17,5,7,5),(2,3,11,9,5,3),(3,7,11,9,5,2),$
$(3,1,23,3,8,5),$
$(5,7,17,5,7,2),(5,3,23,3,8,3)$
$$(p,k_1)=(85\tau^2+49\tau+7,17\tau+5),(99\tau^2+53\tau+7,11\tau+3)$$
$$(p,k_1)=(99\tau^2+35\tau+3,11\tau+2),(69\tau^2+29\tau+3,23\tau+5)$$
$$(p,k_1)=(85\tau^2+19\tau+1,17\tau+2),(69\tau^2+17\tau+1,23\tau+3).$$

If $i\ge 6$ holds,
$\gamma'<\frac{3a-2}{2a+3}< \frac{3}{2}$ holds.
Since the assumption is $\gamma'\ge 2$, we can get to the assertion.
\hfill$\Box$

\begin{figure}[htbp]
\begin{center}
\input{height4.tex}
\caption{}
\label{hei4}
\end{center}
\end{figure}

\subsubsection{The $\epsilon_1\epsilon_2=-1$ case.}
We assume that $0< \gamma'-n+n'<n$ and $\epsilon_1\epsilon_2=-1$.

We consider $\gamma'-n+n'=1$.
From integral condition $2n'-1=0\ (n)$ holds.
Then we have $n'=\frac{n+1}{2}$.
Thus we have $\gamma'=\frac{n+1}{2}=n'$.
We have $k_1=n\tau+\frac{n+1}{2}$.
Taking $2a(n-n')=2nn_1+2\epsilon_2$ as ${\mathbb Z}/p{\mathbb Z}$,
we have $a=-2\epsilon_2\ (n)$.
For some integer $u$ we have $a=un-2\epsilon_2$.
Then we have $(p,k_1)=(uk_1^2-\epsilon_2(2\tau+1)(k_1+1),n\tau+\frac{n+1}{2})$.
This case is realized by type ($\text{IV}_+$).

We may assume that $\gamma'-n+n'\ge2$ holds.
Hence $\lfloor\frac{[-\epsilon_2n_1q_1]_p}{k_1}\rfloor=\tau-1+\lfloor\frac{n'+\gamma'}{n}\rfloor=\tau$ holds.

Now we assume that $\lfloor\frac{(n\tau+2)p}{|q_2|}\rfloor-\lfloor\frac{(n\tau+1)p}{|q_2|}\rfloor=\tau$.
Equivalently $\alpha\le[\alpha(n\tau+2)]_n<n$ holds.
Here we have $[\alpha(n\tau+2)]_n=[2(\gamma'-1)+\frac{an\tau+2a}{|q_2|}]_n$ and
$$\frac{an\tau+2a}{|q_2|}=1-\frac{a\gamma'-a+\epsilon_1}{|q_2|}<1.$$
Thus $\alpha<[2\gamma'-2+\frac{an\tau+2a}{|q_2|}]_n$ holds and $2\gamma'-2<n$ holds.
As a result $\gamma'<\frac{n+2}{2}$ holds.
We put $\gamma''=n-\gamma'$ and $\alpha'=\gamma''+1-\frac{a}{|q_2|}=n-\alpha$.

\begin{lem}
\label{ee=-1sepr}
Let $(p,k_1)$ be a data with $\epsilon_1\epsilon_2=-1$, $0<n'-\gamma''<n$, and $n'\ge \gamma''+2$.
Then the possibilities are either of the following.
\begin{enumerate}
\item If $\epsilon_1=-1$, then $\lfloor\frac{n}{\alpha'}\rfloor=n'-\gamma''$.
\item If $\epsilon_1=1$, then $\lfloor\frac{n}{\alpha'}\rfloor=n'-\gamma''$ or $\lfloor\frac{2n}{\alpha'}\rfloor=n'-\gamma''$.
\end{enumerate}
\end{lem}
{\bf Proof.} \
Suppose that $(\epsilon_1,\epsilon_2)=(-1,1)$.
If $\lfloor\frac{n}{\alpha'}\rfloor<n'-\gamma''$, the global view of $A(i+jk_1)$ is as Figure~\ref{ee-1}.
Here we must note that $B_{\lfloor\frac{n}{\alpha'}\rfloor-1,0}$ is the leftmost block with the top width $\tau$ in the blocks of $B_{\ast,0}\ (\ast\ge 0)$.
\begin{figure}[htbp]
\begin{center}
\input{ee=-1gen.tex}
\caption{A global view of $(\epsilon_1,\epsilon_2)=(-1,1)$}
\label{ee-1}
\end{center}
\end{figure}
Therefore we can find a forbidden subsequence as indicated by the box.
Thus we have $\lfloor\frac{n}{\alpha'}\rfloor\ge n'-\gamma''$.
Since the block $B_{\lfloor\frac{n}{\alpha'}\rfloor,0}$ is the leftmost block the top width $\tau$ we have
$\lfloor\frac{n}{\alpha'}\rfloor=n'-\gamma''$.

Suppose that $(\epsilon_1,\epsilon_2)=(1,-1)$.
If $\lfloor\frac{2n}{\alpha'}\rfloor<n'-\gamma''$, the global view of $A(i+jk_1)$ is as Figure~\ref{ee-12}.
\begin{figure}[htbp]
\begin{center}
\input{ee=-1gen2.tex}
\caption{A global view of $(\epsilon_1,\epsilon_2)=(1,-1)$}
\label{ee-12}
\end{center}
\end{figure}
In this time there are two admitted subsequences in this view.
This is contradiction.
Thus $\lfloor\frac{2n}{\alpha'}\rfloor\ge n'-\gamma''$ holds.
Since the block $B_{n'-\gamma'',0}$ has the bottom width $\tau$, $n'-\gamma''=\lfloor\frac{n}{\alpha'}\rfloor$ or $\lfloor\frac{2n}{\alpha'}\rfloor$.
\hfill$\Box$\\

We first classify the case of $\lfloor\frac{2n}{\alpha'}\rfloor=n'-\gamma''$.
\begin{lem}
\label{cde2}
Suppose that $\epsilon_1=1$, $\lfloor\frac{2n}{\alpha'}\rfloor=n'-\gamma''$ and $a\ge 3$.
Then the data is realized by one of $\text{C}_i$ and $\text{D}_i$ ($i=1,2$).
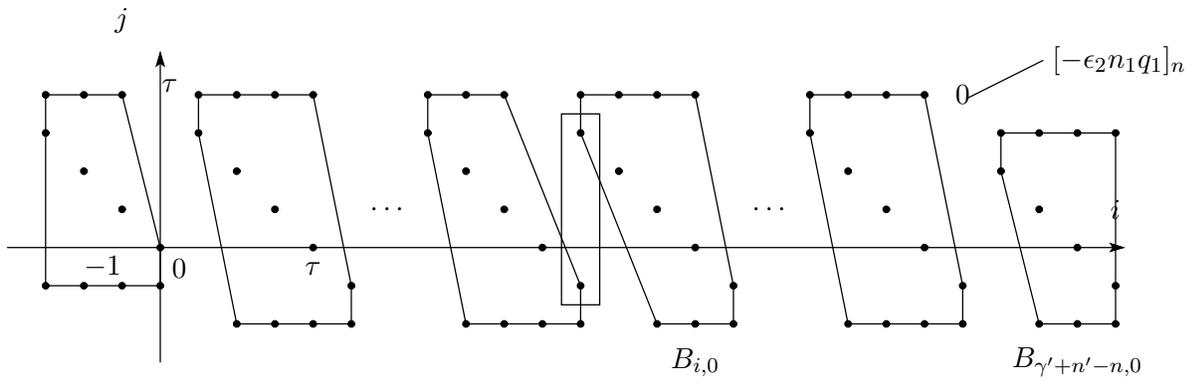
\begin{figure}[htbp]
\begin{center}
\input{ntau+2taugen.tex}
\caption{A global view of Lemma~\ref{cde2}}
\label{ntau+2taugen}
\end{center}
\end{figure}
\end{lem}
{\bf Proof.} \
From the assumption the global view of $(p,k_1)$ is as in Figure~\ref{ntau+2taugen}.
First for a positive integer $j$ the equation $\frac{a(n\tau+1+j)}{|q_2|}=1-\frac{a\gamma'-aj+1}{|q_2|}$ holds.
Furthermore we have $$0<\frac{a\gamma'-aj+1}{|q_2|}<\frac{a(\gamma'-1)+1}{2an+a\gamma'-a+1}=1-\frac{2an}{2an+a\gamma'-a+1}<1.$$
As a result we have $0<\frac{a(n\tau+1+j)}{|q_2|}<1$.

If 
$$\lfloor\frac{(n\tau+1+j)p}{|q_2|}\rfloor-\lfloor\frac{(n\tau+j)p}{|q_2|}\rfloor=\tau+1$$
holds, then from Lemma~\ref{ifif} $0\le [\alpha(n\tau+j+1)]_n=[(j+1)(\gamma'-1)+\frac{a(n\tau+1+j)}{|q_2|}]_n<\alpha$ holds.
This condition is equivalent to the following:
\begin{eqnarray}
&\Leftrightarrow& 0\le\left[(j+1)(\gamma'-1)+\frac{a(n\tau+i)}{|q_2|}\right]_n< \gamma'-1\nonumber\\
&\Leftrightarrow& 0\le[(j+1)(\gamma'-1)]_n< \gamma'-1.\nonumber
\end{eqnarray}

Now suppose that the global view is as Figure~\ref{ntau+2taugen}.
Then we have 
$0\le[(j+1)(\gamma'-1)]_n< \gamma'-1$ for some integer with $0<j<i$ and we have $\gamma'-1\le[(i+1)(\gamma'-1)]_n< n$.
Therefore we have 
\begin{equation}
(i-1)n\le i(\gamma'-1)<(i-1)n+\gamma'-1<(i+1)(\gamma'-1)<in.
\label{gadivi}
\end{equation}
Hence 
$$\frac{(i-1)n+i}{i}\le\gamma'<\frac{in+i+1}{i+1}$$ holds.

In the Figure~\ref{ntau+2taugen} we can find an admitted subsequence as indicated by a box in the figure.
Furthermore from Lemma~\ref{liformula} the block $B_{\star,0}$ that the point $(i,j)=\partial^2_{lt}B_{\star,0}$ satisfies
$i+jk_1=[-\epsilon_2n_1q_1]_p$ and is the closest to the origin is $B_{\gamma'+n'-n,0}$.
Now from the symmetry of Alexander polynomial we put $\gamma'+n'-n=2i$.
Hence (\ref{gadivi}) implies
$$\frac{n+2i^2+i-1}{i+1}<n'\le \frac{n+2i^2-i}{i}.$$

From the integral condition $(2i+1)n'=4i^2\ (n)$ and $4i^2a=2i+1\ (n)$ hold.
Therefore $n'=\frac{2n+4i^2}{2i+1}$ and $a=\frac{un+2i+1}{4i^2}$ hold.
Hence
\begin{equation}
\begin{cases}
a=\frac{un+2i+1}{4i^2}& n_1=\frac{(2i-1)a-u}{2i+1}\\
n'=\frac{2n+4i^2}{2i+1} & \gamma'=\frac{(2i-1)n+2i}{2i+1}.
\end{cases}
\label{2n}
\end{equation}
Furthermore from the symmetry of Alexander polynomial
$$-(-\bar{g}+(-a-1)k_1+1)=-\bar{g}+(-a-1)k_1+(-\epsilon_2n_1q_1+1-2k_1)k_2+1\ (p).$$
Hence we have $n_1=1$.
Using (\ref{2n}) and $n_1=1$
$$((2i-1)n-4i^2)u=8i^3+1.$$
In the case of $i=1$ we get
$$(p,k_1)=(52\tau^2+41\tau+8,13\tau+5),(42\tau^2+37\tau+8,7\tau+3).$$
In the case of $i=2$ we get
$$(p,k_1)=(54\tau^2+69\tau+22,27\tau+17),(42\tau^2+61\tau+22,7\tau+5).$$
The former case is inconsistent with $a\ge 3$.
In the case of $i=3$ there does not exist any suitable pair $(u,n)$.
In the case of $i=4$ we get
$$(p,k_1)=(52\tau^2+89\tau+38,13\tau+11).$$
In the case of $i=5$ we get
$$(p,k_1)=(54\tau^2+93\tau+40,27\tau+23).$$
This is inconsistent with $a\ge 3$.

If $i\ge 6$, then $\gamma'-n+n'=2i\ge 12$ holds.
From the form of the blocks, $\lfloor\frac{n}{n-\alpha}\rfloor=\frac{\gamma'-n+n'}{2}$ holds.

Thus if $\gamma''\neq 2$ or $(\gamma'',a)\neq(3,3)$ we have
$$n'<\frac{n_{a-1}(\gamma'')^2-2}{n_{a-1}\gamma''-2a},$$
where $\gamma''=n-\gamma'$.

Hence we have
$$\gamma''+12<\frac{n_{a-1}(\gamma'')^2-2}{n_{a-1}\gamma''-2a}\Leftrightarrow\gamma''<\frac{12a-1}{5n_{a-1}-n_1}=\frac{12a-1}{5a-6}\le \frac{35}{9}.$$
Thus we have $\gamma''=3$ and $a\neq 3$.

If $\gamma''=3$, then from the integral condition we have $9-4n'=0(n)$.
Then we have $n'=\frac{n+9}{4}$ since $\frac{n'-3}{2}=\lfloor\frac{n}{4-\frac{a}{|q_2|}}\rfloor<\frac{2n}{7}$.
Thus we have $n=8i+3$ and $n'=2i+3$.
Here since
$$i=\lfloor\frac{8i+3}{4-\frac{a}{|q_2|}}\rfloor\ge\lfloor\frac{8i+3}{4}\rfloor=2i$$
this is a contradiction.

If $\gamma''=2$, then from the integral condition we have $n'=\frac{2n+4}{3}$.
Then we have $n=3i+1$, $(i-1)u=4i^2-2i+1$ and $u=\frac{4i^2-2i+1}{i-1}=4i+2+\frac{3}{i-1}$ hold.
Thus $i-1=1,3$ holds these are inconsistent with $i\le 6$.
\hfill$\Box$

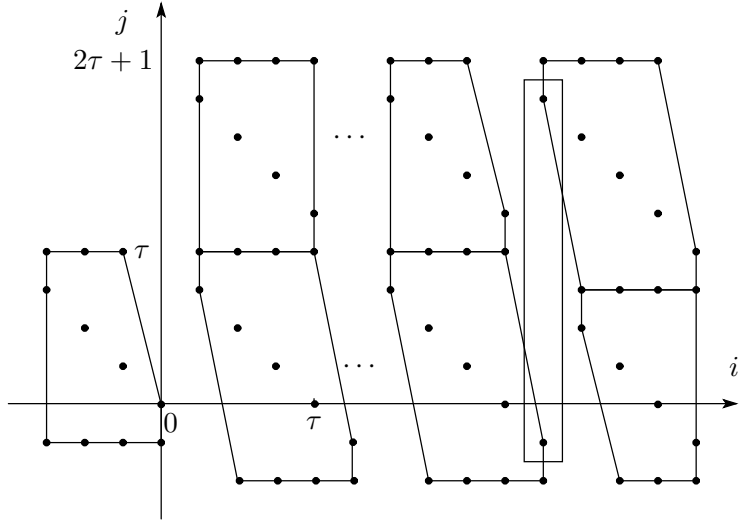
\begin{figure}[htbp]
\begin{center}
\input{ntau+2tau.tex}
\caption{A forbidden or admitted subsequences with $\tau=4$}
\label{ntau+2tau}
\end{center}
\end{figure}

\begin{lem}
\label{fgh2}
If $0<n'<\frac{n}{2}$ and $\lfloor\frac{n}{\alpha'}\rfloor=n'-\gamma''$ hold,
then $\epsilon_1=1$ and the global view of a data $(p,k_1)$ is Figure~\ref{ntau+2tau}, the data is realized by
one of type $F_i,G_i$ and $H_i$ for $i=1,2$.
\end{lem}
{\bf Proof} \
From the assumption the global view is as Figure~\ref{ntau+2tau}.
For since from Lemma~\ref{taufloorfunc} 
$\lfloor\frac{[-\epsilon_2n_1q_1]_p}{k_1}\rfloor=\tau+\lfloor\frac{n'-\gamma''}{n}\rfloor=\tau$ and
$\lfloor\frac{[-\epsilon_2n_2q_1]_p}{k_1}\rfloor=2\tau-1+\lfloor\frac{2\gamma'+[2n']_n}{n}\rfloor=2\tau+1+\lfloor\frac{2n'-2\gamma''}{n}\rfloor=2\tau+1$,
we have $H(B_{0,0})=\tau+3$ and $H(B_{0,1})=\tau+2$.
The leftmost block $B_{\ast,0}$ in which $(i,j)=\partial^2_{lt}B_{\ast,0}$
satisfies $i+jk_1=[-\epsilon_2n_1q_1]_p$ for positive integer $\ast$ is $B_{n'-\gamma'',0}$.
Namely $H(B_{j,0})=\tau+3\ (0\le j<n'-\gamma'')$ and $H(B_{n'-\gamma'',0})=\tau+2$.
This is due to Lemma~\ref{liformula}.

From the global view we can find a subsequence as indicated by the box in Figure~\ref{ntau+2tau}.
This subsequence is the forbidden subsequence when $\epsilon_1=-1$, and the admitted subsequence when $\epsilon_1=1$.
To realize lens space surgery over $S^3$ or $\Sigma(2,3,5)$ we must take $\epsilon_1=1$.
Therefore 1. part in Lemma~\ref{ee=-1sepr} does not occur as long as lens space surgery over $S^3$ or $\Sigma(2,3,5)$.

From Lemma~\ref{cde2} 
$$\frac{n-i-1}{i+1}<\gamma''\le\frac{n-i}{i}$$ holds.
We denote $n'-\gamma''$ by $i$.

From the integrality condition we have for some integer $u$, we get
$$\begin{cases}
a=\frac{un+i+1}{i^2}& n_1=\frac{ia-u}{i+1}\\
n'=\frac{n+i^2}{i+1}& \gamma''=\frac{n-i}{i+1}.
\end{cases}$$

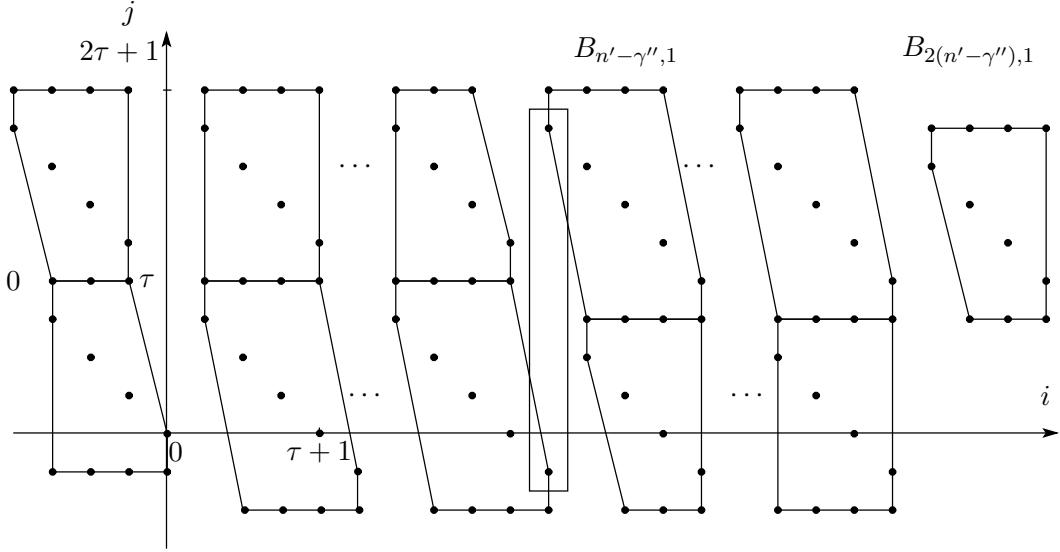
\begin{figure}[htbp]
\begin{center}
\input{ntau+2tausym.tex}
\caption{A forbidden or admitted subsequences with $\tau=4$}
\label{ntau+2tausym}
\end{center}
\end{figure}

Furthermore from the symmetry of Alexander polynomial we have
$$-(\bar{g}+(-a-1)k_1+1)=-\bar{g}+(-a-1)k_1+(-\epsilon_2n_2q_1+1-2k_1)k_2+1$$
and $n_2=1$.
Figure~\ref{ntau+2tausym} illustrates the symmetry.

Thus we have $n_1=\frac{a+1}{2}$.
Hence
\begin{equation}
((i-1)n-2i^2)u=i^3+1\label{eq2}.
\end{equation}
In the case of $i=2$, the possibilities are $(n,u)=(9,9),(11,3),(17,1)$.
If $(n,u)=(9,9)$, then $n'=\frac{13}{3}\not\in{\mathbb Z}$, hence this is inconsistent.
If $(n,u)=(11,3),(17,1)$, then
$$(p,k_1)=(99\tau^2+97\tau+53,11\tau+8),(85\tau^2+121\tau+43,17\tau+12).$$

In the case of $i=3$, the possibilities are $(n,u)=(11,7),(23,1)$.
Hence we have
$$(p,k_1)=(99\tau^2+163\tau+67,11\tau+9),(69\tau^2+109\tau+43,23\tau+18).$$

In the case of $i=4$ the possibilities are $(n,u)=(11,65),(15,5)$.
If $(n,u)=(11,65)$, then $n'=\frac{27}{5}\not\in{\mathbb Z}$, hence this is inconsistent.
If $(n,u)=(15,5)$, then $n'=\frac{31}{5}\not\in{\mathbb Z}$, hence this is inconsistent as well.

In the case of $i=5$, the possibilities are $(n,u)=(13,63),(17,7),(23,3)$.
If $(n,u)=(13,63)$, then $n'=\frac{13}{2}\not\in{\mathbb Z}$, hence this is inconsistent.
If $(n,u)=(17,7),(23,3)$, then
$$(p,k_1)=(85\tau^2+151\tau+67,17\tau+15),(69\tau^2+121\tau+53,23\tau+20)$$

If $i\ge 6$, then using $\lfloor\frac{n}{n-\alpha}\rfloor=n'-\gamma''$ we get
$$n'-\gamma''<\frac{n}{\gamma''+\frac{1}{2}},$$
so that we get
$$n'<\frac{n_{a-1}\gamma''(\gamma''+\frac{1}{2})-1}{n_{a-1}(\gamma''+\frac{1}{2})-a}.$$
Hence we have $\gamma''<\frac{9a+1}{2(2a-3)}=\frac{9}{4}+\frac{29}{4(2a-3)}\le\frac{14}{3}$.
Thus $\gamma''\le 4$ holds.

If $\gamma''=4$, then from the integrality condition $5n'=16(n)$, then
from $n'=\frac{n+16}{5}=\frac{n+i^2}{i+1}$ we have $n=5i+4$.
Using (\ref{eq2}), we have $u=\frac{i^2-i+1}{3i-4}$.
If $i=3i_1$ for some positive integer $i_1$, then $u=\frac{9i_1^2-3i_1+1}{9i_1-4}=i_1+\frac{i_1+1}{9i_1-4}\not\in{\mathbb Z}$.
If $i=3i_1+1$ for some positive integer $i_1$, then $u=\frac{9i_1^2+3i_1+1}{9i_1-1}=i_1+\frac{4i_1+1}{9i_1-1}\not\in{\mathbb Z}$.
Thus we cannot find any integer solution in (\ref{eq2}), this implies that there is no lens space surgery on $S^3$ or $\Sigma(2,3,5)$ in this case.

If $\gamma''=3$, then from the integrality condition $4n'=9(n)$, then
from $n'=\frac{n+9}{4}=\frac{n+i^2}{i+1}$ we have $n=4i+3$.
Using (\ref{eq2}), we have $u=\frac{i^2-i+1}{2i-3}$.
If $i=2i_1$ for some positive integer $i_1$, then $u=\frac{4i_1^2-2i_1+1}{2i_1-3}=i_1+\frac{i_1+1}{4i_1-3}\not\in{\mathbb Z}$.
If $i=2i_1+1$ for some positive integer $i_1$, then $u=\frac{4i_1^2+2i_1+1}{4i_1-1}=i_1+\frac{3i_1+1}{4i_1-1}$.
Thus we have $i_1=2$, $u=3$, and $i=5$ and this is inconsistent with $i\ge 6$.
There is no lens space surgery on $S^3$ or $\Sigma(2,3,5)$ in this case as well.

If $\gamma''=2$ holds, then from the integrality condition $3n'=4(n)$ and $n'=\frac{n+4}{3}=\frac{n+i^2}{i+1}$ we have $n=3i+2$.
Using (\ref{eq2}), we have $u=\frac{i^2-i+1}{i-2}=i+1+\frac{3}{i-2}$.
Thus we have $i-2=1,3$ these cases are inconsistent with $i\ge 6$.
Therefore in the case of $i\ge 6$ we cannot find any lens space surgery over $S^3$ or $\Sigma(2,3,5)$.\hfill$\Box$

In this point we can classify all stable data $(p,k_1)$ with lens space surgery $S^3$ or $\Sigma(2,3,5)$.
To prove that Berge's list is fully complete, we must argue unstable data: namely $\tau=1$ or $0$.
This argument will be done in sequent papers.

\section*{Acknowledgments}
This research is partially supported by the Grant-in-Aid for JSPS Fellow for Young Scientist (21-1458).
The author is grateful for Professor M.Ue and T.Ohtsuki affiliated in Mathematics department
and RIMS, Kyoto University respectively,
who give me a lot of useful comments in a series of arguments.
Furthermore a part of the main theorem was proved while my visit of Michigan State University consequently
the author thanks for the hospitality of the institute and the acceptance of my visit by Professor Akbulut.
Finally, the author deeply appreciates Professors T.Kadokami and Y.Yamada giving me the motivation
for this research.

 \noindent
 Motoo Tange\\
 Research Institute for Mathematical Sciences, \\
 Kyoto University, \\
 Kyoto 606-8502, Japan. \\
 tange@kurims.kyoto-u.ac.jp

\end{document}

%% file: ee=1not.tex
\unitlength 0.1in
\begin{picture}( 40.6500, 12.0000)( 11.3500,-26.0000)
\put(28.0000,-16.0000){\makebox(0,0){$\lfloor\frac{p\ell}{|q_2|}\rfloor-1$}}%
\put(36.0000,-16.0000){\makebox(0,0){$\lfloor\frac{p\ell}{|q_2|}\rfloor$}}%
\put(44.0000,-16.0000){\makebox(0,0){$\lfloor\frac{p\ell}{|q_2|}\rfloor+1$}}%
\put(28.0000,-19.0000){\makebox(0,0){$0$}}%
\put(36.0000,-19.0000){\makebox(0,0){$1$}}%
\put(44.0000,-19.0000){\makebox(0,0){$0$}}%
\put(28.0000,-21.0000){\makebox(0,0){$-1$}}%
\put(36.0000,-21.0000){\makebox(0,0){$0$}}%
\put(44.0000,-21.0000){\makebox(0,0){$-1$}}%
%
\special{pn 8}%
\special{pa 2400 1800}%
\special{pa 4800 1800}%
\special{fp}%
\special{pa 4800 1800}%
\special{pa 4800 2200}%
\special{fp}%
\special{pa 4800 2200}%
\special{pa 2400 2200}%
\special{fp}%
\special{pa 2400 2200}%
\special{pa 2400 1800}%
\special{fp}%
\special{pa 2400 2000}%
\special{pa 4800 2000}%
\special{fp}%
\special{pa 3200 1800}%
\special{pa 3200 2200}%
\special{fp}%
\special{pa 4000 1800}%
\special{pa 4000 2200}%
\special{fp}%
%
\special{pn 8}%
\special{pa 2400 1800}%
\special{pa 2400 1400}%
\special{fp}%
\special{pa 2400 1400}%
\special{pa 4800 1400}%
\special{fp}%
\special{pa 4800 1800}%
\special{pa 4800 1800}%
\special{fp}%
\special{pa 4800 1400}%
\special{pa 4800 1800}%
\special{fp}%
\special{pa 4000 1400}%
\special{pa 4000 1800}%
\special{fp}%
\special{pa 3200 1400}%
\special{pa 3200 1800}%
\special{fp}%
%
\special{pn 8}%
\special{pa 4800 1400}%
\special{pa 5200 1400}%
\special{fp}%
\special{pa 5200 1400}%
\special{pa 5200 2200}%
\special{fp}%
\special{pa 5200 2200}%
\special{pa 4800 2200}%
\special{fp}%
%
\special{pn 8}%
\special{pa 4800 1800}%
\special{pa 5200 1800}%
\special{fp}%
\put(50.0000,-19.0000){\makebox(0,0){(+,i)}}%
\put(50.0000,-21.0000){\makebox(0,0){(+,ii)}}%
\put(50.0000,-16.0000){\makebox(0,0){$x$}}%
%
\special{pn 8}%
\special{pa 4800 2000}%
\special{pa 5200 2000}%
\special{fp}%
%
\special{pn 8}%
\special{ar 4000 2400 200 200  6.2831853 6.2831853}%
\special{ar 4000 2400 200 200  0.0000000 3.1415927}%
%
\special{pn 8}%
\special{pa 4200 2416}%
\special{pa 4200 2400}%
\special{fp}%
\special{sh 1}%
\special{pa 4200 2400}%
\special{pa 4176 2466}%
\special{pa 4196 2454}%
\special{pa 4216 2468}%
\special{pa 4200 2400}%
\special{fp}%
%
\special{pn 8}%
\special{ar 3200 2400 200 200  6.2831853 6.2831853}%
\special{ar 3200 2400 200 200  0.0000000 3.1415927}%
%
\special{pn 8}%
\special{pa 3400 2416}%
\special{pa 3400 2400}%
\special{fp}%
\special{sh 1}%
\special{pa 3400 2400}%
\special{pa 3376 2466}%
\special{pa 3396 2454}%
\special{pa 3416 2468}%
\special{pa 3400 2400}%
\special{fp}%
\put(35.0000,-26.0500){\makebox(0,0){$+1$}}%
\put(43.0000,-26.0500){\makebox(0,0){$-1$}}%
%
\special{pn 8}%
\special{pa 2400 2200}%
\special{pa 2400 2400}%
\special{fp}%
\special{pa 2400 2400}%
\special{pa 5200 2400}%
\special{fp}%
\special{pa 5200 2400}%
\special{pa 5200 2200}%
\special{fp}%
\special{pa 4800 2200}%
\special{pa 4800 2400}%
\special{fp}%
\special{pa 4000 2200}%
\special{pa 4000 2400}%
\special{fp}%
\special{pa 3200 2200}%
\special{pa 3200 2400}%
\special{fp}%
\put(28.0000,-23.0000){\makebox(0,0){$1$}}%
\put(36.0000,-23.0000){\makebox(0,0){$2$}}%
\put(44.0000,-23.0000){\makebox(0,0){$1$}}%
\put(50.0000,-23.0000){\makebox(0,0){(+,iii)}}%
\end{picture}%

%% file: ee=1eq.tex
\unitlength 0.1in
\begin{picture}( 40.5500, 10.0000)( 11.4500,-32.0000)
\put(28.0000,-24.0000){\makebox(0,0){$-\epsilon_2q_1j-1$}}%
\put(36.0000,-24.0000){\makebox(0,0){$-\epsilon_2q_1j$}}%
\put(44.0000,-24.0000){\makebox(0,0){$-\epsilon_2q_1j+1$}}%
\put(28.0500,-27.0000){\makebox(0,0){$1$}}%
\put(36.0500,-27.0000){\makebox(0,0){$0$}}%
\put(44.0500,-27.0000){\makebox(0,0){$-1$}}%
%
\special{pn 8}%
\special{pa 2400 2600}%
\special{pa 2400 2200}%
\special{fp}%
\special{pa 2400 2200}%
\special{pa 4800 2200}%
\special{fp}%
\special{pa 4800 2600}%
\special{pa 4800 2600}%
\special{fp}%
\special{pa 4800 2200}%
\special{pa 4800 2600}%
\special{fp}%
\special{pa 4000 2200}%
\special{pa 4000 2600}%
\special{fp}%
\special{pa 3200 2200}%
\special{pa 3200 2600}%
\special{fp}%
%
\special{pn 8}%
\special{pa 2400 2600}%
\special{pa 2400 2800}%
\special{fp}%
\special{pa 2400 2800}%
\special{pa 4800 2800}%
\special{fp}%
\special{pa 4800 2800}%
\special{pa 4800 2600}%
\special{fp}%
\special{pa 4800 2600}%
\special{pa 2400 2600}%
\special{fp}%
\special{pa 3200 2600}%
\special{pa 3200 2800}%
\special{fp}%
\special{pa 4000 2600}%
\special{pa 4000 2800}%
\special{fp}%
%
\special{pn 8}%
\special{pa 2400 2600}%
\special{pa 1600 2600}%
\special{fp}%
\special{pa 2400 2800}%
\special{pa 1600 2800}%
\special{fp}%
\special{pa 1600 2800}%
\special{pa 1600 2200}%
\special{fp}%
\special{pa 1600 2200}%
\special{pa 2400 2200}%
\special{fp}%
\put(20.0000,-27.0000){\makebox(0,0){$0$}}%
\put(20.0000,-24.0000){\makebox(0,0){$-\epsilon_2q_1j-2$}}%
%
\special{pn 8}%
\special{pa 4800 2200}%
\special{pa 5200 2200}%
\special{fp}%
\special{pa 5200 2200}%
\special{pa 5200 2800}%
\special{fp}%
\special{pa 5200 2800}%
\special{pa 4800 2800}%
\special{fp}%
\special{pa 4800 2600}%
\special{pa 5200 2600}%
\special{fp}%
\put(50.0000,-24.0000){\makebox(0,0){$x$}}%
\put(50.0000,-27.0000){\makebox(0,0){(+,iv)}}%
%
\special{pn 4}%
\special{ar 4000 3000 200 200  6.2831853 6.2831853}%
\special{ar 4000 3000 200 200  0.0000000 3.1415927}%
%
\special{pn 4}%
\special{pa 4200 3016}%
\special{pa 4200 3000}%
\special{fp}%
\special{sh 1}%
\special{pa 4200 3000}%
\special{pa 4176 3066}%
\special{pa 4196 3054}%
\special{pa 4216 3068}%
\special{pa 4200 3000}%
\special{fp}%
%
\special{pn 4}%
\special{ar 3200 3000 200 200  6.2831853 6.2831853}%
\special{ar 3200 3000 200 200  0.0000000 3.1415927}%
%
\special{pn 4}%
\special{pa 3400 3016}%
\special{pa 3400 3000}%
\special{fp}%
\special{sh 1}%
\special{pa 3400 3000}%
\special{pa 3376 3066}%
\special{pa 3396 3054}%
\special{pa 3416 3068}%
\special{pa 3400 3000}%
\special{fp}%
%
\special{pn 4}%
\special{ar 2400 3000 200 200  6.2831853 6.2831853}%
\special{ar 2400 3000 200 200  0.0000000 3.1415927}%
%
\special{pn 4}%
\special{pa 2600 3016}%
\special{pa 2600 3000}%
\special{fp}%
\special{sh 1}%
\special{pa 2600 3000}%
\special{pa 2576 3066}%
\special{pa 2596 3054}%
\special{pa 2616 3068}%
\special{pa 2600 3000}%
\special{fp}%
\put(43.3000,-31.8500){\makebox(0,0){$-1$}}%
\put(35.3000,-31.8500){\makebox(0,0){$-1$}}%
\put(27.3000,-31.8500){\makebox(0,0){$+1$}}%
%
\special{pn 8}%
\special{pa 1600 2800}%
\special{pa 1600 3000}%
\special{fp}%
\special{pa 5200 3000}%
\special{pa 5200 3000}%
\special{fp}%
\special{pa 5200 2800}%
\special{pa 5200 3000}%
\special{fp}%
\special{pa 5200 3000}%
\special{pa 1600 3000}%
\special{fp}%
\special{pa 2400 2800}%
\special{pa 2400 3000}%
\special{fp}%
\special{pa 3200 2800}%
\special{pa 3200 3000}%
\special{fp}%
\special{pa 4000 2800}%
\special{pa 4000 3000}%
\special{fp}%
\special{pa 4800 2800}%
\special{pa 4800 3000}%
\special{fp}%
\put(28.0000,-29.0000){\makebox(0,0){$2$}}%
\put(36.0000,-29.0000){\makebox(0,0){$1$}}%
\put(44.0000,-29.0000){\makebox(0,0){$0$}}%
\put(20.0000,-29.0000){\makebox(0,0){$1$}}%
\put(50.0000,-29.0000){\makebox(0,0){(+,v)}}%
\end{picture}%

%% file: ee=1k1.tex
\unitlength 0.1in
\begin{picture}( 39.6500, 10.0000)( 12.3500,-32.0000)
\put(28.0000,-23.0000){\makebox(0,0){$-\epsilon_2q_1j$}}%
\put(36.0000,-23.0000){\makebox(0,0){$-\epsilon_2q_1j$}}%
\put(44.0000,-23.0000){\makebox(0,0){$-\epsilon_2q_1j$}}%
\put(28.0500,-27.0000){\makebox(0,0){$0$}}%
\put(36.0500,-27.0000){\makebox(0,0){$1$}}%
\put(44.0500,-27.0000){\makebox(0,0){$0$}}%
%
\special{pn 8}%
\special{pa 2400 2600}%
\special{pa 2400 2200}%
\special{fp}%
\special{pa 2400 2200}%
\special{pa 4800 2200}%
\special{fp}%
\special{pa 4800 2600}%
\special{pa 4800 2600}%
\special{fp}%
\special{pa 4800 2200}%
\special{pa 4800 2600}%
\special{fp}%
\special{pa 4000 2200}%
\special{pa 4000 2600}%
\special{fp}%
\special{pa 3200 2200}%
\special{pa 3200 2600}%
\special{fp}%
%
\special{pn 8}%
\special{pa 2400 2600}%
\special{pa 2400 2800}%
\special{fp}%
\special{pa 2400 2800}%
\special{pa 4800 2800}%
\special{fp}%
\special{pa 4800 2800}%
\special{pa 4800 2600}%
\special{fp}%
\special{pa 4800 2600}%
\special{pa 2400 2600}%
\special{fp}%
\special{pa 3200 2600}%
\special{pa 3200 2800}%
\special{fp}%
\special{pa 4000 2600}%
\special{pa 4000 2800}%
\special{fp}%
%
\special{pn 8}%
\special{pa 2400 2600}%
\special{pa 1600 2600}%
\special{fp}%
\special{pa 2400 2800}%
\special{pa 1600 2800}%
\special{fp}%
\special{pa 1600 2800}%
\special{pa 1600 2200}%
\special{fp}%
\special{pa 1600 2200}%
\special{pa 2400 2200}%
\special{fp}%
\put(20.0000,-23.0000){\makebox(0,0){$-\epsilon_2q_1j$}}%
\put(36.0000,-25.0000){\makebox(0,0){$+k_1+1$}}%
\put(44.0000,-25.0000){\makebox(0,0){$+k_1+2$}}%
\put(20.0000,-25.0000){\makebox(0,0){$+k_1-1$}}%
\put(28.0000,-25.0000){\makebox(0,0){$+k_1$}}%
\put(20.1000,-27.0000){\makebox(0,0){$-1$}}%
%
\special{pn 8}%
\special{pa 4800 2200}%
\special{pa 5200 2200}%
\special{fp}%
\special{pa 5200 2200}%
\special{pa 5200 2800}%
\special{fp}%
\special{pa 5200 2800}%
\special{pa 4800 2800}%
\special{fp}%
\special{pa 4800 2600}%
\special{pa 5200 2600}%
\special{fp}%
\put(50.0000,-24.0000){\makebox(0,0){$x$}}%
\put(50.0000,-27.0000){\makebox(0,0){(+,vi)}}%
%
\special{pn 8}%
\special{pa 1600 2800}%
\special{pa 1600 3000}%
\special{fp}%
\special{pa 1600 3000}%
\special{pa 5200 3000}%
\special{fp}%
\special{pa 5200 3000}%
\special{pa 5200 2800}%
\special{fp}%
\special{pa 4800 2800}%
\special{pa 4800 3000}%
\special{fp}%
\special{pa 4000 2800}%
\special{pa 4000 3000}%
\special{fp}%
\special{pa 3200 2800}%
\special{pa 3200 3000}%
\special{fp}%
\special{pa 2400 2800}%
\special{pa 2400 3000}%
\special{fp}%
\put(36.0000,-29.0000){\makebox(0,0){$2$}}%
%
\special{pn 8}%
\special{ar 2400 3000 200 200  6.2831853 6.2831853}%
\special{ar 2400 3000 200 200  0.0000000 3.1415927}%
%
\special{pn 8}%
\special{pa 2600 3016}%
\special{pa 2600 3000}%
\special{fp}%
\special{sh 1}%
\special{pa 2600 3000}%
\special{pa 2576 3066}%
\special{pa 2596 3054}%
\special{pa 2616 3068}%
\special{pa 2600 3000}%
\special{fp}%
%
\special{pn 8}%
\special{ar 3200 3000 200 200  6.2831853 6.2831853}%
\special{ar 3200 3000 200 200  0.0000000 3.1415927}%
%
\special{pn 8}%
\special{pa 3400 3016}%
\special{pa 3400 3000}%
\special{fp}%
\special{sh 1}%
\special{pa 3400 3000}%
\special{pa 3376 3066}%
\special{pa 3396 3054}%
\special{pa 3416 3068}%
\special{pa 3400 3000}%
\special{fp}%
%
\special{pn 8}%
\special{ar 4000 3000 200 200  6.2831853 6.2831853}%
\special{ar 4000 3000 200 200  0.0000000 3.1415927}%
%
\special{pn 8}%
\special{pa 4200 3016}%
\special{pa 4200 3000}%
\special{fp}%
\special{sh 1}%
\special{pa 4200 3000}%
\special{pa 4176 3066}%
\special{pa 4196 3054}%
\special{pa 4216 3068}%
\special{pa 4200 3000}%
\special{fp}%
\put(28.0000,-29.0000){\makebox(0,0){$1$}}%
\put(20.0000,-29.0000){\makebox(0,0){$0$}}%
\put(44.0000,-29.0000){\makebox(0,0){$1$}}%
\put(50.0000,-29.0000){\makebox(0,0){(+,vii)}}%
\put(26.8000,-31.8500){\makebox(0,0){$+1$}}%
\put(34.8000,-31.8500){\makebox(0,0){$+1$}}%
\put(42.8000,-31.8500){\makebox(0,0){$-1$}}%
\end{picture}%

%% file: ee=1notp.tex
\unitlength 0.1in
\begin{picture}( 39.7500, 10.0000)(  2.3000,-24.0000)
\put(18.0500,-16.0000){\makebox(0,0){$\lfloor\frac{p\ell}{|q_2|}\rfloor$}}%
\put(26.0500,-16.0000){\makebox(0,0){$\lfloor\frac{p\ell}{|q_2|}\rfloor+1$}}%
\put(34.0500,-16.0000){\makebox(0,0){$\lfloor\frac{p\ell}{|q_2|}\rfloor+2$}}%
\put(34.1000,-19.0000){\makebox(0,0){$1$}}%
\put(34.1000,-21.0000){\makebox(0,0){$0$}}%
\put(26.1000,-21.0000){\makebox(0,0){$-1$}}%
\put(26.1000,-19.0000){\makebox(0,0){$0$}}%
%
\special{pn 8}%
\special{pa 1406 1800}%
\special{pa 3806 1800}%
\special{fp}%
\special{pa 3806 1800}%
\special{pa 3806 2200}%
\special{fp}%
\special{pa 3806 2200}%
\special{pa 1406 2200}%
\special{fp}%
\special{pa 1406 2200}%
\special{pa 1406 1800}%
\special{fp}%
\special{pa 1406 2000}%
\special{pa 3806 2000}%
\special{fp}%
\special{pa 2206 1800}%
\special{pa 2206 2200}%
\special{fp}%
\special{pa 3006 1800}%
\special{pa 3006 2200}%
\special{fp}%
%
\special{pn 8}%
\special{pa 1406 1800}%
\special{pa 1406 1400}%
\special{fp}%
\special{pa 1406 1400}%
\special{pa 3806 1400}%
\special{fp}%
\special{pa 3806 1800}%
\special{pa 3806 1800}%
\special{fp}%
\special{pa 3806 1400}%
\special{pa 3806 1800}%
\special{fp}%
\special{pa 3006 1400}%
\special{pa 3006 1800}%
\special{fp}%
\special{pa 2206 1400}%
\special{pa 2206 1800}%
\special{fp}%
%
\special{pn 8}%
\special{pa 3806 1400}%
\special{pa 4206 1400}%
\special{fp}%
\special{pa 4206 1400}%
\special{pa 4206 2200}%
\special{fp}%
\special{pa 4206 2200}%
\special{pa 3806 2200}%
\special{fp}%
%
\special{pn 8}%
\special{pa 3806 1800}%
\special{pa 4206 1800}%
\special{fp}%
\put(40.0500,-19.0000){\makebox(0,0){(+,iix)}}%
\put(40.0500,-21.0000){\makebox(0,0){(+,ix)}}%
\put(40.0500,-16.0000){\makebox(0,0){$x$}}%
%
\special{pn 8}%
\special{pa 3806 2000}%
\special{pa 4206 2000}%
\special{fp}%
%
\special{pn 8}%
\special{ar 3010 2200 200 200  6.2831853 6.2831853}%
\special{ar 3010 2200 200 200  0.0000000 3.1415927}%
%
\special{pn 8}%
\special{pa 2812 2216}%
\special{pa 2810 2200}%
\special{fp}%
\special{sh 1}%
\special{pa 2810 2200}%
\special{pa 2794 2268}%
\special{pa 2814 2254}%
\special{pa 2834 2266}%
\special{pa 2810 2200}%
\special{fp}%
%
\special{pn 8}%
\special{ar 2210 2200 200 200  6.2831853 6.2831853}%
\special{ar 2210 2200 200 200  0.0000000 3.1415927}%
%
\special{pn 8}%
\special{pa 2012 2216}%
\special{pa 2010 2200}%
\special{fp}%
\special{sh 1}%
\special{pa 2010 2200}%
\special{pa 1994 2268}%
\special{pa 2014 2254}%
\special{pa 2034 2266}%
\special{pa 2010 2200}%
\special{fp}%
\put(25.1000,-24.0500){\makebox(0,0){$+1$}}%
\put(33.1000,-24.0500){\makebox(0,0){$-1$}}%
\put(18.1000,-21.0000){\makebox(0,0){$0$}}%
\put(18.1000,-19.0000){\makebox(0,0){$1$}}%
\end{picture}%

%% file: ee=1eqp.tex
\unitlength 0.1in
\begin{picture}( 40.5500, 10.0000)( 11.4500,-32.0000)
\put(28.0000,-24.0000){\makebox(0,0){$-\epsilon_2q_1j$}}%
\put(36.0000,-24.0000){\makebox(0,0){$-\epsilon_2q_1j+1$}}%
\put(44.0000,-24.0000){\makebox(0,0){$-\epsilon_2q_1j+2$}}%
\put(28.0500,-27.0000){\makebox(0,0){$-1$}}%
\put(36.0500,-27.0000){\makebox(0,0){$0$}}%
\put(44.0500,-27.0000){\makebox(0,0){$1$}}%
%
\special{pn 8}%
\special{pa 2400 2600}%
\special{pa 2400 2200}%
\special{fp}%
\special{pa 2400 2200}%
\special{pa 4800 2200}%
\special{fp}%
\special{pa 4800 2600}%
\special{pa 4800 2600}%
\special{fp}%
\special{pa 4800 2200}%
\special{pa 4800 2600}%
\special{fp}%
\special{pa 4000 2200}%
\special{pa 4000 2600}%
\special{fp}%
\special{pa 3200 2200}%
\special{pa 3200 2600}%
\special{fp}%
%
\special{pn 8}%
\special{pa 2400 2600}%
\special{pa 2400 2800}%
\special{fp}%
\special{pa 2400 2800}%
\special{pa 4800 2800}%
\special{fp}%
\special{pa 4800 2800}%
\special{pa 4800 2600}%
\special{fp}%
\special{pa 4800 2600}%
\special{pa 2400 2600}%
\special{fp}%
\special{pa 3200 2600}%
\special{pa 3200 2800}%
\special{fp}%
\special{pa 4000 2600}%
\special{pa 4000 2800}%
\special{fp}%
%
\special{pn 8}%
\special{pa 2400 2600}%
\special{pa 1600 2600}%
\special{fp}%
\special{pa 2400 2800}%
\special{pa 1600 2800}%
\special{fp}%
\special{pa 1600 2800}%
\special{pa 1600 2200}%
\special{fp}%
\special{pa 1600 2200}%
\special{pa 2400 2200}%
\special{fp}%
\put(20.0000,-27.0000){\makebox(0,0){$0$}}%
\put(20.0000,-24.0000){\makebox(0,0){$-\epsilon_2q_1j-1$}}%
%
\special{pn 8}%
\special{pa 4800 2200}%
\special{pa 5200 2200}%
\special{fp}%
\special{pa 5200 2200}%
\special{pa 5200 2800}%
\special{fp}%
\special{pa 5200 2800}%
\special{pa 4800 2800}%
\special{fp}%
\special{pa 4800 2600}%
\special{pa 5200 2600}%
\special{fp}%
\put(50.0000,-24.0000){\makebox(0,0){$x$}}%
\put(50.0000,-27.0000){\makebox(0,0){(+,x)}}%
%
\special{pn 4}%
\special{ar 4000 3000 200 200  6.2831853 6.2831853}%
\special{ar 4000 3000 200 200  0.0000000 3.1415927}%
%
\special{pn 4}%
\special{pa 3802 3016}%
\special{pa 3800 3000}%
\special{fp}%
\special{sh 1}%
\special{pa 3800 3000}%
\special{pa 3784 3068}%
\special{pa 3804 3054}%
\special{pa 3824 3066}%
\special{pa 3800 3000}%
\special{fp}%
%
\special{pn 4}%
\special{ar 3200 3000 200 200  6.2831853 6.2831853}%
\special{ar 3200 3000 200 200  0.0000000 3.1415927}%
%
\special{pn 4}%
\special{pa 3002 3016}%
\special{pa 3000 3000}%
\special{fp}%
\special{sh 1}%
\special{pa 3000 3000}%
\special{pa 2984 3068}%
\special{pa 3004 3054}%
\special{pa 3024 3066}%
\special{pa 3000 3000}%
\special{fp}%
%
\special{pn 4}%
\special{ar 2400 3000 200 200  6.2831853 6.2831853}%
\special{ar 2400 3000 200 200  0.0000000 3.1415927}%
%
\special{pn 4}%
\special{pa 2202 3016}%
\special{pa 2200 3000}%
\special{fp}%
\special{sh 1}%
\special{pa 2200 3000}%
\special{pa 2184 3068}%
\special{pa 2204 3054}%
\special{pa 2224 3066}%
\special{pa 2200 3000}%
\special{fp}%
\put(43.3000,-31.8500){\makebox(0,0){$-1$}}%
\put(35.3000,-31.8500){\makebox(0,0){$-1$}}%
\put(27.3000,-31.8500){\makebox(0,0){$+1$}}%
%
\special{pn 8}%
\special{pa 1600 2800}%
\special{pa 1600 3000}%
\special{fp}%
\special{pa 5200 3000}%
\special{pa 5200 3000}%
\special{fp}%
\special{pa 5200 2800}%
\special{pa 5200 3000}%
\special{fp}%
\special{pa 5200 3000}%
\special{pa 1600 3000}%
\special{fp}%
\special{pa 2400 2800}%
\special{pa 2400 3000}%
\special{fp}%
\special{pa 3200 2800}%
\special{pa 3200 3000}%
\special{fp}%
\special{pa 4000 2800}%
\special{pa 4000 3000}%
\special{fp}%
\special{pa 4800 2800}%
\special{pa 4800 3000}%
\special{fp}%
\put(28.0000,-29.0000){\makebox(0,0){$0$}}%
\put(36.0000,-29.0000){\makebox(0,0){$1$}}%
\put(44.0000,-29.0000){\makebox(0,0){$2$}}%
\put(20.0000,-29.0000){\makebox(0,0){$1$}}%
\put(50.0000,-29.0000){\makebox(0,0){(+,xi)}}%
\end{picture}%

%% file: ee=1k1p.tex
\unitlength 0.1in
\begin{picture}( 39.6500, 10.0000)( 12.3500,-32.0000)
\put(28.0000,-23.0000){\makebox(0,0){$-\epsilon_2q_1j$}}%
\put(36.0000,-23.0000){\makebox(0,0){$-\epsilon_2q_1j$}}%
\put(44.0000,-23.0000){\makebox(0,0){$-\epsilon_2q_1j$}}%
\put(28.0500,-27.0000){\makebox(0,0){$0$}}%
\put(36.0500,-27.0000){\makebox(0,0){$-1$}}%
\put(44.0500,-27.0000){\makebox(0,0){$0$}}%
%
\special{pn 8}%
\special{pa 2400 2600}%
\special{pa 2400 2200}%
\special{fp}%
\special{pa 2400 2200}%
\special{pa 4800 2200}%
\special{fp}%
\special{pa 4800 2600}%
\special{pa 4800 2600}%
\special{fp}%
\special{pa 4800 2200}%
\special{pa 4800 2600}%
\special{fp}%
\special{pa 4000 2200}%
\special{pa 4000 2600}%
\special{fp}%
\special{pa 3200 2200}%
\special{pa 3200 2600}%
\special{fp}%
%
\special{pn 8}%
\special{pa 2400 2600}%
\special{pa 2400 2800}%
\special{fp}%
\special{pa 2400 2800}%
\special{pa 4800 2800}%
\special{fp}%
\special{pa 4800 2800}%
\special{pa 4800 2600}%
\special{fp}%
\special{pa 4800 2600}%
\special{pa 2400 2600}%
\special{fp}%
\special{pa 3200 2600}%
\special{pa 3200 2800}%
\special{fp}%
\special{pa 4000 2600}%
\special{pa 4000 2800}%
\special{fp}%
%
\special{pn 8}%
\special{pa 2400 2600}%
\special{pa 1600 2600}%
\special{fp}%
\special{pa 2400 2800}%
\special{pa 1600 2800}%
\special{fp}%
\special{pa 1600 2800}%
\special{pa 1600 2200}%
\special{fp}%
\special{pa 1600 2200}%
\special{pa 2400 2200}%
\special{fp}%
\put(20.0000,-23.0000){\makebox(0,0){$-\epsilon_2q_1j$}}%
\put(36.0000,-25.0000){\makebox(0,0){$+k_1+2$}}%
\put(44.0000,-25.0000){\makebox(0,0){$+k_1+3$}}%
\put(20.0000,-25.0000){\makebox(0,0){$+k_1$}}%
\put(28.0000,-25.0000){\makebox(0,0){$+k_1+1$}}%
\put(20.1000,-27.0000){\makebox(0,0){$1$}}%
%
\special{pn 8}%
\special{pa 4800 2200}%
\special{pa 5200 2200}%
\special{fp}%
\special{pa 5200 2200}%
\special{pa 5200 2800}%
\special{fp}%
\special{pa 5200 2800}%
\special{pa 4800 2800}%
\special{fp}%
\special{pa 4800 2600}%
\special{pa 5200 2600}%
\special{fp}%
\put(50.0000,-24.0000){\makebox(0,0){$x$}}%
\put(50.0000,-27.0000){\makebox(0,0){(+,xii)}}%
%
\special{pn 8}%
\special{pa 1600 2800}%
\special{pa 1600 3000}%
\special{fp}%
\special{pa 1600 3000}%
\special{pa 5200 3000}%
\special{fp}%
\special{pa 5200 3000}%
\special{pa 5200 2800}%
\special{fp}%
\special{pa 4800 2800}%
\special{pa 4800 3000}%
\special{fp}%
\special{pa 4000 2800}%
\special{pa 4000 3000}%
\special{fp}%
\special{pa 3200 2800}%
\special{pa 3200 3000}%
\special{fp}%
\special{pa 2400 2800}%
\special{pa 2400 3000}%
\special{fp}%
\put(36.0000,-29.0000){\makebox(0,0){$0$}}%
%
\special{pn 8}%
\special{ar 2400 3000 200 200  6.2831853 6.2831853}%
\special{ar 2400 3000 200 200  0.0000000 3.1415927}%
%
\special{pn 8}%
\special{pa 2202 3016}%
\special{pa 2200 3000}%
\special{fp}%
\special{sh 1}%
\special{pa 2200 3000}%
\special{pa 2184 3068}%
\special{pa 2204 3054}%
\special{pa 2224 3066}%
\special{pa 2200 3000}%
\special{fp}%
%
\special{pn 8}%
\special{ar 3200 3000 200 200  6.2831853 6.2831853}%
\special{ar 3200 3000 200 200  0.0000000 3.1415927}%
%
\special{pn 8}%
\special{pa 3002 3016}%
\special{pa 3000 3000}%
\special{fp}%
\special{sh 1}%
\special{pa 3000 3000}%
\special{pa 2984 3068}%
\special{pa 3004 3054}%
\special{pa 3024 3066}%
\special{pa 3000 3000}%
\special{fp}%
%
\special{pn 8}%
\special{ar 4000 3000 200 200  6.2831853 6.2831853}%
\special{ar 4000 3000 200 200  0.0000000 3.1415927}%
%
\special{pn 8}%
\special{pa 3802 3016}%
\special{pa 3800 3000}%
\special{fp}%
\special{sh 1}%
\special{pa 3800 3000}%
\special{pa 3784 3068}%
\special{pa 3804 3054}%
\special{pa 3824 3066}%
\special{pa 3800 3000}%
\special{fp}%
\put(28.0000,-29.0000){\makebox(0,0){$1$}}%
\put(20.0000,-29.0000){\makebox(0,0){$2$}}%
\put(44.0000,-29.0000){\makebox(0,0){$1$}}%
\put(50.0000,-29.0000){\makebox(0,0){(+,xiii)}}%
\put(26.8000,-31.8500){\makebox(0,0){$+1$}}%
\put(34.8000,-31.8500){\makebox(0,0){$+1$}}%
\put(42.8000,-31.8500){\makebox(0,0){$-1$}}%
\end{picture}%

%% file: stair.tex
\unitlength 0.1in
\begin{picture}( 37.4000, 10.0500)(  9.2500,-21.2500)
\put(36.0000,-20.0500){\makebox(0,0){$+$}}%
\put(38.0000,-18.0500){\makebox(0,0){$+$}}%
\put(40.0000,-16.0500){\makebox(0,0){$+$}}%
\put(36.0000,-22.0500){\makebox(0,0){$-$}}%
\put(38.0000,-22.0500){\makebox(0,0){$-$}}%
\put(42.0000,-12.0500){\makebox(0,0){$-$}}%
\put(42.0000,-14.0500){\makebox(0,0){$+$}}%
\put(40.0000,-12.0500){\makebox(0,0){$-$}}%
\put(38.0000,-12.0500){\makebox(0,0){$-$}}%
\put(36.0000,-12.0500){\makebox(0,0){$-$}}%
\put(34.0000,-12.0500){\makebox(0,0){$\cdots$}}%
\put(34.0000,-22.0500){\makebox(0,0){$0$}}%
\put(32.0000,-22.0500){\makebox(0,0){$-$}}%
\put(30.0000,-22.0500){\makebox(0,0){$\cdots$}}%
\put(44.0000,-12.0500){\makebox(0,0){$0$}}%
\put(46.0000,-12.0500){\makebox(0,0){$-$}}%
\put(48.0000,-12.0500){\makebox(0,0){$\cdots$}}%
\put(42.0000,-22.1000){\makebox(0,0){$\cdots$}}%
\put(44.0000,-22.1000){\makebox(0,0){$-$}}%
\put(46.0000,-22.1000){\makebox(0,0){$0$}}%
\put(48.0000,-22.1000){\makebox(0,0){$-$}}%
\put(48.0000,-20.1000){\makebox(0,0){$+$}}%
\put(40.0000,-22.1000){\makebox(0,0){$\cdots$}}%
\put(50.0000,-22.0500){\makebox(0,0){$\cdots$}}%
\put(32.0000,-12.0500){\makebox(0,0){$-$}}%
\put(30.0000,-12.0500){\makebox(0,0){$0$}}%
\put(28.0000,-12.0500){\makebox(0,0){$-$}}%
\put(28.0000,-14.0500){\makebox(0,0){$+$}}%
\put(16.0000,-16.0000){\makebox(0,0){positive region}}%
\put(26.0000,-12.1000){\makebox(0,0){$\cdots$}}%
\put(26.0000,-16.1000){\makebox(0,0){$+$}}%
\put(24.0000,-18.1000){\makebox(0,0){$+$}}%
\put(22.0000,-20.1000){\makebox(0,0){$+$}}%
\put(22.0000,-22.1000){\makebox(0,0){$-$}}%
\put(24.0000,-22.1000){\makebox(0,0){$-$}}%
\put(26.0000,-22.1000){\makebox(0,0){$\cdots$}}%
\put(20.0000,-22.1000){\makebox(0,0){$0$}}%
\put(18.0000,-22.1000){\makebox(0,0){$-$}}%
\put(16.0000,-22.1000){\makebox(0,0){$\cdots$}}%
\put(49.6000,-18.9000){\makebox(0,0){$\odots$}}%
%
\special{pn 4}%
\special{pa 2540 1300}%
\special{pa 1750 2090}%
\special{fp}%
\special{pa 2480 1300}%
\special{pa 1690 2090}%
\special{fp}%
\special{pa 2420 1300}%
\special{pa 1630 2090}%
\special{fp}%
\special{pa 2360 1300}%
\special{pa 1570 2090}%
\special{fp}%
\special{pa 2300 1300}%
\special{pa 1510 2090}%
\special{fp}%
\special{pa 2240 1300}%
\special{pa 1450 2090}%
\special{fp}%
\special{pa 2180 1300}%
\special{pa 1390 2090}%
\special{fp}%
\special{pa 2120 1300}%
\special{pa 1330 2090}%
\special{fp}%
\special{pa 2060 1300}%
\special{pa 1270 2090}%
\special{fp}%
\special{pa 2000 1300}%
\special{pa 1210 2090}%
\special{fp}%
\special{pa 1940 1300}%
\special{pa 1150 2090}%
\special{fp}%
\special{pa 1880 1300}%
\special{pa 1090 2090}%
\special{fp}%
\special{pa 1820 1300}%
\special{pa 1030 2090}%
\special{fp}%
\special{pa 1760 1300}%
\special{pa 970 2090}%
\special{fp}%
\special{pa 1700 1300}%
\special{pa 950 2050}%
\special{fp}%
\special{pa 1640 1300}%
\special{pa 950 1990}%
\special{fp}%
\special{pa 1580 1300}%
\special{pa 950 1930}%
\special{fp}%
\special{pa 1520 1300}%
\special{pa 950 1870}%
\special{fp}%
\special{pa 1460 1300}%
\special{pa 950 1810}%
\special{fp}%
\special{pa 1400 1300}%
\special{pa 950 1750}%
\special{fp}%
\special{pa 1340 1300}%
\special{pa 950 1690}%
\special{fp}%
\special{pa 1280 1300}%
\special{pa 950 1630}%
\special{fp}%
\special{pa 1220 1300}%
\special{pa 950 1570}%
\special{fp}%
\special{pa 1160 1300}%
\special{pa 950 1510}%
\special{fp}%
\special{pa 1100 1300}%
\special{pa 950 1450}%
\special{fp}%
\special{pa 1040 1300}%
\special{pa 950 1390}%
\special{fp}%
\special{pa 980 1300}%
\special{pa 950 1330}%
\special{fp}%
\special{pa 2600 1300}%
\special{pa 1810 2090}%
\special{fp}%
\special{pa 2660 1300}%
\special{pa 1870 2090}%
\special{fp}%
\special{pa 2720 1300}%
\special{pa 1930 2090}%
\special{fp}%
%
\special{pn 4}%
\special{pa 2780 1300}%
\special{pa 1990 2090}%
\special{fp}%
\special{pa 2840 1300}%
\special{pa 2050 2090}%
\special{fp}%
\special{pa 2900 1300}%
\special{pa 2110 2090}%
\special{fp}%
\special{pa 2960 1300}%
\special{pa 2170 2090}%
\special{fp}%
\special{pa 3020 1300}%
\special{pa 2230 2090}%
\special{fp}%
\special{pa 3080 1300}%
\special{pa 2290 2090}%
\special{fp}%
\special{pa 3140 1300}%
\special{pa 2350 2090}%
\special{fp}%
\special{pa 3200 1300}%
\special{pa 2410 2090}%
\special{fp}%
\special{pa 3260 1300}%
\special{pa 2470 2090}%
\special{fp}%
\special{pa 3320 1300}%
\special{pa 2530 2090}%
\special{fp}%
\special{pa 3380 1300}%
\special{pa 2590 2090}%
\special{fp}%
\special{pa 3440 1300}%
\special{pa 2650 2090}%
\special{fp}%
\special{pa 3500 1300}%
\special{pa 2710 2090}%
\special{fp}%
\special{pa 3560 1300}%
\special{pa 2770 2090}%
\special{fp}%
\special{pa 3620 1300}%
\special{pa 2830 2090}%
\special{fp}%
\special{pa 3680 1300}%
\special{pa 2890 2090}%
\special{fp}%
\special{pa 3740 1300}%
\special{pa 2950 2090}%
\special{fp}%
\special{pa 3800 1300}%
\special{pa 3010 2090}%
\special{fp}%
\special{pa 3860 1300}%
\special{pa 3070 2090}%
\special{fp}%
\special{pa 3920 1300}%
\special{pa 3130 2090}%
\special{fp}%
\special{pa 3980 1300}%
\special{pa 3190 2090}%
\special{fp}%
\special{pa 4040 1300}%
\special{pa 3250 2090}%
\special{fp}%
\special{pa 4100 1300}%
\special{pa 3310 2090}%
\special{fp}%
\special{pa 4160 1300}%
\special{pa 3370 2090}%
\special{fp}%
\special{pa 4220 1300}%
\special{pa 3430 2090}%
\special{fp}%
\special{pa 4280 1300}%
\special{pa 3490 2090}%
\special{fp}%
\special{pa 4340 1300}%
\special{pa 3550 2090}%
\special{fp}%
\special{pa 4400 1300}%
\special{pa 3610 2090}%
\special{fp}%
\special{pa 4460 1300}%
\special{pa 3670 2090}%
\special{fp}%
\special{pa 4520 1300}%
\special{pa 3730 2090}%
\special{fp}%
%
\special{pn 4}%
\special{pa 4580 1300}%
\special{pa 3790 2090}%
\special{fp}%
\special{pa 4640 1300}%
\special{pa 3850 2090}%
\special{fp}%
\special{pa 4700 1300}%
\special{pa 3910 2090}%
\special{fp}%
\special{pa 4760 1300}%
\special{pa 3970 2090}%
\special{fp}%
\special{pa 4820 1300}%
\special{pa 4030 2090}%
\special{fp}%
\special{pa 4880 1300}%
\special{pa 4090 2090}%
\special{fp}%
\special{pa 4940 1300}%
\special{pa 4150 2090}%
\special{fp}%
\special{pa 5000 1300}%
\special{pa 4210 2090}%
\special{fp}%
\special{pa 5060 1300}%
\special{pa 4270 2090}%
\special{fp}%
\special{pa 5120 1300}%
\special{pa 4330 2090}%
\special{fp}%
\special{pa 5180 1300}%
\special{pa 4390 2090}%
\special{fp}%
\special{pa 5240 1300}%
\special{pa 4450 2090}%
\special{fp}%
\special{pa 5300 1300}%
\special{pa 4510 2090}%
\special{fp}%
\special{pa 5360 1300}%
\special{pa 4570 2090}%
\special{fp}%
\special{pa 5420 1300}%
\special{pa 4630 2090}%
\special{fp}%
\special{pa 5480 1300}%
\special{pa 4690 2090}%
\special{fp}%
\special{pa 5520 1320}%
\special{pa 4750 2090}%
\special{fp}%
\special{pa 5520 1380}%
\special{pa 4810 2090}%
\special{fp}%
\special{pa 5520 1440}%
\special{pa 4870 2090}%
\special{fp}%
\special{pa 5520 1500}%
\special{pa 4930 2090}%
\special{fp}%
\special{pa 5520 1560}%
\special{pa 4990 2090}%
\special{fp}%
\special{pa 5520 1620}%
\special{pa 5050 2090}%
\special{fp}%
\special{pa 5520 1680}%
\special{pa 5110 2090}%
\special{fp}%
\special{pa 5520 1740}%
\special{pa 5170 2090}%
\special{fp}%
\special{pa 5520 1800}%
\special{pa 5230 2090}%
\special{fp}%
\special{pa 5520 1860}%
\special{pa 5290 2090}%
\special{fp}%
\special{pa 5520 1920}%
\special{pa 5350 2090}%
\special{fp}%
\special{pa 5520 1980}%
\special{pa 5410 2090}%
\special{fp}%
\special{pa 5520 2040}%
\special{pa 5470 2090}%
\special{fp}%
\end{picture}%

%% file: tableee=1.tex
\unitlength 0.1in
\begin{picture}( 49.3500, 28.3500)(  6.0000,-34.3500)
\put(15.3500,-31.1700){\makebox(0,0){$0$}}%
\put(21.3500,-31.1700){\makebox(0,0){$-$}}%
\put(18.3500,-31.1700){\makebox(0,0){$-$}}%
\put(18.3500,-35.1700){\makebox(0,0){$1$}}%
\put(21.3500,-35.1700){\makebox(0,0){$2$}}%
\put(7.3500,-29.1700){\makebox(0,0){$1$}}%
\put(7.3500,-27.1700){\makebox(0,0){$2$}}%
\put(7.3500,-25.1700){\makebox(0,0){$3$}}%
\put(12.3500,-31.1700){\makebox(0,0){$+$}}%
\put(15.3500,-35.1700){\makebox(0,0){$0$}}%
\put(12.3500,-35.1700){\makebox(0,0){$-1$}}%
\put(7.3500,-31.1700){\makebox(0,0){$0$}}%
%
\special{pn 8}%
\special{pa 936 3318}%
\special{pa 5326 3318}%
\special{fp}%
\special{sh 1}%
\special{pa 5326 3318}%
\special{pa 5258 3298}%
\special{pa 5272 3318}%
\special{pa 5258 3338}%
\special{pa 5326 3318}%
\special{fp}%
\special{pa 936 3318}%
\special{pa 936 600}%
\special{fp}%
\special{sh 1}%
\special{pa 936 600}%
\special{pa 916 668}%
\special{pa 936 654}%
\special{pa 956 668}%
\special{pa 936 600}%
\special{fp}%
\put(7.3500,-7.1700){\makebox(0,0){$j$}}%
\put(53.3500,-35.1700){\makebox(0,0){$i$}}%
\put(33.3500,-35.1700){\makebox(0,0){$\lfloor\frac{p}{|q_2|}\rfloor$}}%
\put(24.3500,-31.1700){\makebox(0,0){$\cdots$}}%
\put(33.3500,-31.1700){\makebox(0,0){$0$}}%
\put(30.3500,-31.1700){\makebox(0,0){$-$}}%
\put(36.3500,-31.1700){\makebox(0,0){$-$}}%
\put(12.3500,-29.1700){\makebox(0,0){$-$}}%
\put(18.3500,-29.1700){\makebox(0,0){$+$}}%
\put(15.3500,-29.1700){\makebox(0,0){$0$}}%
\put(21.3500,-29.1700){\makebox(0,0){$0$}}%
\put(24.3500,-29.1700){\makebox(0,0){$\cdots$}}%
\put(21.3500,-27.1700){\makebox(0,0){$+$}}%
\put(18.3500,-27.1700){\makebox(0,0){$0$}}%
\put(15.3500,-27.1700){\makebox(0,0){$0$}}%
\put(12.3500,-27.1700){\makebox(0,0){$0$}}%
\put(39.3500,-31.1700){\makebox(0,0){$\cdots$}}%
\put(24.3500,-27.1700){\makebox(0,0){$0$}}%
\put(24.3500,-25.1700){\makebox(0,0){\odots}}%
\put(27.3500,-23.1700){\makebox(0,0){\odots}}%
\put(30.3500,-21.1400){\makebox(0,0){$+$}}%
\put(30.3500,-19.1400){\makebox(0,0){$-$}}%
\put(27.3500,-19.1400){\makebox(0,0){$-$}}%
\put(21.3500,-19.1400){\makebox(0,0){$-$}}%
\put(24.3500,-19.1400){\makebox(0,0){$-$}}%
\put(41.5500,-13.2700){\makebox(0,0){$0$}}%
\put(50.5500,-11.2700){\makebox(0,0){$\cdots$}}%
\put(47.5500,-11.2700){\makebox(0,0){$0$}}%
\put(32.5500,-11.2700){\makebox(0,0){$-$}}%
\put(35.5500,-11.2700){\makebox(0,0){$\cdots$}}%
\put(32.5500,-13.2700){\makebox(0,0){$\cdots$}}%
\put(35.5500,-13.2700){\makebox(0,0){$0$}}%
\put(50.5500,-13.2700){\makebox(0,0){$\cdots$}}%
\put(47.5500,-13.2700){\makebox(0,0){$-$}}%
\put(41.5500,-11.2700){\makebox(0,0){$0$}}%
\put(44.5500,-11.2700){\makebox(0,0){$+$}}%
\put(38.5500,-11.2700){\makebox(0,0){$-$}}%
\put(38.5500,-13.2700){\makebox(0,0){$+$}}%
\put(44.5500,-13.2700){\makebox(0,0){$-$}}%
\put(33.3500,-19.1400){\makebox(0,0){$0$}}%
\put(33.3500,-21.1400){\makebox(0,0){$0$}}%
\put(36.3500,-21.1400){\makebox(0,0){$\cdots$}}%
\put(36.3500,-19.1400){\makebox(0,0){$\cdots$}}%
\put(27.3500,-31.1700){\makebox(0,0){$\cdots$}}%
\put(42.3500,-31.1700){\makebox(0,0){$-$}}%
\put(45.3500,-31.1700){\makebox(0,0){$0$}}%
\put(44.4500,-17.2700){\makebox(0,0){$\alpha(\ast)$}}%
\put(43.4500,-7.2700){\makebox(0,0){$\beta(\ast)$}}%
%
\special{pn 8}%
\special{pa 4346 828}%
\special{pa 4146 1028}%
\special{fp}%
%
\special{pn 8}%
\special{pa 4446 1628}%
\special{pa 4146 1428}%
\special{fp}%
%
\special{pn 4}%
\special{pa 1396 1218}%
\special{pa 1196 1018}%
\special{fp}%
\special{pa 1336 1218}%
\special{pa 1136 1018}%
\special{fp}%
\special{pa 1276 1218}%
\special{pa 1076 1018}%
\special{fp}%
\special{pa 1216 1218}%
\special{pa 1016 1018}%
\special{fp}%
\special{pa 1156 1218}%
\special{pa 1006 1068}%
\special{fp}%
\special{pa 1096 1218}%
\special{pa 1006 1128}%
\special{fp}%
\special{pa 1036 1218}%
\special{pa 1006 1188}%
\special{fp}%
\special{pa 1456 1218}%
\special{pa 1256 1018}%
\special{fp}%
\special{pa 1516 1218}%
\special{pa 1316 1018}%
\special{fp}%
\special{pa 1576 1218}%
\special{pa 1376 1018}%
\special{fp}%
\special{pa 1636 1218}%
\special{pa 1436 1018}%
\special{fp}%
\special{pa 1696 1218}%
\special{pa 1496 1018}%
\special{fp}%
\special{pa 1756 1218}%
\special{pa 1556 1018}%
\special{fp}%
\special{pa 1816 1218}%
\special{pa 1616 1018}%
\special{fp}%
\special{pa 1876 1218}%
\special{pa 1676 1018}%
\special{fp}%
\special{pa 1936 1218}%
\special{pa 1736 1018}%
\special{fp}%
\special{pa 1996 1218}%
\special{pa 1796 1018}%
\special{fp}%
\special{pa 2056 1218}%
\special{pa 1856 1018}%
\special{fp}%
\special{pa 2116 1218}%
\special{pa 1916 1018}%
\special{fp}%
\special{pa 2176 1218}%
\special{pa 1976 1018}%
\special{fp}%
\special{pa 2236 1218}%
\special{pa 2036 1018}%
\special{fp}%
\special{pa 2296 1218}%
\special{pa 2096 1018}%
\special{fp}%
\special{pa 2356 1218}%
\special{pa 2156 1018}%
\special{fp}%
\special{pa 2416 1218}%
\special{pa 2216 1018}%
\special{fp}%
\special{pa 2476 1218}%
\special{pa 2276 1018}%
\special{fp}%
\special{pa 2536 1218}%
\special{pa 2336 1018}%
\special{fp}%
\special{pa 2596 1218}%
\special{pa 2396 1018}%
\special{fp}%
\special{pa 2656 1218}%
\special{pa 2456 1018}%
\special{fp}%
\special{pa 2716 1218}%
\special{pa 2516 1018}%
\special{fp}%
\special{pa 2776 1218}%
\special{pa 2576 1018}%
\special{fp}%
%
\special{pn 4}%
\special{pa 2836 1218}%
\special{pa 2636 1018}%
\special{fp}%
\special{pa 2896 1218}%
\special{pa 2696 1018}%
\special{fp}%
\special{pa 2956 1218}%
\special{pa 2756 1018}%
\special{fp}%
\special{pa 3016 1218}%
\special{pa 2816 1018}%
\special{fp}%
\special{pa 3076 1218}%
\special{pa 2876 1018}%
\special{fp}%
\special{pa 3136 1218}%
\special{pa 2936 1018}%
\special{fp}%
\special{pa 3196 1218}%
\special{pa 2996 1018}%
\special{fp}%
\special{pa 3256 1218}%
\special{pa 3056 1018}%
\special{fp}%
\special{pa 3316 1218}%
\special{pa 3116 1018}%
\special{fp}%
\special{pa 3376 1218}%
\special{pa 3176 1018}%
\special{fp}%
\special{pa 3436 1218}%
\special{pa 3236 1018}%
\special{fp}%
\special{pa 3496 1218}%
\special{pa 3296 1018}%
\special{fp}%
\special{pa 3556 1218}%
\special{pa 3356 1018}%
\special{fp}%
\special{pa 3616 1218}%
\special{pa 3416 1018}%
\special{fp}%
\special{pa 3676 1218}%
\special{pa 3476 1018}%
\special{fp}%
\special{pa 3736 1218}%
\special{pa 3536 1018}%
\special{fp}%
\special{pa 3796 1218}%
\special{pa 3596 1018}%
\special{fp}%
\special{pa 3856 1218}%
\special{pa 3656 1018}%
\special{fp}%
\special{pa 3916 1218}%
\special{pa 3716 1018}%
\special{fp}%
\special{pa 3976 1218}%
\special{pa 3776 1018}%
\special{fp}%
\special{pa 4006 1188}%
\special{pa 3836 1018}%
\special{fp}%
\special{pa 4006 1128}%
\special{pa 3896 1018}%
\special{fp}%
\special{pa 4006 1068}%
\special{pa 3956 1018}%
\special{fp}%
%
\special{pn 4}%
\special{pa 5076 3218}%
\special{pa 4876 3018}%
\special{fp}%
\special{pa 5016 3218}%
\special{pa 4816 3018}%
\special{fp}%
\special{pa 4956 3218}%
\special{pa 4756 3018}%
\special{fp}%
\special{pa 4896 3218}%
\special{pa 4696 3018}%
\special{fp}%
\special{pa 4836 3218}%
\special{pa 4636 3018}%
\special{fp}%
\special{pa 4776 3218}%
\special{pa 4576 3018}%
\special{fp}%
\special{pa 4716 3218}%
\special{pa 4516 3018}%
\special{fp}%
\special{pa 4656 3218}%
\special{pa 4456 3018}%
\special{fp}%
\special{pa 4596 3218}%
\special{pa 4396 3018}%
\special{fp}%
\special{pa 4536 3218}%
\special{pa 4336 3018}%
\special{fp}%
\special{pa 4476 3218}%
\special{pa 4276 3018}%
\special{fp}%
\special{pa 4416 3218}%
\special{pa 4216 3018}%
\special{fp}%
\special{pa 4356 3218}%
\special{pa 4156 3018}%
\special{fp}%
\special{pa 4296 3218}%
\special{pa 4096 3018}%
\special{fp}%
\special{pa 4236 3218}%
\special{pa 4036 3018}%
\special{fp}%
\special{pa 4176 3218}%
\special{pa 3976 3018}%
\special{fp}%
\special{pa 4116 3218}%
\special{pa 3916 3018}%
\special{fp}%
\special{pa 4056 3218}%
\special{pa 3856 3018}%
\special{fp}%
\special{pa 3996 3218}%
\special{pa 3796 3018}%
\special{fp}%
\special{pa 3936 3218}%
\special{pa 3736 3018}%
\special{fp}%
\special{pa 3876 3218}%
\special{pa 3676 3018}%
\special{fp}%
\special{pa 3816 3218}%
\special{pa 3616 3018}%
\special{fp}%
\special{pa 3756 3218}%
\special{pa 3556 3018}%
\special{fp}%
\special{pa 3696 3218}%
\special{pa 3496 3018}%
\special{fp}%
\special{pa 3636 3218}%
\special{pa 3436 3018}%
\special{fp}%
\special{pa 3576 3218}%
\special{pa 3376 3018}%
\special{fp}%
\special{pa 3516 3218}%
\special{pa 3316 3018}%
\special{fp}%
\special{pa 3456 3218}%
\special{pa 3256 3018}%
\special{fp}%
\special{pa 3396 3218}%
\special{pa 3196 3018}%
\special{fp}%
\special{pa 3336 3218}%
\special{pa 3136 3018}%
\special{fp}%
%
\special{pn 4}%
\special{pa 3276 3218}%
\special{pa 3076 3018}%
\special{fp}%
\special{pa 3216 3218}%
\special{pa 3016 3018}%
\special{fp}%
\special{pa 3156 3218}%
\special{pa 2956 3018}%
\special{fp}%
\special{pa 3096 3218}%
\special{pa 2896 3018}%
\special{fp}%
\special{pa 3036 3218}%
\special{pa 2836 3018}%
\special{fp}%
\special{pa 2976 3218}%
\special{pa 2776 3018}%
\special{fp}%
\special{pa 2916 3218}%
\special{pa 2716 3018}%
\special{fp}%
\special{pa 2856 3218}%
\special{pa 2656 3018}%
\special{fp}%
\special{pa 2796 3218}%
\special{pa 2596 3018}%
\special{fp}%
\special{pa 2736 3218}%
\special{pa 2536 3018}%
\special{fp}%
\special{pa 2676 3218}%
\special{pa 2476 3018}%
\special{fp}%
\special{pa 2616 3218}%
\special{pa 2416 3018}%
\special{fp}%
\special{pa 2556 3218}%
\special{pa 2356 3018}%
\special{fp}%
\special{pa 2496 3218}%
\special{pa 2296 3018}%
\special{fp}%
\special{pa 2436 3218}%
\special{pa 2236 3018}%
\special{fp}%
\special{pa 2376 3218}%
\special{pa 2176 3018}%
\special{fp}%
\special{pa 2316 3218}%
\special{pa 2116 3018}%
\special{fp}%
\special{pa 2256 3218}%
\special{pa 2056 3018}%
\special{fp}%
\special{pa 2196 3218}%
\special{pa 1996 3018}%
\special{fp}%
\special{pa 2136 3218}%
\special{pa 1936 3018}%
\special{fp}%
\special{pa 2076 3218}%
\special{pa 1876 3018}%
\special{fp}%
\special{pa 2016 3218}%
\special{pa 1816 3018}%
\special{fp}%
\special{pa 1956 3218}%
\special{pa 1756 3018}%
\special{fp}%
\special{pa 1896 3218}%
\special{pa 1696 3018}%
\special{fp}%
\special{pa 1836 3218}%
\special{pa 1686 3068}%
\special{fp}%
\special{pa 1776 3218}%
\special{pa 1686 3128}%
\special{fp}%
\special{pa 1716 3218}%
\special{pa 1686 3188}%
\special{fp}%
\special{pa 5136 3218}%
\special{pa 4936 3018}%
\special{fp}%
\special{pa 5196 3218}%
\special{pa 4996 3018}%
\special{fp}%
\special{pa 5256 3218}%
\special{pa 5056 3018}%
\special{fp}%
%
\special{pn 4}%
\special{pa 5316 3218}%
\special{pa 5116 3018}%
\special{fp}%
\special{pa 5376 3218}%
\special{pa 5176 3018}%
\special{fp}%
\special{pa 5436 3218}%
\special{pa 5236 3018}%
\special{fp}%
\special{pa 5486 3208}%
\special{pa 5296 3018}%
\special{fp}%
\special{pa 5486 3148}%
\special{pa 5356 3018}%
\special{fp}%
\special{pa 5486 3088}%
\special{pa 5416 3018}%
\special{fp}%
%
\special{pn 4}%
\special{pa 1276 3018}%
\special{pa 1086 2828}%
\special{fp}%
\special{pa 1336 3018}%
\special{pa 1136 2818}%
\special{fp}%
\special{pa 1386 3008}%
\special{pa 1196 2818}%
\special{fp}%
\special{pa 1386 2948}%
\special{pa 1256 2818}%
\special{fp}%
\special{pa 1386 2888}%
\special{pa 1316 2818}%
\special{fp}%
\special{pa 1216 3018}%
\special{pa 1086 2888}%
\special{fp}%
\special{pa 1156 3018}%
\special{pa 1086 2948}%
\special{fp}%
%
\special{pn 4}%
\special{pa 1536 2018}%
\special{pa 1336 1818}%
\special{fp}%
\special{pa 1476 2018}%
\special{pa 1276 1818}%
\special{fp}%
\special{pa 1416 2018}%
\special{pa 1216 1818}%
\special{fp}%
\special{pa 1356 2018}%
\special{pa 1156 1818}%
\special{fp}%
\special{pa 1296 2018}%
\special{pa 1096 1818}%
\special{fp}%
\special{pa 1236 2018}%
\special{pa 1046 1828}%
\special{fp}%
\special{pa 1176 2018}%
\special{pa 1036 1878}%
\special{fp}%
\special{pa 1116 2018}%
\special{pa 1036 1938}%
\special{fp}%
\special{pa 1596 2018}%
\special{pa 1396 1818}%
\special{fp}%
\special{pa 1656 2018}%
\special{pa 1456 1818}%
\special{fp}%
\special{pa 1716 2018}%
\special{pa 1516 1818}%
\special{fp}%
\special{pa 1776 2018}%
\special{pa 1576 1818}%
\special{fp}%
\special{pa 1836 2018}%
\special{pa 1636 1818}%
\special{fp}%
\special{pa 1896 2018}%
\special{pa 1696 1818}%
\special{fp}%
\special{pa 1956 2018}%
\special{pa 1756 1818}%
\special{fp}%
\special{pa 2016 2018}%
\special{pa 1816 1818}%
\special{fp}%
\special{pa 2076 2018}%
\special{pa 1876 1818}%
\special{fp}%
\special{pa 2136 2018}%
\special{pa 1936 1818}%
\special{fp}%
\special{pa 2196 2018}%
\special{pa 1996 1818}%
\special{fp}%
\special{pa 2256 2018}%
\special{pa 2056 1818}%
\special{fp}%
\special{pa 2316 2018}%
\special{pa 2116 1818}%
\special{fp}%
\special{pa 2376 2018}%
\special{pa 2176 1818}%
\special{fp}%
\special{pa 2436 2018}%
\special{pa 2236 1818}%
\special{fp}%
\special{pa 2496 2018}%
\special{pa 2296 1818}%
\special{fp}%
\special{pa 2556 2018}%
\special{pa 2356 1818}%
\special{fp}%
\special{pa 2616 2018}%
\special{pa 2416 1818}%
\special{fp}%
\special{pa 2676 2018}%
\special{pa 2476 1818}%
\special{fp}%
\special{pa 2736 2018}%
\special{pa 2536 1818}%
\special{fp}%
\special{pa 2796 2018}%
\special{pa 2596 1818}%
\special{fp}%
\special{pa 2856 2018}%
\special{pa 2656 1818}%
\special{fp}%
%
\special{pn 4}%
\special{pa 2916 2018}%
\special{pa 2716 1818}%
\special{fp}%
\special{pa 2976 2018}%
\special{pa 2776 1818}%
\special{fp}%
\special{pa 3036 2018}%
\special{pa 2836 1818}%
\special{fp}%
\special{pa 3096 2018}%
\special{pa 2896 1818}%
\special{fp}%
\special{pa 3156 2018}%
\special{pa 2956 1818}%
\special{fp}%
\special{pa 3216 2018}%
\special{pa 3016 1818}%
\special{fp}%
\special{pa 3276 2018}%
\special{pa 3076 1818}%
\special{fp}%
\special{pa 3336 2018}%
\special{pa 3136 1818}%
\special{fp}%
\special{pa 3396 2018}%
\special{pa 3196 1818}%
\special{fp}%
\special{pa 3456 2018}%
\special{pa 3256 1818}%
\special{fp}%
\special{pa 3516 2018}%
\special{pa 3316 1818}%
\special{fp}%
\special{pa 3576 2018}%
\special{pa 3376 1818}%
\special{fp}%
\special{pa 3636 2018}%
\special{pa 3436 1818}%
\special{fp}%
\special{pa 3696 2018}%
\special{pa 3496 1818}%
\special{fp}%
\special{pa 3756 2018}%
\special{pa 3556 1818}%
\special{fp}%
\special{pa 3816 2018}%
\special{pa 3616 1818}%
\special{fp}%
\special{pa 3876 2018}%
\special{pa 3676 1818}%
\special{fp}%
\special{pa 3936 2018}%
\special{pa 3736 1818}%
\special{fp}%
\special{pa 3996 2018}%
\special{pa 3796 1818}%
\special{fp}%
\special{pa 4056 2018}%
\special{pa 3856 1818}%
\special{fp}%
\special{pa 4116 2018}%
\special{pa 3916 1818}%
\special{fp}%
\special{pa 4176 2018}%
\special{pa 3976 1818}%
\special{fp}%
\special{pa 4236 2018}%
\special{pa 4036 1818}%
\special{fp}%
\special{pa 4296 2018}%
\special{pa 4096 1818}%
\special{fp}%
\special{pa 4356 2018}%
\special{pa 4156 1818}%
\special{fp}%
\special{pa 4416 2018}%
\special{pa 4216 1818}%
\special{fp}%
\special{pa 4476 2018}%
\special{pa 4276 1818}%
\special{fp}%
\special{pa 4536 2018}%
\special{pa 4336 1818}%
\special{fp}%
\special{pa 4596 2018}%
\special{pa 4396 1818}%
\special{fp}%
\special{pa 4636 1998}%
\special{pa 4456 1818}%
\special{fp}%
%
\special{pn 4}%
\special{pa 4636 1938}%
\special{pa 4516 1818}%
\special{fp}%
\special{pa 4636 1878}%
\special{pa 4576 1818}%
\special{fp}%
%
\special{pn 4}%
\special{pa 5336 2218}%
\special{pa 5136 2018}%
\special{fp}%
\special{pa 5276 2218}%
\special{pa 5076 2018}%
\special{fp}%
\special{pa 5216 2218}%
\special{pa 5016 2018}%
\special{fp}%
\special{pa 5156 2218}%
\special{pa 4956 2018}%
\special{fp}%
\special{pa 5096 2218}%
\special{pa 4936 2058}%
\special{fp}%
\special{pa 5036 2218}%
\special{pa 4936 2118}%
\special{fp}%
\special{pa 4976 2218}%
\special{pa 4936 2178}%
\special{fp}%
\special{pa 5396 2218}%
\special{pa 5196 2018}%
\special{fp}%
\special{pa 5456 2218}%
\special{pa 5256 2018}%
\special{fp}%
\special{pa 5516 2218}%
\special{pa 5316 2018}%
\special{fp}%
\special{pa 5536 2178}%
\special{pa 5376 2018}%
\special{fp}%
\special{pa 5536 2118}%
\special{pa 5436 2018}%
\special{fp}%
\special{pa 5536 2058}%
\special{pa 5496 2018}%
\special{fp}%
\put(47.8500,-21.1700){\makebox(0,0){$0$}}%
\put(47.8500,-19.1700){\makebox(0,0){$0$}}%
\put(27.3500,-27.1700){\makebox(0,0){$\cdots$}}%
\put(19.3500,-17.1700){\makebox(0,0)[lb]{negative region}}%
%
\special{pn 8}%
\special{pa 1740 1920}%
\special{pa 1940 1720}%
\special{fp}%
%
\special{pn 8}%
\special{pa 2336 1568}%
\special{pa 2536 1168}%
\special{fp}%
\put(33.2000,-26.6000){\makebox(0,0)[lb]{positive region}}%
%
\special{pn 4}%
\special{pa 5500 1420}%
\special{pa 5300 1220}%
\special{fp}%
\special{pa 5440 1420}%
\special{pa 5240 1220}%
\special{fp}%
\special{pa 5380 1420}%
\special{pa 5180 1220}%
\special{fp}%
\special{pa 5320 1420}%
\special{pa 5120 1220}%
\special{fp}%
\special{pa 5260 1420}%
\special{pa 5060 1220}%
\special{fp}%
\special{pa 5200 1420}%
\special{pa 5000 1220}%
\special{fp}%
\special{pa 5140 1420}%
\special{pa 4940 1220}%
\special{fp}%
\special{pa 5080 1420}%
\special{pa 4880 1220}%
\special{fp}%
\special{pa 5020 1420}%
\special{pa 4820 1220}%
\special{fp}%
\special{pa 4960 1420}%
\special{pa 4760 1220}%
\special{fp}%
\special{pa 4900 1420}%
\special{pa 4700 1220}%
\special{fp}%
\special{pa 4840 1420}%
\special{pa 4640 1220}%
\special{fp}%
\special{pa 4780 1420}%
\special{pa 4580 1220}%
\special{fp}%
\special{pa 4720 1420}%
\special{pa 4520 1220}%
\special{fp}%
\special{pa 4660 1420}%
\special{pa 4460 1220}%
\special{fp}%
\special{pa 4600 1420}%
\special{pa 4400 1220}%
\special{fp}%
\special{pa 4540 1420}%
\special{pa 4340 1220}%
\special{fp}%
\special{pa 4480 1420}%
\special{pa 4320 1260}%
\special{fp}%
\special{pa 4420 1420}%
\special{pa 4320 1320}%
\special{fp}%
\special{pa 4360 1420}%
\special{pa 4320 1380}%
\special{fp}%
\special{pa 5520 1380}%
\special{pa 5360 1220}%
\special{fp}%
\special{pa 5520 1320}%
\special{pa 5420 1220}%
\special{fp}%
\special{pa 5520 1260}%
\special{pa 5480 1220}%
\special{fp}%
\put(49.5000,-15.3000){\makebox(0,0){$\beta(\ast)$}}%
%
\special{pn 8}%
\special{pa 5050 1630}%
\special{pa 4850 1830}%
\special{fp}%
%
\special{pn 8}%
\special{pa 5050 2430}%
\special{pa 4750 2230}%
\special{fp}%
\put(50.5000,-25.3000){\makebox(0,0){$\alpha(\ast)$}}%
\put(45.4000,-35.2000){\makebox(0,0){$\lfloor\frac{2p}{|q_2|}\rfloor$}}%
\put(45.2000,-21.2000){\makebox(0,0){$+$}}%
\put(50.6000,-19.1000){\makebox(0,0){$+$}}%
\put(45.2000,-19.2000){\makebox(0,0){$-$}}%
\put(50.6000,-21.1000){\makebox(0,0){$-$}}%
\end{picture}%

%% file: tableee=1n.tex
\unitlength 0.1in
\begin{picture}( 49.3500, 28.3500)(  6.0000,-34.3500)
\put(15.3500,-31.1700){\makebox(0,0){$0$}}%
\put(21.3500,-31.1700){\makebox(0,0){$+$}}%
\put(18.3500,-31.1700){\makebox(0,0){$+$}}%
\put(18.3500,-35.1700){\makebox(0,0){$2$}}%
\put(21.3500,-35.1700){\makebox(0,0){$3$}}%
\put(7.3500,-29.1700){\makebox(0,0){$1$}}%
\put(7.3500,-27.1700){\makebox(0,0){$2$}}%
\put(7.3500,-25.1700){\makebox(0,0){$3$}}%
\put(12.3500,-31.1700){\makebox(0,0){$-$}}%
\put(15.3500,-35.1700){\makebox(0,0){$1$}}%
\put(12.3500,-35.1700){\makebox(0,0){$0$}}%
\put(7.3500,-31.1700){\makebox(0,0){$0$}}%
%
\special{pn 8}%
\special{pa 936 3318}%
\special{pa 5326 3318}%
\special{fp}%
\special{sh 1}%
\special{pa 5326 3318}%
\special{pa 5258 3298}%
\special{pa 5272 3318}%
\special{pa 5258 3338}%
\special{pa 5326 3318}%
\special{fp}%
\special{pa 936 3318}%
\special{pa 936 600}%
\special{fp}%
\special{sh 1}%
\special{pa 936 600}%
\special{pa 916 668}%
\special{pa 936 654}%
\special{pa 956 668}%
\special{pa 936 600}%
\special{fp}%
\put(7.3500,-7.1700){\makebox(0,0){$j$}}%
\put(53.3500,-35.1700){\makebox(0,0){$i$}}%
\put(33.3500,-35.1700){\makebox(0,0){$\lfloor\frac{p}{|q_2|}\rfloor+1$}}%
\put(24.3500,-31.1700){\makebox(0,0){$\cdots$}}%
\put(33.3500,-31.1700){\makebox(0,0){$0$}}%
\put(30.3500,-31.1700){\makebox(0,0){$+$}}%
\put(36.3500,-31.1700){\makebox(0,0){$+$}}%
\put(12.3500,-29.1700){\makebox(0,0){$+$}}%
\put(18.3500,-29.1700){\makebox(0,0){$-$}}%
\put(15.3500,-29.1700){\makebox(0,0){$0$}}%
\put(21.3500,-29.1700){\makebox(0,0){$0$}}%
\put(24.3500,-29.1700){\makebox(0,0){$\cdots$}}%
\put(21.3500,-27.1700){\makebox(0,0){$-$}}%
\put(18.3500,-27.1700){\makebox(0,0){$0$}}%
\put(15.3500,-27.1700){\makebox(0,0){$0$}}%
\put(12.3500,-27.1700){\makebox(0,0){$0$}}%
\put(39.3500,-31.1700){\makebox(0,0){$\cdots$}}%
\put(24.3500,-27.1700){\makebox(0,0){$0$}}%
\put(24.3500,-25.1700){\makebox(0,0){\odots}}%
\put(27.3500,-23.1700){\makebox(0,0){\odots}}%
\put(30.3500,-21.1400){\makebox(0,0){$-$}}%
\put(30.3500,-19.1400){\makebox(0,0){$+$}}%
\put(27.3500,-19.1400){\makebox(0,0){$+$}}%
\put(21.3500,-19.1400){\makebox(0,0){$+$}}%
\put(24.3500,-19.1400){\makebox(0,0){$+$}}%
\put(41.5500,-13.2700){\makebox(0,0){$0$}}%
\put(50.5500,-11.2700){\makebox(0,0){$\cdots$}}%
\put(47.5500,-11.2700){\makebox(0,0){$0$}}%
\put(32.5500,-11.2700){\makebox(0,0){$+$}}%
\put(35.5500,-11.2700){\makebox(0,0){$\cdots$}}%
\put(32.5500,-13.2700){\makebox(0,0){$\cdots$}}%
\put(35.5500,-13.2700){\makebox(0,0){$0$}}%
\put(50.5500,-13.2700){\makebox(0,0){$\cdots$}}%
\put(47.5500,-13.2700){\makebox(0,0){$+$}}%
\put(41.5500,-11.2700){\makebox(0,0){$0$}}%
\put(44.5500,-11.2700){\makebox(0,0){$-$}}%
\put(38.5500,-11.2700){\makebox(0,0){$+$}}%
\put(38.5500,-13.2700){\makebox(0,0){$-$}}%
\put(44.5500,-13.2700){\makebox(0,0){$+$}}%
\put(33.3500,-19.1400){\makebox(0,0){$0$}}%
\put(33.3500,-21.1400){\makebox(0,0){$0$}}%
\put(36.3500,-21.1400){\makebox(0,0){$\cdots$}}%
\put(36.3500,-19.1400){\makebox(0,0){$\cdots$}}%
\put(27.3500,-31.1700){\makebox(0,0){$\cdots$}}%
\put(42.3500,-31.1700){\makebox(0,0){$+$}}%
\put(45.3500,-31.1700){\makebox(0,0){$0$}}%
\put(44.4500,-17.2700){\makebox(0,0){$\alpha(\ast)+1$}}%
\put(43.4500,-7.2700){\makebox(0,0){$\beta(\ast)+1$}}%
%
\special{pn 8}%
\special{pa 4346 828}%
\special{pa 4146 1028}%
\special{fp}%
%
\special{pn 8}%
\special{pa 4446 1628}%
\special{pa 4146 1428}%
\special{fp}%
%
\special{pn 4}%
\special{pa 1396 1218}%
\special{pa 1196 1018}%
\special{fp}%
\special{pa 1336 1218}%
\special{pa 1136 1018}%
\special{fp}%
\special{pa 1276 1218}%
\special{pa 1076 1018}%
\special{fp}%
\special{pa 1216 1218}%
\special{pa 1016 1018}%
\special{fp}%
\special{pa 1156 1218}%
\special{pa 1006 1068}%
\special{fp}%
\special{pa 1096 1218}%
\special{pa 1006 1128}%
\special{fp}%
\special{pa 1036 1218}%
\special{pa 1006 1188}%
\special{fp}%
\special{pa 1456 1218}%
\special{pa 1256 1018}%
\special{fp}%
\special{pa 1516 1218}%
\special{pa 1316 1018}%
\special{fp}%
\special{pa 1576 1218}%
\special{pa 1376 1018}%
\special{fp}%
\special{pa 1636 1218}%
\special{pa 1436 1018}%
\special{fp}%
\special{pa 1696 1218}%
\special{pa 1496 1018}%
\special{fp}%
\special{pa 1756 1218}%
\special{pa 1556 1018}%
\special{fp}%
\special{pa 1816 1218}%
\special{pa 1616 1018}%
\special{fp}%
\special{pa 1876 1218}%
\special{pa 1676 1018}%
\special{fp}%
\special{pa 1936 1218}%
\special{pa 1736 1018}%
\special{fp}%
\special{pa 1996 1218}%
\special{pa 1796 1018}%
\special{fp}%
\special{pa 2056 1218}%
\special{pa 1856 1018}%
\special{fp}%
\special{pa 2116 1218}%
\special{pa 1916 1018}%
\special{fp}%
\special{pa 2176 1218}%
\special{pa 1976 1018}%
\special{fp}%
\special{pa 2236 1218}%
\special{pa 2036 1018}%
\special{fp}%
\special{pa 2296 1218}%
\special{pa 2096 1018}%
\special{fp}%
\special{pa 2356 1218}%
\special{pa 2156 1018}%
\special{fp}%
\special{pa 2416 1218}%
\special{pa 2216 1018}%
\special{fp}%
\special{pa 2476 1218}%
\special{pa 2276 1018}%
\special{fp}%
\special{pa 2536 1218}%
\special{pa 2336 1018}%
\special{fp}%
\special{pa 2596 1218}%
\special{pa 2396 1018}%
\special{fp}%
\special{pa 2656 1218}%
\special{pa 2456 1018}%
\special{fp}%
\special{pa 2716 1218}%
\special{pa 2516 1018}%
\special{fp}%
\special{pa 2776 1218}%
\special{pa 2576 1018}%
\special{fp}%
%
\special{pn 4}%
\special{pa 2836 1218}%
\special{pa 2636 1018}%
\special{fp}%
\special{pa 2896 1218}%
\special{pa 2696 1018}%
\special{fp}%
\special{pa 2956 1218}%
\special{pa 2756 1018}%
\special{fp}%
\special{pa 3016 1218}%
\special{pa 2816 1018}%
\special{fp}%
\special{pa 3076 1218}%
\special{pa 2876 1018}%
\special{fp}%
\special{pa 3136 1218}%
\special{pa 2936 1018}%
\special{fp}%
\special{pa 3196 1218}%
\special{pa 2996 1018}%
\special{fp}%
\special{pa 3256 1218}%
\special{pa 3056 1018}%
\special{fp}%
\special{pa 3316 1218}%
\special{pa 3116 1018}%
\special{fp}%
\special{pa 3376 1218}%
\special{pa 3176 1018}%
\special{fp}%
\special{pa 3436 1218}%
\special{pa 3236 1018}%
\special{fp}%
\special{pa 3496 1218}%
\special{pa 3296 1018}%
\special{fp}%
\special{pa 3556 1218}%
\special{pa 3356 1018}%
\special{fp}%
\special{pa 3616 1218}%
\special{pa 3416 1018}%
\special{fp}%
\special{pa 3676 1218}%
\special{pa 3476 1018}%
\special{fp}%
\special{pa 3736 1218}%
\special{pa 3536 1018}%
\special{fp}%
\special{pa 3796 1218}%
\special{pa 3596 1018}%
\special{fp}%
\special{pa 3856 1218}%
\special{pa 3656 1018}%
\special{fp}%
\special{pa 3916 1218}%
\special{pa 3716 1018}%
\special{fp}%
\special{pa 3976 1218}%
\special{pa 3776 1018}%
\special{fp}%
\special{pa 4006 1188}%
\special{pa 3836 1018}%
\special{fp}%
\special{pa 4006 1128}%
\special{pa 3896 1018}%
\special{fp}%
\special{pa 4006 1068}%
\special{pa 3956 1018}%
\special{fp}%
%
\special{pn 4}%
\special{pa 5076 3218}%
\special{pa 4876 3018}%
\special{fp}%
\special{pa 5016 3218}%
\special{pa 4816 3018}%
\special{fp}%
\special{pa 4956 3218}%
\special{pa 4756 3018}%
\special{fp}%
\special{pa 4896 3218}%
\special{pa 4696 3018}%
\special{fp}%
\special{pa 4836 3218}%
\special{pa 4636 3018}%
\special{fp}%
\special{pa 4776 3218}%
\special{pa 4576 3018}%
\special{fp}%
\special{pa 4716 3218}%
\special{pa 4516 3018}%
\special{fp}%
\special{pa 4656 3218}%
\special{pa 4456 3018}%
\special{fp}%
\special{pa 4596 3218}%
\special{pa 4396 3018}%
\special{fp}%
\special{pa 4536 3218}%
\special{pa 4336 3018}%
\special{fp}%
\special{pa 4476 3218}%
\special{pa 4276 3018}%
\special{fp}%
\special{pa 4416 3218}%
\special{pa 4216 3018}%
\special{fp}%
\special{pa 4356 3218}%
\special{pa 4156 3018}%
\special{fp}%
\special{pa 4296 3218}%
\special{pa 4096 3018}%
\special{fp}%
\special{pa 4236 3218}%
\special{pa 4036 3018}%
\special{fp}%
\special{pa 4176 3218}%
\special{pa 3976 3018}%
\special{fp}%
\special{pa 4116 3218}%
\special{pa 3916 3018}%
\special{fp}%
\special{pa 4056 3218}%
\special{pa 3856 3018}%
\special{fp}%
\special{pa 3996 3218}%
\special{pa 3796 3018}%
\special{fp}%
\special{pa 3936 3218}%
\special{pa 3736 3018}%
\special{fp}%
\special{pa 3876 3218}%
\special{pa 3676 3018}%
\special{fp}%
\special{pa 3816 3218}%
\special{pa 3616 3018}%
\special{fp}%
\special{pa 3756 3218}%
\special{pa 3556 3018}%
\special{fp}%
\special{pa 3696 3218}%
\special{pa 3496 3018}%
\special{fp}%
\special{pa 3636 3218}%
\special{pa 3436 3018}%
\special{fp}%
\special{pa 3576 3218}%
\special{pa 3376 3018}%
\special{fp}%
\special{pa 3516 3218}%
\special{pa 3316 3018}%
\special{fp}%
\special{pa 3456 3218}%
\special{pa 3256 3018}%
\special{fp}%
\special{pa 3396 3218}%
\special{pa 3196 3018}%
\special{fp}%
\special{pa 3336 3218}%
\special{pa 3136 3018}%
\special{fp}%
%
\special{pn 4}%
\special{pa 3276 3218}%
\special{pa 3076 3018}%
\special{fp}%
\special{pa 3216 3218}%
\special{pa 3016 3018}%
\special{fp}%
\special{pa 3156 3218}%
\special{pa 2956 3018}%
\special{fp}%
\special{pa 3096 3218}%
\special{pa 2896 3018}%
\special{fp}%
\special{pa 3036 3218}%
\special{pa 2836 3018}%
\special{fp}%
\special{pa 2976 3218}%
\special{pa 2776 3018}%
\special{fp}%
\special{pa 2916 3218}%
\special{pa 2716 3018}%
\special{fp}%
\special{pa 2856 3218}%
\special{pa 2656 3018}%
\special{fp}%
\special{pa 2796 3218}%
\special{pa 2596 3018}%
\special{fp}%
\special{pa 2736 3218}%
\special{pa 2536 3018}%
\special{fp}%
\special{pa 2676 3218}%
\special{pa 2476 3018}%
\special{fp}%
\special{pa 2616 3218}%
\special{pa 2416 3018}%
\special{fp}%
\special{pa 2556 3218}%
\special{pa 2356 3018}%
\special{fp}%
\special{pa 2496 3218}%
\special{pa 2296 3018}%
\special{fp}%
\special{pa 2436 3218}%
\special{pa 2236 3018}%
\special{fp}%
\special{pa 2376 3218}%
\special{pa 2176 3018}%
\special{fp}%
\special{pa 2316 3218}%
\special{pa 2116 3018}%
\special{fp}%
\special{pa 2256 3218}%
\special{pa 2056 3018}%
\special{fp}%
\special{pa 2196 3218}%
\special{pa 1996 3018}%
\special{fp}%
\special{pa 2136 3218}%
\special{pa 1936 3018}%
\special{fp}%
\special{pa 2076 3218}%
\special{pa 1876 3018}%
\special{fp}%
\special{pa 2016 3218}%
\special{pa 1816 3018}%
\special{fp}%
\special{pa 1956 3218}%
\special{pa 1756 3018}%
\special{fp}%
\special{pa 1896 3218}%
\special{pa 1696 3018}%
\special{fp}%
\special{pa 1836 3218}%
\special{pa 1686 3068}%
\special{fp}%
\special{pa 1776 3218}%
\special{pa 1686 3128}%
\special{fp}%
\special{pa 1716 3218}%
\special{pa 1686 3188}%
\special{fp}%
\special{pa 5136 3218}%
\special{pa 4936 3018}%
\special{fp}%
\special{pa 5196 3218}%
\special{pa 4996 3018}%
\special{fp}%
\special{pa 5256 3218}%
\special{pa 5056 3018}%
\special{fp}%
%
\special{pn 4}%
\special{pa 5316 3218}%
\special{pa 5116 3018}%
\special{fp}%
\special{pa 5376 3218}%
\special{pa 5176 3018}%
\special{fp}%
\special{pa 5436 3218}%
\special{pa 5236 3018}%
\special{fp}%
\special{pa 5486 3208}%
\special{pa 5296 3018}%
\special{fp}%
\special{pa 5486 3148}%
\special{pa 5356 3018}%
\special{fp}%
\special{pa 5486 3088}%
\special{pa 5416 3018}%
\special{fp}%
%
\special{pn 4}%
\special{pa 1276 3018}%
\special{pa 1086 2828}%
\special{fp}%
\special{pa 1336 3018}%
\special{pa 1136 2818}%
\special{fp}%
\special{pa 1386 3008}%
\special{pa 1196 2818}%
\special{fp}%
\special{pa 1386 2948}%
\special{pa 1256 2818}%
\special{fp}%
\special{pa 1386 2888}%
\special{pa 1316 2818}%
\special{fp}%
\special{pa 1216 3018}%
\special{pa 1086 2888}%
\special{fp}%
\special{pa 1156 3018}%
\special{pa 1086 2948}%
\special{fp}%
%
\special{pn 4}%
\special{pa 1536 2018}%
\special{pa 1336 1818}%
\special{fp}%
\special{pa 1476 2018}%
\special{pa 1276 1818}%
\special{fp}%
\special{pa 1416 2018}%
\special{pa 1216 1818}%
\special{fp}%
\special{pa 1356 2018}%
\special{pa 1156 1818}%
\special{fp}%
\special{pa 1296 2018}%
\special{pa 1096 1818}%
\special{fp}%
\special{pa 1236 2018}%
\special{pa 1046 1828}%
\special{fp}%
\special{pa 1176 2018}%
\special{pa 1036 1878}%
\special{fp}%
\special{pa 1116 2018}%
\special{pa 1036 1938}%
\special{fp}%
\special{pa 1596 2018}%
\special{pa 1396 1818}%
\special{fp}%
\special{pa 1656 2018}%
\special{pa 1456 1818}%
\special{fp}%
\special{pa 1716 2018}%
\special{pa 1516 1818}%
\special{fp}%
\special{pa 1776 2018}%
\special{pa 1576 1818}%
\special{fp}%
\special{pa 1836 2018}%
\special{pa 1636 1818}%
\special{fp}%
\special{pa 1896 2018}%
\special{pa 1696 1818}%
\special{fp}%
\special{pa 1956 2018}%
\special{pa 1756 1818}%
\special{fp}%
\special{pa 2016 2018}%
\special{pa 1816 1818}%
\special{fp}%
\special{pa 2076 2018}%
\special{pa 1876 1818}%
\special{fp}%
\special{pa 2136 2018}%
\special{pa 1936 1818}%
\special{fp}%
\special{pa 2196 2018}%
\special{pa 1996 1818}%
\special{fp}%
\special{pa 2256 2018}%
\special{pa 2056 1818}%
\special{fp}%
\special{pa 2316 2018}%
\special{pa 2116 1818}%
\special{fp}%
\special{pa 2376 2018}%
\special{pa 2176 1818}%
\special{fp}%
\special{pa 2436 2018}%
\special{pa 2236 1818}%
\special{fp}%
\special{pa 2496 2018}%
\special{pa 2296 1818}%
\special{fp}%
\special{pa 2556 2018}%
\special{pa 2356 1818}%
\special{fp}%
\special{pa 2616 2018}%
\special{pa 2416 1818}%
\special{fp}%
\special{pa 2676 2018}%
\special{pa 2476 1818}%
\special{fp}%
\special{pa 2736 2018}%
\special{pa 2536 1818}%
\special{fp}%
\special{pa 2796 2018}%
\special{pa 2596 1818}%
\special{fp}%
\special{pa 2856 2018}%
\special{pa 2656 1818}%
\special{fp}%
%
\special{pn 4}%
\special{pa 2916 2018}%
\special{pa 2716 1818}%
\special{fp}%
\special{pa 2976 2018}%
\special{pa 2776 1818}%
\special{fp}%
\special{pa 3036 2018}%
\special{pa 2836 1818}%
\special{fp}%
\special{pa 3096 2018}%
\special{pa 2896 1818}%
\special{fp}%
\special{pa 3156 2018}%
\special{pa 2956 1818}%
\special{fp}%
\special{pa 3216 2018}%
\special{pa 3016 1818}%
\special{fp}%
\special{pa 3276 2018}%
\special{pa 3076 1818}%
\special{fp}%
\special{pa 3336 2018}%
\special{pa 3136 1818}%
\special{fp}%
\special{pa 3396 2018}%
\special{pa 3196 1818}%
\special{fp}%
\special{pa 3456 2018}%
\special{pa 3256 1818}%
\special{fp}%
\special{pa 3516 2018}%
\special{pa 3316 1818}%
\special{fp}%
\special{pa 3576 2018}%
\special{pa 3376 1818}%
\special{fp}%
\special{pa 3636 2018}%
\special{pa 3436 1818}%
\special{fp}%
\special{pa 3696 2018}%
\special{pa 3496 1818}%
\special{fp}%
\special{pa 3756 2018}%
\special{pa 3556 1818}%
\special{fp}%
\special{pa 3816 2018}%
\special{pa 3616 1818}%
\special{fp}%
\special{pa 3876 2018}%
\special{pa 3676 1818}%
\special{fp}%
\special{pa 3936 2018}%
\special{pa 3736 1818}%
\special{fp}%
\special{pa 3996 2018}%
\special{pa 3796 1818}%
\special{fp}%
\special{pa 4056 2018}%
\special{pa 3856 1818}%
\special{fp}%
\special{pa 4116 2018}%
\special{pa 3916 1818}%
\special{fp}%
\special{pa 4176 2018}%
\special{pa 3976 1818}%
\special{fp}%
\special{pa 4236 2018}%
\special{pa 4036 1818}%
\special{fp}%
\special{pa 4296 2018}%
\special{pa 4096 1818}%
\special{fp}%
\special{pa 4356 2018}%
\special{pa 4156 1818}%
\special{fp}%
\special{pa 4416 2018}%
\special{pa 4216 1818}%
\special{fp}%
\special{pa 4476 2018}%
\special{pa 4276 1818}%
\special{fp}%
\special{pa 4536 2018}%
\special{pa 4336 1818}%
\special{fp}%
\special{pa 4596 2018}%
\special{pa 4396 1818}%
\special{fp}%
\special{pa 4636 1998}%
\special{pa 4456 1818}%
\special{fp}%
%
\special{pn 4}%
\special{pa 4636 1938}%
\special{pa 4516 1818}%
\special{fp}%
\special{pa 4636 1878}%
\special{pa 4576 1818}%
\special{fp}%
%
\special{pn 4}%
\special{pa 5336 2218}%
\special{pa 5136 2018}%
\special{fp}%
\special{pa 5276 2218}%
\special{pa 5076 2018}%
\special{fp}%
\special{pa 5216 2218}%
\special{pa 5016 2018}%
\special{fp}%
\special{pa 5156 2218}%
\special{pa 4956 2018}%
\special{fp}%
\special{pa 5096 2218}%
\special{pa 4936 2058}%
\special{fp}%
\special{pa 5036 2218}%
\special{pa 4936 2118}%
\special{fp}%
\special{pa 4976 2218}%
\special{pa 4936 2178}%
\special{fp}%
\special{pa 5396 2218}%
\special{pa 5196 2018}%
\special{fp}%
\special{pa 5456 2218}%
\special{pa 5256 2018}%
\special{fp}%
\special{pa 5516 2218}%
\special{pa 5316 2018}%
\special{fp}%
\special{pa 5536 2178}%
\special{pa 5376 2018}%
\special{fp}%
\special{pa 5536 2118}%
\special{pa 5436 2018}%
\special{fp}%
\special{pa 5536 2058}%
\special{pa 5496 2018}%
\special{fp}%
\put(47.8500,-21.1700){\makebox(0,0){$0$}}%
\put(47.8500,-19.1700){\makebox(0,0){$0$}}%
\put(27.3500,-27.1700){\makebox(0,0){$\cdots$}}%
\put(19.3500,-17.1700){\makebox(0,0)[lb]{positive region}}%
%
\special{pn 8}%
\special{pa 1740 1920}%
\special{pa 1940 1720}%
\special{fp}%
%
\special{pn 8}%
\special{pa 2336 1568}%
\special{pa 2536 1168}%
\special{fp}%
\put(33.8000,-27.2000){\makebox(0,0)[lb]{negative region}}%
%
\special{pn 4}%
\special{pa 5500 1420}%
\special{pa 5300 1220}%
\special{fp}%
\special{pa 5440 1420}%
\special{pa 5240 1220}%
\special{fp}%
\special{pa 5380 1420}%
\special{pa 5180 1220}%
\special{fp}%
\special{pa 5320 1420}%
\special{pa 5120 1220}%
\special{fp}%
\special{pa 5260 1420}%
\special{pa 5060 1220}%
\special{fp}%
\special{pa 5200 1420}%
\special{pa 5000 1220}%
\special{fp}%
\special{pa 5140 1420}%
\special{pa 4940 1220}%
\special{fp}%
\special{pa 5080 1420}%
\special{pa 4880 1220}%
\special{fp}%
\special{pa 5020 1420}%
\special{pa 4820 1220}%
\special{fp}%
\special{pa 4960 1420}%
\special{pa 4760 1220}%
\special{fp}%
\special{pa 4900 1420}%
\special{pa 4700 1220}%
\special{fp}%
\special{pa 4840 1420}%
\special{pa 4640 1220}%
\special{fp}%
\special{pa 4780 1420}%
\special{pa 4580 1220}%
\special{fp}%
\special{pa 4720 1420}%
\special{pa 4520 1220}%
\special{fp}%
\special{pa 4660 1420}%
\special{pa 4460 1220}%
\special{fp}%
\special{pa 4600 1420}%
\special{pa 4400 1220}%
\special{fp}%
\special{pa 4540 1420}%
\special{pa 4340 1220}%
\special{fp}%
\special{pa 4480 1420}%
\special{pa 4320 1260}%
\special{fp}%
\special{pa 4420 1420}%
\special{pa 4320 1320}%
\special{fp}%
\special{pa 4360 1420}%
\special{pa 4320 1380}%
\special{fp}%
\special{pa 5520 1380}%
\special{pa 5360 1220}%
\special{fp}%
\special{pa 5520 1320}%
\special{pa 5420 1220}%
\special{fp}%
\special{pa 5520 1260}%
\special{pa 5480 1220}%
\special{fp}%
%
\special{pn 8}%
\special{pa 5050 1630}%
\special{pa 4850 1830}%
\special{fp}%
\put(51.5000,-15.3000){\makebox(0,0){$\beta(\ast)+1$}}%
\put(50.5000,-25.3000){\makebox(0,0){$\alpha(\ast)+1$}}%
%
\special{pn 8}%
\special{pa 5100 2430}%
\special{pa 4800 2230}%
\special{fp}%
\put(45.4000,-35.2000){\makebox(0,0){$\lfloor\frac{2p}{|q_2|}\rfloor+1$}}%
\put(45.4000,-21.2000){\makebox(0,0){$-$}}%
\put(51.4000,-19.2000){\makebox(0,0){$-$}}%
\put(51.4000,-21.2000){\makebox(0,0){$+$}}%
\put(45.4000,-19.2000){\makebox(0,0){$+$}}%
\put(33.0000,-24.0000){\makebox(0,0){{\tiny $\lfloor\frac{\ell p}{|q_2|}\rfloor+1$}}}%
%
\special{pn 8}%
\special{pa 3300 2200}%
\special{pa 3300 2300}%
\special{fp}%
\special{pa 1770 2310}%
\special{pa 1770 2310}%
\special{fp}%
\special{pa 1770 2310}%
\special{pa 1770 2310}%
\special{fp}%
\put(27.4000,-21.2000){\makebox(0,0){$0$}}%
\put(42.4000,-21.2000){\makebox(0,0){$0$}}%
\put(24.4000,-21.2000){\makebox(0,0){$\cdots$}}%
\put(40.4000,-21.2000){\makebox(0,0){$\cdots$}}%
\end{picture}%

%% file: ee=-1notp.tex
\unitlength 0.1in
\begin{picture}( 40.6500, 12.0000)(  3.1000,-38.0000)
\put(19.7500,-28.0000){\makebox(0,0){$\lfloor\frac{p\ell}{|q_2|}\rfloor-1$}}%
\put(27.7500,-28.0000){\makebox(0,0){$\lfloor\frac{p\ell}{|q_2|}\rfloor$}}%
\put(35.7500,-28.0000){\makebox(0,0){$\lfloor\frac{p\ell}{|q_2|}\rfloor+1$}}%
\put(25.7500,-38.6000){\makebox(0,0)[lb]{$+1$}}%
\put(33.7500,-38.6000){\makebox(0,0)[lb]{$-1$}}%
\put(19.7500,-31.0000){\makebox(0,0){$0$}}%
\put(27.7500,-31.0000){\makebox(0,0){$1$}}%
\put(35.7500,-31.0000){\makebox(0,0){$0$}}%
\put(19.7500,-33.0000){\makebox(0,0){$-1$}}%
\put(27.7500,-33.0000){\makebox(0,0){$0$}}%
\put(35.7500,-33.0000){\makebox(0,0){$-1$}}%
%
\special{pn 8}%
\special{pa 1576 3000}%
\special{pa 3976 3000}%
\special{fp}%
\special{pa 3976 3000}%
\special{pa 3976 3400}%
\special{fp}%
\special{pa 3976 3400}%
\special{pa 1576 3400}%
\special{fp}%
\special{pa 1576 3400}%
\special{pa 1576 3000}%
\special{fp}%
\special{pa 1576 3200}%
\special{pa 3976 3200}%
\special{fp}%
\special{pa 2376 3000}%
\special{pa 2376 3400}%
\special{fp}%
\special{pa 3176 3000}%
\special{pa 3176 3400}%
\special{fp}%
%
\special{pn 8}%
\special{pa 1576 3000}%
\special{pa 1576 2600}%
\special{fp}%
\special{pa 1576 2600}%
\special{pa 3976 2600}%
\special{fp}%
\special{pa 3976 3000}%
\special{pa 3976 3000}%
\special{fp}%
\special{pa 3976 2600}%
\special{pa 3976 3000}%
\special{fp}%
\special{pa 3176 2600}%
\special{pa 3176 3000}%
\special{fp}%
\special{pa 2376 2600}%
\special{pa 2376 3000}%
\special{fp}%
%
\special{pn 8}%
\special{pa 3976 2600}%
\special{pa 4376 2600}%
\special{fp}%
\special{pa 4376 2600}%
\special{pa 4376 3400}%
\special{fp}%
\special{pa 4376 3400}%
\special{pa 3976 3400}%
\special{fp}%
%
\special{pn 8}%
\special{pa 3976 3000}%
\special{pa 4376 3000}%
\special{fp}%
\special{pa 3976 3200}%
\special{pa 4376 3200}%
\special{fp}%
\put(41.7500,-28.0000){\makebox(0,0)[lb]{$x$}}%
\put(41.7500,-31.0000){\makebox(0,0){(-,i)}}%
\put(41.7500,-33.0000){\makebox(0,0){(-,ii)}}%
%
\special{pn 8}%
\special{pa 1576 3400}%
\special{pa 1576 3600}%
\special{fp}%
\special{pa 1576 3600}%
\special{pa 4376 3600}%
\special{fp}%
\special{pa 4376 3600}%
\special{pa 4376 3400}%
\special{fp}%
\special{pa 3976 3400}%
\special{pa 3976 3600}%
\special{fp}%
\special{pa 3176 3400}%
\special{pa 3176 3600}%
\special{fp}%
\special{pa 2376 3400}%
\special{pa 2376 3600}%
\special{fp}%
\put(19.7500,-35.0000){\makebox(0,0){$1$}}%
\put(35.7500,-35.0000){\makebox(0,0){$1$}}%
\put(27.7500,-35.0000){\makebox(0,0){$2$}}%
\put(41.7500,-35.0000){\makebox(0,0){(-,iii)}}%
%
\special{pn 8}%
\special{ar 2380 3600 200 200  6.2831853 6.2831853}%
\special{ar 2380 3600 200 200  0.0000000 3.1415927}%
%
\special{pn 8}%
\special{pa 2580 3616}%
\special{pa 2580 3600}%
\special{fp}%
\special{sh 1}%
\special{pa 2580 3600}%
\special{pa 2556 3666}%
\special{pa 2576 3654}%
\special{pa 2596 3668}%
\special{pa 2580 3600}%
\special{fp}%
%
\special{pn 8}%
\special{ar 3180 3600 200 200  6.2831853 6.2831853}%
\special{ar 3180 3600 200 200  0.0000000 3.1415927}%
%
\special{pn 8}%
\special{pa 3380 3616}%
\special{pa 3380 3600}%
\special{fp}%
\special{sh 1}%
\special{pa 3380 3600}%
\special{pa 3356 3666}%
\special{pa 3376 3654}%
\special{pa 3396 3668}%
\special{pa 3380 3600}%
\special{fp}%
\end{picture}%

%% file: ee=-1eqp.tex
\unitlength 0.1in
\begin{picture}( 32.5500, 12.0000)( 19.4500,-34.0000)
%
\special{pn 8}%
\special{ar 3200 3200 200 200  6.2831853 6.2831853}%
\special{ar 3200 3200 200 200  0.0000000 3.1415927}%
%
\special{pn 8}%
\special{pa 3400 3216}%
\special{pa 3400 3200}%
\special{fp}%
\special{sh 1}%
\special{pa 3400 3200}%
\special{pa 3376 3266}%
\special{pa 3396 3254}%
\special{pa 3416 3268}%
\special{pa 3400 3200}%
\special{fp}%
%
\special{pn 8}%
\special{ar 4000 3200 200 200  6.2831853 6.2831853}%
\special{ar 4000 3200 200 200  0.0000000 3.1415927}%
%
\special{pn 8}%
\special{pa 4200 3216}%
\special{pa 4200 3200}%
\special{fp}%
\special{sh 1}%
\special{pa 4200 3200}%
\special{pa 4176 3266}%
\special{pa 4196 3254}%
\special{pa 4216 3268}%
\special{pa 4200 3200}%
\special{fp}%
\put(28.0000,-24.0000){\makebox(0,0){$-\epsilon_2q_1j-2$}}%
\put(36.0000,-24.0000){\makebox(0,0){$-\epsilon_2q_1j-1$}}%
\put(44.0000,-24.0000){\makebox(0,0){$-\epsilon_2q_1j$}}%
\put(42.0000,-34.7000){\makebox(0,0)[lb]{$+1$}}%
\put(34.0000,-34.7000){\makebox(0,0)[lb]{$0$}}%
\put(28.0500,-27.0000){\makebox(0,0){$-1$}}%
\put(36.0500,-27.0000){\makebox(0,0){$-1$}}%
\put(44.0500,-27.0000){\makebox(0,0){$0$}}%
%
\special{pn 8}%
\special{pa 2400 2600}%
\special{pa 2400 2200}%
\special{fp}%
\special{pa 2400 2200}%
\special{pa 4800 2200}%
\special{fp}%
\special{pa 4800 2600}%
\special{pa 4800 2600}%
\special{fp}%
\special{pa 4800 2200}%
\special{pa 4800 2600}%
\special{fp}%
\special{pa 4000 2200}%
\special{pa 4000 2600}%
\special{fp}%
\special{pa 3200 2200}%
\special{pa 3200 2600}%
\special{fp}%
%
\special{pn 8}%
\special{pa 2400 2600}%
\special{pa 2400 2800}%
\special{fp}%
\special{pa 2400 2800}%
\special{pa 4800 2800}%
\special{fp}%
\special{pa 4800 2800}%
\special{pa 4800 2600}%
\special{fp}%
\special{pa 4800 2600}%
\special{pa 2400 2600}%
\special{fp}%
\special{pa 3200 2600}%
\special{pa 3200 2800}%
\special{fp}%
\special{pa 4000 2600}%
\special{pa 4000 2800}%
\special{fp}%
\put(44.0000,-29.0000){\makebox(0,0){$1$}}%
\put(36.1000,-29.0000){\makebox(0,0){$0$}}%
\put(28.1000,-29.0000){\makebox(0,0){$0$}}%
%
\special{pn 8}%
\special{pa 2400 2800}%
\special{pa 2400 3000}%
\special{fp}%
\special{pa 2400 3000}%
\special{pa 4800 3000}%
\special{fp}%
\special{pa 4800 3000}%
\special{pa 4800 2800}%
\special{fp}%
\special{pa 4000 2800}%
\special{pa 4000 3000}%
\special{fp}%
\special{pa 3200 2800}%
\special{pa 3200 3000}%
\special{fp}%
%
\special{pn 8}%
\special{pa 4800 2200}%
\special{pa 5200 2200}%
\special{fp}%
\special{pa 5200 2200}%
\special{pa 5200 3000}%
\special{fp}%
\special{pa 5200 3000}%
\special{pa 4800 3000}%
\special{fp}%
\special{pa 4800 2600}%
\special{pa 5200 2600}%
\special{fp}%
\special{pa 4800 2800}%
\special{pa 5200 2800}%
\special{fp}%
\put(50.0000,-24.0000){\makebox(0,0){$x$}}%
\put(50.0000,-27.0000){\makebox(0,0){(-,iv)}}%
\put(50.0000,-29.0000){\makebox(0,0){(-,v)}}%
%
\special{pn 8}%
\special{pa 2400 3000}%
\special{pa 2400 3200}%
\special{fp}%
\special{pa 2400 3200}%
\special{pa 5200 3200}%
\special{fp}%
\special{pa 5200 3200}%
\special{pa 5200 3200}%
\special{fp}%
\special{pa 5200 3200}%
\special{pa 5200 3000}%
\special{fp}%
\special{pa 4800 3200}%
\special{pa 4800 3000}%
\special{fp}%
\special{pa 4000 3200}%
\special{pa 4000 3000}%
\special{fp}%
\special{pa 3200 3200}%
\special{pa 3200 3000}%
\special{fp}%
\put(44.0000,-31.0000){\makebox(0,0){$2$}}%
\put(36.0000,-31.0000){\makebox(0,0){$1$}}%
\put(28.0000,-31.0000){\makebox(0,0){$1$}}%
\put(50.0000,-31.0000){\makebox(0,0){(-,vi)}}%
\end{picture}%

%% file: ee=-1k1p.tex
\unitlength 0.1in
\begin{picture}( 31.6500, 10.0000)( 20.3500,-32.0000)
\put(28.0000,-23.0000){\makebox(0,0){$-\epsilon_2q_1j$}}%
\put(36.0000,-23.0000){\makebox(0,0){$-\epsilon_2q_1j$}}%
\put(44.0000,-23.0000){\makebox(0,0){$-\epsilon_2q_1j$}}%
%
\special{pn 8}%
\special{ar 3200 3000 200 200  6.2831853 6.2831853}%
\special{ar 3200 3000 200 200  0.0000000 3.1415927}%
%
\special{pn 8}%
\special{pa 3400 3016}%
\special{pa 3400 3000}%
\special{fp}%
\special{sh 1}%
\special{pa 3400 3000}%
\special{pa 3376 3066}%
\special{pa 3396 3054}%
\special{pa 3416 3068}%
\special{pa 3400 3000}%
\special{fp}%
%
\special{pn 8}%
\special{ar 4000 3000 200 200  6.2831853 6.2831853}%
\special{ar 4000 3000 200 200  0.0000000 3.1415927}%
%
\special{pn 8}%
\special{pa 4200 3016}%
\special{pa 4200 3000}%
\special{fp}%
\special{sh 1}%
\special{pa 4200 3000}%
\special{pa 4176 3066}%
\special{pa 4196 3054}%
\special{pa 4216 3068}%
\special{pa 4200 3000}%
\special{fp}%
\put(42.0000,-32.7000){\makebox(0,0)[lb]{$-1$}}%
\put(34.0000,-32.7000){\makebox(0,0)[lb]{$0$}}%
\put(28.0500,-27.0000){\makebox(0,0){$1$}}%
\put(36.0500,-27.0000){\makebox(0,0){$1$}}%
\put(44.0500,-27.0000){\makebox(0,0){$0$}}%
%
\special{pn 8}%
\special{pa 2400 2600}%
\special{pa 2400 2200}%
\special{fp}%
\special{pa 2400 2200}%
\special{pa 4800 2200}%
\special{fp}%
\special{pa 4800 2600}%
\special{pa 4800 2600}%
\special{fp}%
\special{pa 4800 2200}%
\special{pa 4800 2600}%
\special{fp}%
\special{pa 4000 2200}%
\special{pa 4000 2600}%
\special{fp}%
\special{pa 3200 2200}%
\special{pa 3200 2600}%
\special{fp}%
%
\special{pn 8}%
\special{pa 2400 2600}%
\special{pa 2400 2800}%
\special{fp}%
\special{pa 2400 2800}%
\special{pa 4800 2800}%
\special{fp}%
\special{pa 4800 2800}%
\special{pa 4800 2600}%
\special{fp}%
\special{pa 4800 2600}%
\special{pa 2400 2600}%
\special{fp}%
\special{pa 3200 2600}%
\special{pa 3200 2800}%
\special{fp}%
\special{pa 4000 2600}%
\special{pa 4000 2800}%
\special{fp}%
\put(36.0000,-25.0000){\makebox(0,0){$-k_1$}}%
\put(44.0000,-25.0000){\makebox(0,0){$-k_1+1$}}%
\put(28.0000,-25.0000){\makebox(0,0){$-k_1-1$}}%
\put(44.1000,-29.0000){\makebox(0,0){$-1$}}%
\put(36.1000,-29.0000){\makebox(0,0){$0$}}%
\put(28.1000,-29.0000){\makebox(0,0){$0$}}%
%
\special{pn 8}%
\special{pa 2400 2800}%
\special{pa 2400 3000}%
\special{fp}%
\special{pa 2400 3000}%
\special{pa 4800 3000}%
\special{fp}%
\special{pa 4800 3000}%
\special{pa 4800 2800}%
\special{fp}%
\special{pa 4000 2800}%
\special{pa 4000 3000}%
\special{fp}%
\special{pa 3200 2800}%
\special{pa 3200 3000}%
\special{fp}%
%
\special{pn 8}%
\special{pa 4800 2200}%
\special{pa 5200 2200}%
\special{fp}%
\special{pa 5200 2200}%
\special{pa 5200 3000}%
\special{fp}%
\special{pa 5200 3000}%
\special{pa 4800 3000}%
\special{fp}%
%
\special{pn 8}%
\special{pa 4800 2600}%
\special{pa 5200 2600}%
\special{fp}%
\special{pa 4800 2800}%
\special{pa 5200 2800}%
\special{fp}%
\put(50.0000,-24.0000){\makebox(0,0){$x$}}%
\put(50.0000,-27.0000){\makebox(0,0){(-,vii)}}%
\put(50.0000,-29.0000){\makebox(0,0){(-,iix)}}%
\end{picture}%

%% file: an2g=p.tex
\unitlength 0.1in
\begin{picture}( 24.0000, 15.5300)( 14.0000,-21.5500)
\put(24.0000,-14.0000){\makebox(0,0){$2$}}%
\put(26.0000,-14.0000){\makebox(0,0){$1$}}%
\put(28.0000,-14.0000){\makebox(0,0){$0$}}%
\put(30.0000,-14.0000){\makebox(0,0){$1$}}%
\put(24.0000,-12.0000){\makebox(0,0){$-1$}}%
\put(24.0000,-16.0000){\makebox(0,0){$-1$}}%
\put(26.0000,-12.0000){\makebox(0,0){$0$}}%
\put(22.0000,-12.0000){\makebox(0,0){$0$}}%
\put(22.0000,-16.0000){\makebox(0,0){$0$}}%
\put(26.0000,-16.0000){\makebox(0,0){$0$}}%
\put(22.0000,-14.0000){\makebox(0,0){$1$}}%
\put(20.0000,-14.0000){\makebox(0,0){$0$}}%
\put(18.0000,-14.0000){\makebox(0,0){$1$}}%
\put(20.0000,-12.0000){\makebox(0,0){$1$}}%
\put(28.0000,-16.0000){\makebox(0,0){$1$}}%
%
\special{pn 8}%
\special{pa 1400 2000}%
\special{pa 3730 2000}%
\special{fp}%
\special{sh 1}%
\special{pa 3730 2000}%
\special{pa 3664 1980}%
\special{pa 3678 2000}%
\special{pa 3664 2020}%
\special{pa 3730 2000}%
\special{fp}%
\special{pa 1556 2156}%
\special{pa 1556 602}%
\special{fp}%
\special{sh 1}%
\special{pa 1556 602}%
\special{pa 1536 670}%
\special{pa 1556 656}%
\special{pa 1576 670}%
\special{pa 1556 602}%
\special{fp}%
%
\special{pn 8}%
\special{pa 3100 1300}%
\special{pa 2300 1300}%
\special{pa 2300 1500}%
\special{pa 3100 1500}%
\special{pa 3100 1300}%
\special{fp}%
%
\special{pn 8}%
\special{pa 1700 1270}%
\special{pa 2500 1270}%
\special{pa 2500 1530}%
\special{pa 1700 1530}%
\special{pa 1700 1270}%
\special{fp}%
%
\special{pn 8}%
\special{pa 2100 1700}%
\special{pa 2900 1700}%
\special{pa 2900 1510}%
\special{pa 2100 1510}%
\special{pa 2100 1700}%
\special{fp}%
%
\special{pn 8}%
\special{pa 1900 1300}%
\special{pa 2700 1300}%
\special{pa 2700 1110}%
\special{pa 1900 1110}%
\special{pa 1900 1300}%
\special{fp}%
\put(38.0000,-20.0000){\makebox(0,0)[lb]{$i$}}%
\put(16.0000,-8.0000){\makebox(0,0)[lb]{$j$}}%
\end{picture}%

%% file: int.tex
\unitlength 0.1in
\begin{picture}( 18.0000,  3.7000)( 12.0000, -8.0000)
%
\special{pn 8}%
\special{pa 1200 600}%
\special{pa 1200 800}%
\special{fp}%
\special{pa 1200 700}%
\special{pa 2200 700}%
\special{fp}%
\special{pa 2200 600}%
\special{pa 2200 800}%
\special{fp}%
%
\special{pn 8}%
\special{pa 2000 600}%
\special{pa 2000 800}%
\special{fp}%
\special{pa 2000 700}%
\special{pa 3000 700}%
\special{fp}%
\special{pa 3000 600}%
\special{pa 3000 800}%
\special{fp}%
\put(12.0000,-6.0000){\makebox(0,0)[lb]{$I_{0,0}$}}%
\put(30.0000,-6.0000){\makebox(0,0)[lb]{$I_{1,0}$}}%
\end{picture}%

%% file: interval.tex
\unitlength 0.1in
\begin{picture}( 57.9500, 14.0000)(  0.3500,-18.0000)
%
\special{pn 8}%
\special{pa 800 800}%
\special{pa 2200 800}%
\special{fp}%
%
\special{pn 8}%
\special{pa 800 700}%
\special{pa 800 900}%
\special{fp}%
%
\special{pn 8}%
\special{pa 2200 700}%
\special{pa 2200 900}%
\special{fp}%
\put(22.0000,-5.0000){\makebox(0,0){$\alpha(a_1+j_2)$}}%
\put(8.0000,-5.0000){\makebox(0,0){$\gamma(a_1+j_2)$}}%
\put(30.0000,-13.0000){\makebox(0,0){$\alpha(j_2)$}}%
\put(16.0000,-13.0000){\makebox(0,0){$\gamma(j_2)$}}%
%
\special{pn 8}%
\special{pa 1600 1000}%
\special{pa 3000 1000}%
\special{fp}%
%
\special{pn 8}%
\special{pa 1600 900}%
\special{pa 1600 1100}%
\special{fp}%
%
\special{pn 8}%
\special{pa 3000 900}%
\special{pa 3000 1100}%
\special{fp}%
%
\special{pn 8}%
\special{pa 800 1000}%
\special{pa 800 1800}%
\special{dt 0.045}%
%
\special{pn 8}%
\special{pa 2200 1000}%
\special{pa 2200 1800}%
\special{dt 0.045}%
\put(14.5000,-18.0000){\makebox(0,0)[lb]{$k_1$}}%
%
\special{pn 8}%
\special{ar 1500 1600 700 140  1.8622531 3.1415927}%
%
\special{pn 8}%
\special{ar 1500 1600 700 140  6.2831853 6.2831853}%
\special{ar 1500 1600 700 140  0.0000000 1.2793395}%
%
\special{pn 8}%
\special{pa 3626 786}%
\special{pa 5026 786}%
\special{fp}%
%
\special{pn 8}%
\special{pa 3610 700}%
\special{pa 3610 900}%
\special{fp}%
%
\special{pn 8}%
\special{pa 5026 686}%
\special{pa 5026 886}%
\special{fp}%
\put(50.2500,-4.8500){\makebox(0,0){$\alpha(a_1+j_2+1)$}}%
\put(36.2500,-4.8500){\makebox(0,0){$\gamma(a_1+j_2+1)$}}%
\put(58.2500,-12.8500){\makebox(0,0){$\alpha(j_2+1)$}}%
\put(44.2500,-12.8500){\makebox(0,0){$\gamma(j_2+1)$}}%
%
\special{pn 8}%
\special{pa 4426 986}%
\special{pa 5826 986}%
\special{fp}%
%
\special{pn 8}%
\special{pa 4420 900}%
\special{pa 4420 1100}%
\special{fp}%
%
\special{pn 8}%
\special{pa 5830 890}%
\special{pa 5830 1090}%
\special{fp}%
\end{picture}%

%% file: ee=-1notm.tex
\unitlength 0.1in
\begin{picture}( 39.7500, 10.0000)( 12.2500,-24.0000)
%
\special{pn 8}%
\special{ar 3200 2200 200 200  6.2831853 6.2831853}%
\special{ar 3200 2200 200 200  0.0000000 3.1415927}%
%
\special{pn 8}%
\special{pa 3002 2216}%
\special{pa 3000 2200}%
\special{fp}%
\special{sh 1}%
\special{pa 3000 2200}%
\special{pa 2984 2268}%
\special{pa 3004 2254}%
\special{pa 3024 2266}%
\special{pa 3000 2200}%
\special{fp}%
%
\special{pn 8}%
\special{ar 4000 2200 200 200  6.2831853 6.2831853}%
\special{ar 4000 2200 200 200  0.0000000 3.1415927}%
%
\special{pn 8}%
\special{pa 3802 2216}%
\special{pa 3800 2200}%
\special{fp}%
\special{sh 1}%
\special{pa 3800 2200}%
\special{pa 3784 2268}%
\special{pa 3804 2254}%
\special{pa 3824 2266}%
\special{pa 3800 2200}%
\special{fp}%
\put(28.0000,-16.0000){\makebox(0,0){$\lfloor\frac{p\ell}{|q_2|}\rfloor$}}%
\put(36.0000,-16.0000){\makebox(0,0){$\lfloor\frac{p\ell}{|q_2|}\rfloor+1$}}%
\put(44.0000,-16.0000){\makebox(0,0){$\lfloor\frac{p\ell}{|q_2|}\rfloor+2$}}%
\put(34.0000,-24.6000){\makebox(0,0)[lb]{$+1$}}%
\put(42.0000,-24.6000){\makebox(0,0)[lb]{$-1$}}%
\put(36.0000,-19.0000){\makebox(0,0){$0$}}%
\put(28.0000,-19.0000){\makebox(0,0){$1$}}%
\put(36.0000,-21.0000){\makebox(0,0){$-1$}}%
\put(28.0000,-21.0000){\makebox(0,0){$0$}}%
%
\special{pn 8}%
\special{pa 2400 1800}%
\special{pa 4800 1800}%
\special{fp}%
\special{pa 4800 1800}%
\special{pa 4800 2200}%
\special{fp}%
\special{pa 4800 2200}%
\special{pa 2400 2200}%
\special{fp}%
\special{pa 2400 2200}%
\special{pa 2400 1800}%
\special{fp}%
\special{pa 2400 2000}%
\special{pa 4800 2000}%
\special{fp}%
\special{pa 3200 1800}%
\special{pa 3200 2200}%
\special{fp}%
\special{pa 4000 1800}%
\special{pa 4000 2200}%
\special{fp}%
%
\special{pn 8}%
\special{pa 2400 1800}%
\special{pa 2400 1400}%
\special{fp}%
\special{pa 2400 1400}%
\special{pa 4800 1400}%
\special{fp}%
\special{pa 4800 1800}%
\special{pa 4800 1800}%
\special{fp}%
\special{pa 4800 1400}%
\special{pa 4800 1800}%
\special{fp}%
\special{pa 4000 1400}%
\special{pa 4000 1800}%
\special{fp}%
\special{pa 3200 1400}%
\special{pa 3200 1800}%
\special{fp}%
\put(44.0000,-19.0000){\makebox(0,0){$1$}}%
\put(44.0000,-21.0000){\makebox(0,0){$0$}}%
%
\special{pn 8}%
\special{pa 4800 1400}%
\special{pa 5200 1400}%
\special{fp}%
\special{pa 5200 1400}%
\special{pa 5200 2200}%
\special{fp}%
\special{pa 5200 2200}%
\special{pa 4800 2200}%
\special{fp}%
\special{pa 4800 1800}%
\special{pa 5200 1800}%
\special{fp}%
\special{pa 4800 2000}%
\special{pa 5200 2000}%
\special{fp}%
\put(50.0000,-16.0000){\makebox(0,0){$x$}}%
\put(50.0000,-19.0000){\makebox(0,0){(-,ix)}}%
\put(50.0000,-21.0000){\makebox(0,0){(-,x)}}%
\end{picture}%

%% file: ee=-1eqm.tex
\unitlength 0.1in
\begin{picture}( 32.5500, 10.0000)( 19.4500,-32.0000)
%
\special{pn 8}%
\special{ar 3200 3000 200 200  6.2831853 6.2831853}%
\special{ar 3200 3000 200 200  0.0000000 3.1415927}%
%
\special{pn 8}%
\special{pa 3002 3016}%
\special{pa 3000 3000}%
\special{fp}%
\special{sh 1}%
\special{pa 3000 3000}%
\special{pa 2984 3068}%
\special{pa 3004 3054}%
\special{pa 3024 3066}%
\special{pa 3000 3000}%
\special{fp}%
%
\special{pn 8}%
\special{ar 4000 3000 200 200  6.2831853 6.2831853}%
\special{ar 4000 3000 200 200  0.0000000 3.1415927}%
%
\special{pn 8}%
\special{pa 3802 3016}%
\special{pa 3800 3000}%
\special{fp}%
\special{sh 1}%
\special{pa 3800 3000}%
\special{pa 3784 3068}%
\special{pa 3804 3054}%
\special{pa 3824 3066}%
\special{pa 3800 3000}%
\special{fp}%
\put(28.0000,-24.0000){\makebox(0,0){$-\epsilon_2q_1j-1$}}%
\put(36.0000,-24.0000){\makebox(0,0){$-\epsilon_2q_1j$}}%
\put(44.0000,-24.0000){\makebox(0,0){$-\epsilon_2q_1j+1$}}%
\put(42.0000,-32.7000){\makebox(0,0)[lb]{$+1$}}%
\put(34.0000,-32.7000){\makebox(0,0)[lb]{$0$}}%
\put(28.0500,-27.0000){\makebox(0,0){$1$}}%
\put(36.0500,-27.0000){\makebox(0,0){$1$}}%
\put(44.0500,-27.0000){\makebox(0,0){$0$}}%
%
\special{pn 8}%
\special{pa 2400 2600}%
\special{pa 2400 2200}%
\special{fp}%
\special{pa 2400 2200}%
\special{pa 4800 2200}%
\special{fp}%
\special{pa 4800 2600}%
\special{pa 4800 2600}%
\special{fp}%
\special{pa 4800 2200}%
\special{pa 4800 2600}%
\special{fp}%
\special{pa 4000 2200}%
\special{pa 4000 2600}%
\special{fp}%
\special{pa 3200 2200}%
\special{pa 3200 2600}%
\special{fp}%
%
\special{pn 8}%
\special{pa 2400 2600}%
\special{pa 2400 2800}%
\special{fp}%
\special{pa 2400 2800}%
\special{pa 4800 2800}%
\special{fp}%
\special{pa 4800 2800}%
\special{pa 4800 2600}%
\special{fp}%
\special{pa 4800 2600}%
\special{pa 2400 2600}%
\special{fp}%
\special{pa 3200 2600}%
\special{pa 3200 2800}%
\special{fp}%
\special{pa 4000 2600}%
\special{pa 4000 2800}%
\special{fp}%
\put(44.0000,-29.0000){\makebox(0,0){$-1$}}%
\put(36.1000,-29.0000){\makebox(0,0){$0$}}%
\put(28.1000,-29.0000){\makebox(0,0){$0$}}%
%
\special{pn 8}%
\special{pa 2400 2800}%
\special{pa 2400 3000}%
\special{fp}%
\special{pa 2400 3000}%
\special{pa 4800 3000}%
\special{fp}%
\special{pa 4800 3000}%
\special{pa 4800 2800}%
\special{fp}%
\special{pa 4000 2800}%
\special{pa 4000 3000}%
\special{fp}%
\special{pa 3200 2800}%
\special{pa 3200 3000}%
\special{fp}%
%
\special{pn 8}%
\special{pa 4800 2200}%
\special{pa 5200 2200}%
\special{fp}%
\special{pa 5200 2200}%
\special{pa 5200 3000}%
\special{fp}%
\special{pa 5200 3000}%
\special{pa 4800 3000}%
\special{fp}%
\special{pa 4800 2600}%
\special{pa 5200 2600}%
\special{fp}%
\special{pa 4800 2800}%
\special{pa 5200 2800}%
\special{fp}%
\put(50.0000,-27.0000){\makebox(0,0){(-,xi)}}%
\put(50.0000,-29.0000){\makebox(0,0){(-,xii)}}%
\put(50.0000,-24.0000){\makebox(0,0){$x$}}%
\end{picture}%

%% file: ee=-1k1m.tex
\unitlength 0.1in
\begin{picture}( 31.6500, 10.0000)( 20.3500,-32.0000)
\put(28.0000,-23.0000){\makebox(0,0){$-\epsilon_2q_1j$}}%
\put(36.0000,-23.0000){\makebox(0,0){$-\epsilon_2q_1j$}}%
\put(44.0000,-23.0000){\makebox(0,0){$-\epsilon_2q_1j$}}%
%
\special{pn 8}%
\special{ar 3200 3000 200 200  6.2831853 6.2831853}%
\special{ar 3200 3000 200 200  0.0000000 3.1415927}%
%
\special{pn 8}%
\special{pa 3002 3016}%
\special{pa 3000 3000}%
\special{fp}%
\special{sh 1}%
\special{pa 3000 3000}%
\special{pa 2984 3068}%
\special{pa 3004 3054}%
\special{pa 3024 3066}%
\special{pa 3000 3000}%
\special{fp}%
%
\special{pn 8}%
\special{ar 4000 3000 200 200  6.2831853 6.2831853}%
\special{ar 4000 3000 200 200  0.0000000 3.1415927}%
%
\special{pn 8}%
\special{pa 3802 3016}%
\special{pa 3800 3000}%
\special{fp}%
\special{sh 1}%
\special{pa 3800 3000}%
\special{pa 3784 3068}%
\special{pa 3804 3054}%
\special{pa 3824 3066}%
\special{pa 3800 3000}%
\special{fp}%
\put(42.0000,-32.7000){\makebox(0,0)[lb]{$-1$}}%
\put(34.0000,-32.7000){\makebox(0,0)[lb]{$0$}}%
\put(28.0500,-27.0000){\makebox(0,0){$-1$}}%
\put(36.0500,-27.0000){\makebox(0,0){$-1$}}%
\put(44.0500,-27.0000){\makebox(0,0){$0$}}%
%
\special{pn 8}%
\special{pa 2400 2600}%
\special{pa 2400 2200}%
\special{fp}%
\special{pa 2400 2200}%
\special{pa 4800 2200}%
\special{fp}%
\special{pa 4800 2600}%
\special{pa 4800 2600}%
\special{fp}%
\special{pa 4800 2200}%
\special{pa 4800 2600}%
\special{fp}%
\special{pa 4000 2200}%
\special{pa 4000 2600}%
\special{fp}%
\special{pa 3200 2200}%
\special{pa 3200 2600}%
\special{fp}%
%
\special{pn 8}%
\special{pa 2400 2600}%
\special{pa 2400 2800}%
\special{fp}%
\special{pa 2400 2800}%
\special{pa 4800 2800}%
\special{fp}%
\special{pa 4800 2800}%
\special{pa 4800 2600}%
\special{fp}%
\special{pa 4800 2600}%
\special{pa 2400 2600}%
\special{fp}%
\special{pa 3200 2600}%
\special{pa 3200 2800}%
\special{fp}%
\special{pa 4000 2600}%
\special{pa 4000 2800}%
\special{fp}%
\put(36.0000,-25.0000){\makebox(0,0){$-k_1+1$}}%
\put(44.0000,-25.0000){\makebox(0,0){$-k_1+2$}}%
\put(28.0000,-25.0000){\makebox(0,0){$-k_1$}}%
\put(44.1000,-29.0000){\makebox(0,0){$1$}}%
\put(36.1000,-29.0000){\makebox(0,0){$0$}}%
\put(28.1000,-29.0000){\makebox(0,0){$0$}}%
%
\special{pn 8}%
\special{pa 2400 2800}%
\special{pa 2400 3000}%
\special{fp}%
\special{pa 2400 3000}%
\special{pa 4800 3000}%
\special{fp}%
\special{pa 4800 3000}%
\special{pa 4800 2800}%
\special{fp}%
\special{pa 4000 2800}%
\special{pa 4000 3000}%
\special{fp}%
\special{pa 3200 2800}%
\special{pa 3200 3000}%
\special{fp}%
%
\special{pn 8}%
\special{pa 4800 2200}%
\special{pa 5200 2200}%
\special{fp}%
\special{pa 5200 2200}%
\special{pa 5200 3000}%
\special{fp}%
\special{pa 5200 3000}%
\special{pa 4800 3000}%
\special{fp}%
\special{pa 4800 2600}%
\special{pa 5200 2600}%
\special{fp}%
\special{pa 4800 2800}%
\special{pa 5200 2800}%
\special{fp}%
\put(50.0000,-24.0000){\makebox(0,0){$x$}}%
\put(50.0000,-27.0000){\makebox(0,0){(-,xiii)}}%
\put(50.0000,-29.0000){\makebox(0,0){(-,xiv)}}%
\end{picture}%

%% file: tableee=-1.tex
\unitlength 0.1in
\begin{picture}( 48.3000, 28.2800)(  5.6000,-34.7000)
\put(15.4000,-29.1700){\makebox(0,0){$+$}}%
\put(21.4000,-31.1700){\makebox(0,0){$\cdots$}}%
\put(18.4000,-31.1700){\makebox(0,0){$0$}}%
\put(18.4500,-35.5500){\makebox(0,0){$1$}}%
\put(21.4500,-35.5500){\makebox(0,0){$2$}}%
\put(7.3500,-27.1400){\makebox(0,0){$1$}}%
\put(7.4000,-33.1700){\makebox(0,0){$-2$}}%
\put(7.4000,-31.1700){\makebox(0,0){$-1$}}%
\put(12.3500,-29.1400){\makebox(0,0){$0$}}%
\put(15.4000,-35.5200){\makebox(0,0){$0$}}%
\put(12.4000,-35.5200){\makebox(0,0){$-1$}}%
\put(7.3500,-29.1400){\makebox(0,0){$0$}}%
%
\special{pn 8}%
\special{pa 1000 3438}%
\special{pa 5390 3438}%
\special{fp}%
\special{sh 1}%
\special{pa 5390 3438}%
\special{pa 5324 3418}%
\special{pa 5338 3438}%
\special{pa 5324 3458}%
\special{pa 5390 3438}%
\special{fp}%
\special{pa 1000 3438}%
\special{pa 1000 720}%
\special{fp}%
\special{sh 1}%
\special{pa 1000 720}%
\special{pa 980 788}%
\special{pa 1000 774}%
\special{pa 1020 788}%
\special{pa 1000 720}%
\special{fp}%
\put(7.4000,-9.2000){\makebox(0,0){$j$}}%
\put(33.4000,-31.1700){\makebox(0,0){$0$}}%
\put(30.4000,-31.1700){\makebox(0,0){$+$}}%
\put(36.4000,-33.1700){\makebox(0,0){$-$}}%
\put(12.4000,-27.1700){\makebox(0,0){$+$}}%
\put(15.3500,-27.1400){\makebox(0,0){$0$}}%
\put(24.3500,-27.1400){\makebox(0,0){$+$}}%
\put(21.3500,-25.1400){\makebox(0,0){$+$}}%
\put(12.3500,-25.1400){\makebox(0,0){$0$}}%
\put(39.4000,-33.1700){\makebox(0,0){$\cdots$}}%
\put(30.4000,-21.0700){\makebox(0,0){$-$}}%
\put(27.4000,-21.0700){\makebox(0,0){$-$}}%
\put(21.4000,-21.0700){\makebox(0,0){$-$}}%
\put(24.4000,-21.0700){\makebox(0,0){$-$}}%
\put(41.5500,-13.2700){\makebox(0,0){$0$}}%
\put(50.5500,-11.2700){\makebox(0,0){$\cdots$}}%
\put(47.5500,-11.2700){\makebox(0,0){$0$}}%
\put(32.5500,-11.2700){\makebox(0,0){$-$}}%
\put(35.5500,-11.2700){\makebox(0,0){$\cdots$}}%
\put(32.5500,-13.2700){\makebox(0,0){$\cdots$}}%
\put(35.5500,-13.2700){\makebox(0,0){$0$}}%
\put(50.5500,-13.2700){\makebox(0,0){$\cdots$}}%
\put(47.5500,-13.2700){\makebox(0,0){$-$}}%
\put(41.5500,-11.2700){\makebox(0,0){$0$}}%
\put(44.5500,-11.2700){\makebox(0,0){$0$}}%
\put(38.5500,-11.2700){\makebox(0,0){$-$}}%
\put(38.5500,-13.2700){\makebox(0,0){$0$}}%
\put(44.5500,-13.2700){\makebox(0,0){$-$}}%
\put(33.4000,-21.0700){\makebox(0,0){$0$}}%
\put(36.4000,-21.0700){\makebox(0,0){$\cdots$}}%
\put(27.3500,-31.1400){\makebox(0,0){$0$}}%
\put(42.4000,-33.1700){\makebox(0,0){$-$}}%
\put(45.4000,-33.1700){\makebox(0,0){$0$}}%
\put(44.5000,-17.3000){\makebox(0,0){$\gamma(\ast)$}}%
\put(43.4500,-7.2700){\makebox(0,0){$\alpha(\ast)$}}%
%
\special{pn 8}%
\special{pa 4346 828}%
\special{pa 4146 1028}%
\special{fp}%
%
\special{pn 8}%
\special{pa 4446 1628}%
\special{pa 4146 1428}%
\special{fp}%
%
\special{pn 4}%
\special{pa 1540 2208}%
\special{pa 1340 2008}%
\special{fp}%
\special{pa 1480 2208}%
\special{pa 1280 2008}%
\special{fp}%
\special{pa 1420 2208}%
\special{pa 1220 2008}%
\special{fp}%
\special{pa 1360 2208}%
\special{pa 1160 2008}%
\special{fp}%
\special{pa 1300 2208}%
\special{pa 1100 2008}%
\special{fp}%
\special{pa 1240 2208}%
\special{pa 1050 2018}%
\special{fp}%
\special{pa 1180 2208}%
\special{pa 1040 2068}%
\special{fp}%
\special{pa 1120 2208}%
\special{pa 1040 2128}%
\special{fp}%
\special{pa 1600 2208}%
\special{pa 1400 2008}%
\special{fp}%
\special{pa 1660 2208}%
\special{pa 1460 2008}%
\special{fp}%
\special{pa 1720 2208}%
\special{pa 1520 2008}%
\special{fp}%
\special{pa 1780 2208}%
\special{pa 1580 2008}%
\special{fp}%
\special{pa 1840 2208}%
\special{pa 1640 2008}%
\special{fp}%
\special{pa 1900 2208}%
\special{pa 1700 2008}%
\special{fp}%
\special{pa 1960 2208}%
\special{pa 1760 2008}%
\special{fp}%
\special{pa 2020 2208}%
\special{pa 1820 2008}%
\special{fp}%
\special{pa 2080 2208}%
\special{pa 1880 2008}%
\special{fp}%
\special{pa 2140 2208}%
\special{pa 1940 2008}%
\special{fp}%
\special{pa 2200 2208}%
\special{pa 2000 2008}%
\special{fp}%
\special{pa 2260 2208}%
\special{pa 2060 2008}%
\special{fp}%
\special{pa 2320 2208}%
\special{pa 2120 2008}%
\special{fp}%
\special{pa 2380 2208}%
\special{pa 2180 2008}%
\special{fp}%
\special{pa 2440 2208}%
\special{pa 2240 2008}%
\special{fp}%
\special{pa 2500 2208}%
\special{pa 2300 2008}%
\special{fp}%
\special{pa 2560 2208}%
\special{pa 2360 2008}%
\special{fp}%
\special{pa 2620 2208}%
\special{pa 2420 2008}%
\special{fp}%
\special{pa 2680 2208}%
\special{pa 2480 2008}%
\special{fp}%
\special{pa 2740 2208}%
\special{pa 2540 2008}%
\special{fp}%
\special{pa 2800 2208}%
\special{pa 2600 2008}%
\special{fp}%
\special{pa 2860 2208}%
\special{pa 2660 2008}%
\special{fp}%
%
\special{pn 4}%
\special{pa 2920 2208}%
\special{pa 2720 2008}%
\special{fp}%
\special{pa 2980 2208}%
\special{pa 2780 2008}%
\special{fp}%
\special{pa 3040 2208}%
\special{pa 2840 2008}%
\special{fp}%
\special{pa 3100 2208}%
\special{pa 2900 2008}%
\special{fp}%
\special{pa 3160 2208}%
\special{pa 2960 2008}%
\special{fp}%
\special{pa 3220 2208}%
\special{pa 3020 2008}%
\special{fp}%
\special{pa 3280 2208}%
\special{pa 3080 2008}%
\special{fp}%
\special{pa 3340 2208}%
\special{pa 3140 2008}%
\special{fp}%
\special{pa 3400 2208}%
\special{pa 3200 2008}%
\special{fp}%
\special{pa 3460 2208}%
\special{pa 3260 2008}%
\special{fp}%
\special{pa 3520 2208}%
\special{pa 3320 2008}%
\special{fp}%
\special{pa 3580 2208}%
\special{pa 3380 2008}%
\special{fp}%
\special{pa 3640 2208}%
\special{pa 3440 2008}%
\special{fp}%
\special{pa 3700 2208}%
\special{pa 3500 2008}%
\special{fp}%
\special{pa 3760 2208}%
\special{pa 3560 2008}%
\special{fp}%
\special{pa 3820 2208}%
\special{pa 3620 2008}%
\special{fp}%
\special{pa 3880 2208}%
\special{pa 3680 2008}%
\special{fp}%
\special{pa 3940 2208}%
\special{pa 3740 2008}%
\special{fp}%
\special{pa 4000 2208}%
\special{pa 3800 2008}%
\special{fp}%
\special{pa 4060 2208}%
\special{pa 3860 2008}%
\special{fp}%
\special{pa 4120 2208}%
\special{pa 3920 2008}%
\special{fp}%
\special{pa 4180 2208}%
\special{pa 3980 2008}%
\special{fp}%
\special{pa 4240 2208}%
\special{pa 4040 2008}%
\special{fp}%
\special{pa 4300 2208}%
\special{pa 4100 2008}%
\special{fp}%
\special{pa 4360 2208}%
\special{pa 4160 2008}%
\special{fp}%
\special{pa 4420 2208}%
\special{pa 4220 2008}%
\special{fp}%
\special{pa 4480 2208}%
\special{pa 4280 2008}%
\special{fp}%
\special{pa 4540 2208}%
\special{pa 4340 2008}%
\special{fp}%
\special{pa 4600 2208}%
\special{pa 4400 2008}%
\special{fp}%
\special{pa 4640 2188}%
\special{pa 4460 2008}%
\special{fp}%
%
\special{pn 4}%
\special{pa 4640 2128}%
\special{pa 4520 2008}%
\special{fp}%
\special{pa 4640 2068}%
\special{pa 4580 2008}%
\special{fp}%
\put(47.9000,-21.1700){\makebox(0,0){$0$}}%
\put(47.9000,-23.0700){\makebox(0,0){$0$}}%
\put(13.4000,-17.1700){\makebox(0,0)[lb]{negative region}}%
%
\special{pn 8}%
\special{pa 1740 1570}%
\special{pa 1940 1170}%
\special{fp}%
\put(29.2000,-26.6000){\makebox(0,0)[lb]{positive region}}%
\put(51.5000,-17.3000){\makebox(0,0){$\alpha(\ast)$}}%
%
\special{pn 8}%
\special{pa 4950 1828}%
\special{pa 4750 2028}%
\special{fp}%
%
\special{pn 8}%
\special{pa 5050 2628}%
\special{pa 4750 2428}%
\special{fp}%
\put(50.5000,-27.2700){\makebox(0,0){$\gamma(\ast)$}}%
\put(45.2000,-21.1700){\makebox(0,0){$-$}}%
\put(50.6000,-23.0700){\makebox(0,0){$-$}}%
\put(45.2000,-23.1700){\makebox(0,0){$0$}}%
\put(50.6000,-21.0700){\makebox(0,0){$0$}}%
\put(15.4000,-31.1700){\makebox(0,0){$-$}}%
\put(12.4000,-31.1700){\makebox(0,0){$-$}}%
\put(18.4000,-29.1700){\makebox(0,0){$0$}}%
\put(18.4000,-33.1700){\makebox(0,0){$0$}}%
\put(21.4000,-33.1700){\makebox(0,0){$-$}}%
\put(24.4000,-33.1700){\makebox(0,0){$-$}}%
\put(27.4000,-33.1700){\makebox(0,0){$\cdots$}}%
\put(33.4000,-33.1700){\makebox(0,0){$0$}}%
\put(30.4000,-33.1700){\makebox(0,0){$-$}}%
\put(27.4000,-29.1700){\makebox(0,0){$+$}}%
\put(18.4000,-21.0700){\makebox(0,0){$-$}}%
\put(18.4000,-23.1700){\makebox(0,0){$+$}}%
\put(42.4000,-31.1700){\makebox(0,0){$+$}}%
%
\special{pn 8}%
\special{pa 5100 3398}%
\special{pa 4920 3218}%
\special{fp}%
\special{pa 5040 3398}%
\special{pa 4860 3218}%
\special{fp}%
\special{pa 4980 3398}%
\special{pa 4800 3218}%
\special{fp}%
\special{pa 4920 3398}%
\special{pa 4740 3218}%
\special{fp}%
\special{pa 4860 3398}%
\special{pa 4680 3218}%
\special{fp}%
\special{pa 4800 3398}%
\special{pa 4620 3218}%
\special{fp}%
\special{pa 4740 3398}%
\special{pa 4560 3218}%
\special{fp}%
\special{pa 4680 3398}%
\special{pa 4500 3218}%
\special{fp}%
\special{pa 4620 3398}%
\special{pa 4440 3218}%
\special{fp}%
\special{pa 4560 3398}%
\special{pa 4380 3218}%
\special{fp}%
\special{pa 4500 3398}%
\special{pa 4320 3218}%
\special{fp}%
\special{pa 4440 3398}%
\special{pa 4260 3218}%
\special{fp}%
\special{pa 4380 3398}%
\special{pa 4200 3218}%
\special{fp}%
\special{pa 4320 3398}%
\special{pa 4140 3218}%
\special{fp}%
\special{pa 4260 3398}%
\special{pa 4080 3218}%
\special{fp}%
\special{pa 4200 3398}%
\special{pa 4020 3218}%
\special{fp}%
\special{pa 4140 3398}%
\special{pa 3960 3218}%
\special{fp}%
\special{pa 4080 3398}%
\special{pa 3900 3218}%
\special{fp}%
\special{pa 4020 3398}%
\special{pa 3840 3218}%
\special{fp}%
\special{pa 3960 3398}%
\special{pa 3780 3218}%
\special{fp}%
\special{pa 3900 3398}%
\special{pa 3720 3218}%
\special{fp}%
\special{pa 3840 3398}%
\special{pa 3660 3218}%
\special{fp}%
\special{pa 3780 3398}%
\special{pa 3600 3218}%
\special{fp}%
\special{pa 3720 3398}%
\special{pa 3540 3218}%
\special{fp}%
\special{pa 3660 3398}%
\special{pa 3480 3218}%
\special{fp}%
\special{pa 3600 3398}%
\special{pa 3420 3218}%
\special{fp}%
\special{pa 3540 3398}%
\special{pa 3360 3218}%
\special{fp}%
\special{pa 3480 3398}%
\special{pa 3300 3218}%
\special{fp}%
\special{pa 3420 3398}%
\special{pa 3240 3218}%
\special{fp}%
\special{pa 3360 3398}%
\special{pa 3180 3218}%
\special{fp}%
%
\special{pn 8}%
\special{pa 3300 3398}%
\special{pa 3120 3218}%
\special{fp}%
\special{pa 3240 3398}%
\special{pa 3060 3218}%
\special{fp}%
\special{pa 3180 3398}%
\special{pa 3000 3218}%
\special{fp}%
\special{pa 3120 3398}%
\special{pa 2940 3218}%
\special{fp}%
\special{pa 3060 3398}%
\special{pa 2880 3218}%
\special{fp}%
\special{pa 3000 3398}%
\special{pa 2820 3218}%
\special{fp}%
\special{pa 2940 3398}%
\special{pa 2760 3218}%
\special{fp}%
\special{pa 2880 3398}%
\special{pa 2700 3218}%
\special{fp}%
\special{pa 2820 3398}%
\special{pa 2640 3218}%
\special{fp}%
\special{pa 2760 3398}%
\special{pa 2580 3218}%
\special{fp}%
\special{pa 2700 3398}%
\special{pa 2520 3218}%
\special{fp}%
\special{pa 2640 3398}%
\special{pa 2460 3218}%
\special{fp}%
\special{pa 2580 3398}%
\special{pa 2400 3218}%
\special{fp}%
\special{pa 2520 3398}%
\special{pa 2340 3218}%
\special{fp}%
\special{pa 2460 3398}%
\special{pa 2280 3218}%
\special{fp}%
\special{pa 2400 3398}%
\special{pa 2220 3218}%
\special{fp}%
\special{pa 2340 3398}%
\special{pa 2160 3218}%
\special{fp}%
\special{pa 2280 3398}%
\special{pa 2100 3218}%
\special{fp}%
\special{pa 2220 3398}%
\special{pa 2040 3218}%
\special{fp}%
\special{pa 2160 3398}%
\special{pa 1990 3228}%
\special{fp}%
\special{pa 2100 3398}%
\special{pa 1980 3278}%
\special{fp}%
\special{pa 2040 3398}%
\special{pa 1980 3338}%
\special{fp}%
\special{pa 5160 3398}%
\special{pa 4980 3218}%
\special{fp}%
\special{pa 5220 3398}%
\special{pa 5040 3218}%
\special{fp}%
\special{pa 5270 3388}%
\special{pa 5100 3218}%
\special{fp}%
\special{pa 5270 3328}%
\special{pa 5160 3218}%
\special{fp}%
\special{pa 5270 3268}%
\special{pa 5220 3218}%
\special{fp}%
%
\special{pn 4}%
\special{pa 1450 3228}%
\special{pa 1240 3018}%
\special{fp}%
\special{pa 1390 3228}%
\special{pa 1180 3018}%
\special{fp}%
\special{pa 1330 3228}%
\special{pa 1120 3018}%
\special{fp}%
\special{pa 1270 3228}%
\special{pa 1070 3028}%
\special{fp}%
\special{pa 1210 3228}%
\special{pa 1060 3078}%
\special{fp}%
\special{pa 1150 3228}%
\special{pa 1060 3138}%
\special{fp}%
\special{pa 1090 3228}%
\special{pa 1060 3198}%
\special{fp}%
\special{pa 1510 3228}%
\special{pa 1300 3018}%
\special{fp}%
\special{pa 1570 3228}%
\special{pa 1360 3018}%
\special{fp}%
\special{pa 1630 3228}%
\special{pa 1420 3018}%
\special{fp}%
\special{pa 1680 3218}%
\special{pa 1480 3018}%
\special{fp}%
\special{pa 1680 3158}%
\special{pa 1540 3018}%
\special{fp}%
\special{pa 1680 3098}%
\special{pa 1600 3018}%
\special{fp}%
%
\special{pn 4}%
\special{pa 5110 2408}%
\special{pa 4910 2208}%
\special{fp}%
\special{pa 5170 2408}%
\special{pa 4960 2198}%
\special{fp}%
\special{pa 5220 2398}%
\special{pa 5020 2198}%
\special{fp}%
\special{pa 5220 2338}%
\special{pa 5080 2198}%
\special{fp}%
\special{pa 5220 2278}%
\special{pa 5140 2198}%
\special{fp}%
\special{pa 5050 2408}%
\special{pa 4900 2258}%
\special{fp}%
\special{pa 4990 2408}%
\special{pa 4900 2318}%
\special{fp}%
\special{pa 4930 2408}%
\special{pa 4900 2378}%
\special{fp}%
%
\special{pn 4}%
\special{pa 2850 1230}%
\special{pa 2650 1030}%
\special{fp}%
\special{pa 2790 1230}%
\special{pa 2590 1030}%
\special{fp}%
\special{pa 2730 1230}%
\special{pa 2530 1030}%
\special{fp}%
\special{pa 2670 1230}%
\special{pa 2470 1030}%
\special{fp}%
\special{pa 2610 1230}%
\special{pa 2410 1030}%
\special{fp}%
\special{pa 2550 1230}%
\special{pa 2350 1030}%
\special{fp}%
\special{pa 2490 1230}%
\special{pa 2290 1030}%
\special{fp}%
\special{pa 2430 1230}%
\special{pa 2230 1030}%
\special{fp}%
\special{pa 2370 1230}%
\special{pa 2170 1030}%
\special{fp}%
\special{pa 2310 1230}%
\special{pa 2110 1030}%
\special{fp}%
\special{pa 2250 1230}%
\special{pa 2050 1030}%
\special{fp}%
\special{pa 2190 1230}%
\special{pa 1990 1030}%
\special{fp}%
\special{pa 2130 1230}%
\special{pa 1930 1030}%
\special{fp}%
\special{pa 2070 1230}%
\special{pa 1870 1030}%
\special{fp}%
\special{pa 2010 1230}%
\special{pa 1810 1030}%
\special{fp}%
\special{pa 1950 1230}%
\special{pa 1750 1030}%
\special{fp}%
\special{pa 1890 1230}%
\special{pa 1690 1030}%
\special{fp}%
\special{pa 1830 1230}%
\special{pa 1630 1030}%
\special{fp}%
\special{pa 1770 1230}%
\special{pa 1570 1030}%
\special{fp}%
\special{pa 1710 1230}%
\special{pa 1510 1030}%
\special{fp}%
\special{pa 1650 1230}%
\special{pa 1450 1030}%
\special{fp}%
\special{pa 1590 1230}%
\special{pa 1390 1030}%
\special{fp}%
\special{pa 1530 1230}%
\special{pa 1330 1030}%
\special{fp}%
\special{pa 1470 1230}%
\special{pa 1270 1030}%
\special{fp}%
\special{pa 1410 1230}%
\special{pa 1210 1030}%
\special{fp}%
\special{pa 1350 1230}%
\special{pa 1150 1030}%
\special{fp}%
\special{pa 1290 1230}%
\special{pa 1090 1030}%
\special{fp}%
\special{pa 1230 1230}%
\special{pa 1050 1050}%
\special{fp}%
\special{pa 1170 1230}%
\special{pa 1050 1110}%
\special{fp}%
\special{pa 1110 1230}%
\special{pa 1050 1170}%
\special{fp}%
%
\special{pn 4}%
\special{pa 2910 1230}%
\special{pa 2710 1030}%
\special{fp}%
\special{pa 2970 1230}%
\special{pa 2770 1030}%
\special{fp}%
\special{pa 3030 1230}%
\special{pa 2830 1030}%
\special{fp}%
\special{pa 3090 1230}%
\special{pa 2890 1030}%
\special{fp}%
\special{pa 3150 1230}%
\special{pa 2950 1030}%
\special{fp}%
\special{pa 3210 1230}%
\special{pa 3010 1030}%
\special{fp}%
\special{pa 3270 1230}%
\special{pa 3070 1030}%
\special{fp}%
\special{pa 3330 1230}%
\special{pa 3130 1030}%
\special{fp}%
\special{pa 3390 1230}%
\special{pa 3190 1030}%
\special{fp}%
\special{pa 3450 1230}%
\special{pa 3250 1030}%
\special{fp}%
\special{pa 3510 1230}%
\special{pa 3310 1030}%
\special{fp}%
\special{pa 3570 1230}%
\special{pa 3370 1030}%
\special{fp}%
\special{pa 3630 1230}%
\special{pa 3430 1030}%
\special{fp}%
\special{pa 3690 1230}%
\special{pa 3490 1030}%
\special{fp}%
\special{pa 3750 1230}%
\special{pa 3550 1030}%
\special{fp}%
\special{pa 3810 1230}%
\special{pa 3610 1030}%
\special{fp}%
\special{pa 3870 1230}%
\special{pa 3670 1030}%
\special{fp}%
\special{pa 3930 1230}%
\special{pa 3730 1030}%
\special{fp}%
\special{pa 3980 1220}%
\special{pa 3790 1030}%
\special{fp}%
\special{pa 3990 1170}%
\special{pa 3850 1030}%
\special{fp}%
\special{pa 3990 1110}%
\special{pa 3910 1030}%
\special{fp}%
%
\special{pn 4}%
\special{pa 5190 1410}%
\special{pa 5010 1230}%
\special{fp}%
\special{pa 5130 1410}%
\special{pa 4950 1230}%
\special{fp}%
\special{pa 5070 1410}%
\special{pa 4890 1230}%
\special{fp}%
\special{pa 5010 1410}%
\special{pa 4830 1230}%
\special{fp}%
\special{pa 4950 1410}%
\special{pa 4770 1230}%
\special{fp}%
\special{pa 4890 1410}%
\special{pa 4710 1230}%
\special{fp}%
\special{pa 4830 1410}%
\special{pa 4650 1230}%
\special{fp}%
\special{pa 4770 1410}%
\special{pa 4590 1230}%
\special{fp}%
\special{pa 4710 1410}%
\special{pa 4530 1230}%
\special{fp}%
\special{pa 4650 1410}%
\special{pa 4470 1230}%
\special{fp}%
\special{pa 4590 1410}%
\special{pa 4410 1230}%
\special{fp}%
\special{pa 4530 1410}%
\special{pa 4350 1230}%
\special{fp}%
\special{pa 4470 1410}%
\special{pa 4320 1260}%
\special{fp}%
\special{pa 4410 1410}%
\special{pa 4320 1320}%
\special{fp}%
\special{pa 4350 1410}%
\special{pa 4320 1380}%
\special{fp}%
\special{pa 5250 1410}%
\special{pa 5070 1230}%
\special{fp}%
\special{pa 5280 1380}%
\special{pa 5130 1230}%
\special{fp}%
\special{pa 5280 1320}%
\special{pa 5190 1230}%
\special{fp}%
\special{pa 5280 1260}%
\special{pa 5250 1230}%
\special{fp}%
%
\special{pn 8}%
\special{pa 1650 1750}%
\special{pa 1450 2150}%
\special{fp}%
\put(15.4000,-21.1000){\makebox(0,0){$0$}}%
\end{picture}%

%% file: tableee=-1n.tex
\unitlength 0.1in
\begin{picture}( 47.7500, 27.9500)(  5.5500,-36.3500)
\put(15.3500,-31.1700){\makebox(0,0){$+$}}%
\put(21.3500,-31.1700){\makebox(0,0){$0$}}%
\put(18.3500,-31.1700){\makebox(0,0){$0$}}%
\put(15.4000,-37.1500){\makebox(0,0){$1$}}%
\put(18.4000,-37.2000){\makebox(0,0){$2$}}%
\put(7.4000,-33.2000){\makebox(0,0){$-2$}}%
\put(12.3500,-31.1700){\makebox(0,0){$+$}}%
\put(12.4000,-37.1500){\makebox(0,0){$0$}}%
\put(7.3500,-31.1700){\makebox(0,0){$-1$}}%
\put(7.4000,-11.2000){\makebox(0,0){$j$}}%
\put(53.4000,-37.1700){\makebox(0,0){$i$}}%
\put(30.4000,-37.2000){\makebox(0,0){$\lfloor\frac{p}{|q_2|}\rfloor+1$}}%
\put(24.3500,-31.1700){\makebox(0,0){$\cdots$}}%
\put(33.3500,-31.1700){\makebox(0,0){$-$}}%
\put(30.3500,-31.1700){\makebox(0,0){$0$}}%
\put(36.3500,-31.1700){\makebox(0,0){$0$}}%
\put(39.3500,-31.1700){\makebox(0,0){$\cdots$}}%
\put(27.4000,-27.1000){\makebox(0,0){$-$}}%
\put(41.6000,-15.2500){\makebox(0,0){$0$}}%
\put(50.6000,-13.2500){\makebox(0,0){$\cdots$}}%
\put(47.6000,-13.2500){\makebox(0,0){$0$}}%
\put(32.6000,-13.2500){\makebox(0,0){$+$}}%
\put(35.6000,-13.2500){\makebox(0,0){$\cdots$}}%
\put(32.6000,-15.2500){\makebox(0,0){$\cdots$}}%
\put(35.6000,-15.2500){\makebox(0,0){$0$}}%
\put(50.6000,-15.2500){\makebox(0,0){$\cdots$}}%
\put(47.6000,-15.2500){\makebox(0,0){$+$}}%
\put(41.6000,-13.2500){\makebox(0,0){$0$}}%
\put(44.6000,-13.2500){\makebox(0,0){$0$}}%
\put(38.6000,-13.2500){\makebox(0,0){$+$}}%
\put(38.6000,-15.2500){\makebox(0,0){$0$}}%
\put(44.6000,-15.2500){\makebox(0,0){$+$}}%
\put(27.3500,-31.1700){\makebox(0,0){$\cdots$}}%
\put(34.5500,-19.2800){\makebox(0,0){$\gamma(\ast)$}}%
\put(43.5000,-9.2500){\makebox(0,0){$\alpha(\ast)$}}%
%
\special{pn 8}%
\special{pa 4150 1030}%
\special{pa 3950 1230}%
\special{fp}%
%
\special{pn 4}%
\special{pa 1400 1416}%
\special{pa 1200 1216}%
\special{fp}%
\special{pa 1340 1416}%
\special{pa 1140 1216}%
\special{fp}%
\special{pa 1280 1416}%
\special{pa 1080 1216}%
\special{fp}%
\special{pa 1220 1416}%
\special{pa 1020 1216}%
\special{fp}%
\special{pa 1160 1416}%
\special{pa 1010 1266}%
\special{fp}%
\special{pa 1100 1416}%
\special{pa 1010 1326}%
\special{fp}%
\special{pa 1040 1416}%
\special{pa 1010 1386}%
\special{fp}%
\special{pa 1460 1416}%
\special{pa 1260 1216}%
\special{fp}%
\special{pa 1520 1416}%
\special{pa 1320 1216}%
\special{fp}%
\special{pa 1580 1416}%
\special{pa 1380 1216}%
\special{fp}%
\special{pa 1640 1416}%
\special{pa 1440 1216}%
\special{fp}%
\special{pa 1700 1416}%
\special{pa 1500 1216}%
\special{fp}%
\special{pa 1760 1416}%
\special{pa 1560 1216}%
\special{fp}%
\special{pa 1820 1416}%
\special{pa 1620 1216}%
\special{fp}%
\special{pa 1880 1416}%
\special{pa 1680 1216}%
\special{fp}%
\special{pa 1940 1416}%
\special{pa 1740 1216}%
\special{fp}%
\special{pa 2000 1416}%
\special{pa 1800 1216}%
\special{fp}%
\special{pa 2060 1416}%
\special{pa 1860 1216}%
\special{fp}%
\special{pa 2120 1416}%
\special{pa 1920 1216}%
\special{fp}%
\special{pa 2180 1416}%
\special{pa 1980 1216}%
\special{fp}%
\special{pa 2240 1416}%
\special{pa 2040 1216}%
\special{fp}%
\special{pa 2300 1416}%
\special{pa 2100 1216}%
\special{fp}%
\special{pa 2360 1416}%
\special{pa 2160 1216}%
\special{fp}%
\special{pa 2420 1416}%
\special{pa 2220 1216}%
\special{fp}%
\special{pa 2480 1416}%
\special{pa 2280 1216}%
\special{fp}%
\special{pa 2540 1416}%
\special{pa 2340 1216}%
\special{fp}%
\special{pa 2600 1416}%
\special{pa 2400 1216}%
\special{fp}%
\special{pa 2660 1416}%
\special{pa 2460 1216}%
\special{fp}%
\special{pa 2720 1416}%
\special{pa 2520 1216}%
\special{fp}%
\special{pa 2780 1416}%
\special{pa 2580 1216}%
\special{fp}%
%
\special{pn 4}%
\special{pa 2840 1416}%
\special{pa 2640 1216}%
\special{fp}%
\special{pa 2900 1416}%
\special{pa 2700 1216}%
\special{fp}%
\special{pa 2960 1416}%
\special{pa 2760 1216}%
\special{fp}%
\special{pa 3020 1416}%
\special{pa 2820 1216}%
\special{fp}%
\special{pa 3080 1416}%
\special{pa 2880 1216}%
\special{fp}%
\special{pa 3140 1416}%
\special{pa 2940 1216}%
\special{fp}%
\special{pa 3200 1416}%
\special{pa 3000 1216}%
\special{fp}%
\special{pa 3260 1416}%
\special{pa 3060 1216}%
\special{fp}%
\special{pa 3320 1416}%
\special{pa 3120 1216}%
\special{fp}%
\special{pa 3380 1416}%
\special{pa 3180 1216}%
\special{fp}%
\special{pa 3440 1416}%
\special{pa 3240 1216}%
\special{fp}%
\special{pa 3500 1416}%
\special{pa 3300 1216}%
\special{fp}%
\special{pa 3560 1416}%
\special{pa 3360 1216}%
\special{fp}%
\special{pa 3620 1416}%
\special{pa 3420 1216}%
\special{fp}%
\special{pa 3680 1416}%
\special{pa 3480 1216}%
\special{fp}%
\special{pa 3740 1416}%
\special{pa 3540 1216}%
\special{fp}%
\special{pa 3800 1416}%
\special{pa 3600 1216}%
\special{fp}%
\special{pa 3860 1416}%
\special{pa 3660 1216}%
\special{fp}%
\special{pa 3920 1416}%
\special{pa 3720 1216}%
\special{fp}%
\special{pa 3980 1416}%
\special{pa 3780 1216}%
\special{fp}%
\special{pa 4010 1386}%
\special{pa 3840 1216}%
\special{fp}%
\special{pa 4010 1326}%
\special{pa 3900 1216}%
\special{fp}%
\special{pa 4010 1266}%
\special{pa 3960 1216}%
\special{fp}%
\put(37.9000,-25.2000){\makebox(0,0){$0$}}%
\put(37.9000,-23.2000){\makebox(0,0){$0$}}%
\put(19.4000,-19.1500){\makebox(0,0)[lb]{positive region}}%
%
\special{pn 8}%
\special{pa 2340 1766}%
\special{pa 2540 1366}%
\special{fp}%
\put(11.0000,-28.0000){\makebox(0,0)[lb]{negative region}}%
\put(43.5500,-19.2800){\makebox(0,0){$\alpha(\ast)$}}%
\put(42.7000,-28.0000){\makebox(0,0){$\gamma(\ast)$}}%
%
\special{pn 8}%
\special{pa 3920 2800}%
\special{pa 3620 2600}%
\special{fp}%
\put(35.4000,-25.2000){\makebox(0,0){$0$}}%
\put(40.6000,-23.2000){\makebox(0,0){$0$}}%
\put(40.5000,-25.2000){\makebox(0,0){$+$}}%
\put(35.4000,-23.2000){\makebox(0,0){$+$}}%
\put(30.4000,-29.1000){\makebox(0,0){$-$}}%
\put(24.4000,-25.1000){\makebox(0,0){$-$}}%
\put(18.4500,-23.1500){\makebox(0,0){$+$}}%
%
\special{pn 8}%
\special{pa 940 3518}%
\special{pa 5330 3518}%
\special{fp}%
\special{sh 1}%
\special{pa 5330 3518}%
\special{pa 5264 3498}%
\special{pa 5278 3518}%
\special{pa 5264 3538}%
\special{pa 5330 3518}%
\special{fp}%
%
\special{pn 8}%
\special{pa 940 3518}%
\special{pa 940 1110}%
\special{fp}%
\special{sh 1}%
\special{pa 940 1110}%
\special{pa 920 1178}%
\special{pa 940 1164}%
\special{pa 960 1178}%
\special{pa 940 1110}%
\special{fp}%
\put(43.1000,-23.2000){\makebox(0,0){$\cdots$}}%
\put(33.1000,-25.2000){\makebox(0,0){$\cdots$}}%
%
\special{pn 4}%
\special{pa 3050 2410}%
\special{pa 2850 2210}%
\special{fp}%
\special{pa 2990 2410}%
\special{pa 2790 2210}%
\special{fp}%
\special{pa 2930 2410}%
\special{pa 2730 2210}%
\special{fp}%
\special{pa 2870 2410}%
\special{pa 2670 2210}%
\special{fp}%
\special{pa 2810 2410}%
\special{pa 2610 2210}%
\special{fp}%
\special{pa 2750 2410}%
\special{pa 2550 2210}%
\special{fp}%
\special{pa 2690 2410}%
\special{pa 2490 2210}%
\special{fp}%
\special{pa 2630 2410}%
\special{pa 2430 2210}%
\special{fp}%
\special{pa 2570 2410}%
\special{pa 2370 2210}%
\special{fp}%
\special{pa 2510 2410}%
\special{pa 2310 2210}%
\special{fp}%
\special{pa 2450 2410}%
\special{pa 2250 2210}%
\special{fp}%
\special{pa 2390 2410}%
\special{pa 2190 2210}%
\special{fp}%
\special{pa 2330 2410}%
\special{pa 2130 2210}%
\special{fp}%
\special{pa 2270 2410}%
\special{pa 2070 2210}%
\special{fp}%
\special{pa 2210 2410}%
\special{pa 2010 2210}%
\special{fp}%
\special{pa 2150 2410}%
\special{pa 1950 2210}%
\special{fp}%
\special{pa 2090 2410}%
\special{pa 1890 2210}%
\special{fp}%
\special{pa 2030 2410}%
\special{pa 1830 2210}%
\special{fp}%
\special{pa 1970 2410}%
\special{pa 1770 2210}%
\special{fp}%
\special{pa 1910 2410}%
\special{pa 1710 2210}%
\special{fp}%
\special{pa 1850 2410}%
\special{pa 1650 2210}%
\special{fp}%
\special{pa 1790 2410}%
\special{pa 1590 2210}%
\special{fp}%
\special{pa 1730 2410}%
\special{pa 1530 2210}%
\special{fp}%
\special{pa 1670 2410}%
\special{pa 1470 2210}%
\special{fp}%
\special{pa 1610 2410}%
\special{pa 1410 2210}%
\special{fp}%
\special{pa 1550 2410}%
\special{pa 1350 2210}%
\special{fp}%
\special{pa 1490 2410}%
\special{pa 1290 2210}%
\special{fp}%
\special{pa 1430 2410}%
\special{pa 1230 2210}%
\special{fp}%
\special{pa 1370 2410}%
\special{pa 1170 2210}%
\special{fp}%
\special{pa 1310 2410}%
\special{pa 1110 2210}%
\special{fp}%
%
\special{pn 4}%
\special{pa 1250 2410}%
\special{pa 1050 2210}%
\special{fp}%
\special{pa 1190 2410}%
\special{pa 1020 2240}%
\special{fp}%
\special{pa 1130 2410}%
\special{pa 1020 2300}%
\special{fp}%
\special{pa 1070 2410}%
\special{pa 1020 2360}%
\special{fp}%
\special{pa 3110 2410}%
\special{pa 2910 2210}%
\special{fp}%
\special{pa 3170 2410}%
\special{pa 2970 2210}%
\special{fp}%
\special{pa 3230 2410}%
\special{pa 3030 2210}%
\special{fp}%
\special{pa 3290 2410}%
\special{pa 3090 2210}%
\special{fp}%
\special{pa 3350 2410}%
\special{pa 3150 2210}%
\special{fp}%
\special{pa 3410 2410}%
\special{pa 3210 2210}%
\special{fp}%
\special{pa 3470 2410}%
\special{pa 3270 2210}%
\special{fp}%
\special{pa 3530 2410}%
\special{pa 3330 2210}%
\special{fp}%
\special{pa 3590 2410}%
\special{pa 3390 2210}%
\special{fp}%
\special{pa 3650 2410}%
\special{pa 3450 2210}%
\special{fp}%
\special{pa 3670 2370}%
\special{pa 3510 2210}%
\special{fp}%
\special{pa 3670 2310}%
\special{pa 3570 2210}%
\special{fp}%
\special{pa 3670 2250}%
\special{pa 3630 2210}%
\special{fp}%
%
\special{pn 4}%
\special{pa 4930 2610}%
\special{pa 4730 2410}%
\special{fp}%
\special{pa 4870 2610}%
\special{pa 4670 2410}%
\special{fp}%
\special{pa 4810 2610}%
\special{pa 4610 2410}%
\special{fp}%
\special{pa 4750 2610}%
\special{pa 4550 2410}%
\special{fp}%
\special{pa 4690 2610}%
\special{pa 4490 2410}%
\special{fp}%
\special{pa 4630 2610}%
\special{pa 4430 2410}%
\special{fp}%
\special{pa 4570 2610}%
\special{pa 4370 2410}%
\special{fp}%
\special{pa 4510 2610}%
\special{pa 4310 2410}%
\special{fp}%
\special{pa 4450 2610}%
\special{pa 4250 2410}%
\special{fp}%
\special{pa 4390 2610}%
\special{pa 4190 2410}%
\special{fp}%
\special{pa 4330 2610}%
\special{pa 4130 2410}%
\special{fp}%
\special{pa 4270 2610}%
\special{pa 4070 2410}%
\special{fp}%
\special{pa 4210 2610}%
\special{pa 4010 2410}%
\special{fp}%
\special{pa 4150 2610}%
\special{pa 3950 2410}%
\special{fp}%
\special{pa 4090 2610}%
\special{pa 3920 2440}%
\special{fp}%
\special{pa 4030 2610}%
\special{pa 3920 2500}%
\special{fp}%
\special{pa 3970 2610}%
\special{pa 3920 2560}%
\special{fp}%
\special{pa 4990 2610}%
\special{pa 4790 2410}%
\special{fp}%
\special{pa 5050 2610}%
\special{pa 4850 2410}%
\special{fp}%
\special{pa 5110 2610}%
\special{pa 4910 2410}%
\special{fp}%
\special{pa 5170 2610}%
\special{pa 4970 2410}%
\special{fp}%
\special{pa 5230 2610}%
\special{pa 5030 2410}%
\special{fp}%
\special{pa 5290 2610}%
\special{pa 5090 2410}%
\special{fp}%
\special{pa 5320 2580}%
\special{pa 5150 2410}%
\special{fp}%
\special{pa 5320 2520}%
\special{pa 5210 2410}%
\special{fp}%
\special{pa 5320 2460}%
\special{pa 5270 2410}%
\special{fp}%
%
\special{pn 4}%
\special{pa 5316 1608}%
\special{pa 5126 1418}%
\special{fp}%
\special{pa 5266 1618}%
\special{pa 5066 1418}%
\special{fp}%
\special{pa 5206 1618}%
\special{pa 5006 1418}%
\special{fp}%
\special{pa 5146 1618}%
\special{pa 4946 1418}%
\special{fp}%
\special{pa 5086 1618}%
\special{pa 4886 1418}%
\special{fp}%
\special{pa 5026 1618}%
\special{pa 4826 1418}%
\special{fp}%
\special{pa 4966 1618}%
\special{pa 4766 1418}%
\special{fp}%
\special{pa 4906 1618}%
\special{pa 4706 1418}%
\special{fp}%
\special{pa 4846 1618}%
\special{pa 4646 1418}%
\special{fp}%
\special{pa 4786 1618}%
\special{pa 4586 1418}%
\special{fp}%
\special{pa 4726 1618}%
\special{pa 4526 1418}%
\special{fp}%
\special{pa 4666 1618}%
\special{pa 4466 1418}%
\special{fp}%
\special{pa 4606 1618}%
\special{pa 4406 1418}%
\special{fp}%
\special{pa 4546 1618}%
\special{pa 4346 1418}%
\special{fp}%
\special{pa 4486 1618}%
\special{pa 4316 1448}%
\special{fp}%
\special{pa 4426 1618}%
\special{pa 4316 1508}%
\special{fp}%
\special{pa 4366 1618}%
\special{pa 4316 1568}%
\special{fp}%
\special{pa 5316 1548}%
\special{pa 5186 1418}%
\special{fp}%
\special{pa 5316 1488}%
\special{pa 5246 1418}%
\special{fp}%
\put(43.1000,-25.2000){\makebox(0,0){$\cdots$}}%
\put(33.1000,-23.2000){\makebox(0,0){$\cdots$}}%
%
\special{pn 8}%
\special{pa 3590 2210}%
\special{pa 4090 2010}%
\special{fp}%
\special{pa 3690 2110}%
\special{pa 3690 2110}%
\special{fp}%
\put(36.4000,-33.2000){\makebox(0,0){$0$}}%
\put(33.4000,-33.1000){\makebox(0,0){$+$}}%
\put(30.4000,-33.1000){\makebox(0,0){$+$}}%
\put(27.4000,-33.1000){\makebox(0,0){$\cdots$}}%
\put(21.4000,-33.2000){\makebox(0,0){$+$}}%
\put(18.4000,-33.2000){\makebox(0,0){$0$}}%
\put(15.4000,-33.2000){\makebox(0,0){$0$}}%
\put(21.4000,-23.2000){\makebox(0,0){$0$}}%
\put(24.4000,-23.2000){\makebox(0,0){$+$}}%
%
\special{pn 8}%
\special{pa 4580 3440}%
\special{pa 4370 3230}%
\special{fp}%
\special{pa 4520 3440}%
\special{pa 4310 3230}%
\special{fp}%
\special{pa 4460 3440}%
\special{pa 4250 3230}%
\special{fp}%
\special{pa 4400 3440}%
\special{pa 4190 3230}%
\special{fp}%
\special{pa 4340 3440}%
\special{pa 4130 3230}%
\special{fp}%
\special{pa 4280 3440}%
\special{pa 4070 3230}%
\special{fp}%
\special{pa 4220 3440}%
\special{pa 4010 3230}%
\special{fp}%
\special{pa 4160 3440}%
\special{pa 3950 3230}%
\special{fp}%
\special{pa 4100 3440}%
\special{pa 3890 3230}%
\special{fp}%
\special{pa 4040 3440}%
\special{pa 3830 3230}%
\special{fp}%
\special{pa 3980 3440}%
\special{pa 3770 3230}%
\special{fp}%
\special{pa 3920 3440}%
\special{pa 3710 3230}%
\special{fp}%
\special{pa 3860 3440}%
\special{pa 3650 3230}%
\special{fp}%
\special{pa 3800 3440}%
\special{pa 3590 3230}%
\special{fp}%
\special{pa 3740 3440}%
\special{pa 3530 3230}%
\special{fp}%
\special{pa 3680 3440}%
\special{pa 3470 3230}%
\special{fp}%
\special{pa 3620 3440}%
\special{pa 3410 3230}%
\special{fp}%
\special{pa 3560 3440}%
\special{pa 3350 3230}%
\special{fp}%
\special{pa 3500 3440}%
\special{pa 3290 3230}%
\special{fp}%
\special{pa 3440 3440}%
\special{pa 3230 3230}%
\special{fp}%
\special{pa 3380 3440}%
\special{pa 3170 3230}%
\special{fp}%
\special{pa 3320 3440}%
\special{pa 3110 3230}%
\special{fp}%
\special{pa 3260 3440}%
\special{pa 3050 3230}%
\special{fp}%
\special{pa 3200 3440}%
\special{pa 2990 3230}%
\special{fp}%
\special{pa 3140 3440}%
\special{pa 2930 3230}%
\special{fp}%
\special{pa 3080 3440}%
\special{pa 2870 3230}%
\special{fp}%
\special{pa 3020 3440}%
\special{pa 2810 3230}%
\special{fp}%
\special{pa 2960 3440}%
\special{pa 2750 3230}%
\special{fp}%
\special{pa 2900 3440}%
\special{pa 2690 3230}%
\special{fp}%
\special{pa 2840 3440}%
\special{pa 2630 3230}%
\special{fp}%
%
\special{pn 8}%
\special{pa 2780 3440}%
\special{pa 2570 3230}%
\special{fp}%
\special{pa 2720 3440}%
\special{pa 2510 3230}%
\special{fp}%
\special{pa 2660 3440}%
\special{pa 2450 3230}%
\special{fp}%
\special{pa 2600 3440}%
\special{pa 2390 3230}%
\special{fp}%
\special{pa 2540 3440}%
\special{pa 2330 3230}%
\special{fp}%
\special{pa 2480 3440}%
\special{pa 2270 3230}%
\special{fp}%
\special{pa 2420 3440}%
\special{pa 2210 3230}%
\special{fp}%
\special{pa 2360 3440}%
\special{pa 2150 3230}%
\special{fp}%
\special{pa 2300 3440}%
\special{pa 2090 3230}%
\special{fp}%
\special{pa 2240 3440}%
\special{pa 2030 3230}%
\special{fp}%
\special{pa 2180 3440}%
\special{pa 1980 3240}%
\special{fp}%
\special{pa 2120 3440}%
\special{pa 1980 3300}%
\special{fp}%
\special{pa 2060 3440}%
\special{pa 1980 3360}%
\special{fp}%
\special{pa 4640 3440}%
\special{pa 4430 3230}%
\special{fp}%
\special{pa 4700 3440}%
\special{pa 4490 3230}%
\special{fp}%
\special{pa 4760 3440}%
\special{pa 4550 3230}%
\special{fp}%
\special{pa 4820 3440}%
\special{pa 4610 3230}%
\special{fp}%
\special{pa 4880 3440}%
\special{pa 4670 3230}%
\special{fp}%
\special{pa 4940 3440}%
\special{pa 4730 3230}%
\special{fp}%
\special{pa 5000 3440}%
\special{pa 4790 3230}%
\special{fp}%
\special{pa 5060 3440}%
\special{pa 4850 3230}%
\special{fp}%
\special{pa 5120 3440}%
\special{pa 4910 3230}%
\special{fp}%
\special{pa 5180 3440}%
\special{pa 4970 3230}%
\special{fp}%
\special{pa 5240 3440}%
\special{pa 5030 3230}%
\special{fp}%
\special{pa 5300 3440}%
\special{pa 5090 3230}%
\special{fp}%
\special{pa 5320 3400}%
\special{pa 5150 3230}%
\special{fp}%
\special{pa 5320 3340}%
\special{pa 5210 3230}%
\special{fp}%
\special{pa 5320 3280}%
\special{pa 5270 3230}%
\special{fp}%
%
\special{pn 4}%
\special{pa 1130 3230}%
\special{pa 1030 3130}%
\special{fp}%
\special{pa 1190 3230}%
\special{pa 1030 3070}%
\special{fp}%
\special{pa 1250 3230}%
\special{pa 1040 3020}%
\special{fp}%
\special{pa 1310 3230}%
\special{pa 1090 3010}%
\special{fp}%
\special{pa 1370 3230}%
\special{pa 1150 3010}%
\special{fp}%
\special{pa 1430 3230}%
\special{pa 1210 3010}%
\special{fp}%
\special{pa 1490 3230}%
\special{pa 1270 3010}%
\special{fp}%
\special{pa 1550 3230}%
\special{pa 1330 3010}%
\special{fp}%
\special{pa 1610 3230}%
\special{pa 1390 3010}%
\special{fp}%
\special{pa 1670 3230}%
\special{pa 1450 3010}%
\special{fp}%
\special{pa 1700 3200}%
\special{pa 1510 3010}%
\special{fp}%
\special{pa 1700 3140}%
\special{pa 1570 3010}%
\special{fp}%
\special{pa 1700 3080}%
\special{pa 1630 3010}%
\special{fp}%
\special{pa 1070 3230}%
\special{pa 1030 3190}%
\special{fp}%
%
\special{pn 8}%
\special{pa 3850 1600}%
\special{pa 3650 1800}%
\special{fp}%
\end{picture}%

%% file: block.tex
\unitlength 0.1in
\begin{picture}( 49.4500, 14.5500)(  7.7500,-19.5000)
\put(18.0000,-10.0000){\makebox(0,0){$\mp$}}%
\put(21.0000,-10.0000){\makebox(0,0){$\mp$}}%
\put(23.9500,-9.9500){\makebox(0,0){$\cdots$}}%
\put(26.9500,-9.9500){\makebox(0,0){$\mp$}}%
\put(27.0000,-18.0000){\makebox(0,0){$\cdots$}}%
\put(30.1000,-18.1000){\makebox(0,0){$\mp$}}%
\put(21.0000,-17.9500){\makebox(0,0){$\mp$}}%
\put(26.9500,-12.0500){\makebox(0,0){$\pm$}}%
\put(21.0500,-16.0500){\makebox(0,0){$\pm$}}%
\put(16.0000,-10.0000){\makebox(0,0){$0$}}%
\put(30.0000,-10.0000){\makebox(0,0){$0$}}%
\put(18.0000,-12.1000){\makebox(0,0){$0$}}%
\put(18.1000,-18.0000){\makebox(0,0){$0$}}%
\put(13.0000,-10.0000){\makebox(0,0){$\cdots$}}%
\put(33.0000,-10.0000){\makebox(0,0){$\cdots$}}%
\put(54.0000,-10.0000){\makebox(0,0){$\mp$}}%
\put(51.8000,-10.0500){\makebox(0,0){$\mp$}}%
\put(48.0000,-10.0000){\makebox(0,0){$\cdots$}}%
\put(44.9000,-10.0000){\makebox(0,0){$\mp$}}%
\put(42.0000,-18.0000){\makebox(0,0){$\mp$}}%
\put(51.0000,-18.0000){\makebox(0,0){$\mp$}}%
\put(45.0000,-18.0000){\makebox(0,0){$\cdots$}}%
\put(48.0000,-18.0000){\makebox(0,0){$\mp$}}%
\put(45.0000,-12.0000){\makebox(0,0){$\pm$}}%
\put(51.0000,-16.0000){\makebox(0,0){$\pm$}}%
\put(57.0000,-10.0000){\makebox(0,0){$0$}}%
\put(42.0000,-10.0000){\makebox(0,0){$0$}}%
\put(54.0000,-18.0000){\makebox(0,0){$0$}}%
\put(60.8000,-10.0500){\makebox(0,0){$\cdots$}}%
\put(39.8000,-10.0500){\makebox(0,0){$\cdots$}}%
\put(24.0000,-14.1000){\makebox(0,0){$\odots$}}%
\put(33.1000,-18.1000){\makebox(0,0){$0$}}%
\put(30.0000,-16.1000){\makebox(0,0){$0$}}%
\put(24.0000,-18.0000){\makebox(0,0){$\mp$}}%
\put(48.0000,-14.0000){\makebox(0,0){$\ddots$}}%
\put(18.1000,-16.0000){\makebox(0,0){$\vdots$}}%
\put(18.1000,-14.0000){\makebox(0,0){$\vdots$}}%
\put(54.0000,-16.0000){\makebox(0,0){$\vdots$}}%
\put(54.0000,-14.0000){\makebox(0,0){$\vdots$}}%
\put(54.0000,-12.0000){\makebox(0,0){$0$}}%
\put(42.0000,-16.0000){\makebox(0,0){$0$}}%
\put(42.0000,-14.0000){\makebox(0,0){$\vdots$}}%
\put(42.0000,-12.0000){\makebox(0,0){$\vdots$}}%
\put(30.0000,-14.0000){\makebox(0,0){$\vdots$}}%
\put(30.0000,-12.0000){\makebox(0,0){$\vdots$}}%
\put(15.0000,-18.0000){\makebox(0,0){$\cdots$}}%
\put(57.0000,-18.0000){\makebox(0,0){$\cdots$}}%
\put(18.0000,-8.0000){\makebox(0,0){$t_0$}}%
\put(27.0000,-8.0000){\makebox(0,0){$t_1$}}%
%
\special{pn 8}%
\special{pa 1000 900}%
\special{pa 3530 900}%
\special{fp}%
\special{sh 1}%
\special{pa 3530 900}%
\special{pa 3464 880}%
\special{pa 3478 900}%
\special{pa 3464 920}%
\special{pa 3530 900}%
\special{fp}%
%
\special{pn 8}%
\special{pa 1200 1950}%
\special{pa 1200 690}%
\special{fp}%
\special{sh 1}%
\special{pa 1200 690}%
\special{pa 1180 758}%
\special{pa 1200 744}%
\special{pa 1220 758}%
\special{pa 1200 690}%
\special{fp}%
\put(21.0000,-8.0000){\makebox(0,0){$b_0$}}%
\put(30.0000,-8.0000){\makebox(0,0){$b_1$}}%
\put(10.0000,-10.0000){\makebox(0,0){$s_1$}}%
\put(10.0000,-18.0000){\makebox(0,0){$s_0$}}%
\put(34.0000,-8.0000){\makebox(0,0){$i$}}%
\put(12.0000,-5.8000){\makebox(0,0){$j$}}%
%
\special{pn 20}%
\special{pa 1800 1000}%
\special{pa 2100 1600}%
\special{pa 2100 1800}%
\special{pa 3000 1800}%
\special{pa 2700 1200}%
\special{pa 2700 1000}%
\special{pa 1800 1000}%
\special{dt 0.054}%
%
\special{pn 20}%
\special{pa 5400 1000}%
\special{pa 5100 1600}%
\special{pa 5100 1800}%
\special{pa 4200 1800}%
\special{pa 4500 1200}%
\special{pa 4500 1000}%
\special{pa 5400 1000}%
\special{dt 0.054}%
\end{picture}%

%% file: Bij.tex
\unitlength 0.1in
\begin{picture}( 28.9500, 15.8500)(  1.0500,-19.0000)
%
\special{pn 8}%
\special{pa 200 1400}%
\special{pa 3000 1400}%
\special{fp}%
\special{sh 1}%
\special{pa 3000 1400}%
\special{pa 2934 1380}%
\special{pa 2948 1400}%
\special{pa 2934 1420}%
\special{pa 3000 1400}%
\special{fp}%
%
\special{pn 8}%
\special{pa 1600 1900}%
\special{pa 1600 410}%
\special{fp}%
\special{sh 1}%
\special{pa 1600 410}%
\special{pa 1580 478}%
\special{pa 1600 464}%
\special{pa 1620 478}%
\special{pa 1600 410}%
\special{fp}%
\put(30.0000,-16.0000){\makebox(0,0)[lb]{$i$}}%
\put(14.0000,-4.0000){\makebox(0,0){$j$}}%
\put(14.7000,-16.0000){\makebox(0,0)[lb]{O}}%
\put(20.0000,-12.0000){\makebox(0,0){$B_{0,0}$}}%
\put(26.0000,-12.0000){\makebox(0,0){$B_{1,0}$}}%
\put(26.0000,-16.0000){\makebox(0,0){$B_{1,-1}$}}%
\put(20.0000,-16.0000){\makebox(0,0){$B_{0,-1}$}}%
\put(12.0000,-12.0000){\makebox(0,0){$B_{-1,0}$}}%
\put(12.0000,-16.0000){\makebox(0,0){$B_{-1,-1}$}}%
\put(6.0000,-16.0000){\makebox(0,0){$B_{-2,-1}$}}%
\put(6.0000,-12.0000){\makebox(0,0){$B_{-2,0}$}}%
%
\special{pn 8}%
\special{pa 1000 1000}%
\special{pa 1400 1000}%
\special{pa 1400 1400}%
\special{pa 1000 1400}%
\special{pa 1000 1000}%
\special{fp}%
%
\special{pn 8}%
\special{pa 1800 1000}%
\special{pa 2200 1000}%
\special{pa 2200 1400}%
\special{pa 1800 1400}%
\special{pa 1800 1000}%
\special{fp}%
%
\special{pn 8}%
\special{pa 2400 1000}%
\special{pa 2800 1000}%
\special{pa 2800 1400}%
\special{pa 2400 1400}%
\special{pa 2400 1000}%
\special{fp}%
%
\special{pn 8}%
\special{pa 2400 1400}%
\special{pa 2800 1400}%
\special{pa 2800 1800}%
\special{pa 2400 1800}%
\special{pa 2400 1400}%
\special{fp}%
%
\special{pn 8}%
\special{pa 1800 1400}%
\special{pa 2200 1400}%
\special{pa 2200 1800}%
\special{pa 1800 1800}%
\special{pa 1800 1400}%
\special{fp}%
%
\special{pn 8}%
\special{pa 1000 1400}%
\special{pa 1400 1400}%
\special{pa 1400 1800}%
\special{pa 1000 1800}%
\special{pa 1000 1400}%
\special{fp}%
%
\special{pn 8}%
\special{pa 400 1400}%
\special{pa 800 1400}%
\special{pa 800 1800}%
\special{pa 400 1800}%
\special{pa 400 1400}%
\special{fp}%
%
\special{pn 8}%
\special{pa 400 1000}%
\special{pa 800 1000}%
\special{pa 800 1400}%
\special{pa 400 1400}%
\special{pa 400 1000}%
\special{fp}%
%
\special{pn 8}%
\special{pa 1800 600}%
\special{pa 2200 600}%
\special{pa 2200 1000}%
\special{pa 1800 1000}%
\special{pa 1800 600}%
\special{fp}%
%
\special{pn 8}%
\special{pa 2400 600}%
\special{pa 2800 600}%
\special{pa 2800 1000}%
\special{pa 2400 1000}%
\special{pa 2400 600}%
\special{fp}%
%
\special{pn 8}%
\special{pa 1000 600}%
\special{pa 1400 600}%
\special{pa 1400 1000}%
\special{pa 1000 1000}%
\special{pa 1000 600}%
\special{fp}%
%
\special{pn 8}%
\special{pa 400 600}%
\special{pa 800 600}%
\special{pa 800 1000}%
\special{pa 400 1000}%
\special{pa 400 600}%
\special{fp}%
\put(20.0000,-8.0000){\makebox(0,0){$B_{0,1}$}}%
\put(26.0000,-8.0000){\makebox(0,0){$B_{1,1}$}}%
\put(12.0000,-8.0000){\makebox(0,0){$B_{-1,1}$}}%
\put(6.0000,-8.0000){\makebox(0,0){$B_{-2,1}$}}%
\end{picture}%

%% file: blocks.tex
\unitlength 0.1in
\begin{picture}( 50.5000, 34.0000)( 11.5000,-40.0000)
%
\special{pn 20}%
\special{sh 1}%
\special{ar 2000 600 10 10 0  6.28318530717959E+0000}%
\special{sh 1}%
\special{ar 2200 600 10 10 0  6.28318530717959E+0000}%
\special{sh 1}%
\special{ar 2400 600 10 10 0  6.28318530717959E+0000}%
\special{sh 1}%
\special{ar 2600 600 10 10 0  6.28318530717959E+0000}%
\special{sh 1}%
\special{ar 2600 800 10 10 0  6.28318530717959E+0000}%
\special{sh 1}%
\special{ar 2400 1000 10 10 0  6.28318530717959E+0000}%
\special{sh 1}%
\special{ar 2200 1200 10 10 0  6.28318530717959E+0000}%
\special{sh 1}%
\special{ar 2200 1400 10 10 0  6.28318530717959E+0000}%
\special{sh 1}%
\special{ar 2400 1400 10 10 0  6.28318530717959E+0000}%
\special{sh 1}%
\special{ar 2600 1400 10 10 0  6.28318530717959E+0000}%
\special{sh 1}%
\special{ar 2800 1400 10 10 0  6.28318530717959E+0000}%
%
\special{pn 8}%
\special{pa 2000 600}%
\special{pa 2600 600}%
\special{pa 2600 800}%
\special{pa 2800 1400}%
\special{pa 2200 1400}%
\special{pa 2200 1200}%
\special{pa 2000 600}%
\special{fp}%
%
\special{pn 20}%
\special{sh 1}%
\special{ar 3000 600 10 10 0  6.28318530717959E+0000}%
\special{sh 1}%
\special{ar 3200 600 10 10 0  6.28318530717959E+0000}%
\special{sh 1}%
\special{ar 3400 600 10 10 0  6.28318530717959E+0000}%
\special{sh 1}%
\special{ar 3600 600 10 10 0  6.28318530717959E+0000}%
\special{sh 1}%
\special{ar 3600 800 10 10 0  6.28318530717959E+0000}%
\special{sh 1}%
\special{ar 3400 1000 10 10 0  6.28318530717959E+0000}%
\special{sh 1}%
\special{ar 3200 1200 10 10 0  6.28318530717959E+0000}%
\special{sh 1}%
\special{ar 3200 1400 10 10 0  6.28318530717959E+0000}%
\special{sh 1}%
\special{ar 3400 1400 10 10 0  6.28318530717959E+0000}%
\special{sh 1}%
\special{ar 3600 1400 10 10 0  6.28318530717959E+0000}%
\special{sh 1}%
\special{ar 3800 1400 10 10 0  6.28318530717959E+0000}%
\special{sh 1}%
\special{ar 4000 1400 10 10 0  6.28318530717959E+0000}%
%
\special{pn 8}%
\special{pa 3000 600}%
\special{pa 3600 600}%
\special{pa 3600 800}%
\special{pa 4000 1400}%
\special{pa 3200 1400}%
\special{pa 3200 1200}%
\special{pa 3000 600}%
\special{fp}%
\put(16.0000,-10.0000){\makebox(0,0){$H=\tau-1$}}%
%
\special{pn 20}%
\special{sh 1}%
\special{ar 4000 600 10 10 0  6.28318530717959E+0000}%
\special{sh 1}%
\special{ar 4200 600 10 10 0  6.28318530717959E+0000}%
\special{sh 1}%
\special{ar 4400 600 10 10 0  6.28318530717959E+0000}%
\special{sh 1}%
\special{ar 4600 600 10 10 0  6.28318530717959E+0000}%
\special{sh 1}%
\special{ar 4800 600 10 10 0  6.28318530717959E+0000}%
\special{sh 1}%
\special{ar 4800 800 10 10 0  6.28318530717959E+0000}%
\special{sh 1}%
\special{ar 4600 1000 10 10 0  6.28318530717959E+0000}%
\special{sh 1}%
\special{ar 4400 1200 10 10 0  6.28318530717959E+0000}%
\special{sh 1}%
\special{ar 4400 1400 10 10 0  6.28318530717959E+0000}%
\special{sh 1}%
\special{ar 4600 1400 10 10 0  6.28318530717959E+0000}%
\special{sh 1}%
\special{ar 4800 1400 10 10 0  6.28318530717959E+0000}%
\special{sh 1}%
\special{ar 5000 1400 10 10 0  6.28318530717959E+0000}%
%
\special{pn 8}%
\special{pa 4000 600}%
\special{pa 4400 1200}%
\special{pa 4400 1400}%
\special{pa 5000 1400}%
\special{pa 4800 800}%
\special{pa 4800 600}%
\special{pa 4000 600}%
\special{fp}%
\put(16.0000,-20.0000){\makebox(0,0){$H=\tau$}}%
%
\special{pn 20}%
\special{sh 1}%
\special{ar 2200 1600 10 10 0  6.28318530717959E+0000}%
\special{sh 1}%
\special{ar 2400 1600 10 10 0  6.28318530717959E+0000}%
\special{sh 1}%
\special{ar 2600 1600 10 10 0  6.28318530717959E+0000}%
\special{sh 1}%
\special{ar 2800 1600 10 10 0  6.28318530717959E+0000}%
\special{sh 1}%
\special{ar 2800 1800 10 10 0  6.28318530717959E+0000}%
\special{sh 1}%
\special{ar 2600 2000 10 10 0  6.28318530717959E+0000}%
\special{sh 1}%
\special{ar 2400 2200 10 10 0  6.28318530717959E+0000}%
\special{sh 1}%
\special{ar 2200 2400 10 10 0  6.28318530717959E+0000}%
\special{sh 1}%
\special{ar 2200 2600 10 10 0  6.28318530717959E+0000}%
\special{sh 1}%
\special{ar 2400 2600 10 10 0  6.28318530717959E+0000}%
\special{sh 1}%
\special{ar 2600 2600 10 10 0  6.28318530717959E+0000}%
\special{sh 1}%
\special{ar 2800 2600 10 10 0  6.28318530717959E+0000}%
\special{sh 1}%
\special{ar 2800 2600 10 10 0  6.28318530717959E+0000}%
%
\special{pn 8}%
\special{pa 2200 1600}%
\special{pa 2200 2600}%
\special{pa 2800 2600}%
\special{pa 2800 1600}%
\special{pa 2200 1600}%
\special{fp}%
%
\special{pn 20}%
\special{sh 1}%
\special{ar 3200 1600 10 10 0  6.28318530717959E+0000}%
\special{sh 1}%
\special{ar 3400 1600 10 10 0  6.28318530717959E+0000}%
\special{sh 1}%
\special{ar 3600 1600 10 10 0  6.28318530717959E+0000}%
\special{sh 1}%
\special{ar 3800 1600 10 10 0  6.28318530717959E+0000}%
\special{sh 1}%
\special{ar 3800 1800 10 10 0  6.28318530717959E+0000}%
\special{sh 1}%
\special{ar 3600 2000 10 10 0  6.28318530717959E+0000}%
\special{sh 1}%
\special{ar 3400 2200 10 10 0  6.28318530717959E+0000}%
\special{sh 1}%
\special{ar 3200 2400 10 10 0  6.28318530717959E+0000}%
\special{sh 1}%
\special{ar 3200 2600 10 10 0  6.28318530717959E+0000}%
\special{sh 1}%
\special{ar 3400 2600 10 10 0  6.28318530717959E+0000}%
\special{sh 1}%
\special{ar 3600 2600 10 10 0  6.28318530717959E+0000}%
\special{sh 1}%
\special{ar 3800 2600 10 10 0  6.28318530717959E+0000}%
\special{sh 1}%
\special{ar 4000 2600 10 10 0  6.28318530717959E+0000}%
%
\special{pn 20}%
\special{sh 1}%
\special{ar 4200 1600 10 10 0  6.28318530717959E+0000}%
\special{sh 1}%
\special{ar 4400 1600 10 10 0  6.28318530717959E+0000}%
\special{sh 1}%
\special{ar 4600 1600 10 10 0  6.28318530717959E+0000}%
\special{sh 1}%
\special{ar 4800 1600 10 10 0  6.28318530717959E+0000}%
\special{sh 1}%
\special{ar 5000 1600 10 10 0  6.28318530717959E+0000}%
\special{sh 1}%
\special{ar 5000 1800 10 10 0  6.28318530717959E+0000}%
\special{sh 1}%
\special{ar 4800 2000 10 10 0  6.28318530717959E+0000}%
\special{sh 1}%
\special{ar 4600 2200 10 10 0  6.28318530717959E+0000}%
\special{sh 1}%
\special{ar 4400 2400 10 10 0  6.28318530717959E+0000}%
\special{sh 1}%
\special{ar 4400 2600 10 10 0  6.28318530717959E+0000}%
\special{sh 1}%
\special{ar 4600 2600 10 10 0  6.28318530717959E+0000}%
\special{sh 1}%
\special{ar 4800 2600 10 10 0  6.28318530717959E+0000}%
\special{sh 1}%
\special{ar 5000 2600 10 10 0  6.28318530717959E+0000}%
%
\special{pn 8}%
\special{pa 3200 1600}%
\special{pa 3800 1600}%
\special{pa 3800 1800}%
\special{pa 4000 2600}%
\special{pa 3200 2600}%
\special{pa 3200 1600}%
\special{fp}%
%
\special{pn 8}%
\special{pa 4400 2600}%
\special{pa 5000 2600}%
\special{pa 5000 1600}%
\special{pa 4200 1600}%
\special{pa 4400 2400}%
\special{pa 4400 2600}%
\special{fp}%
\put(16.0000,-34.0000){\makebox(0,0){$H=\tau+1$}}%
%
\special{pn 20}%
\special{sh 1}%
\special{ar 2200 2800 10 10 0  6.28318530717959E+0000}%
\special{sh 1}%
\special{ar 2400 2800 10 10 0  6.28318530717959E+0000}%
\special{sh 1}%
\special{ar 2600 2800 10 10 0  6.28318530717959E+0000}%
\special{sh 1}%
\special{ar 2800 2800 10 10 0  6.28318530717959E+0000}%
\special{sh 1}%
\special{ar 2800 3000 10 10 0  6.28318530717959E+0000}%
\special{sh 1}%
\special{ar 2600 3200 10 10 0  6.28318530717959E+0000}%
\special{sh 1}%
\special{ar 2400 3400 10 10 0  6.28318530717959E+0000}%
\special{sh 1}%
\special{ar 2200 3600 10 10 0  6.28318530717959E+0000}%
\special{sh 1}%
\special{ar 2000 3800 10 10 0  6.28318530717959E+0000}%
\special{sh 1}%
\special{ar 2000 4000 10 10 0  6.28318530717959E+0000}%
\special{sh 1}%
\special{ar 2200 4000 10 10 0  6.28318530717959E+0000}%
\special{sh 1}%
\special{ar 2400 4000 10 10 0  6.28318530717959E+0000}%
\special{sh 1}%
\special{ar 2600 4000 10 10 0  6.28318530717959E+0000}%
%
\special{pn 8}%
\special{pa 2200 2800}%
\special{pa 2800 2800}%
\special{pa 2800 3000}%
\special{pa 2600 4000}%
\special{pa 2000 4000}%
\special{pa 2000 3800}%
\special{pa 2200 2800}%
\special{fp}%
%
\special{pn 20}%
\special{sh 1}%
\special{ar 3200 2800 10 10 0  6.28318530717959E+0000}%
\special{sh 1}%
\special{ar 3400 2800 10 10 0  6.28318530717959E+0000}%
\special{sh 1}%
\special{ar 3600 2800 10 10 0  6.28318530717959E+0000}%
\special{sh 1}%
\special{ar 3800 2800 10 10 0  6.28318530717959E+0000}%
\special{sh 1}%
\special{ar 3800 3000 10 10 0  6.28318530717959E+0000}%
\special{sh 1}%
\special{ar 3600 3200 10 10 0  6.28318530717959E+0000}%
\special{sh 1}%
\special{ar 3400 3400 10 10 0  6.28318530717959E+0000}%
\special{sh 1}%
\special{ar 3200 3600 10 10 0  6.28318530717959E+0000}%
\special{sh 1}%
\special{ar 3000 3800 10 10 0  6.28318530717959E+0000}%
\special{sh 1}%
\special{ar 3000 4000 10 10 0  6.28318530717959E+0000}%
\special{sh 1}%
\special{ar 3200 4000 10 10 0  6.28318530717959E+0000}%
\special{sh 1}%
\special{ar 3400 4000 10 10 0  6.28318530717959E+0000}%
\special{sh 1}%
\special{ar 3600 4000 10 10 0  6.28318530717959E+0000}%
\special{sh 1}%
\special{ar 3800 4000 10 10 0  6.28318530717959E+0000}%
%
\special{pn 8}%
\special{pa 3000 4000}%
\special{pa 3800 4000}%
\special{pa 3800 2800}%
\special{pa 3200 2800}%
\special{pa 3000 3800}%
\special{pa 3000 4000}%
\special{fp}%
%
\special{pn 20}%
\special{sh 1}%
\special{ar 4200 2800 10 10 0  6.28318530717959E+0000}%
\special{sh 1}%
\special{ar 4400 2800 10 10 0  6.28318530717959E+0000}%
\special{sh 1}%
\special{ar 4600 2800 10 10 0  6.28318530717959E+0000}%
\special{sh 1}%
\special{ar 4800 2800 10 10 0  6.28318530717959E+0000}%
\special{sh 1}%
\special{ar 5000 2800 10 10 0  6.28318530717959E+0000}%
\special{sh 1}%
\special{ar 5000 3000 10 10 0  6.28318530717959E+0000}%
\special{sh 1}%
\special{ar 4800 3200 10 10 0  6.28318530717959E+0000}%
\special{sh 1}%
\special{ar 4600 3400 10 10 0  6.28318530717959E+0000}%
\special{sh 1}%
\special{ar 4400 3600 10 10 0  6.28318530717959E+0000}%
\special{sh 1}%
\special{ar 4200 3800 10 10 0  6.28318530717959E+0000}%
\special{sh 1}%
\special{ar 4200 4000 10 10 0  6.28318530717959E+0000}%
\special{sh 1}%
\special{ar 4400 4000 10 10 0  6.28318530717959E+0000}%
\special{sh 1}%
\special{ar 4600 4000 10 10 0  6.28318530717959E+0000}%
\special{sh 1}%
\special{ar 4800 4000 10 10 0  6.28318530717959E+0000}%
%
\special{pn 8}%
\special{pa 4200 2800}%
\special{pa 4200 4000}%
\special{pa 4800 4000}%
\special{pa 5000 3000}%
\special{pa 5000 2800}%
\special{pa 4200 2800}%
\special{fp}%
%
\special{pn 20}%
\special{sh 1}%
\special{ar 5000 600 10 10 0  6.28318530717959E+0000}%
\special{sh 1}%
\special{ar 5200 600 10 10 0  6.28318530717959E+0000}%
\special{sh 1}%
\special{ar 5400 600 10 10 0  6.28318530717959E+0000}%
\special{sh 1}%
\special{ar 5600 600 10 10 0  6.28318530717959E+0000}%
\special{sh 1}%
\special{ar 5800 600 10 10 0  6.28318530717959E+0000}%
\special{sh 1}%
\special{ar 5800 800 10 10 0  6.28318530717959E+0000}%
\special{sh 1}%
\special{ar 5600 1000 10 10 0  6.28318530717959E+0000}%
\special{sh 1}%
\special{ar 5400 1200 10 10 0  6.28318530717959E+0000}%
\special{sh 1}%
\special{ar 5400 1400 10 10 0  6.28318530717959E+0000}%
\special{sh 1}%
\special{ar 5600 1400 10 10 0  6.28318530717959E+0000}%
\special{sh 1}%
\special{ar 5800 1400 10 10 0  6.28318530717959E+0000}%
\special{sh 1}%
\special{ar 6000 1400 10 10 0  6.28318530717959E+0000}%
\special{sh 1}%
\special{ar 6200 1400 10 10 0  6.28318530717959E+0000}%
%
\special{pn 8}%
\special{pa 5000 600}%
\special{pa 5400 1200}%
\special{pa 5400 1400}%
\special{pa 6200 1400}%
\special{pa 5800 800}%
\special{pa 5800 600}%
\special{pa 5000 600}%
\special{fp}%
%
\special{pn 20}%
\special{sh 1}%
\special{ar 5200 1600 10 10 0  6.28318530717959E+0000}%
\special{sh 1}%
\special{ar 5400 1600 10 10 0  6.28318530717959E+0000}%
\special{sh 1}%
\special{ar 5600 1600 10 10 0  6.28318530717959E+0000}%
\special{sh 1}%
\special{ar 5800 1600 10 10 0  6.28318530717959E+0000}%
\special{sh 1}%
\special{ar 6000 1600 10 10 0  6.28318530717959E+0000}%
\special{sh 1}%
\special{ar 6000 1800 10 10 0  6.28318530717959E+0000}%
\special{sh 1}%
\special{ar 5800 2000 10 10 0  6.28318530717959E+0000}%
\special{sh 1}%
\special{ar 5600 2200 10 10 0  6.28318530717959E+0000}%
\special{sh 1}%
\special{ar 5400 2400 10 10 0  6.28318530717959E+0000}%
\special{sh 1}%
\special{ar 5400 2600 10 10 0  6.28318530717959E+0000}%
\special{sh 1}%
\special{ar 5600 2600 10 10 0  6.28318530717959E+0000}%
\special{sh 1}%
\special{ar 5800 2600 10 10 0  6.28318530717959E+0000}%
\special{sh 1}%
\special{ar 6000 2600 10 10 0  6.28318530717959E+0000}%
\special{sh 1}%
\special{ar 6200 2600 10 10 0  6.28318530717959E+0000}%
%
\special{pn 8}%
\special{pa 6200 2600}%
\special{pa 6000 1800}%
\special{pa 6000 1600}%
\special{pa 5200 1600}%
\special{pa 5400 2400}%
\special{pa 5400 2600}%
\special{pa 6200 2600}%
\special{fp}%
%
\special{pn 20}%
\special{sh 1}%
\special{ar 5400 2800 10 10 0  6.28318530717959E+0000}%
\special{sh 1}%
\special{ar 5600 2800 10 10 0  6.28318530717959E+0000}%
\special{sh 1}%
\special{ar 5800 2800 10 10 0  6.28318530717959E+0000}%
\special{sh 1}%
\special{ar 6000 2800 10 10 0  6.28318530717959E+0000}%
\special{sh 1}%
\special{ar 6200 2800 10 10 0  6.28318530717959E+0000}%
\special{sh 1}%
\special{ar 6200 3000 10 10 0  6.28318530717959E+0000}%
\special{sh 1}%
\special{ar 6000 3200 10 10 0  6.28318530717959E+0000}%
\special{sh 1}%
\special{ar 5800 3400 10 10 0  6.28318530717959E+0000}%
\special{sh 1}%
\special{ar 5600 3600 10 10 0  6.28318530717959E+0000}%
\special{sh 1}%
\special{ar 5400 3800 10 10 0  6.28318530717959E+0000}%
\special{sh 1}%
\special{ar 5400 4000 10 10 0  6.28318530717959E+0000}%
\special{sh 1}%
\special{ar 5600 4000 10 10 0  6.28318530717959E+0000}%
\special{sh 1}%
\special{ar 5800 4000 10 10 0  6.28318530717959E+0000}%
\special{sh 1}%
\special{ar 6000 4000 10 10 0  6.28318530717959E+0000}%
\special{sh 1}%
\special{ar 6200 4000 10 10 0  6.28318530717959E+0000}%
%
\special{pn 8}%
\special{pa 5400 2800}%
\special{pa 5400 4000}%
\special{pa 6200 4000}%
\special{pa 6200 2800}%
\special{pa 5400 2800}%
\special{fp}%
\end{picture}%

%% file: blockboun.tex
\unitlength 0.1in
\begin{picture}( 28.8000, 18.0000)(  1.2000,-23.1500)
%
\special{pn 8}%
\special{pa 1800 2000}%
\special{pa 2600 2000}%
\special{pa 2400 1200}%
\special{pa 2400 1000}%
\special{pa 1400 1000}%
\special{pa 1800 1800}%
\special{pa 1800 2000}%
\special{pa 1800 2000}%
\special{fp}%
%
\special{pn 20}%
\special{sh 1}%
\special{ar 1400 1000 10 10 0  6.28318530717959E+0000}%
\special{sh 1}%
\special{ar 1600 1000 10 10 0  6.28318530717959E+0000}%
\special{sh 1}%
\special{ar 1800 1000 10 10 0  6.28318530717959E+0000}%
\special{sh 1}%
\special{ar 2000 1000 10 10 0  6.28318530717959E+0000}%
\special{sh 1}%
\special{ar 2200 1000 10 10 0  6.28318530717959E+0000}%
\special{sh 1}%
\special{ar 2400 1000 10 10 0  6.28318530717959E+0000}%
\special{sh 1}%
\special{ar 2400 1200 10 10 0  6.28318530717959E+0000}%
\special{sh 1}%
\special{ar 2200 1400 10 10 0  6.28318530717959E+0000}%
\special{sh 1}%
\special{ar 2000 1600 10 10 0  6.28318530717959E+0000}%
\special{sh 1}%
\special{ar 1800 1800 10 10 0  6.28318530717959E+0000}%
\special{sh 1}%
\special{ar 1800 2000 10 10 0  6.28318530717959E+0000}%
\special{sh 1}%
\special{ar 2000 2000 10 10 0  6.28318530717959E+0000}%
\special{sh 1}%
\special{ar 2200 2000 10 10 0  6.28318530717959E+0000}%
\special{sh 1}%
\special{ar 2400 2000 10 10 0  6.28318530717959E+0000}%
\special{sh 1}%
\special{ar 2600 2000 10 10 0  6.28318530717959E+0000}%
\put(12.0000,-8.0000){\makebox(0,0){$\partial^2_{lt}B_{i,j}$}}%
\put(26.0000,-8.0000){\makebox(0,0){$\partial^2_{rt}B_{i,j}$}}%
\put(28.0000,-22.0000){\makebox(0,0){$\partial^2_{rb}B_{i,j}$}}%
\put(16.0000,-22.0000){\makebox(0,0){$\partial^2_{lb}B_{i,j}$}}%
\put(19.0000,-6.0000){\makebox(0,0){$\partial_tB_{i,j}$}}%
\put(22.0000,-24.0000){\makebox(0,0){$\partial_bB_{i,j}$}}%
\put(18.0000,-14.0000){\makebox(0,0){$B_{i,j}$}}%
%
\special{pn 20}%
\special{sh 1}%
\special{ar 1200 1000 10 10 0  6.28318530717959E+0000}%
\special{sh 1}%
\special{ar 1000 1000 10 10 0  6.28318530717959E+0000}%
\special{sh 1}%
\special{ar 2600 1000 10 10 0  6.28318530717959E+0000}%
\special{sh 1}%
\special{ar 2800 1000 10 10 0  6.28318530717959E+0000}%
\special{sh 1}%
\special{ar 2800 2000 10 10 0  6.28318530717959E+0000}%
\special{sh 1}%
\special{ar 3000 2000 10 10 0  6.28318530717959E+0000}%
\special{sh 1}%
\special{ar 1600 2000 10 10 0  6.28318530717959E+0000}%
\special{sh 1}%
\special{ar 1400 2000 10 10 0  6.28318530717959E+0000}%
\special{sh 1}%
\special{ar 1400 2000 10 10 0  6.28318530717959E+0000}%
%
\special{pn 8}%
\special{ar 2400 1000 100 100  0.0000000 6.2831853}%
%
\special{pn 8}%
\special{ar 1000 1000 100 100  0.0000000 6.2831853}%
%
\special{pn 8}%
\special{ar 2600 2000 100 100  0.0000000 6.2831853}%
%
\special{pn 8}%
\special{ar 1400 2000 100 100  0.0000000 6.2831853}%
\end{picture}%

%% file: tau+1.tex
\unitlength 0.1in
\begin{picture}( 61.7000, 20.0000)(  4.4000,-22.0000)
\put(16.0000,-22.0000){\makebox(0,0){$2$}}%
\put(18.0000,-22.0000){\makebox(0,0){$3$}}%
\put(8.0000,-18.0000){\makebox(0,0){$0$}}%
\put(8.0000,-16.0000){\makebox(0,0){$1$}}%
\put(8.0000,-14.0000){\makebox(0,0){$2$}}%
\put(14.0000,-22.0000){\makebox(0,0){$1$}}%
\put(12.0000,-22.0000){\makebox(0,0){$0$}}%
\put(9.6000,-2.9000){\makebox(0,0){$j$}}%
\put(34.0000,-22.0000){\makebox(0,0){$i$}}%
\put(22.0000,-22.0000){\makebox(0,0){$\tau+1$}}%
%
\special{pn 8}%
\special{pa 1000 2000}%
\special{pa 3410 2000}%
\special{fp}%
\special{sh 1}%
\special{pa 3410 2000}%
\special{pa 3344 1980}%
\special{pa 3358 2000}%
\special{pa 3344 2020}%
\special{pa 3410 2000}%
\special{fp}%
%
\special{pn 8}%
\special{pa 4600 2200}%
\special{pa 4600 410}%
\special{fp}%
\special{sh 1}%
\special{pa 4600 410}%
\special{pa 4580 478}%
\special{pa 4600 464}%
\special{pa 4620 478}%
\special{pa 4600 410}%
\special{fp}%
\put(8.0000,-8.0000){\makebox(0,0){$\tau+1$}}%
\put(8.0000,-6.0000){\makebox(0,0){$\tau+2$}}%
%
\special{pn 20}%
\special{sh 1}%
\special{ar 1200 1800 10 10 0  6.28318530717959E+0000}%
\special{sh 1}%
\special{ar 1200 1600 10 10 0  6.28318530717959E+0000}%
%
\special{pn 20}%
\special{sh 1}%
\special{ar 1600 1800 10 10 0  6.28318530717959E+0000}%
\special{sh 1}%
\special{ar 1600 1600 10 10 0  6.28318530717959E+0000}%
\special{sh 1}%
\special{ar 1800 1400 10 10 0  6.28318530717959E+0000}%
\special{sh 1}%
\special{ar 2000 1200 10 10 0  6.28318530717959E+0000}%
\special{sh 1}%
\special{ar 2200 1000 10 10 0  6.28318530717959E+0000}%
\special{sh 1}%
\special{ar 2400 800 10 10 0  6.28318530717959E+0000}%
%
\special{pn 20}%
\special{sh 1}%
\special{ar 1800 1800 10 10 0  6.28318530717959E+0000}%
\special{sh 1}%
\special{ar 2000 1800 10 10 0  6.28318530717959E+0000}%
\put(22.0000,-18.0000){\makebox(0,0){$0$}}%
%
\special{pn 8}%
\special{pa 2300 700}%
\special{pa 2500 700}%
\special{pa 2500 1700}%
\special{pa 2300 1700}%
\special{pa 2300 700}%
\special{fp}%
%
\special{pn 8}%
\special{pa 1200 600}%
\special{pa 1000 600}%
\special{dt 0.045}%
\put(14.0000,-18.0000){\makebox(0,0){$0$}}%
\put(14.0000,-16.0000){\makebox(0,0){$0$}}%
\put(14.0000,-14.0000){\makebox(0,0){$0$}}%
\put(16.0000,-14.0000){\makebox(0,0){$0$}}%
\put(16.0000,-12.0000){\makebox(0,0){$0$}}%
\put(16.0000,-10.0000){\makebox(0,0){$0$}}%
\put(22.0000,-16.0000){\makebox(0,0){$0$}}%
\put(22.0000,-14.0000){\makebox(0,0){$0$}}%
\put(24.0000,-14.0000){\makebox(0,0){$0$}}%
\put(24.0000,-12.0000){\makebox(0,0){$0$}}%
\put(24.0000,-10.0000){\makebox(0,0){$0$}}%
\put(26.0000,-10.0000){\makebox(0,0){$0$}}%
\put(26.0000,-8.0000){\makebox(0,0){$0$}}%
\put(26.0000,-6.0000){\makebox(0,0){$0$}}%
\put(22.0000,-4.0000){\makebox(0,0){$B_{0,1}$}}%
\put(30.0000,-4.0000){\makebox(0,0){$B_{1,1}$}}%
\put(48.0000,-21.9500){\makebox(0,0){$1$}}%
\put(50.0000,-21.9500){\makebox(0,0){$2$}}%
\put(40.0000,-17.9500){\makebox(0,0){$0$}}%
\put(40.0000,-15.9500){\makebox(0,0){$1$}}%
\put(40.0000,-13.9500){\makebox(0,0){$2$}}%
\put(46.0000,-21.9500){\makebox(0,0){$0$}}%
\put(44.0000,-21.9500){\makebox(0,0){$-1$}}%
\put(39.6000,-2.8500){\makebox(0,0){$j$}}%
\put(66.0000,-21.9500){\makebox(0,0){$i$}}%
\put(54.0000,-21.9500){\makebox(0,0){$\tau$}}%
%
\special{pn 8}%
\special{pa 4200 1996}%
\special{pa 6610 1996}%
\special{fp}%
\special{sh 1}%
\special{pa 6610 1996}%
\special{pa 6544 1976}%
\special{pa 6558 1996}%
\special{pa 6544 2016}%
\special{pa 6610 1996}%
\special{fp}%
\put(40.0000,-7.9500){\makebox(0,0){$\tau+1$}}%
\put(40.0000,-5.9500){\makebox(0,0){$\tau+2$}}%
%
\special{pn 20}%
\special{sh 1}%
\special{ar 4400 1796 10 10 0  6.28318530717959E+0000}%
\special{sh 1}%
\special{ar 4400 1596 10 10 0  6.28318530717959E+0000}%
%
\special{pn 20}%
\special{sh 1}%
\special{ar 4800 1796 10 10 0  6.28318530717959E+0000}%
\special{sh 1}%
\special{ar 4800 1596 10 10 0  6.28318530717959E+0000}%
\special{sh 1}%
\special{ar 5000 1396 10 10 0  6.28318530717959E+0000}%
\special{sh 1}%
\special{ar 5200 1196 10 10 0  6.28318530717959E+0000}%
\special{sh 1}%
\special{ar 5400 996 10 10 0  6.28318530717959E+0000}%
\special{sh 1}%
\special{ar 5600 796 10 10 0  6.28318530717959E+0000}%
%
\special{pn 20}%
\special{sh 1}%
\special{ar 5000 1796 10 10 0  6.28318530717959E+0000}%
\special{sh 1}%
\special{ar 5200 1796 10 10 0  6.28318530717959E+0000}%
\put(54.0000,-17.9500){\makebox(0,0){$0$}}%
%
\special{pn 20}%
\special{sh 1}%
\special{ar 5600 596 10 10 0  6.28318530717959E+0000}%
\special{sh 1}%
\special{ar 5400 596 10 10 0  6.28318530717959E+0000}%
\special{sh 1}%
\special{ar 5200 596 10 10 0  6.28318530717959E+0000}%
%
\special{pn 8}%
\special{pa 5500 696}%
\special{pa 5700 696}%
\special{pa 5700 1696}%
\special{pa 5500 1696}%
\special{pa 5500 696}%
\special{fp}%
%
\special{pn 8}%
\special{pa 4400 596}%
\special{pa 4200 596}%
\special{dt 0.045}%
\put(46.0000,-17.9500){\makebox(0,0){$0$}}%
\put(46.0000,-15.9500){\makebox(0,0){$0$}}%
\put(46.0000,-13.9500){\makebox(0,0){$0$}}%
\put(48.0000,-13.9500){\makebox(0,0){$0$}}%
\put(48.0000,-11.9500){\makebox(0,0){$0$}}%
\put(48.0000,-9.9500){\makebox(0,0){$0$}}%
\put(54.0000,-15.9500){\makebox(0,0){$0$}}%
\put(54.0000,-13.9500){\makebox(0,0){$0$}}%
\put(56.0000,-13.9500){\makebox(0,0){$0$}}%
\put(56.0000,-11.9500){\makebox(0,0){$0$}}%
\put(56.0000,-9.9500){\makebox(0,0){$0$}}%
\put(58.0000,-9.9500){\makebox(0,0){$0$}}%
\put(58.0000,-7.9500){\makebox(0,0){$0$}}%
\put(58.0000,-5.9500){\makebox(0,0){$0$}}%
\put(54.0000,-3.9500){\makebox(0,0){$B_{0,1}$}}%
\put(62.0000,-3.9500){\makebox(0,0){$B_{1,1}$}}%
%
\special{pn 20}%
\special{sh 1}%
\special{ar 1000 1600 10 10 0  6.28318530717959E+0000}%
\special{sh 1}%
\special{ar 1000 1200 10 10 0  6.28318530717959E+0000}%
\special{sh 1}%
\special{ar 1200 1000 10 10 0  6.28318530717959E+0000}%
\special{sh 1}%
\special{ar 1400 800 10 10 0  6.28318530717959E+0000}%
%
\special{pn 20}%
\special{sh 1}%
\special{ar 1420 600 10 10 0  6.28318530717959E+0000}%
\special{sh 1}%
\special{ar 1220 600 10 10 0  6.28318530717959E+0000}%
%
\special{pn 20}%
\special{sh 1}%
\special{ar 2400 600 10 10 0  6.28318530717959E+0000}%
\special{sh 1}%
\special{ar 2200 600 10 10 0  6.28318530717959E+0000}%
\special{sh 1}%
\special{ar 2000 600 10 10 0  6.28318530717959E+0000}%
\special{sh 1}%
\special{ar 1800 600 10 10 0  6.28318530717959E+0000}%
%
\special{pn 8}%
\special{pa 2000 1800}%
\special{pa 2400 800}%
\special{pa 2400 600}%
\special{pa 1800 600}%
\special{pa 1600 1600}%
\special{pa 1600 1800}%
\special{pa 2000 1800}%
\special{fp}%
\put(16.0000,-8.0000){\makebox(0,0){$0$}}%
\put(16.0000,-6.0000){\makebox(0,0){$0$}}%
%
\special{pn 20}%
\special{sh 1}%
\special{ar 2800 600 10 10 0  6.28318530717959E+0000}%
\special{sh 1}%
\special{ar 3000 600 10 10 0  6.28318530717959E+0000}%
\special{sh 1}%
\special{ar 3200 600 10 10 0  6.28318530717959E+0000}%
\special{sh 1}%
\special{ar 3200 800 10 10 0  6.28318530717959E+0000}%
\special{sh 1}%
\special{ar 3000 1000 10 10 0  6.28318530717959E+0000}%
\special{sh 1}%
\special{ar 2800 1200 10 10 0  6.28318530717959E+0000}%
\special{sh 1}%
\special{ar 2600 1400 10 10 0  6.28318530717959E+0000}%
\special{sh 1}%
\special{ar 2400 1600 10 10 0  6.28318530717959E+0000}%
\special{sh 1}%
\special{ar 2400 1800 10 10 0  6.28318530717959E+0000}%
\special{sh 1}%
\special{ar 2600 1800 10 10 0  6.28318530717959E+0000}%
\special{sh 1}%
\special{ar 2800 1800 10 10 0  6.28318530717959E+0000}%
\special{sh 1}%
\special{ar 3000 1800 10 10 0  6.28318530717959E+0000}%
%
\special{pn 8}%
\special{pa 3000 1800}%
\special{pa 3200 800}%
\special{pa 3200 600}%
\special{pa 2800 600}%
\special{pa 2400 1600}%
\special{pa 2400 1800}%
\special{pa 3000 1800}%
\special{fp}%
%
\special{pn 20}%
\special{sh 1}%
\special{ar 5600 1800 10 10 0  6.28318530717959E+0000}%
\special{sh 1}%
\special{ar 5600 1600 10 10 0  6.28318530717959E+0000}%
\special{sh 1}%
\special{ar 5800 1400 10 10 0  6.28318530717959E+0000}%
\special{sh 1}%
\special{ar 6000 1200 10 10 0  6.28318530717959E+0000}%
\special{sh 1}%
\special{ar 6200 1000 10 10 0  6.28318530717959E+0000}%
\special{sh 1}%
\special{ar 6400 800 10 10 0  6.28318530717959E+0000}%
\special{sh 1}%
\special{ar 6400 600 10 10 0  6.28318530717959E+0000}%
\special{sh 1}%
\special{ar 6200 600 10 10 0  6.28318530717959E+0000}%
\special{sh 1}%
\special{ar 6000 600 10 10 0  6.28318530717959E+0000}%
\special{sh 1}%
\special{ar 5800 1800 10 10 0  6.28318530717959E+0000}%
\special{sh 1}%
\special{ar 6000 1800 10 10 0  6.28318530717959E+0000}%
\special{sh 1}%
\special{ar 6200 1800 10 10 0  6.28318530717959E+0000}%
%
\special{pn 20}%
\special{sh 1}%
\special{ar 4200 1200 10 10 0  6.28318530717959E+0000}%
\special{sh 1}%
\special{ar 4400 1000 10 10 0  6.28318530717959E+0000}%
\special{sh 1}%
\special{ar 4600 800 10 10 0  6.28318530717959E+0000}%
\special{sh 1}%
\special{ar 4600 600 10 10 0  6.28318530717959E+0000}%
\special{sh 1}%
\special{ar 4400 600 10 10 0  6.28318530717959E+0000}%
\special{sh 1}%
\special{ar 4200 600 10 10 0  6.28318530717959E+0000}%
%
\special{pn 20}%
\special{sh 1}%
\special{ar 5000 600 10 10 0  6.28318530717959E+0000}%
%
\special{pn 8}%
\special{pa 5000 600}%
\special{pa 4800 1600}%
\special{pa 4800 1800}%
\special{pa 5200 1800}%
\special{pa 5600 800}%
\special{pa 5600 600}%
\special{pa 5000 600}%
\special{fp}%
%
\special{pn 8}%
\special{pa 6000 600}%
\special{pa 5600 1600}%
\special{pa 5600 1800}%
\special{pa 6200 1800}%
\special{pa 6400 800}%
\special{pa 6400 600}%
\special{pa 6000 600}%
\special{fp}%
\put(48.0000,-8.0000){\makebox(0,0){$0$}}%
\put(48.0000,-6.0000){\makebox(0,0){$0$}}%
%
\special{pn 8}%
\special{pa 1000 600}%
\special{pa 1400 600}%
\special{fp}%
\special{pa 1400 600}%
\special{pa 1400 800}%
\special{fp}%
\special{pa 1400 800}%
\special{pa 1200 1600}%
\special{fp}%
%
\special{pn 8}%
\special{pa 4200 600}%
\special{pa 4600 600}%
\special{fp}%
\special{pa 4600 600}%
\special{pa 4600 800}%
\special{fp}%
\special{pa 4600 800}%
\special{pa 4400 1600}%
\special{fp}%
\special{pa 4400 1600}%
\special{pa 4200 1600}%
\special{fp}%
%
\special{pn 8}%
\special{pa 1200 2200}%
\special{pa 1200 410}%
\special{fp}%
\special{sh 1}%
\special{pa 1200 410}%
\special{pa 1180 478}%
\special{pa 1200 464}%
\special{pa 1220 478}%
\special{pa 1200 410}%
\special{fp}%
%
\special{pn 8}%
\special{pa 1200 1600}%
\special{pa 1000 1600}%
\special{fp}%
\end{picture}%

%% file: tau+1sym.tex
\unitlength 0.1in
\begin{picture}( 37.6000, 19.1000)(  4.4000,-21.1000)
\put(16.0000,-21.9500){\makebox(0,0){$1$}}%
\put(18.0000,-21.9500){\makebox(0,0){$2$}}%
\put(8.0000,-17.9500){\makebox(0,0){$0$}}%
\put(8.0000,-15.9500){\makebox(0,0){$1$}}%
\put(8.0000,-13.9500){\makebox(0,0){$2$}}%
\put(14.0000,-21.9500){\makebox(0,0){$0$}}%
\put(12.0000,-21.9500){\makebox(0,0){$-1$}}%
\put(7.6000,-2.8500){\makebox(0,0){$j$}}%
\put(34.0000,-21.9500){\makebox(0,0){$i$}}%
\put(22.0000,-21.9500){\makebox(0,0){$\tau$}}%
%
\special{pn 8}%
\special{pa 1000 1996}%
\special{pa 3410 1996}%
\special{fp}%
\special{sh 1}%
\special{pa 3410 1996}%
\special{pa 3344 1976}%
\special{pa 3358 1996}%
\special{pa 3344 2016}%
\special{pa 3410 1996}%
\special{fp}%
%
\special{pn 8}%
\special{pa 1400 2000}%
\special{pa 1400 252}%
\special{fp}%
\special{sh 1}%
\special{pa 1400 252}%
\special{pa 1380 320}%
\special{pa 1400 306}%
\special{pa 1420 320}%
\special{pa 1400 252}%
\special{fp}%
\put(8.0000,-7.9500){\makebox(0,0){$\tau+1$}}%
\put(8.0000,-5.9500){\makebox(0,0){$\tau+2$}}%
%
\special{pn 20}%
\special{sh 1}%
\special{ar 1200 1796 10 10 0  6.28318530717959E+0000}%
\special{sh 1}%
\special{ar 1200 1596 10 10 0  6.28318530717959E+0000}%
%
\special{pn 20}%
\special{sh 1}%
\special{ar 1600 1796 10 10 0  6.28318530717959E+0000}%
\special{sh 1}%
\special{ar 1600 1596 10 10 0  6.28318530717959E+0000}%
\special{sh 1}%
\special{ar 1800 1396 10 10 0  6.28318530717959E+0000}%
\special{sh 1}%
\special{ar 2000 1196 10 10 0  6.28318530717959E+0000}%
\special{sh 1}%
\special{ar 2200 996 10 10 0  6.28318530717959E+0000}%
\special{sh 1}%
\special{ar 2400 796 10 10 0  6.28318530717959E+0000}%
%
\special{pn 20}%
\special{sh 1}%
\special{ar 1800 1796 10 10 0  6.28318530717959E+0000}%
\special{sh 1}%
\special{ar 2000 1796 10 10 0  6.28318530717959E+0000}%
\put(22.0000,-17.9500){\makebox(0,0){$0$}}%
%
\special{pn 20}%
\special{sh 1}%
\special{ar 2400 596 10 10 0  6.28318530717959E+0000}%
\special{sh 1}%
\special{ar 2200 596 10 10 0  6.28318530717959E+0000}%
\special{sh 1}%
\special{ar 2000 596 10 10 0  6.28318530717959E+0000}%
%
\special{pn 8}%
\special{pa 2300 696}%
\special{pa 2500 696}%
\special{pa 2500 1696}%
\special{pa 2300 1696}%
\special{pa 2300 696}%
\special{fp}%
%
\special{pn 8}%
\special{pa 1200 596}%
\special{pa 1000 596}%
\special{dt 0.045}%
\put(14.0000,-17.9500){\makebox(0,0){$0$}}%
\put(14.0000,-15.9500){\makebox(0,0){$0$}}%
\put(14.0000,-13.9500){\makebox(0,0){$0$}}%
\put(16.0000,-13.9500){\makebox(0,0){$0$}}%
\put(16.0000,-11.9500){\makebox(0,0){$0$}}%
\put(16.0000,-9.9500){\makebox(0,0){$0$}}%
\put(22.0000,-15.9500){\makebox(0,0){$0$}}%
\put(22.0000,-13.9500){\makebox(0,0){$0$}}%
\put(24.0000,-13.9500){\makebox(0,0){$0$}}%
\put(24.0000,-11.9500){\makebox(0,0){$0$}}%
\put(24.0000,-9.9500){\makebox(0,0){$0$}}%
\put(26.0000,-9.9500){\makebox(0,0){$0$}}%
\put(26.0000,-7.9500){\makebox(0,0){$0$}}%
\put(26.0000,-5.9500){\makebox(0,0){$0$}}%
\put(22.0000,-3.9500){\makebox(0,0){$B_{0,1}$}}%
\put(30.0000,-3.9500){\makebox(0,0){$B_{1,1}$}}%
%
\special{pn 20}%
\special{sh 1}%
\special{ar 2400 1800 10 10 0  6.28318530717959E+0000}%
\special{sh 1}%
\special{ar 2400 1600 10 10 0  6.28318530717959E+0000}%
\special{sh 1}%
\special{ar 2600 1400 10 10 0  6.28318530717959E+0000}%
\special{sh 1}%
\special{ar 2800 1200 10 10 0  6.28318530717959E+0000}%
\special{sh 1}%
\special{ar 3000 1000 10 10 0  6.28318530717959E+0000}%
\special{sh 1}%
\special{ar 3200 800 10 10 0  6.28318530717959E+0000}%
\special{sh 1}%
\special{ar 3200 600 10 10 0  6.28318530717959E+0000}%
\special{sh 1}%
\special{ar 3000 600 10 10 0  6.28318530717959E+0000}%
\special{sh 1}%
\special{ar 2800 600 10 10 0  6.28318530717959E+0000}%
\special{sh 1}%
\special{ar 2600 1800 10 10 0  6.28318530717959E+0000}%
\special{sh 1}%
\special{ar 2800 1800 10 10 0  6.28318530717959E+0000}%
\special{sh 1}%
\special{ar 3000 1800 10 10 0  6.28318530717959E+0000}%
%
\special{pn 20}%
\special{sh 1}%
\special{ar 1000 1200 10 10 0  6.28318530717959E+0000}%
\special{sh 1}%
\special{ar 1200 1000 10 10 0  6.28318530717959E+0000}%
\special{sh 1}%
\special{ar 1400 800 10 10 0  6.28318530717959E+0000}%
\special{sh 1}%
\special{ar 1400 600 10 10 0  6.28318530717959E+0000}%
\special{sh 1}%
\special{ar 1200 600 10 10 0  6.28318530717959E+0000}%
\special{sh 1}%
\special{ar 1000 600 10 10 0  6.28318530717959E+0000}%
%
\special{pn 20}%
\special{sh 1}%
\special{ar 1800 600 10 10 0  6.28318530717959E+0000}%
%
\special{pn 8}%
\special{pa 1800 600}%
\special{pa 1600 1600}%
\special{pa 1600 1800}%
\special{pa 2000 1800}%
\special{pa 2400 800}%
\special{pa 2400 600}%
\special{pa 1800 600}%
\special{fp}%
%
\special{pn 8}%
\special{pa 2800 600}%
\special{pa 2400 1600}%
\special{pa 2400 1800}%
\special{pa 3000 1800}%
\special{pa 3200 800}%
\special{pa 3200 600}%
\special{pa 2800 600}%
\special{fp}%
\put(16.0000,-8.0000){\makebox(0,0){$0$}}%
\put(16.0000,-6.0000){\makebox(0,0){$0$}}%
\put(32.0000,-18.0000){\makebox(0,0){$0$}}%
\put(32.0000,-16.0000){\makebox(0,0){$0$}}%
\put(32.0000,-14.0000){\makebox(0,0){$0$}}%
%
\special{pn 20}%
\special{sh 1}%
\special{ar 3400 1800 10 10 0  6.28318530717959E+0000}%
\special{sh 1}%
\special{ar 3400 1600 10 10 0  6.28318530717959E+0000}%
\special{sh 1}%
\special{ar 3600 1400 10 10 0  6.28318530717959E+0000}%
\special{sh 1}%
\special{ar 3800 1200 10 10 0  6.28318530717959E+0000}%
\special{sh 1}%
\special{ar 4000 1000 10 10 0  6.28318530717959E+0000}%
\special{sh 1}%
\special{ar 3600 1800 10 10 0  6.28318530717959E+0000}%
\special{sh 1}%
\special{ar 3800 1800 10 10 0  6.28318530717959E+0000}%
\special{sh 1}%
\special{ar 4000 800 10 10 0  6.28318530717959E+0000}%
\special{sh 1}%
\special{ar 3800 800 10 10 0  6.28318530717959E+0000}%
\special{sh 1}%
\special{ar 3600 800 10 10 0  6.28318530717959E+0000}%
%
\special{pn 8}%
\special{pa 3600 800}%
\special{pa 4000 800}%
\special{pa 4000 1000}%
\special{pa 3800 1800}%
\special{pa 3400 1800}%
\special{pa 3400 1600}%
\special{pa 3600 800}%
\special{fp}%
\put(32.0000,-12.0000){\makebox(0,0){$0$}}%
\put(34.0000,-10.0000){\makebox(0,0){$0$}}%
\put(32.0000,-10.0000){\makebox(0,0){$0$}}%
\put(34.0000,-8.0000){\makebox(0,0){$0$}}%
\put(34.0000,-6.0000){\makebox(0,0){$0$}}%
%
\special{pn 20}%
\special{sh 1}%
\special{ar 3600 600 10 10 0  6.28318530717959E+0000}%
\put(46.0000,-4.0000){\makebox(0,0){$\alpha(n_1)$}}%
%
\special{pn 8}%
\special{pa 3400 800}%
\special{pa 4200 400}%
\special{fp}%
%
\special{pn 8}%
\special{pa 1000 600}%
\special{pa 1400 600}%
\special{fp}%
\special{pa 1400 600}%
\special{pa 1400 800}%
\special{fp}%
\special{pa 1400 800}%
\special{pa 1200 1600}%
\special{fp}%
\special{pa 1200 1600}%
\special{pa 1000 1600}%
\special{fp}%
\end{picture}%

%% file: tau+1-1.tex
\unitlength 0.1in
\begin{picture}( 52.9000, 20.8500)(  6.0000,-22.0000)
\put(11.1000,-9.0000){\makebox(0,0){$0$}}%
\put(10.0000,-2.0000){\makebox(0,0){$j$}}%
\put(30.0000,-8.0000){\makebox(0,0){$i$}}%
\put(22.0000,-11.0000){\makebox(0,0){$\tau$}}%
%
\special{pn 8}%
\special{pa 600 1000}%
\special{pa 3090 1000}%
\special{fp}%
\special{sh 1}%
\special{pa 3090 1000}%
\special{pa 3024 980}%
\special{pa 3038 1000}%
\special{pa 3024 1020}%
\special{pa 3090 1000}%
\special{fp}%
%
\special{pn 8}%
\special{pa 1200 2196}%
\special{pa 1200 250}%
\special{fp}%
\special{sh 1}%
\special{pa 1200 250}%
\special{pa 1180 318}%
\special{pa 1200 304}%
\special{pa 1220 318}%
\special{pa 1200 250}%
\special{fp}%
\put(11.0000,-11.0000){\makebox(0,0){$-1$}}%
\put(11.0500,-14.0500){\makebox(0,0){$-2$}}%
\put(14.5000,-11.0000){\makebox(0,0){$0$}}%
%
\special{pn 8}%
\special{pa 2200 970}%
\special{pa 2200 1010}%
\special{fp}%
%
\special{pn 20}%
\special{sh 1}%
\special{ar 1806 1596 10 10 0  6.28318530717959E+0000}%
\special{sh 1}%
\special{ar 2006 1796 10 10 0  6.28318530717959E+0000}%
\special{sh 1}%
\special{ar 2006 1996 10 10 0  6.28318530717959E+0000}%
\special{sh 1}%
\special{ar 1806 1996 10 10 0  6.28318530717959E+0000}%
\special{sh 1}%
\special{ar 1606 1996 10 10 0  6.28318530717959E+0000}%
%
\special{pn 20}%
\special{sh 1}%
\special{ar 2406 1996 10 10 0  6.28318530717959E+0000}%
\special{sh 1}%
\special{ar 2606 1996 10 10 0  6.28318530717959E+0000}%
%
\special{pn 20}%
\special{sh 1}%
\special{ar 1006 1196 10 10 0  6.28318530717959E+0000}%
%
\special{pn 8}%
\special{pa 1806 1396}%
\special{pa 2406 1396}%
\special{pa 2006 1796}%
\special{pa 2006 1996}%
\special{pa 1606 1996}%
\special{pa 1806 1396}%
\special{fp}%
%
\special{pn 20}%
\special{sh 1}%
\special{ar 1200 1996 10 10 0  6.28318530717959E+0000}%
\special{sh 1}%
\special{ar 1200 1796 10 10 0  6.28318530717959E+0000}%
\special{sh 1}%
\special{ar 1000 1596 10 10 0  6.28318530717959E+0000}%
\special{sh 1}%
\special{ar 800 1396 10 10 0  6.28318530717959E+0000}%
\special{sh 1}%
\special{ar 800 1196 10 10 0  6.28318530717959E+0000}%
%
\special{pn 8}%
\special{pa 2306 1296}%
\special{pa 2506 1296}%
\special{pa 2506 2096}%
\special{pa 2306 2096}%
\special{pa 2306 1296}%
\special{fp}%
\put(39.0000,-9.0000){\makebox(0,0){$0$}}%
\put(38.0000,-2.0000){\makebox(0,0){$j$}}%
\put(58.0000,-8.0000){\makebox(0,0){$i$}}%
\put(50.0000,-11.0000){\makebox(0,0){$\tau+1$}}%
\put(39.0000,-11.0000){\makebox(0,0){$-1$}}%
\put(39.0000,-13.9500){\makebox(0,0){$-2$}}%
\put(43.0000,-11.0000){\makebox(0,0){$1$}}%
%
\special{pn 20}%
\special{sh 1}%
\special{ar 4600 1596 10 10 0  6.28318530717959E+0000}%
\special{sh 1}%
\special{ar 4800 1796 10 10 0  6.28318530717959E+0000}%
\special{sh 1}%
\special{ar 4800 1996 10 10 0  6.28318530717959E+0000}%
\special{sh 1}%
\special{ar 4600 1996 10 10 0  6.28318530717959E+0000}%
\special{sh 1}%
\special{ar 4400 1996 10 10 0  6.28318530717959E+0000}%
%
\special{pn 20}%
\special{sh 1}%
\special{ar 5200 1996 10 10 0  6.28318530717959E+0000}%
\special{sh 1}%
\special{ar 5400 1996 10 10 0  6.28318530717959E+0000}%
%
\special{pn 20}%
\special{sh 1}%
\special{ar 3800 1196 10 10 0  6.28318530717959E+0000}%
%
\special{pn 8}%
\special{pa 4600 1396}%
\special{pa 5200 1396}%
\special{pa 4800 1796}%
\special{pa 4800 1996}%
\special{pa 4400 1996}%
\special{pa 4600 1396}%
\special{fp}%
%
\special{pn 20}%
\special{sh 1}%
\special{ar 4000 1996 10 10 0  6.28318530717959E+0000}%
\special{sh 1}%
\special{ar 4000 1796 10 10 0  6.28318530717959E+0000}%
\special{sh 1}%
\special{ar 3800 1596 10 10 0  6.28318530717959E+0000}%
\special{sh 1}%
\special{ar 3600 1396 10 10 0  6.28318530717959E+0000}%
\special{sh 1}%
\special{ar 3600 1196 10 10 0  6.28318530717959E+0000}%
%
\special{pn 8}%
\special{pa 5100 1296}%
\special{pa 5300 1296}%
\special{pa 5300 2096}%
\special{pa 5100 2096}%
\special{pa 5100 1296}%
\special{fp}%
%
\special{pn 20}%
\special{sh 1}%
\special{ar 2406 1396 10 10 0  6.28318530717959E+0000}%
\special{sh 1}%
\special{ar 2206 1396 10 10 0  6.28318530717959E+0000}%
\special{sh 1}%
\special{ar 2006 1396 10 10 0  6.28318530717959E+0000}%
\special{sh 1}%
\special{ar 1806 1396 10 10 0  6.28318530717959E+0000}%
\special{sh 1}%
\special{ar 1406 1196 10 10 0  6.28318530717959E+0000}%
\special{sh 1}%
\special{ar 1206 1196 10 10 0  6.28318530717959E+0000}%
%
\special{pn 20}%
\special{sh 1}%
\special{ar 4000 1196 10 10 0  6.28318530717959E+0000}%
\special{sh 1}%
\special{ar 4200 1196 10 10 0  6.28318530717959E+0000}%
\special{sh 1}%
\special{ar 4600 1396 10 10 0  6.28318530717959E+0000}%
\special{sh 1}%
\special{ar 4800 1396 10 10 0  6.28318530717959E+0000}%
\special{sh 1}%
\special{ar 5000 1396 10 10 0  6.28318530717959E+0000}%
\special{sh 1}%
\special{ar 5200 1396 10 10 0  6.28318530717959E+0000}%
%
\special{pn 20}%
\special{sh 1}%
\special{ar 1006 1996 10 10 0  6.28318530717959E+0000}%
%
\special{pn 20}%
\special{sh 1}%
\special{ar 3600 1996 10 10 0  6.28318530717959E+0000}%
%
\special{pn 20}%
\special{sh 1}%
\special{ar 3800 1996 10 10 0  6.28318530717959E+0000}%
%
\special{pn 20}%
\special{sh 1}%
\special{ar 806 1996 10 10 0  6.28318530717959E+0000}%
%
\special{pn 8}%
\special{pa 3600 1196}%
\special{pa 4200 1196}%
\special{pa 4000 1796}%
\special{pa 4000 1996}%
\special{pa 3600 1996}%
\special{pa 3600 1196}%
\special{fp}%
%
\special{pn 8}%
\special{pa 1400 1196}%
\special{pa 1200 1796}%
\special{pa 1200 1996}%
\special{pa 800 1996}%
\special{pa 800 1196}%
\special{pa 1400 1196}%
\special{fp}%
\put(39.0000,-17.9500){\makebox(0,0){$-\tau$}}%
\put(11.0500,-17.9500){\makebox(0,0){$-\tau$}}%
%
\special{pn 20}%
\special{sh 1}%
\special{ar 2806 1396 10 10 0  6.28318530717959E+0000}%
\special{sh 1}%
\special{ar 2806 1596 10 10 0  6.28318530717959E+0000}%
\special{sh 1}%
\special{ar 3006 1796 10 10 0  6.28318530717959E+0000}%
\special{sh 1}%
\special{ar 3006 1996 10 10 0  6.28318530717959E+0000}%
\special{sh 1}%
\special{ar 2806 1996 10 10 0  6.28318530717959E+0000}%
%
\special{pn 20}%
\special{sh 1}%
\special{ar 5600 1396 10 10 0  6.28318530717959E+0000}%
\special{sh 1}%
\special{ar 5600 1596 10 10 0  6.28318530717959E+0000}%
\special{sh 1}%
\special{ar 5800 1796 10 10 0  6.28318530717959E+0000}%
\special{sh 1}%
\special{ar 5800 1996 10 10 0  6.28318530717959E+0000}%
\special{sh 1}%
\special{ar 5600 1996 10 10 0  6.28318530717959E+0000}%
\put(54.0000,-13.9500){\makebox(0,0){$0$}}%
\put(26.0500,-13.9500){\makebox(0,0){$0$}}%
\put(22.0500,-19.9500){\makebox(0,0){$0$}}%
\put(14.0000,-19.9500){\makebox(0,0){$0$}}%
\put(50.0000,-19.9500){\makebox(0,0){$0$}}%
\put(42.0000,-19.9500){\makebox(0,0){$0$}}%
%
\special{pn 20}%
\special{sh 1}%
\special{ar 3000 1396 10 10 0  6.28318530717959E+0000}%
\special{sh 1}%
\special{ar 5800 1396 10 10 0  6.28318530717959E+0000}%
%
\special{pn 20}%
\special{sh 1}%
\special{ar 5200 1200 10 10 0  6.28318530717959E+0000}%
%
\special{pn 20}%
\special{sh 1}%
\special{ar 5000 1000 10 10 0  6.28318530717959E+0000}%
%
\special{pn 20}%
\special{sh 1}%
\special{ar 4200 1000 10 10 0  6.28318530717959E+0000}%
%
\special{pn 20}%
\special{sh 1}%
\special{ar 1400 1000 10 10 0  6.28318530717959E+0000}%
%
\special{pn 20}%
\special{sh 1}%
\special{ar 2204 1004 10 10 0  6.28318530717959E+0000}%
%
\special{pn 20}%
\special{sh 1}%
\special{ar 2404 1204 10 10 0  6.28318530717959E+0000}%
%
\special{pn 8}%
\special{pa 3000 1396}%
\special{pa 2800 1396}%
\special{fp}%
\special{pa 2800 1396}%
\special{pa 2400 1996}%
\special{fp}%
\special{pa 2400 1996}%
\special{pa 3000 1996}%
\special{fp}%
%
\special{pn 8}%
\special{pa 5800 1396}%
\special{pa 5600 1396}%
\special{fp}%
\special{pa 5200 1996}%
\special{pa 5600 1396}%
\special{fp}%
\special{pa 5200 1996}%
\special{pa 5800 1996}%
\special{fp}%
%
\special{pn 20}%
\special{sh 1}%
\special{ar 2000 800 10 10 0  6.28318530717959E+0000}%
\special{sh 1}%
\special{ar 1800 600 10 10 0  6.28318530717959E+0000}%
\special{sh 1}%
\special{ar 1800 400 10 10 0  6.28318530717959E+0000}%
\special{sh 1}%
\special{ar 1400 400 10 10 0  6.28318530717959E+0000}%
\special{sh 1}%
\special{ar 1200 400 10 10 0  6.28318530717959E+0000}%
\special{sh 1}%
\special{ar 1200 800 10 10 0  6.28318530717959E+0000}%
\special{sh 1}%
\special{ar 1000 600 10 10 0  6.28318530717959E+0000}%
\special{sh 1}%
\special{ar 1000 400 10 10 0  6.28318530717959E+0000}%
\special{sh 1}%
\special{ar 2000 400 10 10 0  6.28318530717959E+0000}%
\special{sh 1}%
\special{ar 2200 400 10 10 0  6.28318530717959E+0000}%
%
\special{pn 8}%
\special{pa 1800 400}%
\special{pa 2200 400}%
\special{pa 2400 1200}%
\special{pa 2400 1400}%
\special{pa 1800 1400}%
\special{pa 1800 400}%
\special{fp}%
%
\special{pn 8}%
\special{pa 800 1200}%
\special{pa 1400 1200}%
\special{pa 1400 400}%
\special{pa 1000 400}%
\special{pa 800 1200}%
\special{fp}%
%
\special{pn 8}%
\special{pa 3600 1200}%
\special{pa 4200 1200}%
\special{pa 4200 400}%
\special{pa 3800 400}%
\special{pa 3600 1200}%
\special{fp}%
%
\special{pn 8}%
\special{pa 4600 400}%
\special{pa 5000 400}%
\special{pa 5200 1200}%
\special{pa 5200 1400}%
\special{pa 4600 1400}%
\special{pa 4600 400}%
\special{fp}%
%
\special{pn 20}%
\special{sh 1}%
\special{ar 4800 800 10 10 0  6.28318530717959E+0000}%
\special{sh 1}%
\special{ar 4600 600 10 10 0  6.28318530717959E+0000}%
\special{sh 1}%
\special{ar 4600 400 10 10 0  6.28318530717959E+0000}%
\special{sh 1}%
\special{ar 4200 400 10 10 0  6.28318530717959E+0000}%
\special{sh 1}%
\special{ar 4000 400 10 10 0  6.28318530717959E+0000}%
\special{sh 1}%
\special{ar 4000 800 10 10 0  6.28318530717959E+0000}%
\special{sh 1}%
\special{ar 3800 600 10 10 0  6.28318530717959E+0000}%
\special{sh 1}%
\special{ar 3800 400 10 10 0  6.28318530717959E+0000}%
\special{sh 1}%
\special{ar 4800 400 10 10 0  6.28318530717959E+0000}%
\special{sh 1}%
\special{ar 5000 400 10 10 0  6.28318530717959E+0000}%
%
\special{pn 8}%
\special{pa 4000 2200}%
\special{pa 4000 256}%
\special{fp}%
\special{sh 1}%
\special{pa 4000 256}%
\special{pa 3980 322}%
\special{pa 4000 308}%
\special{pa 4020 322}%
\special{pa 4000 256}%
\special{fp}%
%
\special{pn 8}%
\special{pa 3400 1000}%
\special{pa 5890 1000}%
\special{fp}%
\special{sh 1}%
\special{pa 5890 1000}%
\special{pa 5824 980}%
\special{pa 5838 1000}%
\special{pa 5824 1020}%
\special{pa 5890 1000}%
\special{fp}%
%
\special{pn 8}%
\special{pa 5000 960}%
\special{pa 5000 1000}%
\special{fp}%
\end{picture}%

%% file: tau+1-1sym.tex
\unitlength 0.1in
\begin{picture}( 35.8000, 22.0500)( 26.0000,-27.1500)
\put(31.0000,-12.9500){\makebox(0,0){$0$}}%
\put(30.0000,-5.9500){\makebox(0,0){$j$}}%
\put(60.0000,-12.0000){\makebox(0,0){$i$}}%
\put(42.0000,-14.9500){\makebox(0,0){$\tau+1$}}%
\put(31.0000,-14.9500){\makebox(0,0){$-1$}}%
\put(31.0000,-17.9000){\makebox(0,0){$-2$}}%
\put(35.0000,-14.9500){\makebox(0,0){$1$}}%
%
\special{pn 20}%
\special{sh 1}%
\special{ar 3800 2000 10 10 0  6.28318530717959E+0000}%
\special{sh 1}%
\special{ar 4000 2200 10 10 0  6.28318530717959E+0000}%
\special{sh 1}%
\special{ar 4000 2400 10 10 0  6.28318530717959E+0000}%
\special{sh 1}%
\special{ar 3800 2400 10 10 0  6.28318530717959E+0000}%
\special{sh 1}%
\special{ar 3600 2400 10 10 0  6.28318530717959E+0000}%
%
\special{pn 20}%
\special{sh 1}%
\special{ar 4400 2400 10 10 0  6.28318530717959E+0000}%
\special{sh 1}%
\special{ar 4600 2400 10 10 0  6.28318530717959E+0000}%
%
\special{pn 20}%
\special{sh 1}%
\special{ar 3000 1600 10 10 0  6.28318530717959E+0000}%
%
\special{pn 20}%
\special{sh 1}%
\special{ar 3200 2400 10 10 0  6.28318530717959E+0000}%
\special{sh 1}%
\special{ar 3200 2200 10 10 0  6.28318530717959E+0000}%
\special{sh 1}%
\special{ar 3000 2000 10 10 0  6.28318530717959E+0000}%
\special{sh 1}%
\special{ar 2800 1800 10 10 0  6.28318530717959E+0000}%
\special{sh 1}%
\special{ar 2800 1600 10 10 0  6.28318530717959E+0000}%
%
\special{pn 8}%
\special{pa 4300 1690}%
\special{pa 4500 1690}%
\special{pa 4500 2490}%
\special{pa 4300 2490}%
\special{pa 4300 1690}%
\special{fp}%
%
\special{pn 20}%
\special{sh 1}%
\special{ar 3200 1600 10 10 0  6.28318530717959E+0000}%
\special{sh 1}%
\special{ar 3400 1600 10 10 0  6.28318530717959E+0000}%
\special{sh 1}%
\special{ar 3800 1800 10 10 0  6.28318530717959E+0000}%
\special{sh 1}%
\special{ar 4000 1800 10 10 0  6.28318530717959E+0000}%
\special{sh 1}%
\special{ar 4200 1800 10 10 0  6.28318530717959E+0000}%
\special{sh 1}%
\special{ar 4400 1800 10 10 0  6.28318530717959E+0000}%
%
\special{pn 20}%
\special{sh 1}%
\special{ar 2800 2400 10 10 0  6.28318530717959E+0000}%
%
\special{pn 20}%
\special{sh 1}%
\special{ar 3000 2400 10 10 0  6.28318530717959E+0000}%
%
\special{pn 8}%
\special{pa 2800 1600}%
\special{pa 3400 1600}%
\special{pa 3200 2200}%
\special{pa 3200 2400}%
\special{pa 2800 2400}%
\special{pa 2800 1600}%
\special{fp}%
\put(31.0000,-21.9000){\makebox(0,0){$-\tau$}}%
%
\special{pn 20}%
\special{sh 1}%
\special{ar 4800 1800 10 10 0  6.28318530717959E+0000}%
\special{sh 1}%
\special{ar 4800 2000 10 10 0  6.28318530717959E+0000}%
\special{sh 1}%
\special{ar 5000 2200 10 10 0  6.28318530717959E+0000}%
\special{sh 1}%
\special{ar 5000 2400 10 10 0  6.28318530717959E+0000}%
\special{sh 1}%
\special{ar 4800 2400 10 10 0  6.28318530717959E+0000}%
\put(46.0000,-17.9000){\makebox(0,0){$0$}}%
\put(42.0000,-23.9000){\makebox(0,0){$0$}}%
\put(34.0000,-23.9000){\makebox(0,0){$0$}}%
%
\special{pn 20}%
\special{sh 1}%
\special{ar 4400 1596 10 10 0  6.28318530717959E+0000}%
%
\special{pn 20}%
\special{sh 1}%
\special{ar 4200 1396 10 10 0  6.28318530717959E+0000}%
%
\special{pn 20}%
\special{sh 1}%
\special{ar 3400 1396 10 10 0  6.28318530717959E+0000}%
%
\special{pn 8}%
\special{pa 2800 1600}%
\special{pa 3400 1600}%
\special{pa 3400 800}%
\special{pa 3000 800}%
\special{pa 2800 1600}%
\special{fp}%
%
\special{pn 8}%
\special{pa 3800 800}%
\special{pa 4200 800}%
\special{pa 4400 1600}%
\special{pa 4400 1800}%
\special{pa 3800 1800}%
\special{pa 3800 800}%
\special{fp}%
%
\special{pn 20}%
\special{sh 1}%
\special{ar 4000 1200 10 10 0  6.28318530717959E+0000}%
\special{sh 1}%
\special{ar 3800 1000 10 10 0  6.28318530717959E+0000}%
\special{sh 1}%
\special{ar 3800 800 10 10 0  6.28318530717959E+0000}%
\special{sh 1}%
\special{ar 3400 800 10 10 0  6.28318530717959E+0000}%
\special{sh 1}%
\special{ar 3200 800 10 10 0  6.28318530717959E+0000}%
\special{sh 1}%
\special{ar 3200 1200 10 10 0  6.28318530717959E+0000}%
\special{sh 1}%
\special{ar 3000 1000 10 10 0  6.28318530717959E+0000}%
\special{sh 1}%
\special{ar 3000 800 10 10 0  6.28318530717959E+0000}%
\special{sh 1}%
\special{ar 4000 800 10 10 0  6.28318530717959E+0000}%
\special{sh 1}%
\special{ar 4200 800 10 10 0  6.28318530717959E+0000}%
%
\special{pn 8}%
\special{pa 3200 2596}%
\special{pa 3200 650}%
\special{fp}%
\special{sh 1}%
\special{pa 3200 650}%
\special{pa 3180 718}%
\special{pa 3200 704}%
\special{pa 3220 718}%
\special{pa 3200 650}%
\special{fp}%
%
\special{pn 8}%
\special{pa 2600 1396}%
\special{pa 6180 1396}%
\special{fp}%
\special{sh 1}%
\special{pa 6180 1396}%
\special{pa 6114 1376}%
\special{pa 6128 1396}%
\special{pa 6114 1416}%
\special{pa 6180 1396}%
\special{fp}%
%
\special{pn 8}%
\special{pa 4200 1356}%
\special{pa 4200 1396}%
\special{fp}%
%
\special{pn 20}%
\special{sh 1}%
\special{ar 5000 1800 10 10 0  6.28318530717959E+0000}%
\special{sh 1}%
\special{ar 5200 1800 10 10 0  6.28318530717959E+0000}%
%
\special{pn 8}%
\special{pa 4400 2400}%
\special{pa 5000 2400}%
\special{pa 5000 2200}%
\special{pa 5200 1800}%
\special{pa 4800 1800}%
\special{pa 4800 2000}%
\special{pa 4400 2400}%
\special{fp}%
%
\special{pn 20}%
\special{sh 1}%
\special{ar 5400 2600 10 10 0  6.28318530717959E+0000}%
\special{sh 1}%
\special{ar 5600 2600 10 10 0  6.28318530717959E+0000}%
\special{sh 1}%
\special{ar 5800 2600 10 10 0  6.28318530717959E+0000}%
\special{sh 1}%
\special{ar 6000 2600 10 10 0  6.28318530717959E+0000}%
\special{sh 1}%
\special{ar 6000 2400 10 10 0  6.28318530717959E+0000}%
\special{sh 1}%
\special{ar 5800 2200 10 10 0  6.28318530717959E+0000}%
\special{sh 1}%
\special{ar 5600 2000 10 10 0  6.28318530717959E+0000}%
\special{sh 1}%
\special{ar 5600 1800 10 10 0  6.28318530717959E+0000}%
\special{sh 1}%
\special{ar 5800 1800 10 10 0  6.28318530717959E+0000}%
\special{sh 1}%
\special{ar 6000 1800 10 10 0  6.28318530717959E+0000}%
%
\special{pn 8}%
\special{pa 5600 1800}%
\special{pa 6000 1800}%
\special{pa 6000 2600}%
\special{pa 5400 2600}%
\special{pa 5600 2000}%
\special{pa 5600 1800}%
\special{fp}%
%
\special{pn 20}%
\special{sh 1}%
\special{ar 6000 1600 10 10 0  6.28318530717959E+0000}%
\special{sh 1}%
\special{ar 5800 1400 10 10 0  6.28318530717959E+0000}%
\special{sh 1}%
\special{ar 5200 1600 10 10 0  6.28318530717959E+0000}%
\special{sh 1}%
\special{ar 5000 1400 10 10 0  6.28318530717959E+0000}%
\special{sh 1}%
\special{ar 4800 1200 10 10 0  6.28318530717959E+0000}%
\special{sh 1}%
\special{ar 4600 1000 10 10 0  6.28318530717959E+0000}%
\special{sh 1}%
\special{ar 4600 800 10 10 0  6.28318530717959E+0000}%
\special{sh 1}%
\special{ar 5600 1200 10 10 0  6.28318530717959E+0000}%
\special{sh 1}%
\special{ar 5400 1000 10 10 0  6.28318530717959E+0000}%
\special{sh 1}%
\special{ar 5400 800 10 10 0  6.28318530717959E+0000}%
\special{sh 1}%
\special{ar 5000 800 10 10 0  6.28318530717959E+0000}%
\special{sh 1}%
\special{ar 4800 800 10 10 0  6.28318530717959E+0000}%
\special{sh 1}%
\special{ar 5600 800 10 10 0  6.28318530717959E+0000}%
\special{sh 1}%
\special{ar 5800 800 10 10 0  6.28318530717959E+0000}%
%
\special{pn 8}%
\special{pa 4400 1800}%
\special{pa 4000 2200}%
\special{pa 4000 2400}%
\special{pa 3600 2400}%
\special{pa 3800 2000}%
\special{pa 3800 1800}%
\special{pa 4400 1800}%
\special{fp}%
\put(54.0000,-18.0000){\makebox(0,0){$0$}}%
\put(52.0000,-24.0000){\makebox(0,0){$0$}}%
%
\special{pn 8}%
\special{pa 4800 1800}%
\special{pa 4600 1000}%
\special{pa 4600 800}%
\special{pa 5000 800}%
\special{pa 5200 1600}%
\special{pa 5200 1800}%
\special{pa 4800 1800}%
\special{fp}%
%
\special{pn 8}%
\special{pa 5600 1800}%
\special{pa 5400 1000}%
\special{pa 5400 800}%
\special{pa 5800 800}%
\special{pa 6000 1600}%
\special{pa 6000 1800}%
\special{pa 5600 1800}%
\special{fp}%
\put(56.0000,-28.0000){\makebox(0,0){$B_{n-n'-\gamma',-1}$}}%
\end{picture}%

%% file: tauIIppad.tex
\unitlength 0.1in
\begin{picture}( 55.7500, 19.0000)(  3.3000,-20.0000)
\put(7.3500,-13.8500){\makebox(0,0){$-\tau+1$}}%
\put(12.0500,-2.0500){\makebox(0,0){$0$}}%
\put(8.5000,-1.9500){\makebox(0,0){$j$}}%
\put(55.2500,-1.8500){\makebox(0,0){$i$}}%
%
\special{pn 8}%
\special{pa 800 400}%
\special{pa 5720 400}%
\special{fp}%
\special{sh 1}%
\special{pa 5720 400}%
\special{pa 5654 380}%
\special{pa 5668 400}%
\special{pa 5654 420}%
\special{pa 5720 400}%
\special{fp}%
\put(8.9000,-5.9000){\makebox(0,0){$0$}}%
\put(7.3500,-11.8500){\makebox(0,0){$-\tau+2$}}%
\put(22.0500,-2.0500){\makebox(0,0){$\tau$}}%
\put(35.0500,-2.0500){\makebox(0,0){$(\lfloor\frac{n}{\alpha}\rfloor-1)\tau$}}%
\put(14.0500,-6.0500){\makebox(0,0){$0$}}%
\put(16.0500,-12.0500){\makebox(0,0){$0$}}%
\put(16.0500,-14.0500){\makebox(0,0){$0$}}%
\put(24.0500,-6.0500){\makebox(0,0){$0$}}%
\put(24.0500,-8.0500){\makebox(0,0){$0$}}%
\put(24.0500,-10.0500){\makebox(0,0){$0$}}%
\put(26.0500,-14.0500){\makebox(0,0){$0$}}%
\put(26.0500,-12.0500){\makebox(0,0){$0$}}%
\put(16.0500,-10.0500){\makebox(0,0){$0$}}%
\put(14.1000,-10.1000){\makebox(0,0){$0$}}%
\put(37.0000,-6.0000){\makebox(0,0){$0$}}%
\put(37.0000,-8.0000){\makebox(0,0){$0$}}%
\put(37.0000,-10.0000){\makebox(0,0){$0$}}%
\put(39.0000,-12.0000){\makebox(0,0){$0$}}%
\put(39.0000,-14.0000){\makebox(0,0){$0$}}%
%
\special{pn 20}%
\special{sh 1}%
\special{ar 4506 600 10 10 0  6.28318530717959E+0000}%
\special{sh 1}%
\special{ar 4306 600 10 10 0  6.28318530717959E+0000}%
\special{sh 1}%
\special{ar 4106 600 10 10 0  6.28318530717959E+0000}%
\special{sh 1}%
\special{ar 3906 600 10 10 0  6.28318530717959E+0000}%
%
\special{pn 20}%
\special{sh 1}%
\special{ar 5300 1400 10 10 0  6.28318530717959E+0000}%
\put(35.0000,-10.0000){\makebox(0,0){$\cdots$}}%
\put(47.0000,-6.0000){\makebox(0,0){$0$}}%
\put(47.0000,-8.0000){\makebox(0,0){$0$}}%
\put(51.0500,-10.0000){\makebox(0,0){$0$}}%
\put(51.0500,-12.0000){\makebox(0,0){$0$}}%
\put(51.0500,-14.0000){\makebox(0,0){$0$}}%
%
\special{pn 8}%
\special{pa 2600 600}%
\special{pa 2600 600}%
\special{pa 2800 1200}%
\special{pa 2800 1400}%
\special{pa 3400 1400}%
\special{pa 3200 800}%
\special{pa 3200 600}%
\special{pa 2600 600}%
\special{fp}%
%
\special{pn 8}%
\special{pa 1600 600}%
\special{pa 1800 1200}%
\special{pa 1800 1400}%
\special{pa 2400 1400}%
\special{pa 2200 800}%
\special{pa 2200 600}%
\special{pa 1600 600}%
\special{fp}%
%
\special{pn 8}%
\special{pa 3906 600}%
\special{pa 4106 1200}%
\special{pa 4106 1400}%
\special{pa 4906 1400}%
\special{pa 4506 800}%
\special{pa 4506 600}%
\special{pa 3906 600}%
\special{fp}%
%
\special{pn 20}%
\special{sh 1}%
\special{ar 5306 1400 10 10 0  6.28318530717959E+0000}%
\special{sh 1}%
\special{ar 5306 1200 10 10 0  6.28318530717959E+0000}%
\special{sh 1}%
\special{ar 5506 1000 10 10 0  6.28318530717959E+0000}%
\special{sh 1}%
\special{ar 5706 800 10 10 0  6.28318530717959E+0000}%
\special{sh 1}%
\special{ar 5706 600 10 10 0  6.28318530717959E+0000}%
\special{sh 1}%
\special{ar 5506 600 10 10 0  6.28318530717959E+0000}%
\special{sh 1}%
\special{ar 5306 600 10 10 0  6.28318530717959E+0000}%
\special{sh 1}%
\special{ar 5106 600 10 10 0  6.28318530717959E+0000}%
\special{sh 1}%
\special{ar 4906 600 10 10 0  6.28318530717959E+0000}%
\special{sh 1}%
\special{ar 5506 1400 10 10 0  6.28318530717959E+0000}%
\special{sh 1}%
\special{ar 5706 1400 10 10 0  6.28318530717959E+0000}%
\special{sh 1}%
\special{ar 5906 1400 10 10 0  6.28318530717959E+0000}%
%
\special{pn 8}%
\special{pa 5006 1500}%
\special{pa 4806 1500}%
\special{pa 4806 500}%
\special{pa 5006 500}%
\special{pa 5006 1500}%
\special{fp}%
%
\special{pn 8}%
\special{pa 1000 2000}%
\special{pa 1000 200}%
\special{fp}%
\special{sh 1}%
\special{pa 1000 200}%
\special{pa 980 268}%
\special{pa 1000 254}%
\special{pa 1020 268}%
\special{pa 1000 200}%
\special{fp}%
%
\special{pn 20}%
\special{sh 1}%
\special{ar 4100 1400 10 10 0  6.28318530717959E+0000}%
\special{sh 1}%
\special{ar 4100 1200 10 10 0  6.28318530717959E+0000}%
\special{sh 1}%
\special{ar 4300 1000 10 10 0  6.28318530717959E+0000}%
\special{sh 1}%
\special{ar 4500 800 10 10 0  6.28318530717959E+0000}%
\special{sh 1}%
\special{ar 4300 1400 10 10 0  6.28318530717959E+0000}%
\special{sh 1}%
\special{ar 4500 1400 10 10 0  6.28318530717959E+0000}%
\special{sh 1}%
\special{ar 4700 1400 10 10 0  6.28318530717959E+0000}%
\special{sh 1}%
\special{ar 4900 1400 10 10 0  6.28318530717959E+0000}%
%
\special{pn 20}%
\special{sh 1}%
\special{ar 2600 600 10 10 0  6.28318530717959E+0000}%
\special{sh 1}%
\special{ar 2800 600 10 10 0  6.28318530717959E+0000}%
\special{sh 1}%
\special{ar 3000 600 10 10 0  6.28318530717959E+0000}%
\special{sh 1}%
\special{ar 3200 600 10 10 0  6.28318530717959E+0000}%
\special{sh 1}%
\special{ar 3200 800 10 10 0  6.28318530717959E+0000}%
\special{sh 1}%
\special{ar 3000 1000 10 10 0  6.28318530717959E+0000}%
\special{sh 1}%
\special{ar 2800 1200 10 10 0  6.28318530717959E+0000}%
\special{sh 1}%
\special{ar 2800 1400 10 10 0  6.28318530717959E+0000}%
\special{sh 1}%
\special{ar 3000 1400 10 10 0  6.28318530717959E+0000}%
\special{sh 1}%
\special{ar 3200 1400 10 10 0  6.28318530717959E+0000}%
\special{sh 1}%
\special{ar 3400 1400 10 10 0  6.28318530717959E+0000}%
%
\special{pn 20}%
\special{sh 1}%
\special{ar 1600 600 10 10 0  6.28318530717959E+0000}%
\special{sh 1}%
\special{ar 1800 600 10 10 0  6.28318530717959E+0000}%
\special{sh 1}%
\special{ar 2000 600 10 10 0  6.28318530717959E+0000}%
\special{sh 1}%
\special{ar 2200 600 10 10 0  6.28318530717959E+0000}%
\special{sh 1}%
\special{ar 2200 800 10 10 0  6.28318530717959E+0000}%
\special{sh 1}%
\special{ar 2000 1000 10 10 0  6.28318530717959E+0000}%
\special{sh 1}%
\special{ar 1800 1200 10 10 0  6.28318530717959E+0000}%
\special{sh 1}%
\special{ar 1800 1400 10 10 0  6.28318530717959E+0000}%
%
\special{pn 20}%
\special{sh 1}%
\special{ar 2400 1400 10 10 0  6.28318530717959E+0000}%
\put(39.0000,-10.0000){\makebox(0,0){$0$}}%
\put(49.1000,-10.0000){\makebox(0,0){$0$}}%
\put(49.1000,-12.0000){\makebox(0,0){$0$}}%
\put(49.1000,-8.0000){\makebox(0,0){$0$}}%
\put(26.1000,-10.1000){\makebox(0,0){$0$}}%
\put(45.0500,-2.0500){\makebox(0,0){$\lfloor\frac{n}{\alpha}\rfloor\tau$}}%
%
\special{pn 8}%
\special{pa 2200 340}%
\special{pa 2200 400}%
\special{fp}%
%
\special{pn 8}%
\special{pa 1200 340}%
\special{pa 1200 400}%
\special{fp}%
%
\special{pn 8}%
\special{pa 3500 340}%
\special{pa 3500 400}%
\special{fp}%
%
\special{pn 8}%
\special{pa 4500 340}%
\special{pa 4500 400}%
\special{fp}%
\put(16.1000,-8.1000){\makebox(0,0){$0$}}%
\put(14.1000,-8.1000){\makebox(0,0){$0$}}%
\put(39.0000,-8.0000){\makebox(0,0){$0$}}%
\put(26.1000,-8.1000){\makebox(0,0){$0$}}%
\put(24.1000,-12.1000){\makebox(0,0){$0$}}%
%
\special{pn 20}%
\special{sh 1}%
\special{ar 2000 1400 10 10 0  6.28318530717959E+0000}%
\special{sh 1}%
\special{ar 2200 1400 10 10 0  6.28318530717959E+0000}%
\put(21.0000,-12.0000){\makebox(0,0){$B_{0,-1}$}}%
\put(31.0000,-12.0000){\makebox(0,0){$B_{1,-1}$}}%
\put(45.0000,-12.0000){\makebox(0,0){$B_{\lfloor\frac{n}{\alpha}\rfloor-1,-1}$}}%
\put(56.0000,-12.0000){\makebox(0,0){$B_{\lfloor\frac{n}{\alpha}\rfloor,-1}$}}%
%
\special{pn 8}%
\special{pa 5500 600}%
\special{pa 4900 600}%
\special{fp}%
\special{pa 4900 600}%
\special{pa 5300 1200}%
\special{fp}%
\special{pa 5300 1200}%
\special{pa 5300 1400}%
\special{fp}%
\special{pa 5700 1400}%
\special{pa 5700 1400}%
\special{fp}%
%
\special{pn 8}%
\special{pa 5700 1400}%
\special{pa 5900 1400}%
\special{dt 0.045}%
%
\special{pn 8}%
\special{pa 5300 1400}%
\special{pa 5700 1400}%
\special{fp}%
%
\special{pn 8}%
\special{pa 5500 600}%
\special{pa 5700 600}%
\special{dt 0.045}%
%
\special{pn 20}%
\special{sh 1}%
\special{ar 1000 1400 10 10 0  6.28318530717959E+0000}%
\special{sh 1}%
\special{ar 1200 1400 10 10 0  6.28318530717959E+0000}%
\special{sh 1}%
\special{ar 1400 1400 10 10 0  6.28318530717959E+0000}%
\special{sh 1}%
\special{ar 1200 600 10 10 0  6.28318530717959E+0000}%
\special{sh 1}%
\special{ar 1000 800 10 10 0  6.28318530717959E+0000}%
\special{sh 1}%
\special{ar 1200 400 10 10 0  6.28318530717959E+0000}%
\special{sh 1}%
\special{ar 1000 400 10 10 0  6.28318530717959E+0000}%
%
\special{pn 8}%
\special{pa 1000 1400}%
\special{pa 1400 1400}%
\special{fp}%
\special{pa 1400 1400}%
\special{pa 1200 600}%
\special{fp}%
\special{pa 1200 600}%
\special{pa 1200 400}%
\special{fp}%
\end{picture}%

%% file: tauIIpp.tex
\unitlength 0.1in
\begin{picture}( 44.9000, 20.5000)(  3.3000,-24.0500)
\put(7.3500,-17.9000){\makebox(0,0){$-\tau+1$}}%
\put(12.0500,-6.1000){\makebox(0,0){$0$}}%
\put(10.0000,-4.4000){\makebox(0,0){$j$}}%
\put(49.3000,-5.9500){\makebox(0,0){$i$}}%
%
\special{pn 8}%
\special{pa 800 806}%
\special{pa 4820 806}%
\special{fp}%
\special{sh 1}%
\special{pa 4820 806}%
\special{pa 4754 786}%
\special{pa 4768 806}%
\special{pa 4754 826}%
\special{pa 4820 806}%
\special{fp}%
\put(8.9000,-9.9500){\makebox(0,0){$0$}}%
\put(7.3500,-15.9000){\makebox(0,0){$-\tau+2$}}%
\put(22.0500,-6.1000){\makebox(0,0){$\tau$}}%
\put(32.0000,-6.0500){\makebox(0,0){$\lfloor\frac{n}{\alpha}\rfloor\tau$}}%
\put(12.0500,-10.1000){\makebox(0,0){$0$}}%
\put(14.0500,-16.1000){\makebox(0,0){$0$}}%
\put(14.0500,-18.1000){\makebox(0,0){$0$}}%
\put(22.0500,-10.1000){\makebox(0,0){$0$}}%
\put(24.0000,-18.0500){\makebox(0,0){$0$}}%
\put(14.0500,-14.1000){\makebox(0,0){$0$}}%
\put(12.1000,-14.1500){\makebox(0,0){$0$}}%
\put(32.0000,-10.0500){\makebox(0,0){$0$}}%
\put(23.0000,-14.0500){\makebox(0,0){$\cdots$}}%
%
\special{pn 8}%
\special{pa 1400 1006}%
\special{pa 1600 1606}%
\special{pa 1600 1806}%
\special{pa 2200 1806}%
\special{pa 2000 1206}%
\special{pa 2000 1006}%
\special{pa 1400 1006}%
\special{fp}%
%
\special{pn 8}%
\special{pa 3500 1906}%
\special{pa 3300 1906}%
\special{pa 3300 906}%
\special{pa 3500 906}%
\special{pa 3500 1906}%
\special{fp}%
%
\special{pn 8}%
\special{pa 1000 2406}%
\special{pa 1000 606}%
\special{fp}%
\special{sh 1}%
\special{pa 1000 606}%
\special{pa 980 672}%
\special{pa 1000 658}%
\special{pa 1020 672}%
\special{pa 1000 606}%
\special{fp}%
%
\special{pn 20}%
\special{sh 1}%
\special{ar 1400 1006 10 10 0  6.28318530717959E+0000}%
\special{sh 1}%
\special{ar 1600 1006 10 10 0  6.28318530717959E+0000}%
\special{sh 1}%
\special{ar 1800 1006 10 10 0  6.28318530717959E+0000}%
\special{sh 1}%
\special{ar 2000 1006 10 10 0  6.28318530717959E+0000}%
\special{sh 1}%
\special{ar 2000 1206 10 10 0  6.28318530717959E+0000}%
\special{sh 1}%
\special{ar 1800 1406 10 10 0  6.28318530717959E+0000}%
\special{sh 1}%
\special{ar 1600 1606 10 10 0  6.28318530717959E+0000}%
\special{sh 1}%
\special{ar 1600 1806 10 10 0  6.28318530717959E+0000}%
%
\special{pn 20}%
\special{sh 1}%
\special{ar 2200 1806 10 10 0  6.28318530717959E+0000}%
\put(44.0000,-10.0500){\makebox(0,0){$0$}}%
\put(42.0000,-6.0500){\makebox(0,0){$(\lfloor\frac{n}{\alpha}\rfloor+1)\tau$}}%
%
\special{pn 8}%
\special{pa 2200 746}%
\special{pa 2200 806}%
\special{fp}%
%
\special{pn 8}%
\special{pa 1200 746}%
\special{pa 1200 806}%
\special{fp}%
%
\special{pn 8}%
\special{pa 3200 756}%
\special{pa 3200 816}%
\special{fp}%
%
\special{pn 8}%
\special{pa 4200 756}%
\special{pa 4200 816}%
\special{fp}%
\put(14.1000,-12.1500){\makebox(0,0){$0$}}%
\put(12.1000,-12.1500){\makebox(0,0){$0$}}%
%
\special{pn 20}%
\special{sh 1}%
\special{ar 1800 1806 10 10 0  6.28318530717959E+0000}%
\special{sh 1}%
\special{ar 2000 1806 10 10 0  6.28318530717959E+0000}%
\put(19.0000,-20.0500){\makebox(0,0){$B_{0,-1}$}}%
%
\special{pn 20}%
\special{sh 1}%
\special{ar 2400 1006 10 10 0  6.28318530717959E+0000}%
\special{sh 1}%
\special{ar 2600 1006 10 10 0  6.28318530717959E+0000}%
\special{sh 1}%
\special{ar 2800 1006 10 10 0  6.28318530717959E+0000}%
\special{sh 1}%
\special{ar 3000 1006 10 10 0  6.28318530717959E+0000}%
\special{sh 1}%
\special{ar 3000 1206 10 10 0  6.28318530717959E+0000}%
\special{sh 1}%
\special{ar 2800 1406 10 10 0  6.28318530717959E+0000}%
\special{sh 1}%
\special{ar 2600 1606 10 10 0  6.28318530717959E+0000}%
\special{sh 1}%
\special{ar 2600 1806 10 10 0  6.28318530717959E+0000}%
\special{sh 1}%
\special{ar 2800 1806 10 10 0  6.28318530717959E+0000}%
\special{sh 1}%
\special{ar 3000 1806 10 10 0  6.28318530717959E+0000}%
\special{sh 1}%
\special{ar 3200 1806 10 10 0  6.28318530717959E+0000}%
\special{sh 1}%
\special{ar 3400 1806 10 10 0  6.28318530717959E+0000}%
%
\special{pn 8}%
\special{pa 2400 1006}%
\special{pa 2600 1606}%
\special{pa 2600 1806}%
\special{pa 3400 1806}%
\special{pa 3000 1206}%
\special{pa 3000 1006}%
\special{pa 2400 1006}%
\special{fp}%
%
\special{pn 20}%
\special{sh 1}%
\special{ar 3400 1006 10 10 0  6.28318530717959E+0000}%
\special{sh 1}%
\special{ar 3600 1006 10 10 0  6.28318530717959E+0000}%
\special{sh 1}%
\special{ar 3800 1006 10 10 0  6.28318530717959E+0000}%
\special{sh 1}%
\special{ar 4000 1006 10 10 0  6.28318530717959E+0000}%
\special{sh 1}%
\special{ar 4200 1006 10 10 0  6.28318530717959E+0000}%
\special{sh 1}%
\special{ar 4200 1206 10 10 0  6.28318530717959E+0000}%
\special{sh 1}%
\special{ar 4000 1406 10 10 0  6.28318530717959E+0000}%
\special{sh 1}%
\special{ar 3800 1606 10 10 0  6.28318530717959E+0000}%
\special{sh 1}%
\special{ar 3800 1806 10 10 0  6.28318530717959E+0000}%
\special{sh 1}%
\special{ar 4000 1806 10 10 0  6.28318530717959E+0000}%
\special{sh 1}%
\special{ar 4400 1806 10 10 0  6.28318530717959E+0000}%
\special{sh 1}%
\special{ar 4200 1806 10 10 0  6.28318530717959E+0000}%
%
\special{pn 8}%
\special{pa 3400 1006}%
\special{pa 3800 1606}%
\special{pa 3800 1806}%
\special{pa 4400 1806}%
\special{pa 4200 1206}%
\special{pa 4200 1006}%
\special{pa 3400 1006}%
\special{fp}%
\put(30.0000,-20.0500){\makebox(0,0){$B_{\lfloor\frac{n}{\alpha}\rfloor-1,-1}$}}%
\put(40.0000,-20.0500){\makebox(0,0){$B_{\lfloor\frac{n}{\alpha}\rfloor,-1}$}}%
\put(36.0000,-18.0500){\makebox(0,0){$0$}}%
\end{picture}%

%% file: tauIIpp2.tex
\unitlength 0.1in
\begin{picture}( 61.7000, 22.0500)(  3.3000,-23.1500)
\put(7.3500,-13.8500){\makebox(0,0){$-\tau+1$}}%
\put(8.5000,-1.9500){\makebox(0,0){$j$}}%
\put(63.3000,-3.9500){\makebox(0,0){$i$}}%
%
\special{pn 8}%
\special{pa 800 600}%
\special{pa 6500 600}%
\special{fp}%
\special{sh 1}%
\special{pa 6500 600}%
\special{pa 6434 580}%
\special{pa 6448 600}%
\special{pa 6434 620}%
\special{pa 6500 600}%
\special{fp}%
\put(8.9000,-7.9000){\makebox(0,0){O}}%
\put(7.3500,-11.8500){\makebox(0,0){$-\tau+2$}}%
\put(22.1000,-4.1000){\makebox(0,0){$\tau$}}%
\put(32.0000,-4.0000){\makebox(0,0){$i\tau$}}%
\put(12.0000,-6.0000){\makebox(0,0){$0$}}%
\put(22.0500,-6.0500){\makebox(0,0){$0$}}%
\put(24.0500,-14.0500){\makebox(0,0){$0$}}%
\put(32.0000,-6.0000){\makebox(0,0){$0$}}%
\put(23.0000,-10.0000){\makebox(0,0){$\cdots$}}%
%
\special{pn 8}%
\special{pa 1400 600}%
\special{pa 1600 1200}%
\special{pa 1600 1400}%
\special{pa 2200 1400}%
\special{pa 2000 800}%
\special{pa 2000 600}%
\special{pa 1400 600}%
\special{fp}%
%
\special{pn 8}%
\special{pa 3500 1500}%
\special{pa 3300 1500}%
\special{pa 3300 500}%
\special{pa 3500 500}%
\special{pa 3500 1500}%
\special{fp}%
%
\special{pn 8}%
\special{pa 1000 2000}%
\special{pa 1000 200}%
\special{fp}%
\special{sh 1}%
\special{pa 1000 200}%
\special{pa 980 268}%
\special{pa 1000 254}%
\special{pa 1020 268}%
\special{pa 1000 200}%
\special{fp}%
%
\special{pn 20}%
\special{sh 1}%
\special{ar 1400 600 10 10 0  6.28318530717959E+0000}%
\special{sh 1}%
\special{ar 1600 600 10 10 0  6.28318530717959E+0000}%
\special{sh 1}%
\special{ar 1800 600 10 10 0  6.28318530717959E+0000}%
\special{sh 1}%
\special{ar 2000 600 10 10 0  6.28318530717959E+0000}%
\special{sh 1}%
\special{ar 2000 800 10 10 0  6.28318530717959E+0000}%
\special{sh 1}%
\special{ar 1800 1000 10 10 0  6.28318530717959E+0000}%
\special{sh 1}%
\special{ar 1600 1200 10 10 0  6.28318530717959E+0000}%
\special{sh 1}%
\special{ar 1600 1400 10 10 0  6.28318530717959E+0000}%
%
\special{pn 20}%
\special{sh 1}%
\special{ar 2200 1400 10 10 0  6.28318530717959E+0000}%
\put(44.0000,-6.0000){\makebox(0,0){$0$}}%
\put(42.0000,-4.0000){\makebox(0,0){$(i+1)\tau$}}%
%
\special{pn 8}%
\special{pa 2200 546}%
\special{pa 2200 606}%
\special{fp}%
%
\special{pn 8}%
\special{pa 3200 560}%
\special{pa 3200 620}%
\special{fp}%
%
\special{pn 8}%
\special{pa 4200 550}%
\special{pa 4200 610}%
\special{fp}%
\put(46.0000,-14.0000){\makebox(0,0){$0$}}%
%
\special{pn 20}%
\special{sh 1}%
\special{ar 1800 1400 10 10 0  6.28318530717959E+0000}%
\special{sh 1}%
\special{ar 2000 1400 10 10 0  6.28318530717959E+0000}%
\put(19.0000,-16.0000){\makebox(0,0){$B_{0,-1}$}}%
\put(61.0000,-18.0000){\makebox(0,0){$B_{n'-\gamma',-1}$}}%
%
\special{pn 20}%
\special{sh 1}%
\special{ar 2400 600 10 10 0  6.28318530717959E+0000}%
\special{sh 1}%
\special{ar 2600 600 10 10 0  6.28318530717959E+0000}%
\special{sh 1}%
\special{ar 2800 600 10 10 0  6.28318530717959E+0000}%
\special{sh 1}%
\special{ar 3000 600 10 10 0  6.28318530717959E+0000}%
\special{sh 1}%
\special{ar 3000 800 10 10 0  6.28318530717959E+0000}%
\special{sh 1}%
\special{ar 2800 1000 10 10 0  6.28318530717959E+0000}%
\special{sh 1}%
\special{ar 2600 1200 10 10 0  6.28318530717959E+0000}%
\special{sh 1}%
\special{ar 2600 1400 10 10 0  6.28318530717959E+0000}%
\special{sh 1}%
\special{ar 2800 1400 10 10 0  6.28318530717959E+0000}%
\special{sh 1}%
\special{ar 3000 1400 10 10 0  6.28318530717959E+0000}%
\special{sh 1}%
\special{ar 3200 1400 10 10 0  6.28318530717959E+0000}%
\special{sh 1}%
\special{ar 3400 1400 10 10 0  6.28318530717959E+0000}%
%
\special{pn 8}%
\special{pa 2400 600}%
\special{pa 2600 1200}%
\special{pa 2600 1400}%
\special{pa 3400 1400}%
\special{pa 3000 800}%
\special{pa 3000 600}%
\special{pa 2400 600}%
\special{fp}%
%
\special{pn 20}%
\special{sh 1}%
\special{ar 3400 600 10 10 0  6.28318530717959E+0000}%
\special{sh 1}%
\special{ar 3600 600 10 10 0  6.28318530717959E+0000}%
\special{sh 1}%
\special{ar 3800 600 10 10 0  6.28318530717959E+0000}%
\special{sh 1}%
\special{ar 4000 600 10 10 0  6.28318530717959E+0000}%
\special{sh 1}%
\special{ar 4200 600 10 10 0  6.28318530717959E+0000}%
\special{sh 1}%
\special{ar 4200 800 10 10 0  6.28318530717959E+0000}%
\special{sh 1}%
\special{ar 4000 1000 10 10 0  6.28318530717959E+0000}%
\special{sh 1}%
\special{ar 3800 1200 10 10 0  6.28318530717959E+0000}%
\special{sh 1}%
\special{ar 3800 1400 10 10 0  6.28318530717959E+0000}%
\special{sh 1}%
\special{ar 4000 1400 10 10 0  6.28318530717959E+0000}%
\special{sh 1}%
\special{ar 4400 1400 10 10 0  6.28318530717959E+0000}%
\special{sh 1}%
\special{ar 4200 1400 10 10 0  6.28318530717959E+0000}%
%
\special{pn 8}%
\special{pa 3400 600}%
\special{pa 3800 1200}%
\special{pa 3800 1400}%
\special{pa 4400 1400}%
\special{pa 4200 800}%
\special{pa 4200 600}%
\special{pa 3400 600}%
\special{fp}%
\put(44.0000,-10.0000){\makebox(0,0){$\cdots$}}%
%
\special{pn 20}%
\special{sh 1}%
\special{ar 5800 1400 10 10 0  6.28318530717959E+0000}%
\special{sh 1}%
\special{ar 5800 1600 10 10 0  6.28318530717959E+0000}%
\special{sh 1}%
\special{ar 6000 1200 10 10 0  6.28318530717959E+0000}%
\special{sh 1}%
\special{ar 6200 1000 10 10 0  6.28318530717959E+0000}%
\special{sh 1}%
\special{ar 6400 800 10 10 0  6.28318530717959E+0000}%
\special{sh 1}%
\special{ar 6400 600 10 10 0  6.28318530717959E+0000}%
\special{sh 1}%
\special{ar 6200 600 10 10 0  6.28318530717959E+0000}%
\special{sh 1}%
\special{ar 6000 600 10 10 0  6.28318530717959E+0000}%
\special{sh 1}%
\special{ar 5800 600 10 10 0  6.28318530717959E+0000}%
\special{sh 1}%
\special{ar 5600 600 10 10 0  6.28318530717959E+0000}%
%
\special{pn 8}%
\special{pa 5600 600}%
\special{pa 6400 600}%
\special{pa 6400 800}%
\special{pa 6400 1600}%
\special{pa 5800 1600}%
\special{pa 5800 1400}%
\special{pa 5600 600}%
\special{fp}%
%
\special{pn 8}%
\special{pa 4600 600}%
\special{pa 4800 1200}%
\special{pa 4800 1400}%
\special{pa 5400 1400}%
\special{pa 5200 800}%
\special{pa 5200 600}%
\special{pa 4600 600}%
\special{fp}%
%
\special{pn 20}%
\special{sh 1}%
\special{ar 4600 600 10 10 0  6.28318530717959E+0000}%
\special{sh 1}%
\special{ar 4800 600 10 10 0  6.28318530717959E+0000}%
\special{sh 1}%
\special{ar 5000 600 10 10 0  6.28318530717959E+0000}%
\special{sh 1}%
\special{ar 5200 600 10 10 0  6.28318530717959E+0000}%
\special{sh 1}%
\special{ar 5200 800 10 10 0  6.28318530717959E+0000}%
\special{sh 1}%
\special{ar 5000 1000 10 10 0  6.28318530717959E+0000}%
\special{sh 1}%
\special{ar 4800 1200 10 10 0  6.28318530717959E+0000}%
\special{sh 1}%
\special{ar 4800 1400 10 10 0  6.28318530717959E+0000}%
%
\special{pn 20}%
\special{sh 1}%
\special{ar 5000 1400 10 10 0  6.28318530717959E+0000}%
\special{sh 1}%
\special{ar 5200 1400 10 10 0  6.28318530717959E+0000}%
%
\special{pn 20}%
\special{sh 1}%
\special{ar 5400 1400 10 10 0  6.28318530717959E+0000}%
%
\special{pn 20}%
\special{sh 1}%
\special{ar 6000 1600 10 10 0  6.28318530717959E+0000}%
\special{sh 1}%
\special{ar 6200 1600 10 10 0  6.28318530717959E+0000}%
\special{sh 1}%
\special{ar 6400 1600 10 10 0  6.28318530717959E+0000}%
\put(36.0000,-14.0000){\makebox(0,0){$0$}}%
\put(56.0000,-14.0000){\makebox(0,0){$0$}}%
\put(54.0000,-6.0000){\makebox(0,0){$0$}}%
\put(14.0000,-14.0000){\makebox(0,0){$0$}}%
%
\special{pn 20}%
\special{sh 1}%
\special{ar 1000 600 10 10 0  6.28318530717959E+0000}%
\special{sh 1}%
\special{ar 1000 400 10 10 0  6.28318530717959E+0000}%
\special{sh 1}%
\special{ar 800 400 10 10 0  6.28318530717959E+0000}%
\special{sh 1}%
\special{ar 1200 1400 10 10 0  6.28318530717959E+0000}%
\special{sh 1}%
\special{ar 1000 1400 10 10 0  6.28318530717959E+0000}%
\special{sh 1}%
\special{ar 800 1400 10 10 0  6.28318530717959E+0000}%
%
\special{pn 8}%
\special{pa 800 400}%
\special{pa 1000 400}%
\special{fp}%
\special{pa 1000 600}%
\special{pa 1000 600}%
\special{fp}%
\special{pa 1200 1400}%
\special{pa 1000 600}%
\special{fp}%
\special{pa 1200 1400}%
\special{pa 800 1400}%
\special{fp}%
\put(42.0000,-16.0000){\makebox(0,0){$B_{i,-1}$}}%
%
\special{pn 8}%
\special{ar 3400 1800 400 400  6.2831853 6.2831853}%
\special{ar 3400 1800 400 400  0.0000000 3.1415927}%
%
\special{pn 8}%
\special{pa 3000 1816}%
\special{pa 3000 1800}%
\special{fp}%
\special{sh 1}%
\special{pa 3000 1800}%
\special{pa 2980 1868}%
\special{pa 3000 1854}%
\special{pa 3020 1868}%
\special{pa 3000 1800}%
\special{fp}%
%
\special{pn 8}%
\special{pa 3800 1816}%
\special{pa 3800 1800}%
\special{fp}%
\special{sh 1}%
\special{pa 3800 1800}%
\special{pa 3780 1868}%
\special{pa 3800 1854}%
\special{pa 3820 1868}%
\special{pa 3800 1800}%
\special{fp}%
\put(34.0000,-24.0000){\makebox(0,0){symmetry}}%
\end{picture}%

%% file: height2.tex
\unitlength 0.1in
\begin{picture}( 56.1000, 26.8500)(  2.9000,-32.0000)
\put(7.3500,-19.8500){\makebox(0,0){$-\tau+1$}}%
\put(12.0500,-8.0500){\makebox(0,0){$0$}}%
\put(8.5000,-7.9500){\makebox(0,0){$j$}}%
\put(55.2500,-7.8500){\makebox(0,0){$i$}}%
%
\special{pn 8}%
\special{pa 800 1000}%
\special{pa 5720 1000}%
\special{fp}%
\special{sh 1}%
\special{pa 5720 1000}%
\special{pa 5654 980}%
\special{pa 5668 1000}%
\special{pa 5654 1020}%
\special{pa 5720 1000}%
\special{fp}%
\put(8.9000,-11.9000){\makebox(0,0){$0$}}%
\put(7.3500,-17.8500){\makebox(0,0){$-\tau+2$}}%
\put(22.0500,-8.0500){\makebox(0,0){$\tau$}}%
\put(35.0500,-8.0500){\makebox(0,0){$(n'-\gamma'-1)\tau$}}%
\put(14.0500,-12.0500){\makebox(0,0){$0$}}%
\put(16.0500,-18.0500){\makebox(0,0){$0$}}%
\put(16.0500,-20.0500){\makebox(0,0){$0$}}%
\put(24.0500,-12.0500){\makebox(0,0){$0$}}%
\put(24.0500,-14.0500){\makebox(0,0){$0$}}%
\put(24.0500,-16.0500){\makebox(0,0){$0$}}%
\put(26.0500,-20.0500){\makebox(0,0){$0$}}%
\put(26.0500,-18.0500){\makebox(0,0){$0$}}%
\put(16.0500,-16.0500){\makebox(0,0){$0$}}%
\put(14.1000,-16.1000){\makebox(0,0){$0$}}%
\put(37.0000,-12.0000){\makebox(0,0){$0$}}%
\put(37.0000,-14.0000){\makebox(0,0){$0$}}%
\put(37.0000,-16.0000){\makebox(0,0){$0$}}%
\put(39.0000,-18.0000){\makebox(0,0){$0$}}%
\put(39.0000,-20.0000){\makebox(0,0){$0$}}%
\put(35.0000,-16.0000){\makebox(0,0){$\cdots$}}%
\put(47.0000,-12.0000){\makebox(0,0){$0$}}%
\put(47.0000,-14.0000){\makebox(0,0){$0$}}%
%
\special{pn 8}%
\special{pa 2600 1200}%
\special{pa 2600 1200}%
\special{pa 2800 1800}%
\special{pa 2800 2000}%
\special{pa 3400 2000}%
\special{pa 3200 1400}%
\special{pa 3200 1200}%
\special{pa 2600 1200}%
\special{fp}%
%
\special{pn 8}%
\special{pa 1600 1200}%
\special{pa 1800 1800}%
\special{pa 1800 2000}%
\special{pa 2400 2000}%
\special{pa 2200 1400}%
\special{pa 2200 1200}%
\special{pa 1600 1200}%
\special{fp}%
%
\special{pn 20}%
\special{sh 1}%
\special{ar 2600 1200 10 10 0  6.28318530717959E+0000}%
\special{sh 1}%
\special{ar 2800 1200 10 10 0  6.28318530717959E+0000}%
\special{sh 1}%
\special{ar 3000 1200 10 10 0  6.28318530717959E+0000}%
\special{sh 1}%
\special{ar 3200 1200 10 10 0  6.28318530717959E+0000}%
\special{sh 1}%
\special{ar 3200 1400 10 10 0  6.28318530717959E+0000}%
\special{sh 1}%
\special{ar 3000 1600 10 10 0  6.28318530717959E+0000}%
\special{sh 1}%
\special{ar 2800 1800 10 10 0  6.28318530717959E+0000}%
\special{sh 1}%
\special{ar 2800 2000 10 10 0  6.28318530717959E+0000}%
\special{sh 1}%
\special{ar 3000 2000 10 10 0  6.28318530717959E+0000}%
\special{sh 1}%
\special{ar 3200 2000 10 10 0  6.28318530717959E+0000}%
\special{sh 1}%
\special{ar 3400 2000 10 10 0  6.28318530717959E+0000}%
%
\special{pn 20}%
\special{sh 1}%
\special{ar 1600 1200 10 10 0  6.28318530717959E+0000}%
\special{sh 1}%
\special{ar 1800 1200 10 10 0  6.28318530717959E+0000}%
\special{sh 1}%
\special{ar 2000 1200 10 10 0  6.28318530717959E+0000}%
\special{sh 1}%
\special{ar 2200 1200 10 10 0  6.28318530717959E+0000}%
\special{sh 1}%
\special{ar 2200 1400 10 10 0  6.28318530717959E+0000}%
\special{sh 1}%
\special{ar 2000 1600 10 10 0  6.28318530717959E+0000}%
\special{sh 1}%
\special{ar 1800 1800 10 10 0  6.28318530717959E+0000}%
\special{sh 1}%
\special{ar 1800 2000 10 10 0  6.28318530717959E+0000}%
%
\special{pn 20}%
\special{sh 1}%
\special{ar 2400 2000 10 10 0  6.28318530717959E+0000}%
\put(39.0000,-16.0000){\makebox(0,0){$0$}}%
\put(49.1000,-16.0000){\makebox(0,0){$0$}}%
\put(49.1000,-18.0000){\makebox(0,0){$0$}}%
\put(49.1000,-14.0000){\makebox(0,0){$0$}}%
\put(26.1000,-16.1000){\makebox(0,0){$0$}}%
\put(45.0500,-8.0500){\makebox(0,0){$(n'-\gamma')\tau$}}%
%
\special{pn 8}%
\special{pa 2200 940}%
\special{pa 2200 1000}%
\special{fp}%
%
\special{pn 8}%
\special{pa 1200 940}%
\special{pa 1200 1000}%
\special{fp}%
%
\special{pn 8}%
\special{pa 3500 940}%
\special{pa 3500 1000}%
\special{fp}%
%
\special{pn 8}%
\special{pa 4500 940}%
\special{pa 4500 1000}%
\special{fp}%
\put(16.1000,-14.1000){\makebox(0,0){$0$}}%
\put(14.1000,-14.1000){\makebox(0,0){$0$}}%
\put(39.0000,-14.0000){\makebox(0,0){$0$}}%
\put(26.1000,-14.1000){\makebox(0,0){$0$}}%
\put(24.1000,-18.1000){\makebox(0,0){$0$}}%
%
\special{pn 20}%
\special{sh 1}%
\special{ar 2000 2000 10 10 0  6.28318530717959E+0000}%
\special{sh 1}%
\special{ar 2200 2000 10 10 0  6.28318530717959E+0000}%
\put(21.0000,-6.0000){\makebox(0,0){$B_{0,-1}$}}%
\put(31.0000,-6.0000){\makebox(0,0){$B_{1,-1}$}}%
%
\special{pn 8}%
\special{pa 3900 1200}%
\special{pa 3900 1200}%
\special{pa 4100 1800}%
\special{pa 4100 2000}%
\special{pa 4700 2000}%
\special{pa 4500 1400}%
\special{pa 4500 1200}%
\special{pa 3900 1200}%
\special{fp}%
%
\special{pn 20}%
\special{sh 1}%
\special{ar 3900 1200 10 10 0  6.28318530717959E+0000}%
\special{sh 1}%
\special{ar 4100 1200 10 10 0  6.28318530717959E+0000}%
\special{sh 1}%
\special{ar 4300 1200 10 10 0  6.28318530717959E+0000}%
\special{sh 1}%
\special{ar 4500 1200 10 10 0  6.28318530717959E+0000}%
\special{sh 1}%
\special{ar 4500 1400 10 10 0  6.28318530717959E+0000}%
\special{sh 1}%
\special{ar 4300 1600 10 10 0  6.28318530717959E+0000}%
\special{sh 1}%
\special{ar 4100 1800 10 10 0  6.28318530717959E+0000}%
\special{sh 1}%
\special{ar 4100 2000 10 10 0  6.28318530717959E+0000}%
\special{sh 1}%
\special{ar 4300 2000 10 10 0  6.28318530717959E+0000}%
\special{sh 1}%
\special{ar 4500 2000 10 10 0  6.28318530717959E+0000}%
\special{sh 1}%
\special{ar 4700 2000 10 10 0  6.28318530717959E+0000}%
\put(47.1000,-16.0000){\makebox(0,0){$0$}}%
%
\special{pn 20}%
\special{sh 1}%
\special{ar 4900 1200 10 10 0  6.28318530717959E+0000}%
\special{sh 1}%
\special{ar 5100 1200 10 10 0  6.28318530717959E+0000}%
\special{sh 1}%
\special{ar 5300 1200 10 10 0  6.28318530717959E+0000}%
\special{sh 1}%
\special{ar 5500 1200 10 10 0  6.28318530717959E+0000}%
\special{sh 1}%
\special{ar 5700 1200 10 10 0  6.28318530717959E+0000}%
\special{sh 1}%
\special{ar 5700 1400 10 10 0  6.28318530717959E+0000}%
\special{sh 1}%
\special{ar 5500 1600 10 10 0  6.28318530717959E+0000}%
\special{sh 1}%
\special{ar 5300 1800 10 10 0  6.28318530717959E+0000}%
\special{sh 1}%
\special{ar 5100 2000 10 10 0  6.28318530717959E+0000}%
\special{sh 1}%
\special{ar 5100 2200 10 10 0  6.28318530717959E+0000}%
\special{sh 1}%
\special{ar 5300 2200 10 10 0  6.28318530717959E+0000}%
\special{sh 1}%
\special{ar 5500 2200 10 10 0  6.28318530717959E+0000}%
\special{sh 1}%
\special{ar 5700 2200 10 10 0  6.28318530717959E+0000}%
\special{sh 1}%
\special{ar 5700 2200 10 10 0  6.28318530717959E+0000}%
%
\special{pn 8}%
\special{pa 5700 2200}%
\special{pa 5100 2200}%
\special{pa 5100 2000}%
\special{pa 4900 1200}%
\special{pa 5700 1200}%
\special{pa 5700 2200}%
\special{fp}%
\put(44.0000,-6.0000){\makebox(0,0){$B_{n'-\gamma'-1,-1}$}}%
\put(53.0000,-6.0000){\makebox(0,0){$B_{n'-\gamma',-1}$}}%
\put(49.1000,-20.0000){\makebox(0,0){$0$}}%
%
\special{pn 20}%
\special{sh 1}%
\special{ar 1400 2000 10 10 0  6.28318530717959E+0000}%
\special{sh 1}%
\special{ar 1200 2000 10 10 0  6.28318530717959E+0000}%
\special{sh 1}%
\special{ar 1000 2000 10 10 0  6.28318530717959E+0000}%
\special{sh 1}%
\special{ar 1200 1200 10 10 0  6.28318530717959E+0000}%
\special{sh 1}%
\special{ar 1000 1400 10 10 0  6.28318530717959E+0000}%
\special{sh 1}%
\special{ar 1200 1000 10 10 0  6.28318530717959E+0000}%
%
\special{pn 8}%
\special{pa 1000 2000}%
\special{pa 1400 2000}%
\special{fp}%
\special{pa 1400 2000}%
\special{pa 1200 1200}%
\special{fp}%
\special{pa 1200 1200}%
\special{pa 1200 1000}%
\special{fp}%
%
\special{pn 8}%
\special{pa 5100 2200}%
\special{pa 5300 2800}%
\special{pa 5300 3000}%
\special{pa 5900 3000}%
\special{pa 5700 2400}%
\special{pa 5700 2200}%
\special{pa 5100 2200}%
\special{fp}%
%
\special{pn 20}%
\special{sh 1}%
\special{ar 5700 2400 10 10 0  6.28318530717959E+0000}%
\special{sh 1}%
\special{ar 5500 2600 10 10 0  6.28318530717959E+0000}%
\special{sh 1}%
\special{ar 5300 2800 10 10 0  6.28318530717959E+0000}%
\special{sh 1}%
\special{ar 5300 3000 10 10 0  6.28318530717959E+0000}%
\special{sh 1}%
\special{ar 5500 3000 10 10 0  6.28318530717959E+0000}%
\special{sh 1}%
\special{ar 5700 3000 10 10 0  6.28318530717959E+0000}%
\special{sh 1}%
\special{ar 5900 3000 10 10 0  6.28318530717959E+0000}%
%
\special{pn 20}%
\special{sh 1}%
\special{ar 4700 2200 10 10 0  6.28318530717959E+0000}%
\special{sh 1}%
\special{ar 4500 2400 10 10 0  6.28318530717959E+0000}%
\special{sh 1}%
\special{ar 4300 2600 10 10 0  6.28318530717959E+0000}%
\special{sh 1}%
\special{ar 4100 2800 10 10 0  6.28318530717959E+0000}%
\special{sh 1}%
\special{ar 4100 3000 10 10 0  6.28318530717959E+0000}%
\special{sh 1}%
\special{ar 4300 3000 10 10 0  6.28318530717959E+0000}%
\special{sh 1}%
\special{ar 4500 3000 10 10 0  6.28318530717959E+0000}%
\special{sh 1}%
\special{ar 4700 3000 10 10 0  6.28318530717959E+0000}%
\special{sh 1}%
\special{ar 4900 3000 10 10 0  6.28318530717959E+0000}%
%
\special{pn 8}%
\special{pa 4900 3000}%
\special{pa 4100 3000}%
\special{pa 4100 2000}%
\special{pa 4700 2000}%
\special{pa 4700 2200}%
\special{pa 4900 3000}%
\special{fp}%
\put(49.1000,-22.0000){\makebox(0,0){$0$}}%
\put(49.1000,-24.0000){\makebox(0,0){$0$}}%
\put(49.1000,-26.0000){\makebox(0,0){$0$}}%
\put(51.1000,-30.0000){\makebox(0,0){$0$}}%
\put(51.1000,-28.0000){\makebox(0,0){$0$}}%
%
\special{pn 20}%
\special{sh 1}%
\special{ar 3400 2200 10 10 0  6.28318530717959E+0000}%
\special{sh 1}%
\special{ar 3200 2400 10 10 0  6.28318530717959E+0000}%
\special{sh 1}%
\special{ar 3000 2600 10 10 0  6.28318530717959E+0000}%
\special{sh 1}%
\special{ar 2800 2800 10 10 0  6.28318530717959E+0000}%
\special{sh 1}%
\special{ar 2800 3000 10 10 0  6.28318530717959E+0000}%
\special{sh 1}%
\special{ar 3000 3000 10 10 0  6.28318530717959E+0000}%
\special{sh 1}%
\special{ar 3200 3000 10 10 0  6.28318530717959E+0000}%
\special{sh 1}%
\special{ar 3400 3000 10 10 0  6.28318530717959E+0000}%
%
\special{pn 20}%
\special{sh 1}%
\special{ar 2400 2200 10 10 0  6.28318530717959E+0000}%
\special{sh 1}%
\special{ar 2200 2400 10 10 0  6.28318530717959E+0000}%
\special{sh 1}%
\special{ar 2000 2600 10 10 0  6.28318530717959E+0000}%
\special{sh 1}%
\special{ar 1800 2800 10 10 0  6.28318530717959E+0000}%
\special{sh 1}%
\special{ar 1800 3000 10 10 0  6.28318530717959E+0000}%
\special{sh 1}%
\special{ar 2000 3000 10 10 0  6.28318530717959E+0000}%
\special{sh 1}%
\special{ar 2200 3000 10 10 0  6.28318530717959E+0000}%
\special{sh 1}%
\special{ar 2400 3000 10 10 0  6.28318530717959E+0000}%
\special{sh 1}%
\special{ar 1400 2200 10 10 0  6.28318530717959E+0000}%
\special{sh 1}%
\special{ar 1200 2400 10 10 0  6.28318530717959E+0000}%
\special{sh 1}%
\special{ar 1000 2600 10 10 0  6.28318530717959E+0000}%
\special{sh 1}%
\special{ar 1400 3000 10 10 0  6.28318530717959E+0000}%
\special{sh 1}%
\special{ar 1200 3000 10 10 0  6.28318530717959E+0000}%
\special{sh 1}%
\special{ar 1000 3000 10 10 0  6.28318530717959E+0000}%
%
\special{pn 8}%
\special{pa 1800 2000}%
\special{pa 1800 2000}%
\special{pa 1800 3000}%
\special{pa 2400 3000}%
\special{pa 2400 2000}%
\special{pa 1800 2000}%
\special{fp}%
%
\special{pn 8}%
\special{pa 2800 2000}%
\special{pa 2800 3000}%
\special{pa 3400 3000}%
\special{pa 3400 2000}%
\special{pa 2800 2000}%
\special{fp}%
%
\special{pn 8}%
\special{pa 1400 2000}%
\special{pa 1400 3000}%
\special{fp}%
\special{pa 1400 3000}%
\special{pa 1000 3000}%
\special{fp}%
%
\special{pn 8}%
\special{pa 1000 3200}%
\special{pa 1000 800}%
\special{fp}%
\special{sh 1}%
\special{pa 1000 800}%
\special{pa 980 868}%
\special{pa 1000 854}%
\special{pa 1020 868}%
\special{pa 1000 800}%
\special{fp}%
%
\special{pn 8}%
\special{pa 4800 1100}%
\special{pa 5000 1100}%
\special{pa 5000 3100}%
\special{pa 4800 3100}%
\special{pa 4800 1100}%
\special{fp}%
\put(39.0000,-22.0000){\makebox(0,0){$0$}}%
\put(39.0000,-24.0000){\makebox(0,0){$0$}}%
\put(39.0000,-26.0000){\makebox(0,0){$0$}}%
\put(39.0000,-28.0000){\makebox(0,0){$0$}}%
\put(39.0000,-30.0000){\makebox(0,0){$0$}}%
\put(7.4000,-29.9000){\makebox(0,0){$-2\tau+1$}}%
\put(21.0000,-32.0000){\makebox(0,0){$B_{0,-2}$}}%
\put(31.0000,-32.0000){\makebox(0,0){$B_{1,-2}$}}%
\put(44.0000,-32.0000){\makebox(0,0){$B_{n'-\gamma'-1,-2}$}}%
\put(56.0000,-32.0000){\makebox(0,0){$B_{n'-\gamma',-2}$}}%
\end{picture}%

%% file: height4.tex
\unitlength 0.1in
\begin{picture}( 60.4000, 26.9500)(  2.9000,-34.0500)
\put(7.3500,-19.8500){\makebox(0,0){$-\tau+1$}}%
\put(12.0500,-8.0500){\makebox(0,0){$0$}}%
\put(8.5000,-7.9500){\makebox(0,0){$j$}}%
\put(63.3000,-11.2000){\makebox(0,0){$i$}}%
%
\special{pn 8}%
\special{pa 800 1000}%
\special{pa 6330 1000}%
\special{fp}%
\special{sh 1}%
\special{pa 6330 1000}%
\special{pa 6264 980}%
\special{pa 6278 1000}%
\special{pa 6264 1020}%
\special{pa 6330 1000}%
\special{fp}%
\put(8.9000,-11.9000){\makebox(0,0){$0$}}%
\put(7.3500,-17.8500){\makebox(0,0){$-\tau+2$}}%
\put(14.0500,-12.0500){\makebox(0,0){$0$}}%
\put(16.0500,-18.0500){\makebox(0,0){$0$}}%
\put(16.0500,-20.0500){\makebox(0,0){$0$}}%
\put(16.0500,-16.0500){\makebox(0,0){$0$}}%
\put(14.1000,-16.1000){\makebox(0,0){$0$}}%
\put(23.0000,-16.0000){\makebox(0,0){$\cdots$}}%
\put(31.2000,-16.0000){\makebox(0,0){$0$}}%
\put(31.2000,-18.0000){\makebox(0,0){$0$}}%
\put(31.2000,-14.0000){\makebox(0,0){$0$}}%
\put(39.1000,-8.1000){\makebox(0,0){$(n'-\gamma'+1)\tau$}}%
%
\special{pn 8}%
\special{pa 1200 940}%
\special{pa 1200 1000}%
\special{fp}%
%
\special{pn 8}%
\special{pa 3900 940}%
\special{pa 3900 1000}%
\special{fp}%
\put(16.1000,-14.1000){\makebox(0,0){$0$}}%
\put(14.1000,-14.1000){\makebox(0,0){$0$}}%
%
\special{pn 20}%
\special{sh 1}%
\special{ar 3106 1200 10 10 0  6.28318530717959E+0000}%
\special{sh 1}%
\special{ar 3306 1200 10 10 0  6.28318530717959E+0000}%
\special{sh 1}%
\special{ar 3506 1200 10 10 0  6.28318530717959E+0000}%
\special{sh 1}%
\special{ar 3706 1200 10 10 0  6.28318530717959E+0000}%
\special{sh 1}%
\special{ar 3906 1200 10 10 0  6.28318530717959E+0000}%
\special{sh 1}%
\special{ar 3906 1400 10 10 0  6.28318530717959E+0000}%
\special{sh 1}%
\special{ar 3706 1600 10 10 0  6.28318530717959E+0000}%
\special{sh 1}%
\special{ar 3506 1800 10 10 0  6.28318530717959E+0000}%
\special{sh 1}%
\special{ar 3306 2000 10 10 0  6.28318530717959E+0000}%
\special{sh 1}%
\special{ar 3306 2200 10 10 0  6.28318530717959E+0000}%
\special{sh 1}%
\special{ar 3506 2200 10 10 0  6.28318530717959E+0000}%
\special{sh 1}%
\special{ar 3706 2200 10 10 0  6.28318530717959E+0000}%
\special{sh 1}%
\special{ar 3906 2200 10 10 0  6.28318530717959E+0000}%
\special{sh 1}%
\special{ar 3906 2200 10 10 0  6.28318530717959E+0000}%
%
\special{pn 8}%
\special{pa 3906 2200}%
\special{pa 3306 2200}%
\special{pa 3306 2000}%
\special{pa 3106 1200}%
\special{pa 3906 1200}%
\special{pa 3906 2200}%
\special{fp}%
\put(31.2000,-20.0000){\makebox(0,0){$0$}}%
%
\special{pn 20}%
\special{sh 1}%
\special{ar 1400 2000 10 10 0  6.28318530717959E+0000}%
\special{sh 1}%
\special{ar 1200 2000 10 10 0  6.28318530717959E+0000}%
\special{sh 1}%
\special{ar 1000 2000 10 10 0  6.28318530717959E+0000}%
\special{sh 1}%
\special{ar 1200 1200 10 10 0  6.28318530717959E+0000}%
\special{sh 1}%
\special{ar 1000 1400 10 10 0  6.28318530717959E+0000}%
\special{sh 1}%
\special{ar 1200 1000 10 10 0  6.28318530717959E+0000}%
%
\special{pn 8}%
\special{pa 1000 2000}%
\special{pa 1400 2000}%
\special{fp}%
\special{pa 1400 2000}%
\special{pa 1200 1200}%
\special{fp}%
\special{pa 1200 1200}%
\special{pa 1200 1000}%
\special{fp}%
%
\special{pn 8}%
\special{pa 3310 2200}%
\special{pa 3510 2800}%
\special{pa 3510 3000}%
\special{pa 4110 3000}%
\special{pa 3910 2400}%
\special{pa 3910 2200}%
\special{pa 3310 2200}%
\special{fp}%
%
\special{pn 20}%
\special{sh 1}%
\special{ar 3910 2400 10 10 0  6.28318530717959E+0000}%
\special{sh 1}%
\special{ar 3710 2600 10 10 0  6.28318530717959E+0000}%
\special{sh 1}%
\special{ar 3510 2800 10 10 0  6.28318530717959E+0000}%
\special{sh 1}%
\special{ar 3510 3000 10 10 0  6.28318530717959E+0000}%
\special{sh 1}%
\special{ar 3710 3000 10 10 0  6.28318530717959E+0000}%
\special{sh 1}%
\special{ar 3910 3000 10 10 0  6.28318530717959E+0000}%
\special{sh 1}%
\special{ar 4110 3000 10 10 0  6.28318530717959E+0000}%
\put(31.1500,-22.0000){\makebox(0,0){$0$}}%
\put(31.2000,-24.0000){\makebox(0,0){$0$}}%
\put(31.2000,-26.0000){\makebox(0,0){$0$}}%
\put(33.1000,-30.0000){\makebox(0,0){$0$}}%
\put(33.2000,-28.0000){\makebox(0,0){$0$}}%
%
\special{pn 8}%
\special{pa 1400 2000}%
\special{pa 1400 3000}%
\special{fp}%
\special{pa 1400 3000}%
\special{pa 1000 3000}%
\special{fp}%
%
\special{pn 8}%
\special{pa 1000 3200}%
\special{pa 1000 800}%
\special{fp}%
\special{sh 1}%
\special{pa 1000 800}%
\special{pa 980 868}%
\special{pa 1000 854}%
\special{pa 1020 868}%
\special{pa 1000 800}%
\special{fp}%
%
\special{pn 8}%
\special{pa 3006 1100}%
\special{pa 3206 1100}%
\special{pa 3206 3100}%
\special{pa 3006 3100}%
\special{pa 3006 1100}%
\special{fp}%
\put(7.4000,-29.9000){\makebox(0,0){$-2\tau+1$}}%
\put(38.0000,-32.0000){\makebox(0,0){$B_{n'-\gamma',-2}$}}%
%
\special{pn 20}%
\special{sh 1}%
\special{ar 3100 3000 10 10 0  6.28318530717959E+0000}%
\special{sh 1}%
\special{ar 2900 3000 10 10 0  6.28318530717959E+0000}%
\special{sh 1}%
\special{ar 2900 2200 10 10 0  6.28318530717959E+0000}%
\special{sh 1}%
\special{ar 2700 2400 10 10 0  6.28318530717959E+0000}%
\special{sh 1}%
\special{ar 2700 3000 10 10 0  6.28318530717959E+0000}%
\special{sh 1}%
\special{ar 2900 2000 10 10 0  6.28318530717959E+0000}%
\special{sh 1}%
\special{ar 2700 2000 10 10 0  6.28318530717959E+0000}%
\special{sh 1}%
\special{ar 2700 1200 10 10 0  6.28318530717959E+0000}%
\special{sh 1}%
\special{ar 2700 1400 10 10 0  6.28318530717959E+0000}%
\special{sh 1}%
\special{ar 2500 1200 10 10 0  6.28318530717959E+0000}%
%
\special{pn 8}%
\special{pa 2700 1200}%
\special{pa 2700 1400}%
\special{fp}%
%
\special{pn 8}%
\special{pa 2700 1400}%
\special{pa 2900 2000}%
\special{fp}%
%
\special{pn 8}%
\special{pa 2700 1200}%
\special{pa 2500 1200}%
\special{dt 0.045}%
%
\special{pn 8}%
\special{pa 2900 2000}%
\special{pa 2700 2000}%
\special{fp}%
\special{pa 2900 2000}%
\special{pa 2900 2200}%
\special{fp}%
\special{pa 2900 2200}%
\special{pa 3100 3000}%
\special{fp}%
%
\special{pn 8}%
\special{pa 3100 3000}%
\special{pa 3100 3000}%
\special{fp}%
\special{pa 2700 3000}%
\special{pa 2700 3000}%
\special{fp}%
%
\special{pn 8}%
\special{pa 3100 3000}%
\special{pa 2700 3000}%
\special{fp}%
%
\special{pn 8}%
\special{pa 2700 2000}%
\special{pa 2700 2000}%
\special{dt 0.045}%
\special{pa 2500 2000}%
\special{pa 2700 2000}%
\special{dt 0.045}%
%
\special{pn 8}%
\special{pa 3100 3000}%
\special{pa 3100 3000}%
\special{fp}%
%
\special{pn 8}%
\special{pa 2700 3000}%
\special{pa 2500 3000}%
\special{dt 0.045}%
%
\special{pn 20}%
\special{sh 1}%
\special{ar 5110 3006 10 10 0  6.28318530717959E+0000}%
\special{sh 1}%
\special{ar 5310 3006 10 10 0  6.28318530717959E+0000}%
\special{sh 1}%
\special{ar 5310 3206 10 10 0  6.28318530717959E+0000}%
\special{sh 1}%
\special{ar 5110 3406 10 10 0  6.28318530717959E+0000}%
\special{sh 1}%
\special{ar 5110 2406 10 10 0  6.28318530717959E+0000}%
\special{sh 1}%
\special{ar 5110 2206 10 10 0  6.28318530717959E+0000}%
\special{sh 1}%
\special{ar 4910 2206 10 10 0  6.28318530717959E+0000}%
\special{sh 1}%
\special{ar 5710 3206 10 10 0  6.28318530717959E+0000}%
\special{sh 1}%
\special{ar 5910 3206 10 10 0  6.28318530717959E+0000}%
\special{sh 1}%
\special{ar 6110 3206 10 10 0  6.28318530717959E+0000}%
\special{sh 1}%
\special{ar 6310 3206 10 10 0  6.28318530717959E+0000}%
\special{sh 1}%
\special{ar 5710 3006 10 10 0  6.28318530717959E+0000}%
\special{sh 1}%
\special{ar 5910 2806 10 10 0  6.28318530717959E+0000}%
\special{sh 1}%
\special{ar 6110 2606 10 10 0  6.28318530717959E+0000}%
\special{sh 1}%
\special{ar 6310 2406 10 10 0  6.28318530717959E+0000}%
\special{sh 1}%
\special{ar 6310 2206 10 10 0  6.28318530717959E+0000}%
\special{sh 1}%
\special{ar 6110 2206 10 10 0  6.28318530717959E+0000}%
\special{sh 1}%
\special{ar 5910 2206 10 10 0  6.28318530717959E+0000}%
\special{sh 1}%
\special{ar 5710 2206 10 10 0  6.28318530717959E+0000}%
\special{sh 1}%
\special{ar 5510 2206 10 10 0  6.28318530717959E+0000}%
%
\special{pn 8}%
\special{pa 5110 2206}%
\special{pa 5110 2206}%
\special{dt 0.045}%
\special{pa 4910 2206}%
\special{pa 5110 2206}%
\special{dt 0.045}%
%
\special{pn 8}%
\special{pa 5110 2206}%
\special{pa 5110 2406}%
\special{fp}%
%
\special{pn 8}%
\special{pa 5310 3000}%
\special{pa 5110 3000}%
\special{dt 0.045}%
%
\special{pn 8}%
\special{pa 5110 2406}%
\special{pa 5310 3006}%
\special{fp}%
%
\special{pn 8}%
\special{pa 5510 2206}%
\special{pa 5710 3006}%
\special{pa 5710 3206}%
\special{pa 6310 3206}%
\special{pa 6310 2206}%
\special{pa 5510 2206}%
\special{fp}%
%
\special{pn 20}%
\special{sh 1}%
\special{ar 5510 2006 10 10 0  6.28318530717959E+0000}%
\special{sh 1}%
\special{ar 5710 1806 10 10 0  6.28318530717959E+0000}%
\special{sh 1}%
\special{ar 5910 1606 10 10 0  6.28318530717959E+0000}%
\special{sh 1}%
\special{ar 4910 1206 10 10 0  6.28318530717959E+0000}%
\special{sh 1}%
\special{ar 5110 1206 10 10 0  6.28318530717959E+0000}%
\special{sh 1}%
\special{ar 5110 1406 10 10 0  6.28318530717959E+0000}%
\special{sh 1}%
\special{ar 4910 1606 10 10 0  6.28318530717959E+0000}%
\special{sh 1}%
\special{ar 5510 1206 10 10 0  6.28318530717959E+0000}%
\special{sh 1}%
\special{ar 5710 1206 10 10 0  6.28318530717959E+0000}%
\special{sh 1}%
\special{ar 5910 1206 10 10 0  6.28318530717959E+0000}%
\special{sh 1}%
\special{ar 6110 1206 10 10 0  6.28318530717959E+0000}%
\special{sh 1}%
\special{ar 6110 1406 10 10 0  6.28318530717959E+0000}%
\special{sh 1}%
\special{ar 6110 1406 10 10 0  6.28318530717959E+0000}%
%
\special{pn 8}%
\special{pa 5510 1206}%
\special{pa 5510 2206}%
\special{pa 6310 2206}%
\special{pa 6110 1406}%
\special{pa 6110 1206}%
\special{pa 5510 1206}%
\special{fp}%
%
\special{pn 8}%
\special{pa 4910 1200}%
\special{pa 5110 1200}%
\special{dt 0.045}%
%
\special{pn 8}%
\special{pa 5110 1406}%
\special{pa 5110 2206}%
\special{fp}%
%
\special{pn 8}%
\special{pa 5110 1206}%
\special{pa 5110 1406}%
\special{fp}%
\put(53.0000,-14.0500){\makebox(0,0){$0$}}%
\put(53.0000,-12.0500){\makebox(0,0){$0$}}%
\put(53.0000,-16.0500){\makebox(0,0){$0$}}%
\put(53.0000,-18.0500){\makebox(0,0){$0$}}%
\put(53.0000,-20.0500){\makebox(0,0){$0$}}%
\put(53.0000,-22.0500){\makebox(0,0){$0$}}%
\put(53.0000,-24.0500){\makebox(0,0){$0$}}%
\put(55.0000,-26.0500){\makebox(0,0){$0$}}%
\put(55.0000,-28.1000){\makebox(0,0){$0$}}%
\put(55.0000,-30.1000){\makebox(0,0){$0$}}%
\put(55.0000,-32.0000){\makebox(0,0){$0$}}%
%
\special{pn 8}%
\special{pa 5300 3000}%
\special{pa 5300 3200}%
\special{dt 0.045}%
%
\special{pn 8}%
\special{pa 1600 1200}%
\special{pa 2000 1200}%
\special{fp}%
\special{pa 1600 1200}%
\special{pa 1800 1800}%
\special{fp}%
\special{pa 1800 1800}%
\special{pa 1800 2000}%
\special{fp}%
\special{pa 1800 2000}%
\special{pa 1800 3000}%
\special{fp}%
\special{pa 1800 3000}%
\special{pa 2000 3000}%
\special{fp}%
\special{pa 1800 2000}%
\special{pa 2000 2000}%
\special{fp}%
%
\special{pn 8}%
\special{pa 2000 1200}%
\special{pa 2000 1200}%
\special{fp}%
\special{pa 2000 1200}%
\special{pa 2200 1200}%
\special{fp}%
%
\special{pn 8}%
\special{pa 2000 2000}%
\special{pa 2200 2000}%
\special{dt 0.045}%
\special{pa 2000 3000}%
\special{pa 2200 3000}%
\special{dt 0.045}%
%
\special{pn 20}%
\special{sh 1}%
\special{ar 1600 1200 10 10 0  6.28318530717959E+0000}%
\special{sh 1}%
\special{ar 1800 1200 10 10 0  6.28318530717959E+0000}%
\special{sh 1}%
\special{ar 2000 1200 10 10 0  6.28318530717959E+0000}%
\special{sh 1}%
\special{ar 1800 1800 10 10 0  6.28318530717959E+0000}%
\special{sh 1}%
\special{ar 2000 1600 10 10 0  6.28318530717959E+0000}%
\special{sh 1}%
\special{ar 1800 2000 10 10 0  6.28318530717959E+0000}%
\special{sh 1}%
\special{ar 2000 2000 10 10 0  6.28318530717959E+0000}%
\special{sh 1}%
\special{ar 1800 3000 10 10 0  6.28318530717959E+0000}%
\special{sh 1}%
\special{ar 2000 3000 10 10 0  6.28318530717959E+0000}%
\special{sh 1}%
\special{ar 1800 2800 10 10 0  6.28318530717959E+0000}%
\special{sh 1}%
\special{ar 2000 2600 10 10 0  6.28318530717959E+0000}%
\special{sh 1}%
\special{ar 1400 2200 10 10 0  6.28318530717959E+0000}%
\special{sh 1}%
\special{ar 1200 2400 10 10 0  6.28318530717959E+0000}%
\special{sh 1}%
\special{ar 1000 2600 10 10 0  6.28318530717959E+0000}%
%
\special{pn 20}%
\special{sh 1}%
\special{ar 4300 1200 10 10 0  6.28318530717959E+0000}%
\special{sh 1}%
\special{ar 4500 1200 10 10 0  6.28318530717959E+0000}%
\special{sh 1}%
\special{ar 4300 2200 10 10 0  6.28318530717959E+0000}%
\special{sh 1}%
\special{ar 4300 2000 10 10 0  6.28318530717959E+0000}%
\special{sh 1}%
\special{ar 4500 1800 10 10 0  6.28318530717959E+0000}%
\special{sh 1}%
\special{ar 4500 2200 10 10 0  6.28318530717959E+0000}%
\special{sh 1}%
\special{ar 4500 3000 10 10 0  6.28318530717959E+0000}%
\special{sh 1}%
\special{ar 4500 2800 10 10 0  6.28318530717959E+0000}%
\special{sh 1}%
\special{ar 4700 3000 10 10 0  6.28318530717959E+0000}%
%
\special{pn 8}%
\special{pa 4300 1200}%
\special{pa 4300 2200}%
\special{fp}%
\special{pa 4300 2200}%
\special{pa 4500 2800}%
\special{fp}%
\special{pa 4500 2800}%
\special{pa 4500 3000}%
\special{fp}%
\special{pa 4500 3000}%
\special{pa 4700 3000}%
\special{fp}%
\special{pa 4300 2200}%
\special{pa 4500 2200}%
\special{fp}%
\special{pa 4300 1200}%
\special{pa 4500 1200}%
\special{fp}%
\put(23.0000,-24.0000){\makebox(0,0){$\cdots$}}%
\put(49.0000,-26.0000){\makebox(0,0){$\cdots$}}%
\put(47.0000,-16.0000){\makebox(0,0){$\cdots$}}%
\put(60.0000,-34.0000){\makebox(0,0){$B_{2n'-2\gamma'-1,-2}$}}%
\put(61.1000,-8.1000){\makebox(0,0){$(2n'-2\gamma')\tau$}}%
%
\special{pn 8}%
\special{pa 6100 940}%
\special{pa 6100 1000}%
\special{fp}%
\put(29.0000,-14.0000){\makebox(0,0){$0$}}%
\put(29.0000,-12.0000){\makebox(0,0){$0$}}%
\put(16.1000,-22.1000){\makebox(0,0){$0$}}%
\put(16.1000,-24.1000){\makebox(0,0){$0$}}%
\put(16.1000,-26.1000){\makebox(0,0){$0$}}%
\put(16.1000,-28.1000){\makebox(0,0){$0$}}%
\put(16.1000,-30.1000){\makebox(0,0){$0$}}%
\put(43.0000,-30.0000){\makebox(0,0){$0$}}%
\put(43.0000,-28.0000){\makebox(0,0){$0$}}%
\put(41.0000,-26.0000){\makebox(0,0){$0$}}%
\put(41.0000,-24.0000){\makebox(0,0){$0$}}%
\put(41.0000,-22.0000){\makebox(0,0){$0$}}%
\put(41.0000,-20.0000){\makebox(0,0){$0$}}%
\put(41.0000,-18.0000){\makebox(0,0){$0$}}%
\put(41.0000,-16.0000){\makebox(0,0){$0$}}%
\put(41.0000,-14.0000){\makebox(0,0){$0$}}%
\put(41.0000,-12.0000){\makebox(0,0){$0$}}%
\end{picture}%

%% file: ee=-1gen.tex
\unitlength 0.1in
\begin{picture}( 38.0500, 16.8500)(  6.6500,-34.0000)
%
\special{pn 8}%
\special{pa 800 2800}%
\special{pa 800 2800}%
\special{fp}%
%
\special{pn 8}%
\special{pa 1000 3400}%
\special{pa 1000 1740}%
\special{fp}%
\special{sh 1}%
\special{pa 1000 1740}%
\special{pa 980 1808}%
\special{pa 1000 1794}%
\special{pa 1020 1808}%
\special{pa 1000 1740}%
\special{fp}%
%
\special{pn 8}%
\special{pa 800 2800}%
\special{pa 4470 2800}%
\special{fp}%
\special{sh 1}%
\special{pa 4470 2800}%
\special{pa 4404 2780}%
\special{pa 4418 2800}%
\special{pa 4404 2820}%
\special{pa 4470 2800}%
\special{fp}%
%
\special{pn 20}%
\special{sh 1}%
\special{ar 1000 2800 10 10 0  6.28318530717959E+0000}%
\special{sh 1}%
\special{ar 800 2600 10 10 0  6.28318530717959E+0000}%
\special{sh 1}%
\special{ar 1000 3000 10 10 0  6.28318530717959E+0000}%
\special{sh 1}%
\special{ar 1800 2800 10 10 0  6.28318530717959E+0000}%
\special{sh 1}%
\special{ar 1600 2600 10 10 0  6.28318530717959E+0000}%
\special{sh 1}%
\special{ar 1400 2400 10 10 0  6.28318530717959E+0000}%
\special{sh 1}%
\special{ar 1200 2200 10 10 0  6.28318530717959E+0000}%
\special{sh 1}%
\special{ar 1200 2000 10 10 0  6.28318530717959E+0000}%
\special{sh 1}%
\special{ar 1400 2000 10 10 0  6.28318530717959E+0000}%
\special{sh 1}%
\special{ar 1600 2000 10 10 0  6.28318530717959E+0000}%
\special{sh 1}%
\special{ar 2000 3000 10 10 0  6.28318530717959E+0000}%
\special{sh 1}%
\special{ar 2000 3200 10 10 0  6.28318530717959E+0000}%
\special{sh 1}%
\special{ar 1800 3200 10 10 0  6.28318530717959E+0000}%
\special{sh 1}%
\special{ar 1600 3200 10 10 0  6.28318530717959E+0000}%
\special{sh 1}%
\special{ar 1400 3200 10 10 0  6.28318530717959E+0000}%
\special{sh 1}%
\special{ar 1800 2000 10 10 0  6.28318530717959E+0000}%
\special{sh 1}%
\special{ar 800 3000 10 10 0  6.28318530717959E+0000}%
%
\special{pn 8}%
\special{pa 1200 2000}%
\special{pa 1800 2000}%
\special{pa 2000 3000}%
\special{pa 2000 3200}%
\special{pa 1400 3200}%
\special{pa 1200 2200}%
\special{pa 1200 2000}%
\special{fp}%
\put(20.0000,-26.0000){\makebox(0,0){$\cdots$}}%
%
\special{pn 20}%
\special{sh 1}%
\special{ar 2200 2000 10 10 0  6.28318530717959E+0000}%
\special{sh 1}%
\special{ar 2400 2000 10 10 0  6.28318530717959E+0000}%
\special{sh 1}%
\special{ar 2600 2000 10 10 0  6.28318530717959E+0000}%
\special{sh 1}%
\special{ar 2200 2200 10 10 0  6.28318530717959E+0000}%
\special{sh 1}%
\special{ar 2400 2400 10 10 0  6.28318530717959E+0000}%
\special{sh 1}%
\special{ar 2600 2600 10 10 0  6.28318530717959E+0000}%
\special{sh 1}%
\special{ar 2800 2800 10 10 0  6.28318530717959E+0000}%
\special{sh 1}%
\special{ar 3000 3000 10 10 0  6.28318530717959E+0000}%
\special{sh 1}%
\special{ar 3000 3200 10 10 0  6.28318530717959E+0000}%
\special{sh 1}%
\special{ar 2800 3200 10 10 0  6.28318530717959E+0000}%
\special{sh 1}%
\special{ar 2600 3200 10 10 0  6.28318530717959E+0000}%
\special{sh 1}%
\special{ar 2400 3200 10 10 0  6.28318530717959E+0000}%
%
\special{pn 8}%
\special{pa 2200 2000}%
\special{pa 2600 2000}%
\special{pa 3000 3000}%
\special{pa 3000 3200}%
\special{pa 2400 3200}%
\special{pa 2200 2200}%
\special{pa 2200 2200}%
\special{pa 2200 2000}%
\special{fp}%
%
\special{pn 20}%
\special{sh 1}%
\special{ar 3000 2000 10 10 0  6.28318530717959E+0000}%
\special{sh 1}%
\special{ar 3200 2000 10 10 0  6.28318530717959E+0000}%
\special{sh 1}%
\special{ar 3400 2000 10 10 0  6.28318530717959E+0000}%
\special{sh 1}%
\special{ar 3600 2000 10 10 0  6.28318530717959E+0000}%
\special{sh 1}%
\special{ar 3000 2200 10 10 0  6.28318530717959E+0000}%
\special{sh 1}%
\special{ar 3200 2400 10 10 0  6.28318530717959E+0000}%
\special{sh 1}%
\special{ar 3400 2600 10 10 0  6.28318530717959E+0000}%
\special{sh 1}%
\special{ar 3600 2800 10 10 0  6.28318530717959E+0000}%
\special{sh 1}%
\special{ar 3800 3000 10 10 0  6.28318530717959E+0000}%
\special{sh 1}%
\special{ar 3800 3200 10 10 0  6.28318530717959E+0000}%
\special{sh 1}%
\special{ar 3600 3200 10 10 0  6.28318530717959E+0000}%
\special{sh 1}%
\special{ar 3400 3200 10 10 0  6.28318530717959E+0000}%
%
\special{pn 8}%
\special{pa 3000 2000}%
\special{pa 3600 2000}%
\special{pa 3800 3000}%
\special{pa 3800 3200}%
\special{pa 3400 3200}%
\special{pa 3000 2200}%
\special{pa 3000 2000}%
\special{fp}%
%
\special{pn 8}%
\special{pa 2900 2100}%
\special{pa 3100 2100}%
\special{pa 3100 3100}%
\special{pa 2900 3100}%
\special{pa 2900 2100}%
\special{fp}%
\put(40.0000,-26.0000){\makebox(0,0){$\cdots$}}%
\put(32.0000,-18.0000){\makebox(0,0){$B_{\lfloor\frac{n}{\alpha'}\rfloor,0}$}}%
\put(14.0000,-18.0000){\makebox(0,0){$B_{0,0}$}}%
\put(8.0000,-18.0000){\makebox(0,0){$j$}}%
\put(46.0000,-28.0000){\makebox(0,0){$i$}}%
%
\special{pn 8}%
\special{pa 800 3000}%
\special{pa 1000 3000}%
\special{fp}%
\special{pa 1000 3000}%
\special{pa 1000 2800}%
\special{fp}%
\special{pa 1000 2800}%
\special{pa 800 2000}%
\special{fp}%
\end{picture}%

%% file: ee=-1gen2.tex
\unitlength 0.1in
\begin{picture}( 56.0000, 16.8500)(  6.6500,-34.0000)
%
\special{pn 8}%
\special{pa 800 2800}%
\special{pa 800 2800}%
\special{fp}%
%
\special{pn 8}%
\special{pa 1000 3400}%
\special{pa 1000 1740}%
\special{fp}%
\special{sh 1}%
\special{pa 1000 1740}%
\special{pa 980 1808}%
\special{pa 1000 1794}%
\special{pa 1020 1808}%
\special{pa 1000 1740}%
\special{fp}%
%
\special{pn 8}%
\special{pa 800 2800}%
\special{pa 6200 2800}%
\special{fp}%
\special{sh 1}%
\special{pa 6200 2800}%
\special{pa 6134 2780}%
\special{pa 6148 2800}%
\special{pa 6134 2820}%
\special{pa 6200 2800}%
\special{fp}%
%
\special{pn 20}%
\special{sh 1}%
\special{ar 1000 2800 10 10 0  6.28318530717959E+0000}%
\special{sh 1}%
\special{ar 800 2600 10 10 0  6.28318530717959E+0000}%
\special{sh 1}%
\special{ar 1000 3000 10 10 0  6.28318530717959E+0000}%
\special{sh 1}%
\special{ar 1800 2800 10 10 0  6.28318530717959E+0000}%
\special{sh 1}%
\special{ar 1600 2600 10 10 0  6.28318530717959E+0000}%
\special{sh 1}%
\special{ar 1400 2400 10 10 0  6.28318530717959E+0000}%
\special{sh 1}%
\special{ar 1200 2200 10 10 0  6.28318530717959E+0000}%
\special{sh 1}%
\special{ar 1200 2000 10 10 0  6.28318530717959E+0000}%
\special{sh 1}%
\special{ar 1400 2000 10 10 0  6.28318530717959E+0000}%
\special{sh 1}%
\special{ar 1600 2000 10 10 0  6.28318530717959E+0000}%
\special{sh 1}%
\special{ar 2000 3000 10 10 0  6.28318530717959E+0000}%
\special{sh 1}%
\special{ar 2000 3200 10 10 0  6.28318530717959E+0000}%
\special{sh 1}%
\special{ar 1800 3200 10 10 0  6.28318530717959E+0000}%
\special{sh 1}%
\special{ar 1600 3200 10 10 0  6.28318530717959E+0000}%
\special{sh 1}%
\special{ar 1400 3200 10 10 0  6.28318530717959E+0000}%
\special{sh 1}%
\special{ar 1800 2000 10 10 0  6.28318530717959E+0000}%
\special{sh 1}%
\special{ar 800 3000 10 10 0  6.28318530717959E+0000}%
%
\special{pn 8}%
\special{pa 1200 2000}%
\special{pa 1800 2000}%
\special{pa 2000 3000}%
\special{pa 2000 3200}%
\special{pa 1400 3200}%
\special{pa 1200 2200}%
\special{pa 1200 2000}%
\special{fp}%
\put(20.0000,-26.0000){\makebox(0,0){$\cdots$}}%
%
\special{pn 20}%
\special{sh 1}%
\special{ar 2200 2000 10 10 0  6.28318530717959E+0000}%
\special{sh 1}%
\special{ar 2400 2000 10 10 0  6.28318530717959E+0000}%
\special{sh 1}%
\special{ar 2600 2000 10 10 0  6.28318530717959E+0000}%
\special{sh 1}%
\special{ar 2200 2200 10 10 0  6.28318530717959E+0000}%
\special{sh 1}%
\special{ar 2400 2400 10 10 0  6.28318530717959E+0000}%
\special{sh 1}%
\special{ar 2600 2600 10 10 0  6.28318530717959E+0000}%
\special{sh 1}%
\special{ar 2800 2800 10 10 0  6.28318530717959E+0000}%
\special{sh 1}%
\special{ar 3000 3000 10 10 0  6.28318530717959E+0000}%
\special{sh 1}%
\special{ar 3000 3200 10 10 0  6.28318530717959E+0000}%
\special{sh 1}%
\special{ar 2800 3200 10 10 0  6.28318530717959E+0000}%
\special{sh 1}%
\special{ar 2600 3200 10 10 0  6.28318530717959E+0000}%
\special{sh 1}%
\special{ar 2400 3200 10 10 0  6.28318530717959E+0000}%
%
\special{pn 8}%
\special{pa 2200 2000}%
\special{pa 2600 2000}%
\special{pa 3000 3000}%
\special{pa 3000 3200}%
\special{pa 2400 3200}%
\special{pa 2200 2200}%
\special{pa 2200 2200}%
\special{pa 2200 2000}%
\special{fp}%
%
\special{pn 20}%
\special{sh 1}%
\special{ar 3000 2000 10 10 0  6.28318530717959E+0000}%
\special{sh 1}%
\special{ar 3200 2000 10 10 0  6.28318530717959E+0000}%
\special{sh 1}%
\special{ar 3400 2000 10 10 0  6.28318530717959E+0000}%
\special{sh 1}%
\special{ar 3600 2000 10 10 0  6.28318530717959E+0000}%
\special{sh 1}%
\special{ar 3000 2200 10 10 0  6.28318530717959E+0000}%
\special{sh 1}%
\special{ar 3200 2400 10 10 0  6.28318530717959E+0000}%
\special{sh 1}%
\special{ar 3400 2600 10 10 0  6.28318530717959E+0000}%
\special{sh 1}%
\special{ar 3600 2800 10 10 0  6.28318530717959E+0000}%
\special{sh 1}%
\special{ar 3800 3000 10 10 0  6.28318530717959E+0000}%
\special{sh 1}%
\special{ar 3800 3200 10 10 0  6.28318530717959E+0000}%
\special{sh 1}%
\special{ar 3600 3200 10 10 0  6.28318530717959E+0000}%
\special{sh 1}%
\special{ar 3400 3200 10 10 0  6.28318530717959E+0000}%
%
\special{pn 8}%
\special{pa 3000 2000}%
\special{pa 3600 2000}%
\special{pa 3800 3000}%
\special{pa 3800 3200}%
\special{pa 3400 3200}%
\special{pa 3000 2200}%
\special{pa 3000 2000}%
\special{fp}%
%
\special{pn 8}%
\special{pa 2900 2100}%
\special{pa 3100 2100}%
\special{pa 3100 3100}%
\special{pa 2900 3100}%
\special{pa 2900 2100}%
\special{fp}%
\put(40.0000,-26.0000){\makebox(0,0){$\cdots$}}%
\put(32.0000,-18.0000){\makebox(0,0){$B_{\lfloor\frac{n}{\alpha'}\rfloor,0}$}}%
\put(14.0000,-18.0000){\makebox(0,0){$B_{0,0}$}}%
\put(8.0000,-18.0000){\makebox(0,0){$j$}}%
\put(64.0000,-28.0000){\makebox(0,0){$i$}}%
%
\special{pn 8}%
\special{pa 800 3000}%
\special{pa 1000 3000}%
\special{fp}%
\special{pa 1000 3000}%
\special{pa 1000 2800}%
\special{fp}%
\special{pa 1000 2800}%
\special{pa 800 2000}%
\special{fp}%
%
\special{pn 8}%
\special{pa 4000 2000}%
\special{pa 4600 2000}%
\special{pa 4800 3000}%
\special{pa 4800 3200}%
\special{pa 4200 3200}%
\special{pa 4000 2200}%
\special{pa 4000 2000}%
\special{fp}%
%
\special{pn 20}%
\special{sh 1}%
\special{ar 4600 2000 10 10 0  6.28318530717959E+0000}%
%
\special{pn 20}%
\special{sh 1}%
\special{ar 4200 3200 10 10 0  6.28318530717959E+0000}%
%
\special{pn 20}%
\special{sh 1}%
\special{ar 4400 3200 10 10 0  6.28318530717959E+0000}%
%
\special{pn 20}%
\special{sh 1}%
\special{ar 4600 3200 10 10 0  6.28318530717959E+0000}%
%
\special{pn 20}%
\special{sh 1}%
\special{ar 4800 3200 10 10 0  6.28318530717959E+0000}%
%
\special{pn 20}%
\special{sh 1}%
\special{ar 4800 3000 10 10 0  6.28318530717959E+0000}%
%
\special{pn 20}%
\special{sh 1}%
\special{ar 4400 2000 10 10 0  6.28318530717959E+0000}%
%
\special{pn 20}%
\special{sh 1}%
\special{ar 4200 2000 10 10 0  6.28318530717959E+0000}%
%
\special{pn 20}%
\special{sh 1}%
\special{ar 4000 2000 10 10 0  6.28318530717959E+0000}%
%
\special{pn 20}%
\special{sh 1}%
\special{ar 4000 2200 10 10 0  6.28318530717959E+0000}%
%
\special{pn 20}%
\special{sh 1}%
\special{ar 4200 2400 10 10 0  6.28318530717959E+0000}%
%
\special{pn 20}%
\special{sh 1}%
\special{ar 4400 2600 10 10 0  6.28318530717959E+0000}%
%
\special{pn 20}%
\special{sh 1}%
\special{ar 4600 2800 10 10 0  6.28318530717959E+0000}%
%
\special{pn 20}%
\special{sh 1}%
\special{ar 5000 2000 10 10 0  6.28318530717959E+0000}%
\special{sh 1}%
\special{ar 5200 2000 10 10 0  6.28318530717959E+0000}%
\special{sh 1}%
\special{ar 5400 2000 10 10 0  6.28318530717959E+0000}%
\special{sh 1}%
\special{ar 5000 2200 10 10 0  6.28318530717959E+0000}%
\special{sh 1}%
\special{ar 5200 2400 10 10 0  6.28318530717959E+0000}%
\special{sh 1}%
\special{ar 5400 2600 10 10 0  6.28318530717959E+0000}%
\special{sh 1}%
\special{ar 5600 2800 10 10 0  6.28318530717959E+0000}%
\special{sh 1}%
\special{ar 5800 3000 10 10 0  6.28318530717959E+0000}%
\special{sh 1}%
\special{ar 5800 3200 10 10 0  6.28318530717959E+0000}%
\special{sh 1}%
\special{ar 5600 3200 10 10 0  6.28318530717959E+0000}%
\special{sh 1}%
\special{ar 5400 3200 10 10 0  6.28318530717959E+0000}%
\special{sh 1}%
\special{ar 5200 3200 10 10 0  6.28318530717959E+0000}%
%
\special{pn 8}%
\special{pa 5000 2000}%
\special{pa 5400 2000}%
\special{pa 5800 3000}%
\special{pa 5800 3200}%
\special{pa 5200 3200}%
\special{pa 5000 2200}%
\special{pa 5000 2200}%
\special{pa 5000 2000}%
\special{fp}%
%
\special{pn 20}%
\special{sh 1}%
\special{ar 5800 2000 10 10 0  6.28318530717959E+0000}%
\special{sh 1}%
\special{ar 5800 2200 10 10 0  6.28318530717959E+0000}%
\special{sh 1}%
\special{ar 6000 2000 10 10 0  6.28318530717959E+0000}%
\special{sh 1}%
\special{ar 6200 3200 10 10 0  6.28318530717959E+0000}%
%
\special{pn 8}%
\special{pa 5800 2000}%
\special{pa 6000 2000}%
\special{fp}%
\special{pa 5800 2000}%
\special{pa 5800 2200}%
\special{fp}%
\special{pa 5800 2200}%
\special{pa 6200 3200}%
\special{fp}%
%
\special{pn 8}%
\special{pa 5700 2100}%
\special{pa 5900 2100}%
\special{pa 5900 3100}%
\special{pa 5700 3100}%
\special{pa 5700 2100}%
\special{fp}%
\put(48.0000,-26.0000){\makebox(0,0){$\cdots$}}%
\put(60.0000,-18.0000){\makebox(0,0){$B_{\lfloor\frac{2n}{\alpha'}\rfloor,0}$}}%
\end{picture}%

%% file: ntau+2taugen.tex
\unitlength 0.1in
\begin{picture}( 58.4000, 18.7500)(  2.0000,-26.0000)
\put(60.0000,-18.0000){\makebox(0,0){$i$}}%
%
\special{pn 8}%
\special{pa 200 2000}%
\special{pa 6040 2000}%
\special{fp}%
\special{sh 1}%
\special{pa 6040 2000}%
\special{pa 5974 1980}%
\special{pa 5988 2000}%
\special{pa 5974 2020}%
\special{pa 6040 2000}%
\special{fp}%
\put(8.0000,-8.1000){\makebox(0,0){$j$}}%
%
\special{pn 8}%
\special{pa 1000 2600}%
\special{pa 1000 980}%
\special{fp}%
\special{sh 1}%
\special{pa 1000 980}%
\special{pa 980 1048}%
\special{pa 1000 1034}%
\special{pa 1020 1048}%
\special{pa 1000 980}%
\special{fp}%
\put(7.0000,-21.0000){\makebox(0,0){$-1$}}%
\put(18.0000,-21.1000){\makebox(0,0){$\tau$}}%
\put(10.5000,-11.5000){\makebox(0,0){$\tau$}}%
%
\special{pn 8}%
\special{pa 3100 1300}%
\special{pa 3300 1300}%
\special{pa 3300 2300}%
\special{pa 3100 2300}%
\special{pa 3100 1300}%
\special{fp}%
\put(11.0000,-21.1000){\makebox(0,0){$0$}}%
%
\special{pn 20}%
\special{sh 1}%
\special{ar 3200 1200 10 10 0  6.28318530717959E+0000}%
\special{sh 1}%
\special{ar 3400 1200 10 10 0  6.28318530717959E+0000}%
\special{sh 1}%
\special{ar 3600 1200 10 10 0  6.28318530717959E+0000}%
\special{sh 1}%
\special{ar 3800 1200 10 10 0  6.28318530717959E+0000}%
\special{sh 1}%
\special{ar 3200 1400 10 10 0  6.28318530717959E+0000}%
\special{sh 1}%
\special{ar 3400 1600 10 10 0  6.28318530717959E+0000}%
\special{sh 1}%
\special{ar 3600 1800 10 10 0  6.28318530717959E+0000}%
\special{sh 1}%
\special{ar 3800 2000 10 10 0  6.28318530717959E+0000}%
\special{sh 1}%
\special{ar 4000 2200 10 10 0  6.28318530717959E+0000}%
\special{sh 1}%
\special{ar 4000 2400 10 10 0  6.28318530717959E+0000}%
\special{sh 1}%
\special{ar 3800 2400 10 10 0  6.28318530717959E+0000}%
\special{sh 1}%
\special{ar 3600 2400 10 10 0  6.28318530717959E+0000}%
%
\special{pn 8}%
\special{pa 3800 1200}%
\special{pa 4000 2200}%
\special{pa 4000 2400}%
\special{pa 3600 2400}%
\special{pa 3200 1400}%
\special{pa 3200 1200}%
\special{pa 3800 1200}%
\special{fp}%
%
\special{pn 8}%
\special{pa 5220 1220}%
\special{pa 5620 1020}%
\special{fp}%
\put(60.2000,-10.2000){\makebox(0,0){$[-\epsilon_2n_1q_1]_n$}}%
\put(52.0000,-12.0000){\makebox(0,0){$0$}}%
\put(58.0000,-26.0000){\makebox(0,0){$B_{\gamma'+n'-n,0}$}}%
%
\special{pn 20}%
\special{sh 1}%
\special{ar 1200 1200 10 10 0  6.28318530717959E+0000}%
\special{sh 1}%
\special{ar 1400 1200 10 10 0  6.28318530717959E+0000}%
\special{sh 1}%
\special{ar 1600 1200 10 10 0  6.28318530717959E+0000}%
\special{sh 1}%
\special{ar 1800 1200 10 10 0  6.28318530717959E+0000}%
\special{sh 1}%
\special{ar 1200 1400 10 10 0  6.28318530717959E+0000}%
\special{sh 1}%
\special{ar 1400 1600 10 10 0  6.28318530717959E+0000}%
\special{sh 1}%
\special{ar 1600 1800 10 10 0  6.28318530717959E+0000}%
\special{sh 1}%
\special{ar 1800 2000 10 10 0  6.28318530717959E+0000}%
\special{sh 1}%
\special{ar 2000 2200 10 10 0  6.28318530717959E+0000}%
\special{sh 1}%
\special{ar 2000 2400 10 10 0  6.28318530717959E+0000}%
\special{sh 1}%
\special{ar 1800 2400 10 10 0  6.28318530717959E+0000}%
\special{sh 1}%
\special{ar 1600 2400 10 10 0  6.28318530717959E+0000}%
\special{sh 1}%
\special{ar 1400 2400 10 10 0  6.28318530717959E+0000}%
%
\special{pn 8}%
\special{pa 1200 1200}%
\special{pa 1800 1200}%
\special{pa 2000 2200}%
\special{pa 2000 2400}%
\special{pa 1400 2400}%
\special{pa 1200 1400}%
\special{pa 1200 1200}%
\special{fp}%
\put(22.0000,-18.0000){\makebox(0,0){$\cdots$}}%
%
\special{pn 20}%
\special{sh 1}%
\special{ar 2400 1200 10 10 0  6.28318530717959E+0000}%
\special{sh 1}%
\special{ar 2400 1400 10 10 0  6.28318530717959E+0000}%
\special{sh 1}%
\special{ar 2600 1600 10 10 0  6.28318530717959E+0000}%
\special{sh 1}%
\special{ar 2800 1800 10 10 0  6.28318530717959E+0000}%
\special{sh 1}%
\special{ar 3000 2000 10 10 0  6.28318530717959E+0000}%
\special{sh 1}%
\special{ar 3200 2200 10 10 0  6.28318530717959E+0000}%
\special{sh 1}%
\special{ar 3200 2400 10 10 0  6.28318530717959E+0000}%
\special{sh 1}%
\special{ar 3000 2400 10 10 0  6.28318530717959E+0000}%
\special{sh 1}%
\special{ar 2800 2400 10 10 0  6.28318530717959E+0000}%
\special{sh 1}%
\special{ar 2600 2400 10 10 0  6.28318530717959E+0000}%
%
\special{pn 20}%
\special{sh 1}%
\special{ar 2600 1200 10 10 0  6.28318530717959E+0000}%
\special{sh 1}%
\special{ar 2800 1200 10 10 0  6.28318530717959E+0000}%
%
\special{pn 8}%
\special{pa 2400 1200}%
\special{pa 2400 1400}%
\special{pa 2600 2400}%
\special{pa 3200 2400}%
\special{pa 3200 2200}%
\special{pa 2800 1200}%
\special{pa 2400 1200}%
\special{fp}%
\put(42.0000,-18.0000){\makebox(0,0){$\cdots$}}%
%
\special{pn 8}%
\special{pa 4400 1200}%
\special{pa 5000 1200}%
\special{pa 5200 2200}%
\special{pa 5200 2400}%
\special{pa 4600 2400}%
\special{pa 4400 1400}%
\special{pa 4400 1200}%
\special{fp}%
%
\special{pn 20}%
\special{sh 1}%
\special{ar 4400 1200 10 10 0  6.28318530717959E+0000}%
\special{sh 1}%
\special{ar 4400 1400 10 10 0  6.28318530717959E+0000}%
\special{sh 1}%
\special{ar 4600 1600 10 10 0  6.28318530717959E+0000}%
\special{sh 1}%
\special{ar 4800 1800 10 10 0  6.28318530717959E+0000}%
\special{sh 1}%
\special{ar 5000 2000 10 10 0  6.28318530717959E+0000}%
\special{sh 1}%
\special{ar 5200 2200 10 10 0  6.28318530717959E+0000}%
\special{sh 1}%
\special{ar 5200 2400 10 10 0  6.28318530717959E+0000}%
\special{sh 1}%
\special{ar 5000 2400 10 10 0  6.28318530717959E+0000}%
\special{sh 1}%
\special{ar 4800 2400 10 10 0  6.28318530717959E+0000}%
\special{sh 1}%
\special{ar 4600 2400 10 10 0  6.28318530717959E+0000}%
%
\special{pn 20}%
\special{sh 1}%
\special{ar 4600 1200 10 10 0  6.28318530717959E+0000}%
\special{sh 1}%
\special{ar 4800 1200 10 10 0  6.28318530717959E+0000}%
\special{sh 1}%
\special{ar 5000 1200 10 10 0  6.28318530717959E+0000}%
\put(38.0000,-26.0000){\makebox(0,0){$B_{i,0}$}}%
%
\special{pn 8}%
\special{pa 5400 1400}%
\special{pa 6000 1400}%
\special{pa 6000 2400}%
\special{pa 5600 2400}%
\special{pa 5400 1600}%
\special{pa 5400 1400}%
\special{fp}%
%
\special{pn 20}%
\special{sh 1}%
\special{ar 5400 1400 10 10 0  6.28318530717959E+0000}%
\special{sh 1}%
\special{ar 5600 1400 10 10 0  6.28318530717959E+0000}%
\special{sh 1}%
\special{ar 5800 1400 10 10 0  6.28318530717959E+0000}%
\special{sh 1}%
\special{ar 6000 1400 10 10 0  6.28318530717959E+0000}%
\special{sh 1}%
\special{ar 5400 1600 10 10 0  6.28318530717959E+0000}%
\special{sh 1}%
\special{ar 5600 1800 10 10 0  6.28318530717959E+0000}%
\special{sh 1}%
\special{ar 5800 2000 10 10 0  6.28318530717959E+0000}%
\special{sh 1}%
\special{ar 6000 2200 10 10 0  6.28318530717959E+0000}%
\special{sh 1}%
\special{ar 6000 2400 10 10 0  6.28318530717959E+0000}%
\special{sh 1}%
\special{ar 5800 2400 10 10 0  6.28318530717959E+0000}%
\special{sh 1}%
\special{ar 5600 2400 10 10 0  6.28318530717959E+0000}%
%
\special{pn 8}%
\special{pa 400 1200}%
\special{pa 800 1200}%
\special{pa 1000 2000}%
\special{pa 1000 2200}%
\special{pa 400 2200}%
\special{pa 400 1200}%
\special{fp}%
%
\special{pn 20}%
\special{sh 1}%
\special{ar 400 1200 10 10 0  6.28318530717959E+0000}%
\special{sh 1}%
\special{ar 600 1200 10 10 0  6.28318530717959E+0000}%
\special{sh 1}%
\special{ar 800 1200 10 10 0  6.28318530717959E+0000}%
\special{sh 1}%
\special{ar 400 1400 10 10 0  6.28318530717959E+0000}%
\special{sh 1}%
\special{ar 600 1600 10 10 0  6.28318530717959E+0000}%
\special{sh 1}%
\special{ar 800 1800 10 10 0  6.28318530717959E+0000}%
\special{sh 1}%
\special{ar 1000 2000 10 10 0  6.28318530717959E+0000}%
\special{sh 1}%
\special{ar 1000 2200 10 10 0  6.28318530717959E+0000}%
\special{sh 1}%
\special{ar 800 2200 10 10 0  6.28318530717959E+0000}%
\special{sh 1}%
\special{ar 600 2200 10 10 0  6.28318530717959E+0000}%
\special{sh 1}%
\special{ar 400 2200 10 10 0  6.28318530717959E+0000}%
\end{picture}%

%% file: ntau+2tau.tex
\unitlength 0.1in
\begin{picture}( 38.1000, 27.0000)(  6.0000,-28.0000)
\put(12.0000,-2.0000){\makebox(0,0){$j$}}%
\put(22.0000,-22.9500){\makebox(0,0){$\tau$}}%
%
\special{pn 8}%
\special{pa 600 2196}%
\special{pa 4410 2196}%
\special{fp}%
\special{sh 1}%
\special{pa 4410 2196}%
\special{pa 4344 2176}%
\special{pa 4358 2196}%
\special{pa 4344 2216}%
\special{pa 4410 2196}%
\special{fp}%
%
\special{pn 8}%
\special{pa 1400 2800}%
\special{pa 1400 100}%
\special{fp}%
\special{sh 1}%
\special{pa 1400 100}%
\special{pa 1380 168}%
\special{pa 1400 154}%
\special{pa 1420 168}%
\special{pa 1400 100}%
\special{fp}%
\put(14.5000,-22.9500){\makebox(0,0){$0$}}%
%
\special{pn 8}%
\special{pa 2200 2170}%
\special{pa 2200 2210}%
\special{fp}%
%
\special{pn 20}%
\special{sh 1}%
\special{ar 2410 2600 10 10 0  6.28318530717959E+0000}%
%
\special{pn 20}%
\special{sh 1}%
\special{ar 1810 2600 10 10 0  6.28318530717959E+0000}%
%
\special{pn 20}%
\special{sh 1}%
\special{ar 2010 2600 10 10 0  6.28318530717959E+0000}%
%
\special{pn 20}%
\special{sh 1}%
\special{ar 2210 2600 10 10 0  6.28318530717959E+0000}%
%
\special{pn 20}%
\special{sh 1}%
\special{ar 2204 2198 10 10 0  6.28318530717959E+0000}%
%
\special{pn 20}%
\special{sh 1}%
\special{ar 2404 2398 10 10 0  6.28318530717959E+0000}%
%
\special{pn 20}%
\special{sh 1}%
\special{ar 2200 1400 10 10 0  6.28318530717959E+0000}%
%
\special{pn 20}%
\special{sh 1}%
\special{ar 1600 1400 10 10 0  6.28318530717959E+0000}%
%
\special{pn 20}%
\special{sh 1}%
\special{ar 2000 1400 10 10 0  6.28318530717959E+0000}%
%
\special{pn 20}%
\special{sh 1}%
\special{ar 1800 1400 10 10 0  6.28318530717959E+0000}%
%
\special{pn 20}%
\special{sh 1}%
\special{ar 2000 2000 10 10 0  6.28318530717959E+0000}%
%
\special{pn 20}%
\special{sh 1}%
\special{ar 1800 1800 10 10 0  6.28318530717959E+0000}%
%
\special{pn 20}%
\special{sh 1}%
\special{ar 1600 1600 10 10 0  6.28318530717959E+0000}%
\put(13.0000,-14.0000){\makebox(0,0){$\tau$}}%
%
\special{pn 8}%
\special{pa 1600 1400}%
\special{pa 1600 1600}%
\special{pa 1800 2600}%
\special{pa 2400 2600}%
\special{pa 2400 2400}%
\special{pa 2200 1400}%
\special{pa 1600 1400}%
\special{fp}%
\put(24.5000,-20.0000){\makebox(0,0){$\cdots$}}%
%
\special{pn 20}%
\special{sh 1}%
\special{ar 2600 1400 10 10 0  6.28318530717959E+0000}%
\special{sh 1}%
\special{ar 2600 1600 10 10 0  6.28318530717959E+0000}%
\special{sh 1}%
\special{ar 2800 1800 10 10 0  6.28318530717959E+0000}%
\special{sh 1}%
\special{ar 3000 2000 10 10 0  6.28318530717959E+0000}%
\special{sh 1}%
\special{ar 3200 2200 10 10 0  6.28318530717959E+0000}%
\special{sh 1}%
\special{ar 3400 2400 10 10 0  6.28318530717959E+0000}%
\special{sh 1}%
\special{ar 3400 2600 10 10 0  6.28318530717959E+0000}%
\special{sh 1}%
\special{ar 3200 2600 10 10 0  6.28318530717959E+0000}%
\special{sh 1}%
\special{ar 3000 2600 10 10 0  6.28318530717959E+0000}%
\special{sh 1}%
\special{ar 2800 2600 10 10 0  6.28318530717959E+0000}%
\special{sh 1}%
\special{ar 2800 1400 10 10 0  6.28318530717959E+0000}%
\special{sh 1}%
\special{ar 3000 1400 10 10 0  6.28318530717959E+0000}%
\special{sh 1}%
\special{ar 3200 1400 10 10 0  6.28318530717959E+0000}%
%
\special{pn 8}%
\special{pa 2600 1400}%
\special{pa 2600 1600}%
\special{pa 2800 2600}%
\special{pa 3400 2600}%
\special{pa 3400 2400}%
\special{pa 3200 1400}%
\special{pa 2600 1400}%
\special{fp}%
%
\special{pn 20}%
\special{sh 1}%
\special{ar 3600 1600 10 10 0  6.28318530717959E+0000}%
\special{sh 1}%
\special{ar 3800 1600 10 10 0  6.28318530717959E+0000}%
\special{sh 1}%
\special{ar 4000 1600 10 10 0  6.28318530717959E+0000}%
\special{sh 1}%
\special{ar 4200 1600 10 10 0  6.28318530717959E+0000}%
\special{sh 1}%
\special{ar 3600 1800 10 10 0  6.28318530717959E+0000}%
\special{sh 1}%
\special{ar 3800 2000 10 10 0  6.28318530717959E+0000}%
\special{sh 1}%
\special{ar 4000 2200 10 10 0  6.28318530717959E+0000}%
\special{sh 1}%
\special{ar 4200 2400 10 10 0  6.28318530717959E+0000}%
\special{sh 1}%
\special{ar 4200 2600 10 10 0  6.28318530717959E+0000}%
\special{sh 1}%
\special{ar 4000 2600 10 10 0  6.28318530717959E+0000}%
\special{sh 1}%
\special{ar 3800 2600 10 10 0  6.28318530717959E+0000}%
%
\special{pn 8}%
\special{pa 3600 1600}%
\special{pa 4200 1600}%
\special{pa 4200 2600}%
\special{pa 3800 2600}%
\special{pa 3600 1800}%
\special{pa 3600 1600}%
\special{fp}%
%
\special{pn 20}%
\special{sh 1}%
\special{ar 2200 1200 10 10 0  6.28318530717959E+0000}%
\special{sh 1}%
\special{ar 2000 1000 10 10 0  6.28318530717959E+0000}%
\special{sh 1}%
\special{ar 1800 800 10 10 0  6.28318530717959E+0000}%
\special{sh 1}%
\special{ar 1600 600 10 10 0  6.28318530717959E+0000}%
\special{sh 1}%
\special{ar 1600 400 10 10 0  6.28318530717959E+0000}%
\special{sh 1}%
\special{ar 1800 400 10 10 0  6.28318530717959E+0000}%
\special{sh 1}%
\special{ar 2000 400 10 10 0  6.28318530717959E+0000}%
\special{sh 1}%
\special{ar 2200 400 10 10 0  6.28318530717959E+0000}%
%
\special{pn 20}%
\special{sh 1}%
\special{ar 4200 1400 10 10 0  6.28318530717959E+0000}%
\special{sh 1}%
\special{ar 4000 1200 10 10 0  6.28318530717959E+0000}%
\special{sh 1}%
\special{ar 3800 1000 10 10 0  6.28318530717959E+0000}%
\special{sh 1}%
\special{ar 3600 800 10 10 0  6.28318530717959E+0000}%
\special{sh 1}%
\special{ar 3400 600 10 10 0  6.28318530717959E+0000}%
\special{sh 1}%
\special{ar 3400 400 10 10 0  6.28318530717959E+0000}%
\special{sh 1}%
\special{ar 3600 400 10 10 0  6.28318530717959E+0000}%
\special{sh 1}%
\special{ar 3800 400 10 10 0  6.28318530717959E+0000}%
\special{sh 1}%
\special{ar 4000 400 10 10 0  6.28318530717959E+0000}%
%
\special{pn 8}%
\special{pa 3400 400}%
\special{pa 4000 400}%
\special{pa 4200 1400}%
\special{pa 4200 1600}%
\special{pa 3600 1600}%
\special{pa 3400 600}%
\special{pa 3400 400}%
\special{fp}%
%
\special{pn 20}%
\special{sh 1}%
\special{ar 3200 1200 10 10 0  6.28318530717959E+0000}%
\special{sh 1}%
\special{ar 3000 1000 10 10 0  6.28318530717959E+0000}%
\special{sh 1}%
\special{ar 2800 800 10 10 0  6.28318530717959E+0000}%
\special{sh 1}%
\special{ar 2600 600 10 10 0  6.28318530717959E+0000}%
\special{sh 1}%
\special{ar 2600 400 10 10 0  6.28318530717959E+0000}%
\special{sh 1}%
\special{ar 2800 400 10 10 0  6.28318530717959E+0000}%
\special{sh 1}%
\special{ar 3000 400 10 10 0  6.28318530717959E+0000}%
\put(24.0000,-8.0000){\makebox(0,0){$\cdots$}}%
%
\special{pn 8}%
\special{pa 2600 400}%
\special{pa 3000 400}%
\special{pa 3200 1200}%
\special{pa 3200 1400}%
\special{pa 2600 1400}%
\special{pa 2600 400}%
\special{fp}%
%
\special{pn 8}%
\special{pa 1600 400}%
\special{pa 2200 400}%
\special{pa 2200 1400}%
\special{pa 1600 1400}%
\special{pa 1600 400}%
\special{fp}%
%
\special{pn 8}%
\special{pa 3300 500}%
\special{pa 3500 500}%
\special{pa 3500 2500}%
\special{pa 3300 2500}%
\special{pa 3300 500}%
\special{fp}%
\put(44.0000,-20.0000){\makebox(0,0){$i$}}%
\put(11.5000,-4.0000){\makebox(0,0){$2\tau+1$}}%
%
\special{pn 8}%
\special{pa 800 1400}%
\special{pa 1200 1400}%
\special{pa 1400 2200}%
\special{pa 1400 2400}%
\special{pa 800 2400}%
\special{pa 800 1400}%
\special{fp}%
%
\special{pn 20}%
\special{sh 1}%
\special{ar 800 1400 10 10 0  6.28318530717959E+0000}%
\special{sh 1}%
\special{ar 1000 1400 10 10 0  6.28318530717959E+0000}%
\special{sh 1}%
\special{ar 1200 1400 10 10 0  6.28318530717959E+0000}%
\special{sh 1}%
\special{ar 800 1600 10 10 0  6.28318530717959E+0000}%
\special{sh 1}%
\special{ar 1000 1800 10 10 0  6.28318530717959E+0000}%
\special{sh 1}%
\special{ar 1200 2000 10 10 0  6.28318530717959E+0000}%
\special{sh 1}%
\special{ar 1400 2200 10 10 0  6.28318530717959E+0000}%
\special{sh 1}%
\special{ar 1400 2400 10 10 0  6.28318530717959E+0000}%
\special{sh 1}%
\special{ar 1200 2400 10 10 0  6.28318530717959E+0000}%
\special{sh 1}%
\special{ar 1000 2400 10 10 0  6.28318530717959E+0000}%
\special{sh 1}%
\special{ar 800 2400 10 10 0  6.28318530717959E+0000}%
\end{picture}%

%% file: ntau+2tausym.tex
\unitlength 0.1in
\begin{picture}( 55.9500, 28.8500)(  8.6000,-32.0000)
\put(16.0000,-4.0000){\makebox(0,0){$j$}}%
\put(25.9500,-26.9500){\makebox(0,0){$\tau+1$}}%
%
\special{pn 8}%
\special{pa 996 2596}%
\special{pa 6456 2596}%
\special{fp}%
\special{sh 1}%
\special{pa 6456 2596}%
\special{pa 6388 2576}%
\special{pa 6402 2596}%
\special{pa 6388 2616}%
\special{pa 6456 2596}%
\special{fp}%
%
\special{pn 8}%
\special{pa 1796 3200}%
\special{pa 1796 500}%
\special{fp}%
\special{sh 1}%
\special{pa 1796 500}%
\special{pa 1776 568}%
\special{pa 1796 554}%
\special{pa 1816 568}%
\special{pa 1796 500}%
\special{fp}%
\put(18.4500,-26.9500){\makebox(0,0){$0$}}%
%
\special{pn 8}%
\special{pa 2596 2570}%
\special{pa 2596 2610}%
\special{fp}%
%
\special{pn 20}%
\special{sh 1}%
\special{ar 2806 3000 10 10 0  6.28318530717959E+0000}%
%
\special{pn 20}%
\special{sh 1}%
\special{ar 2206 3000 10 10 0  6.28318530717959E+0000}%
%
\special{pn 20}%
\special{sh 1}%
\special{ar 2406 3000 10 10 0  6.28318530717959E+0000}%
%
\special{pn 20}%
\special{sh 1}%
\special{ar 2606 3000 10 10 0  6.28318530717959E+0000}%
%
\special{pn 20}%
\special{sh 1}%
\special{ar 2598 2598 10 10 0  6.28318530717959E+0000}%
%
\special{pn 20}%
\special{sh 1}%
\special{ar 2798 2798 10 10 0  6.28318530717959E+0000}%
%
\special{pn 20}%
\special{sh 1}%
\special{ar 2596 1800 10 10 0  6.28318530717959E+0000}%
%
\special{pn 20}%
\special{sh 1}%
\special{ar 1996 1800 10 10 0  6.28318530717959E+0000}%
%
\special{pn 20}%
\special{sh 1}%
\special{ar 2396 1800 10 10 0  6.28318530717959E+0000}%
%
\special{pn 20}%
\special{sh 1}%
\special{ar 2196 1800 10 10 0  6.28318530717959E+0000}%
%
\special{pn 20}%
\special{sh 1}%
\special{ar 2396 2400 10 10 0  6.28318530717959E+0000}%
%
\special{pn 20}%
\special{sh 1}%
\special{ar 2196 2200 10 10 0  6.28318530717959E+0000}%
%
\special{pn 20}%
\special{sh 1}%
\special{ar 1996 2000 10 10 0  6.28318530717959E+0000}%
\put(16.9500,-18.0000){\makebox(0,0){$\tau$}}%
%
\special{pn 8}%
\special{pa 1996 1800}%
\special{pa 1996 2000}%
\special{pa 2196 3000}%
\special{pa 2796 3000}%
\special{pa 2796 2800}%
\special{pa 2596 1800}%
\special{pa 1996 1800}%
\special{fp}%
\put(28.4500,-24.0000){\makebox(0,0){$\cdots$}}%
%
\special{pn 20}%
\special{sh 1}%
\special{ar 2996 1800 10 10 0  6.28318530717959E+0000}%
\special{sh 1}%
\special{ar 2996 2000 10 10 0  6.28318530717959E+0000}%
\special{sh 1}%
\special{ar 3196 2200 10 10 0  6.28318530717959E+0000}%
\special{sh 1}%
\special{ar 3396 2400 10 10 0  6.28318530717959E+0000}%
\special{sh 1}%
\special{ar 3596 2600 10 10 0  6.28318530717959E+0000}%
\special{sh 1}%
\special{ar 3796 2800 10 10 0  6.28318530717959E+0000}%
\special{sh 1}%
\special{ar 3796 3000 10 10 0  6.28318530717959E+0000}%
\special{sh 1}%
\special{ar 3596 3000 10 10 0  6.28318530717959E+0000}%
\special{sh 1}%
\special{ar 3396 3000 10 10 0  6.28318530717959E+0000}%
\special{sh 1}%
\special{ar 3196 3000 10 10 0  6.28318530717959E+0000}%
\special{sh 1}%
\special{ar 3196 1800 10 10 0  6.28318530717959E+0000}%
\special{sh 1}%
\special{ar 3396 1800 10 10 0  6.28318530717959E+0000}%
\special{sh 1}%
\special{ar 3596 1800 10 10 0  6.28318530717959E+0000}%
%
\special{pn 8}%
\special{pa 2996 1800}%
\special{pa 2996 2000}%
\special{pa 3196 3000}%
\special{pa 3796 3000}%
\special{pa 3796 2800}%
\special{pa 3596 1800}%
\special{pa 2996 1800}%
\special{fp}%
%
\special{pn 20}%
\special{sh 1}%
\special{ar 3996 2000 10 10 0  6.28318530717959E+0000}%
\special{sh 1}%
\special{ar 4196 2000 10 10 0  6.28318530717959E+0000}%
\special{sh 1}%
\special{ar 4396 2000 10 10 0  6.28318530717959E+0000}%
\special{sh 1}%
\special{ar 4596 2000 10 10 0  6.28318530717959E+0000}%
\special{sh 1}%
\special{ar 3996 2200 10 10 0  6.28318530717959E+0000}%
\special{sh 1}%
\special{ar 4196 2400 10 10 0  6.28318530717959E+0000}%
\special{sh 1}%
\special{ar 4396 2600 10 10 0  6.28318530717959E+0000}%
\special{sh 1}%
\special{ar 4596 2800 10 10 0  6.28318530717959E+0000}%
\special{sh 1}%
\special{ar 4596 3000 10 10 0  6.28318530717959E+0000}%
\special{sh 1}%
\special{ar 4396 3000 10 10 0  6.28318530717959E+0000}%
\special{sh 1}%
\special{ar 4196 3000 10 10 0  6.28318530717959E+0000}%
%
\special{pn 8}%
\special{pa 3996 2000}%
\special{pa 4596 2000}%
\special{pa 4596 3000}%
\special{pa 4196 3000}%
\special{pa 3996 2200}%
\special{pa 3996 2000}%
\special{fp}%
%
\special{pn 20}%
\special{sh 1}%
\special{ar 2596 1600 10 10 0  6.28318530717959E+0000}%
\special{sh 1}%
\special{ar 2396 1400 10 10 0  6.28318530717959E+0000}%
\special{sh 1}%
\special{ar 2196 1200 10 10 0  6.28318530717959E+0000}%
\special{sh 1}%
\special{ar 1996 1000 10 10 0  6.28318530717959E+0000}%
\special{sh 1}%
\special{ar 1996 800 10 10 0  6.28318530717959E+0000}%
\special{sh 1}%
\special{ar 2196 800 10 10 0  6.28318530717959E+0000}%
\special{sh 1}%
\special{ar 2396 800 10 10 0  6.28318530717959E+0000}%
\special{sh 1}%
\special{ar 2596 800 10 10 0  6.28318530717959E+0000}%
%
\special{pn 20}%
\special{sh 1}%
\special{ar 4596 1800 10 10 0  6.28318530717959E+0000}%
\special{sh 1}%
\special{ar 4396 1600 10 10 0  6.28318530717959E+0000}%
\special{sh 1}%
\special{ar 4196 1400 10 10 0  6.28318530717959E+0000}%
\special{sh 1}%
\special{ar 3996 1200 10 10 0  6.28318530717959E+0000}%
\special{sh 1}%
\special{ar 3796 1000 10 10 0  6.28318530717959E+0000}%
\special{sh 1}%
\special{ar 3796 800 10 10 0  6.28318530717959E+0000}%
\special{sh 1}%
\special{ar 3996 800 10 10 0  6.28318530717959E+0000}%
\special{sh 1}%
\special{ar 4196 800 10 10 0  6.28318530717959E+0000}%
\special{sh 1}%
\special{ar 4396 800 10 10 0  6.28318530717959E+0000}%
%
\special{pn 8}%
\special{pa 3796 800}%
\special{pa 4396 800}%
\special{pa 4596 1800}%
\special{pa 4596 2000}%
\special{pa 3996 2000}%
\special{pa 3796 1000}%
\special{pa 3796 800}%
\special{fp}%
%
\special{pn 20}%
\special{sh 1}%
\special{ar 3596 1600 10 10 0  6.28318530717959E+0000}%
\special{sh 1}%
\special{ar 3396 1400 10 10 0  6.28318530717959E+0000}%
\special{sh 1}%
\special{ar 3196 1200 10 10 0  6.28318530717959E+0000}%
\special{sh 1}%
\special{ar 2996 1000 10 10 0  6.28318530717959E+0000}%
\special{sh 1}%
\special{ar 2996 800 10 10 0  6.28318530717959E+0000}%
\special{sh 1}%
\special{ar 3196 800 10 10 0  6.28318530717959E+0000}%
\special{sh 1}%
\special{ar 3396 800 10 10 0  6.28318530717959E+0000}%
\put(27.9500,-12.0000){\makebox(0,0){$\cdots$}}%
%
\special{pn 8}%
\special{pa 2996 800}%
\special{pa 3396 800}%
\special{pa 3596 1600}%
\special{pa 3596 1800}%
\special{pa 2996 1800}%
\special{pa 2996 800}%
\special{fp}%
%
\special{pn 8}%
\special{pa 1996 800}%
\special{pa 2596 800}%
\special{pa 2596 1800}%
\special{pa 1996 1800}%
\special{pa 1996 800}%
\special{fp}%
%
\special{pn 8}%
\special{pa 3696 900}%
\special{pa 3896 900}%
\special{pa 3896 2900}%
\special{pa 3696 2900}%
\special{pa 3696 900}%
\special{fp}%
\put(63.9500,-24.0000){\makebox(0,0){$i$}}%
\put(15.5000,-6.0000){\makebox(0,0){$2\tau+1$}}%
\put(45.9500,-12.0000){\makebox(0,0){$\cdots$}}%
\put(48.4500,-24.0000){\makebox(0,0){$\cdots$}}%
%
\special{pn 20}%
\special{sh 1}%
\special{ar 4996 3000 10 10 0  6.28318530717959E+0000}%
\special{sh 1}%
\special{ar 5196 3000 10 10 0  6.28318530717959E+0000}%
\special{sh 1}%
\special{ar 5396 3000 10 10 0  6.28318530717959E+0000}%
\special{sh 1}%
\special{ar 5596 3000 10 10 0  6.28318530717959E+0000}%
\special{sh 1}%
\special{ar 5596 2800 10 10 0  6.28318530717959E+0000}%
\special{sh 1}%
\special{ar 5396 2600 10 10 0  6.28318530717959E+0000}%
\special{sh 1}%
\special{ar 5196 2400 10 10 0  6.28318530717959E+0000}%
\special{sh 1}%
\special{ar 4996 2200 10 10 0  6.28318530717959E+0000}%
\special{sh 1}%
\special{ar 4996 2000 10 10 0  6.28318530717959E+0000}%
\special{sh 1}%
\special{ar 5196 2000 10 10 0  6.28318530717959E+0000}%
\special{sh 1}%
\special{ar 5396 2000 10 10 0  6.28318530717959E+0000}%
\special{sh 1}%
\special{ar 5596 2000 10 10 0  6.28318530717959E+0000}%
%
\special{pn 20}%
\special{sh 1}%
\special{ar 5596 1800 10 10 0  6.28318530717959E+0000}%
\special{sh 1}%
\special{ar 5396 1600 10 10 0  6.28318530717959E+0000}%
\special{sh 1}%
\special{ar 5196 1400 10 10 0  6.28318530717959E+0000}%
\special{sh 1}%
\special{ar 4996 1200 10 10 0  6.28318530717959E+0000}%
\special{sh 1}%
\special{ar 4796 1000 10 10 0  6.28318530717959E+0000}%
\special{sh 1}%
\special{ar 4796 800 10 10 0  6.28318530717959E+0000}%
\special{sh 1}%
\special{ar 4996 800 10 10 0  6.28318530717959E+0000}%
\special{sh 1}%
\special{ar 5196 800 10 10 0  6.28318530717959E+0000}%
\special{sh 1}%
\special{ar 5396 800 10 10 0  6.28318530717959E+0000}%
%
\special{pn 8}%
\special{pa 4796 800}%
\special{pa 5396 800}%
\special{pa 5596 1800}%
\special{pa 5596 2000}%
\special{pa 4996 2000}%
\special{pa 4796 1000}%
\special{pa 4796 800}%
\special{fp}%
%
\special{pn 8}%
\special{pa 4996 2000}%
\special{pa 4996 3000}%
\special{pa 5596 3000}%
\special{pa 5596 2000}%
\special{pa 4996 2000}%
\special{fp}%
\put(59.9500,-6.0000){\makebox(0,0){$B_{2(n'-\gamma''),1}$}}%
%
\special{pn 20}%
\special{sh 1}%
\special{ar 1596 1600 10 10 0  6.28318530717959E+0000}%
\special{sh 1}%
\special{ar 1396 1400 10 10 0  6.28318530717959E+0000}%
\special{sh 1}%
\special{ar 1196 1200 10 10 0  6.28318530717959E+0000}%
\special{sh 1}%
\special{ar 996 1000 10 10 0  6.28318530717959E+0000}%
\special{sh 1}%
\special{ar 996 800 10 10 0  6.28318530717959E+0000}%
\special{sh 1}%
\special{ar 1196 800 10 10 0  6.28318530717959E+0000}%
\special{sh 1}%
\special{ar 1396 800 10 10 0  6.28318530717959E+0000}%
\special{sh 1}%
\special{ar 1596 800 10 10 0  6.28318530717959E+0000}%
%
\special{pn 8}%
\special{pa 996 800}%
\special{pa 996 1000}%
\special{pa 1196 1800}%
\special{pa 1596 1800}%
\special{pa 1596 800}%
\special{pa 996 800}%
\special{fp}%
\put(9.9500,-18.0000){\makebox(0,0){$0$}}%
%
\special{pn 20}%
\special{sh 1}%
\special{ar 1600 1800 10 10 0  6.28318530717959E+0000}%
\special{sh 1}%
\special{ar 1400 1800 10 10 0  6.28318530717959E+0000}%
\special{sh 1}%
\special{ar 1200 1800 10 10 0  6.28318530717959E+0000}%
\special{sh 1}%
\special{ar 1200 2000 10 10 0  6.28318530717959E+0000}%
\special{sh 1}%
\special{ar 1400 2200 10 10 0  6.28318530717959E+0000}%
\special{sh 1}%
\special{ar 1600 2400 10 10 0  6.28318530717959E+0000}%
\special{sh 1}%
\special{ar 1800 2600 10 10 0  6.28318530717959E+0000}%
\special{sh 1}%
\special{ar 1800 2800 10 10 0  6.28318530717959E+0000}%
\special{sh 1}%
\special{ar 1600 2800 10 10 0  6.28318530717959E+0000}%
\special{sh 1}%
\special{ar 1400 2800 10 10 0  6.28318530717959E+0000}%
\special{sh 1}%
\special{ar 1200 2800 10 10 0  6.28318530717959E+0000}%
%
\special{pn 8}%
\special{pa 1600 1800}%
\special{pa 1800 2600}%
\special{pa 1800 2800}%
\special{pa 1200 2800}%
\special{pa 1200 1800}%
\special{pa 1600 1800}%
\special{fp}%
%
\special{pn 8}%
\special{sh 1}%
\special{ar 6400 1000 10 10 0  6.28318530717959E+0000}%
%
\special{pn 8}%
\special{pa 5800 1000}%
\special{pa 6400 1000}%
\special{pa 6400 2000}%
\special{pa 6000 2000}%
\special{pa 5800 1200}%
\special{pa 5800 1000}%
\special{fp}%
%
\special{pn 20}%
\special{sh 1}%
\special{ar 5800 1000 10 10 0  6.28318530717959E+0000}%
\special{sh 1}%
\special{ar 6000 1000 10 10 0  6.28318530717959E+0000}%
\special{sh 1}%
\special{ar 6200 1000 10 10 0  6.28318530717959E+0000}%
\special{sh 1}%
\special{ar 6400 1000 10 10 0  6.28318530717959E+0000}%
\special{sh 1}%
\special{ar 5800 1200 10 10 0  6.28318530717959E+0000}%
\special{sh 1}%
\special{ar 6000 1400 10 10 0  6.28318530717959E+0000}%
\special{sh 1}%
\special{ar 6200 1600 10 10 0  6.28318530717959E+0000}%
\special{sh 1}%
\special{ar 6400 1800 10 10 0  6.28318530717959E+0000}%
\special{sh 1}%
\special{ar 6400 2000 10 10 0  6.28318530717959E+0000}%
\special{sh 1}%
\special{ar 6200 2000 10 10 0  6.28318530717959E+0000}%
\special{sh 1}%
\special{ar 6000 2000 10 10 0  6.28318530717959E+0000}%
%
\special{pn 8}%
\special{pa 1780 800}%
\special{pa 1820 800}%
\special{fp}%
\put(42.0000,-6.0000){\makebox(0,0){$B_{n'-\gamma'',1}$}}%
\end{picture}%